\documentclass[11pt,leqno]{article}
\usepackage{a4wide}
\usepackage{amssymb}
\usepackage{amsmath}
\usepackage{amsthm}
\usepackage{mathrsfs}
\usepackage{yfonts}
\usepackage{euscript}
\usepackage{upgreek}
\usepackage[all]{xypic}
\usepackage[usenames]{color}
\usepackage[usenames,dvipsnames]{xcolor} 
\usepackage{mathabx} 
\usepackage[ngerman,shadow,color=blue!30,textsize=footnotesize,colorinlistoftodos]{todonotes}
\usepackage[T1]{fontenc}
\usepackage[utf8]{inputenc}
\usepackage{comment}
\usepackage{xr-hyper}
\usepackage{enotez, translations}
\usepackage[unicode,psdextra]{hyperref}
\usepackage{bookmark}

\newcommand{\com}[1]{{#1}}                
\newcommand{\comb}[1]{{#1}}             
\newcommand{\Footnote}[1]{\footnote{{\color{Sepia} #1}}}  
\renewcommand{\Footnote}[2][]{\relax}  
\theoremstyle{plain}

\newcommand{\Aff}{\operatorname{Aff}}
\newcommand{\res}{\operatorname{res}}
\newcommand{\coker}{\operatorname{coker}}

\newcommand{\id}{\operatorname{id}}
\newcommand{\im}{\operatorname{im}}
\newcommand{\ev}{\operatorname{ev}}

\newcommand{\pr}{\operatorname{pr}}
\newcommand{\divi}{\operatorname{div}}

\newcommand{\Sp}{\operatorname{Sp}}

\newcommand{\Hom}{\operatorname{Hom}}
\newcommand{\End}{\operatorname{End}}

\newcommand{\Fil}{\operatorname{Fil}}
\newcommand{\Mod}{\operatorname{Mod}}
\newcommand{\Lie}{\operatorname{Lie}}

\newcommand{\Rep}{\operatorname{Rep}}

\newcommand{\Aut}{\operatorname{Aut}}

\newcommand{\Gal}{\operatorname{Gal}}
\newcommand{\rk}{\operatorname{rk}}

 \newcommand{\CC}{\mathbb{C}}

 \newcommand{\QQ}{\mathbb{Q}}

 \renewcommand{\cL}{\mathcal{L}}

 \newcommand{\cC}{\mathcal{C}}
 \newcommand{\cO}{\mathcal{O}}

\newcommand{\cE}{\mathcal{E}}
\renewcommand{\cR}{\mathcal{R}}

 \newcommand{\vp}{\varphi}

 \newcommand{\Qp}{\QQ_p}

 \newcommand{\Cp}{\CC_p}

 \newcommand{\be}{\begin{equation}}
\newcommand{\ee}{\end{equation}}

\newcommand{\fpbar}{\ifmmode {\overline{\mathbb{F}_p}}\else$\mathbb{F}_p$\ \fi}

\newcommand{\fp}{\ifmmode {\mathbb{F}_p}\else$\mathbb{F}_p$\ \fi}
\newcommand{\zp}{\ifmmode \mathbb{Z}_p\else$\mathbb{Z}_p$\ \fi}
\newcommand{\zpur}{\ifmmode \widehat{\zp^{ur}}\else $\widehat{\zp^{ur}}$\ \fi}

 \newcommand{\TLT}{T_{\pi}}

 \newcommand{\rig}{\mathrm{rig}}

\renewcommand{\H}{\mathrm{H}}

\newcommand{\kkl}{<}
\newcommand{\kkr}{>}

\newcommand{\KKl}{\{}
\newcommand{\KKr}{\}}

\newcommand{\La}{\ifmmode\Lambda\else$\Lambda$\fi}

\newcommand{\q}{\ifmmode {\mathbb Q}\else${\mathbb Q}$\ \fi}
\newcommand{\qp}{\ifmmode {\mathbb Q}_p\else${\mathbb Q}_p$\ \fi}
\newcommand{\z}{\mathbb{Z}}

\newcommand{\Q}{\ifmmode {\mathbb Q}\else${\mathbb Q}$\ \fi}
\newcommand{\ql}{\ifmmode {{\mathbb Q}_l}\else${\mathbb Q}_l$\ \fi}

\newtheorem{theorem}{Theorem}[subsection]
\newtheorem{corollary}[theorem]{Corollary}
\newtheorem{lemma}[theorem]{Lemma}
\newtheorem{remark}[theorem]{Remark}
\newtheorem{proposition}[theorem]{Proposition}

\newtheorem{definition}[theorem]{Definition}

\newtheorem{question}[theorem]{Question}

\newtheorem{btheorem}{Theorem}

\theoremstyle{remark}

\author{Peter Schneider and Otmar Venjakob}
\begin{document}

\title{\textbf{Reciprocity laws for $(\varphi_L,\Gamma_L)$-modules over Lubin-Tate extensions}}

\maketitle

\begin{abstract}
 In the Lubin-Tate setting we study pairings for analytic $(\varphi_L,\Gamma_L)$-modules and prove an abstract reciprocity law which then implies a relation between the analogue of Perrin-Riou's Big Exponential map as developed by Berger and Fourquaux and a $p$-adic regulator map whose construction relies on the theory of Kisin-Ren modules generalising the concept of Wach modules to the Lubin-Tate situation.
\end{abstract}

\tableofcontents

\section{Introduction}

Classically explicit reciprocity laws or formulas usually mean an explicit computation of Hilbert symbols or (local) cup products using e.g.\ differential forms, (Coleman) power series etc. and    a bunch of manifestations of this idea exists in the literature due to Artin-Hasse, Iwasawa, Wiles, Kolyvagin, Vostokov, Br\"{u}ckner, Coleman, Sen, de Shalit, Fesenko, Bloch-Kato, Benois ... In the same spirit Perrin-Riou's reciprocity law gives an explicit calculation of the Iwasawa cohomology pairing in terms of big exponential and regulator maps for crystalline representations of $G_{\mathbb{Q}_p}$; more precisely, the latter maps are adjoint to each other when also involving the crystalline duality paring after base change to the distribution algebra corresponding to the cyclotomic situation.

The motivating question for this article is to investigate what happens if one replaces the cyclotomic $\mathbb{Z}_p$-extension by a Lubin-Tate extension $L_\infty$ over some finite extension $L$ over $\mathbb{Q}_p$ with Galois group $\Gamma_L=G(L_\infty/L)$ and Lubin-Tate character $\chi_{LT}:G_L\to o_L^\times$ which all arise from a Lubin-Tate formal group attached to a prime $\pi_L\in o_L$  the additive group of the ring of integers $o_L$ of $L$; by $q$ we denote the cardinality of the residue field  $o_L/o_L\pi_L.$ We try to extend the above sketched cyclotomic picture to the Lubin-Tate case at least for {\it $L$-analytic} crystalline representations of the absolute Galois group $G_L$ of $L.$ As pointed out in \cite{SV15} already, the character $\tau := \chi_{cyc}\cdot\chi_{LT}^{-1}$ plays a crucial role.

To this aim we study $(\varphi_L,\Gamma_L)$-modules over different Robba rings with coefficients in suitable complete intermediate fields $L\subseteq K\subseteq\mathbb{C}_p$.  The starting point is the theory of Schneider and Teitelbaum: In  \cite{ST2} they introduced the rigid analytic group variety $\mathfrak{X}$ over $L$, which parametrizes the locally $L$-analytic characters of $o_L$, and similarly $\mathfrak{X}^\times\cong\mathfrak{X}_{\Gamma_L}$ for the locally $L$-analytic groups $o_L^\times\cong\Gamma_L$, the isomorphisms being induced by (the inverse of) $\chi_{LT}$. Under the assumption that   the period $\Omega$ of the dual of the fixed Lubin-Tate group belongs to $K$ they establish an isomorphism $\kappa: \mathbf{B}_K\cong \mathfrak{X}_K$ of rigid analytic varieties over $K,$ called the
Lubin-Tate isomorphism, where $\mathbf{B}$ denotes the rigid analytic open unit disk and the index $K$ indicates base change to $K$. In sections \ref{sec:basics}, \ref{sec:usualRobba}, \ref{sec:groupRobba} we recall or introduce the Robba rings $\cR_K(\mathfrak{Y})$ for all the above varieties $\mathfrak{Y}.$ We call $\cR_K(\Gamma_L):=\cR_K(\mathfrak{X}_{\Gamma_L})$ also the Robba group ring as we can consider it as an extension of the locally $L$-analytic distribution algebra $D(\Gamma_L,K)$ with coefficients in $K$ as follows: The Fourier isomorphism $D(o_L,K)\cong {\mathcal{O}_K(\mathfrak{X})}$ onto the ring of holomorphic functions on $\mathfrak{X} $ induces the Mellin-transform, i.e.,\ a topological isomorphism between $D(\Gamma_L,K)\cong D(o_L^\times,K)$ and the $D(o_L^\times,K)$-submodule $({\mathcal{O}_K(\mathfrak{X})})^{\psi_L^\mathfrak{X}=0}$ of ${\mathcal{O}_K(\mathfrak{X})}$ on which the $\psi_L^\mathfrak{X}$-operator - to be recalled in section \ref{sec:add-Robba} - acts as zero. As a special case of the following theorem we  extend the Mellin transform to an isomorphism of $\cR_K(\Gamma_L)$ and $(\mathcal{R}_K(\mathfrak{X}))^{\psi_L=0}$.
\begin{btheorem}[Theorem \ref{thm:Mpsiequal0-decent}]
If $M$ denotes a $L$-analytic $(\varphi_L,\Gamma_L)$-module over $\mathcal{R}_K(\mathfrak{X})$ for any complete intermediate field $L\subseteq K\subseteq\mathbb{C}_p$, then $M^{\psi_L=0}$ is a free $\cR_K(\Gamma_L)$-module of rank $\rk_{\mathcal{R}_K(\mathfrak{X})}M.$
\end{btheorem}

For $\mathbf{B}$ instead of $\mathfrak{X}$ an analogous statement holds, if $K$ contains $\Omega;$ technically, this is the case we prove first (see Theorem \ref{thm:Mpsiequal0}) and which then, after involving the Lubin-Tate isomorphism,  descends to the Theorem. Under this condition on $K$ we may illustrate that  via Fourier theory and the Lubin-Tate isomorpshism the locally $L$-analytic distribution algebra $D(o_L,K)$  becomes isomorphic to  the subring
  ${\mathcal{O}_K(\mathbf{B})}\subseteq\cR_K(\mathbf{B})$  consisting of those functions which converge on the full open unit disk, while the functions in $\cR_K(\mathbf{B})$ in general only converge on some annulus $r\leq |Z|<1$ for some radius $0<r<1.$ This isomorphism induces the Mellin-transform, i.e.,\ a topological isomorphism between $D(o_L^\times,K)$ and the $D(o_L^\times,K)$-submodule $({\mathcal{O}_K(\mathbf{B})})^{\psi_L=0}$ of ${\mathcal{O}_K(\mathbf{B})}$ on which the $\psi_L$-operator - up to a scalar a left inverse of the Lubin-Tate $\varphi_L$-operator - acts as zero.

A second ingredient is Serre duality on the above rigid analytic varieties $\mathfrak{Y}$, which induces - as developed in this generality in section \ref{sec:duality} -  a residue pairing
\[ \Omega^1_{\mathcal{R}_L(\mathfrak{Y})}\times \mathcal{R}_L(\mathfrak{Y})\to K\] in \ref{Serre2} for the differentials $\Omega^1_{\mathcal{R}_L(\mathfrak{Y})}$ and also a pairing
\begin{align*}
  < \ ,\ >_{\mathfrak{Y}} \, : \mathcal{R}_K(\mathfrak{Y}) \times \mathcal{R}_K(\mathfrak{Y}) & \longrightarrow K  ,
\end{align*}
see \eqref{f:pairing-RobbaX}, \eqref{f:pairing-mult-charvar}. For $\mathfrak{Y}=\mathfrak{X}_{\Gamma_L}$ the latter  induces topological isomorphisms \[\Hom_{K,cts}(\cR_K(\Gamma_L),K) \cong\cR_K(\Gamma_L)\;\;\mbox{ and  } \;\Hom_{K,cts}(\cR_K(\Gamma_L)/D(\Gamma_L,K),K)\cong D(\Gamma_L,K).\] For a $L$-analytic $(\varphi_L,\Gamma_L)$-module $M$ over $\cR:=\cR_K(\mathfrak{Y})$ with $\mathfrak{Y}$ either equal to $\mathfrak{X}$ or $\mathbf{B},$  we finally use these isomorphisms to define on the one hand the two Iwasawa pairings
 \[ \KKl\;,\;\KKr_{M,Iw}^0:\check{M}^{\psi_L=0}\times M^{\psi_L=0} \to
 \cR_K(\Gamma_L)\]
and
\begin{equation*}
\{\;,\;\}_{M,Iw}: \check{M}^{\psi_L=\frac{q}{\pi_L}} \times M^{\psi_L=1} \to D(\Gamma_L,K),
\end{equation*}
where $\check{M}:=\Hom_\cR(M,\Omega^1_\cR).$ They are linked by the commutative diagram
\begin{equation*}
\xymatrix{
     {\KKl\;,\;\KKr_{M,Iw}:  {\check{M}^{\psi_L=\frac{q}{\pi_L}}}}\phantom{ [ , ]: }  \ar[d]_{ \varphi_L -1}\ar@{}[r]|{\times} &  {M^{\psi_L=1}}
     \ar[d]_{ \frac{\pi_L}{q}\varphi_L -1} \ar@{->}[r] & D(\Gamma_L,K)\ar@{^(->}[d]\\
    {\KKl\;,\;\KKr_{M,Iw}^0: \check{M}^{\psi_L=0}\phantom{\KKl\;,\;\KKr_{Iw}^0dd }} \ar@{}[r]|{\times} & M^{\psi_L=0}  \ar[r]^{} &\cR_K(\Gamma_L) . }
\end{equation*}

Now assume that $M$ arises as $D^\dagger_{rig}(W)$ under Berger's equivalence of categories, if $\mathfrak{Y}=\mathbf{B}$, and as $D^\dagger_{rig}(W)_\mathfrak{X}$ under the equivalence from \cite{BSX}, if $\mathfrak{Y}=\mathfrak{X}$, (see Thm.\  \ref{thm:BergerEquiv}) from an $L$-analytic, crystalline representation $W$ of $G_L,$ whence  $\check{M}\cong D^\dagger_{rig}(W^*(\chi_{LT}))$ and  $\check{M}\cong D^\dagger_{rig}(W^*(\chi_{LT}))_\mathfrak{X}$, respectively.
Then, on the other hand  we obtain the pairing
\[[\;,\;]_{D_{cris,L}(W)}:\cR^{\psi_L=0}\otimes_L D_{cris,L}(W^*(\chi_{LT})) \times  \cR^{\psi_L=0}\otimes_L D_{cris,L}(W)  \to \cR_K(\Gamma_L)\]
by base extension of the usual crystalline duality pairing  - if $\mathfrak{Y}=\mathbf{B}$  assuming $\Omega\in K$ -, see \eqref{f:def[]}.
The work of Kisin-Ren and Berger-Schneider-Xie, respectively, provides  comparison isomorphisms
\begin{align*}
\mathrm{comp}_M:M[\frac{1}{t_{\mathfrak{Y}}}]\cong  \cR[\frac{1}{t_{\mathfrak{Y}}}]\otimes_L D_{cris,L}(W)
\end{align*}
and
\begin{align*}
\mathrm{comp}_{\check{M}}:\check{M}[\frac{1}{t_{\mathfrak{Y}}}]\cong  \cR[\frac{1}{t_{\mathfrak{Y}}}]\otimes_L D_{cris,L}(W^*(\chi_{LT})).
\end{align*}
Here $ t_\mathbf{B}:=t_{LT}:=\log_{LT}(Z)\in\cR$ denotes the Lubin-Tate period which stems from the Lubin-Tate logarithm while $t_\mathfrak{X}=\log_\mathfrak{X}$ as defined before Remark \ref{rem:logXlogLT}.  The Lubin-Tate character $\chi_{LT}$ induces  isomorphism $\Gamma_L \xrightarrow{\cong} o_L^\times$ as well as $Lie(\Gamma_L) \xrightarrow{\cong} L$, and we let $\nabla\in Lie(\Gamma_L)$ be the  preimage of $1.$ Then the abstract reciprocity law we prove is the following statement.
\begin{btheorem}[Theorem \ref{thm-recproclawKKK}]
For all $x\in   \check{M}^{\psi_L=0} $ and $y\in   {M}^{\psi_L=0}$, for which the crystalline pairing is defined via the comparison isomorphism, it
holds
\[\frac{q-1}{q}\KKl  {\nabla}x,y\KKr^0_{M,Iw}=[x,y]_{D_{cris,L}(W)},\]
if $\mathfrak{Y}=\mathfrak{X}$, while the analogous statement for $\mathfrak{Y}=\mathbf{B}$ holds upon assuming $\Omega\in K$.
\end{btheorem}
As explained in more detail at the beginning of section \ref{sec:pairings} the proof of this abstract reciprocity law is mainly based on the insight, how the residue maps of $\mathfrak{X}$ and $\mathfrak{X}^\times$ and hence their associated pairings $< \ ,\ >_{\mathfrak{X}}$ and $  < \ ,\ >_{\mathfrak{X}^\times} $ are related to each other by Theorem \ref{thm:residuumidentity} in subsection \ref{subsec:res-alternative}.

As an application for $\mathfrak{Y}=\mathbf{B}$ we show in section \ref{sec:application} the adjointness of big exponential and regulator maps. Recall that Berger and Fourquaux have constructed
for $V$   an $L$-analytic representation of $G_L$ and  an integer $h\geq 1$   such that
\begin{itemize}
  \item $\mathrm{Fil}^{-h}D_{cris,L}(V)=D_{cris,L}(V)$ and
  \item $D_{cris,L}(V)^{\varphi_L=\pi_L^{-h}}=0$
\end{itemize}   a {\it big exponential map} \`{a} la Perrin-Riou
\begin{align*}
\Omega_{V,h}:({\mathcal{O}_K(\mathbf{B})})^{\psi_L=0}\otimes_L D_{cris,L}(V)&\to D_{\mathrm{rig}}^\dagger(V)^{\psi_L=\frac{q}{\pi_L}},
\end{align*}
which up to comparison isomorphism  is for  $h=1$ given by $
f=(1-\varphi_L)x\mapsto \nabla x$ and which interpolates Bloch-Kato exponential maps $exp_{L,V(\chi_{LT}^r)}$.

On the other hand, based on an extension of the work of Kisin and Ren in the first section, we  construct
for a lattice $T\subseteq V,$    such that $V(\tau^{-1})$ is $L$-analytic and crystalline and such that $V$ does not have any quotient isomorphic to
$L(\tau),$     a {\it regulator map} \`{a} la Loeffler and Zerbes
\[\mathcal{L}_V^0: H^1_{Iw}(L_\infty/L,T)\cong D_{LT}(T(\tau^{-1}))^{\psi_L=1}\to   ({\mathcal{O}_K(\mathbf{B})}_{})^{\psi_L=0}\otimes_L
D_{cris,L}(V(\tau^{-1}))\]  as applying the operator
\[{1-\frac{\pi_L}{q}\varphi_L}\]
 up to comparison isomorphism.  Then we derive from the abstract version above with $W=V(\tau^{-1})$ the following reciprocity formula
\begin{btheorem}[Theorem \ref{thm:adjointness}] Assume that $V^*(1)$ is $L$-analytic.  If $\mathrm{Fil}^{-1}D_{cris,L}(V^*(1))=D_{cris,L}(V^*(1))$ and
$D_{cris,L}(V^*(1))^{\varphi_L=\linebreak\pi_L^{-1}}=D_{cris,L}(V^*(1))^{\varphi_L=1}=0$, then the following diagram commutes:
\begin{equation*}\small
\xymatrix@=0.7em@R=1.7em{
    {D_{\mathrm{rig}}^\dagger(V^*(1))^{\psi_L=\frac{q}{\pi_L}}}\ar@{}[r]|{\times} &  {D(V(\tau^{-1}))^{\psi_L=1}}\ar[d]_{\mathcal{L}_V^0}
    \ar[r]^(0.6){\{,\}_{Iw}} & D(\Gamma_L,\Cp) \ar@{=}[d] \\
    ({\mathcal{O}_K(\mathbf{B})})^{\psi_L=0}\otimes_L D_{cris,L}(V^*(1)) \ar[u]_{\Omega_{V^*(1),1}}   \ar@{}[r]|{\times} &  ({\mathcal{O}_K(\mathbf{B})})^{\psi_L=0}\otimes_L
    D_{cris,L}(V(\tau^{-1})) \ar[r]^-{[,]_{}} & D(\Gamma_L,\Cp).}
\end{equation*}
\end{btheorem}
While the crystalline pairing satisfies an interpolation property (Prop.\  \ref{prop:cris}) for trivial reasons, the statement that the second Iwasawa pairing interpolates Tate's cup product pairing is more subtle (Prop.\  \ref{prop:Tate}). Eventually the interpolation property of Berger and Fourquaux for $\Omega_{V,h}$ combined with the adjointness of the latter with $\mathcal{L}_V^0$ implies an
interpolation formula for the regulator map, which interpolates dual Bloch-Kato exponential maps, see Thm.\  \ref{thm:adjointformula}.\\

{\bf Acknowledgements:} We thank Rustam Steingart for discussions about section \ref{subsec:phiGamma}. In his thesis \cite{Stein} he has generalized Theorem \ref{thm:Mpsiequal0} to $(\varphi_L,\Gamma_L)$-modules over families. Both authors are grateful to UBC and PIMS at Vancouver for supporting a fruitful stay.
The project was funded by the Deutsche Forschungsgemeinschaft (DFG, German Research Foundation) under  SFB 1442, {\it Geometry: Deformations and Rigidity}, project-ID 427320536, Germany’s Excellence Strategy EXC 2044 390685587, {\it Mathematics M\"{u}nster: Dynamics–Geometry–Structure}, TRR 326, {\it Geometry and Arithmetic of Uniformized Structures}, project-ID 444845124, and DFG-Forschergruppe award number [1920], {\it Symmetrie, Geometrie und Arithmetik}.

\section{Notation}

Let $\mathbb{Q}_p \subseteq L \subset \mathbb{C}_p$ be a field of finite degree $d$ over $\mathbb{Q}_p$, $o_L$ the ring of integers of $L$, $\pi_L \in o_L$ a fixed prime element, $k_L = o_L/\pi_L o_L$ the residue field,  $q := |k_L|$ and $e$ the absolute ramification index of $L$. We always use the absolute value $|\ |$ on $\mathbb{C}_p$ which is normalized by $|\pi_L| = q^{-1}$. We \textbf{warn} the reader, though, that we will repeatedly use the references \cite{BSX}, \cite{FX}, \cite{La}, \cite{Sc1}, \cite{Sc}, \cite{ST}, and \cite{ST2} in which the absolute value is normalized differently from this paper by $|p| = p^{-1}$. Our absolute value is the $d$th power of the one in these references. The transcription of certain formulas to our convention will usually be done silently.

We fix a Lubin-Tate formal $o_L$-module $LT = LT_{\pi_L}$ over $o_L$ corresponding to the prime element $\pi_L$. We always identify $LT$ with the open unit disk around zero, which gives us a global coordinate $Z$ on $LT$. The $o_L$-action then is given by formal power series $[a](Z) \in o_L[[Z]]$. For simplicity the formal group law will be denoted by $+_{LT}$.

The power series $\frac{\partial (X +_{LT} Y)}{\partial Y}_{|(X,Y) = (Z,0)}$ is a unit in $o_L[[Z]]$ and we let $g_{LT}(Z)$ denote its inverse. Then $g_{LT}(Z) dZ$ is, up to scalars, the unique invariant differential form on $LT$ (\cite[\S5.8]{Haz}). We also let \begin{equation}\label{f:tLT}
  \log_{LT}(Z) = Z + \ldots
\end{equation}
denote the unique formal power series in $L[[Z]]$ whose formal derivative is $g_{LT}$. This $\log_{LT}$ is the logarithm of $LT$ (\cite[8.6]{Lan}). In particular, $g_{LT}dZ = d\log_{LT}$. The invariant derivation $\partial_\mathrm{inv}$ corresponding to the form $d\log_{LT}$ is determined by
\begin{equation*}
  f' dZ = df = \partial_\mathrm{inv} (f) d\log_{LT} = \partial_\mathrm{inv} (f) g_{LT} dZ
\end{equation*}
and hence is given by
\begin{equation}\label{f:inv}
  \partial_\mathrm{inv}(f) = g_{LT}^{-1} f' \ .
\end{equation}
For any $a \in o_L$ we have
\begin{equation}\label{f:dlog}
  \log_{LT} ([a](Z)) = a \cdot \log_{LT} \qquad\text{and hence}\qquad ag_{LT}(Z) = g_{LT}([a](Z))\cdot [a]'(Z)
\end{equation}
(\cite[8.6 Lemma 2]{Lan}).

Let $\TLT$ be the Tate module of $LT$. Then $\TLT$ is
a free $o_L$-module of rank one, say with generator $\eta$, and the action of
$G_L := \Gal(\overline{L}/L)$ on $\TLT$ is given by a continuous character $\chi_{LT} :
 G_L \longrightarrow o_L^\times$. Let $\TLT'$ denote the Tate module of the $p$-divisible group Cartier dual to $LT$ with period $\Omega$ (depending on the choice of a generator of $\TLT'$), which again is a free $o_L$-module of rank one. The Galois action on $\TLT'\cong\TLT^*(1)$ is given by the continuous character $\tau := \chi_{cyc}\cdot\chi_{LT}^{-1}$, where
$\chi_{cyc}$ is the cyclotomic character.

For $n \geq 0$ we let $L_n/L$ denote the extension (in $\mathbb{C}_p$) generated by the $\pi_L^n$-torsion points of $LT$, and we put $L_\infty := \bigcup_n L_n$. The extension $L_\infty/L$ is Galois. We let $\Gamma_L := \Gal(L_\infty/L)$ and $H_L := \Gal(\overline{L}/L_\infty)$. The Lubin-Tate character $\chi_{LT}$ induces an isomorphism $\Gamma_L \xrightarrow{\cong} o_L^\times$.

Henceforth   we use the same notation as in \cite{SV15}.
In particular, the ring endomorphisms induced by sending $Z$ to $[\pi_L](Z)$ are called $\varphi_L$ where applicable; e.g.\  for the ring $\mathscr{A}_L$ defined to be the $\pi_L$-adic completion of $o_L[[Z]][Z^{-1}]$ or $\mathscr{B}_L := \mathscr{A}_L[\pi_L^{-1}]$ which denotes the field of fractions of $\mathscr{A}_L$.  Recall that we also have introduced the unique additive endomorphism $\psi_L$ of $\mathscr{B}_L$ (and then  $\mathscr{A}_L$) which satisfies
\begin{equation*}
  \varphi_L \circ \psi_L = \pi_L^{-1} \cdot trace_{\mathscr{B}_L/\varphi_L(\mathscr{B}_L)} \ .
\end{equation*}
Moreover,  projection formula
\begin{equation*}
  \psi_L(\varphi_L(f_1)f_2) = f_1 \psi_L(f_2) \qquad\text{for any $f_i \in \mathscr{B}_L$}
\end{equation*}
as well as the formula
\begin{equation*}
  \psi_L \circ \varphi_L = \frac{q}{\pi_L} \cdot \id \
\end{equation*}
hold. An  \'{e}tale $(\varphi_L,\Gamma_L)$-module $M$ comes with a Frobenius operator $\varphi_M$ and an induced operator  denoted by $\psi_M$.

Let $\widetilde{\mathbf{E}}^+ := \varprojlim o_{\mathbb{C}_p}/p o_{\mathbb{C}_p}$ with the transition maps being given by the Frobenius $\varphi(a) = a^p$. We may also identify $\widetilde{\mathbf{E}}^+$ with
$\varprojlim o_{\mathbb{C}_p}/\pi_L o_{\mathbb{C}_p}$ with
the transition maps being given by the $q$-Frobenius
$\varphi_q (a) = a^q$. Recall that $\widetilde{\mathbf{E}}^+$ is a complete valuation ring with residue field $\overline{\mathbb{F}_p}$ and its field of fractions $\widetilde{\mathbf{E}} = \varprojlim \mathbb{C}_p$ being algebraically closed of characteristic $p$. Let $\mathfrak{m}_{\widetilde{\mathbf{E}}}$ denote the maximal ideal in $\widetilde{\mathbf{E}}^+$.

The $q$-Frobenius $\varphi_q$ first extends by functoriality to the rings of the Witt vectors $W(\widetilde{\mathbf{E}}^+) \subseteq W(\widetilde{\mathbf{E}})$ and then $o_L$-linearly to $W(\widetilde{\mathbf{E}}^+)_L := W(\widetilde{\mathbf{E}}^+) \otimes_{o_{L_0}} o_L \subseteq W(\widetilde{\mathbf{E}})_L := W(\widetilde{\mathbf{E}}) \otimes_{o_{L_0}} o_L$, where $L_0$ is the maximal unramified subextension of $L$. The Galois group $G_L$ obviously acts on $\widetilde{\mathbf{E}}$ and $W(\widetilde{\mathbf{E}})_L$ by automorphisms commuting with $\varphi_q$.  This $G_L$-action is continuous for the weak topology on $W(\widetilde{\mathbf{E}})_L$ (cf.\ \cite[Lemma 1.5.3]{GAL}).

Sometimes we omit the index $q, L,$ or $M$ from the Frobenius operator, but we always write $\varphi_p$ when dealing with the $p$-Frobenius.
\newpage

\section{Wach-modules \`{a} la Kisin-Ren}

\subsection{Wach-modules}\label{sec:KR}

In this section we recall the theory of Wach-modules \`{a} la Kisin-Ren \cite{KR} (with the simplifying assumption that - in their notation - $K=L$, $m=1$ etc.).

By sending $Z$ to $\omega_{LT}\in W(\widetilde{\mathbf{E}}^+)_L$ (see directly after \cite[Lem.\ 4.1]{SV15}) we obtain an $G_L$-equivariant, Frobenius compatible  embedding of rings
\begin{equation*}
     o_L[[Z]] \longrightarrow W(\widetilde{\mathbf{E}}^+)_L \
\end{equation*}
the image of which we call $\mathbf{A}_L^+$, it is a subring of $\mathbf{A}_L$ (the image of $\mathscr{A}_L$ in $W(\widetilde{\mathbf{E}})_L$).   The latter ring is a complete discrete valuation ring with prime element $\pi_L$ and residue field the image $\mathbf{E}_L$ of $k_L ((Z)) \hookrightarrow \widetilde{\mathbf{E}}$ sending $Z$ to $\omega:=\omega_{LT} \mod \pi_L.$
We form the maximal integral unramified extension ($=$ strict Henselization) $\mathbf{A}_L^{nr}$  of $\mathbf{A}_L$  inside $W(\widetilde{\mathbf{E}})_L$. Its $p$-adic completion $\mathbf{A}$ still is contained in $W(\widetilde{\mathbf{E}})_L$. Note that $\mathbf{A}$ is a complete discrete valuation ring with prime element $\pi_L$ and residue field the separable algebraic closure $\mathbf{E}_L^{sep}$ of $\mathbf{E}_L$ in $\widetilde{\mathbf{E}}$.  By the functoriality properties of strict Henselizations the $q$-Frobenius $\varphi_q$ preserves $\mathbf{A}$. According to \cite[Lemma 1.4]{KR} the $G_L$-action on $W(\widetilde{\mathbf{E}})_L$ respects $\mathbf{A}$ and induces an isomorphism $H_L = \ker(\chi_{LT}) \xrightarrow{\cong} \Aut^{cont}(\mathbf{A}/\mathbf{A}_L)$. We set $\mathbf{A}^+:=\mathbf{A}\cap W(\widetilde{\mathbf{E}}^+)_L.$

Set $Q:= \frac{[\pi_L](\omega_{LT})}{\omega_{LT}}\in \mathbf{A}_L^+$, which satisfies per definitionem $\varphi_L(\omega_{LT})=Q\cdot \omega_{LT}$.

Following \cite{KR} we write $\mathcal{O} = \mathcal{O}_L(\mathbf{B})$ for the ring of rigid analytic functions on the open unit disk $\mathbf{B}$ over $L$, or equivalently the ring of power series in $Z$ over $L$ converging in $\mathbf{B}$. Via sending $\omega_{LT}$ to $Z$ we view $\mathbf{A}_L^+$ as a subring of $\mathcal{O}$.  We denote by $\mathrm{Mod}_{\mathcal{O}}^{\varphi_L,\Gamma_L,an}$ the category consisting of finitely generated free $ \mathcal{O}$-modules $\mathcal{M}$ together with the following data:
\begin{enumerate}
   \item an isomorphism $1 \otimes \varphi_{\mathcal{M}} : (\varphi_L^*\mathcal{M})[\frac{1}{Q}] \cong \mathcal{M}[\frac{1}{Q}].$ \footnote{By $\varphi_L^*\mathcal{M}$ we understand the module $\mathcal{O}\otimes_{\mathcal{O},\varphi_L} \mathcal{M}$.  }
   \item a semi-linear $\Gamma_L$-action on $\mathcal{M}$, commuting with $\varphi_{\mathcal{M}}$ and such that the induced action on $D(\mathcal{M}) := \mathcal{M}/\omega_{LT}\mathcal{M}$ is trivial.
\end{enumerate}
We note that, since $\mathcal{M}/\omega_{LT}\mathcal{M} = \mathcal{M}[\frac{1}{Q}] /\omega_{LT}\mathcal{M}[\frac{1}{Q}]$ the map $\varphi_{\mathcal{M}}$ induces an $L$-linear endomorphism of $D(\mathcal{M})$, which we denote by $\varphi_{D(\mathcal{M})}$. As a consequence of (1) it, in fact, is an automorphism.

The $\Gamma_L$-action on $\mathcal{M}$ is differentiable (\cite[Lemma 3.4.13]{BSX}) \footnote{Note that the statements in (loc.\ cit.) are all over the character variety; but by the introduction to \S 3.4 they are also valid over the open unit ball - with even easier proofs.}, and the corresponding derived action of $\mathrm{Lie}(\Gamma_L)$ is $L$-bilinear (\cite[Remark 3.4.15]{BSX})\footnote{In \cite{KR} being $L$-analytic is an extra condition in the definition of $\mathrm{Mod}_{\mathcal{O}}^{\varphi_L,\Gamma_L,an}$, which by this remark is automatically satisfied, whence the corresponding categories, with and without the superscript 'an' in \cite{KR} coincide!}.

Similarly, we denote by $\mathrm{Mod}_{\mathbf{A}_L^+}^{\varphi_L,\Gamma_L,an}$ the category consisting of finitely generated free $\mathbf{A}_L^+$-modules $N$ together with the following data:
\begin{enumerate}
 \item an isomorphism $1\otimes \varphi_N:(\varphi_L^*N)[\frac{1}{Q}]\cong N[\frac{1}{Q}]$. \footnote{By $\varphi_L^*N$ we understand the module $\mathbf{A}_L^+\otimes_{\mathbf{A}_L^+,\varphi_L} N$, and formally $\varphi_N$ is a map from $N$ to $N[\frac{1}{Q}]$.}
 \item a semi-linear $\Gamma_L$-action on $N$, commuting with $\varphi_N$ and such that the induced action on $N/\omega_{LT}N$ is trivial.
\end{enumerate}
The map $\varphi_N$ induces an $L$-linear automorphism of $D(N) := N[\frac{1}{p}]/\omega_{LT}N[\frac{1}{p}]$ denoted by $\varphi_{D(N)}$.

Obviously we have the base extension functor $\mathcal{O} \otimes_{\mathbf{A}_L^+} - : \mathrm{Mod}_{\mathbf{A}_L^+}^{\varphi_L,\Gamma_L,an} \longrightarrow \mathrm{Mod}_{\mathcal{O}}^{\varphi_L,\Gamma_L,an}$. It satisfies
\begin{equation}\label{f:1}
  D(\mathcal{O} \otimes_{\mathbf{A}_L^+} N) = D(N) \ .
\end{equation}
We write $\mathrm{Mod}_{\mathcal{O}}^{\varphi_L,\Gamma_L,0}$ for the full subcategory of $\mathrm{Mod}_{\mathcal{O}}^{\varphi_L,\Gamma_L, an}$ consisting of all $\mathcal{M}$ such that $\mathcal{R}\otimes_{\mathcal{O}}\mathcal{M}$ is pure of slope $0$. Here $\mathcal{R}$ denotes the Robba ring.

By $\mathrm{Mod}^{F,\varphi_q}_L$ we denote the category of finite dimensional $L$-vector spaces $D$ equipped with an $L$-linear automorphism $\varphi_q : D \xrightarrow{\cong} D$ and a decreasing, separated, and exhaustive filtration, indexed by $\mathbb{Z}$, by $L$-subspaces. In $\mathrm{Mod}^{F,\varphi_q}_L$ we have the full subcategory $\mathrm{Mod}^{F,\varphi_q,wa}_L$ of weakly admissible objects. For $D$ in $\mathrm{Mod}^{F,\varphi_q,wa}_L$ let $V_L(D) = \mathrm{Fil}^0(B_{cris,L}\otimes_L D)^{\varphi_q=1}$ where, as usual, $B_{cris,L} := B_{cris} \otimes_{L_0} L$. In order to formulate the crystalline comparison theorem in this context we also consider the category $\mathrm{Mod}^{F,\varphi_q}_{L_0\otimes_{\mathbb{Q}_p}L}$ of finitely generated free $L_0\otimes_{\mathbb{Q}_p}L$-modules $\mathfrak{D}$ equipped with a $(\varphi_p \otimes \id)$-linear automorphism $\varphi_q : {\mathfrak{D}}\xrightarrow{\cong} {\mathfrak{D}}$ and a decreasing, separated, and exhaustive filtration on $\mathfrak{D}_L := \mathfrak{D}\otimes_{L_0} L$, indexed by $\mathbb{Z}$, by $L\otimes_{\mathbb{Q}_p}L$-submodules. For $\mathfrak{D}$ in $\mathrm{Mod}^{F,\varphi_q}_{L_0\otimes_{\mathbb{Q}_p}L}$ we define, as usual,
\begin{equation*}
  V(\mathfrak{D}) := (B_{cris} \otimes_{L_0} \mathfrak{D})^{\varphi=1} \cap \mathrm{Fil}^0(B_{dR} \otimes_L \mathfrak{D}_L) \ .
\end{equation*}

Let $\Rep_{o_L,f}(G_L)$ denote the category of finitely generated free $o_L$-modules equipped with a continuous linear $G_L$-action and $\Rep_{o_L,f}^{cris,an}(G_L)$ the full subcategory of those $T$ which are free over $o_L$ and such that the representation $V := L \otimes_{o_L} T$ is crystalline and {\em analytic}, i.e., satisfying that, if $D_{dR}(T) := (T \otimes_{\mathbb{Z}_p} B_{dR})^{G_L}$, the filtration on $D_{dR}(T)_\mathfrak{m}$ is trivial for each maximal ideal $\mathfrak{m}$ of $L \otimes_{\mathbb{Q}_p} L$ which does not correspond to the identity $\id : L \to L$. Correspondingly we let $\Rep_L^{cris}(G_L)$, resp.\ $\Rep_L^{cris,an}(G_L)$, denote the category of continuous $G_L$-representations in finite dimensional $L$-vector spaces which are crystalline, resp.\ crystalline and analytic. The base extension functor $L \otimes_{o_L} -$ induces an equivalence of categories $\Rep_{o_L,f}^{cris,an}(G_L) \otimes_{\mathbb{Z}_p} \mathbb{Q}_p \xrightarrow{\simeq} \Rep_L^{cris,an}(G_L)$. Here applying $\otimes_{\mathbb{Z}_p} \mathbb{Q}_p$ to a $\mathbb{Z}_p$-linear category means applying this functor to the $\mathrm{Hom}$-modules. For $V$ in $\Rep_L^{cris,an}(G_L)$ we set $D_{cris,L}(V) := (B_{cris,L} \otimes_{L} V)^{G_L} = (B_{cris} \otimes_{L_0} V)^{G_L}$ and $D_{cris}(V) := (B_{cris} \otimes_{\mathbb{Q}_p} V)^{G_L}$. The usual crystalline comparison theorem says that $D_{cris}$ and $V$ are equivalences of categories between $\Rep_L^{cris}(G_L)$ and the subcategory of weakly admissible objects in $\mathrm{Mod}^{F,\varphi_q}_{L_0\otimes_{\mathbb{Q}_p}L}$.

\begin{lemma}\label{filteredcategories}
(\cite[Lemma 5.3]{ST4} and subsequent discussion, or \cite[Cor.\ 3.3.1]{KR}) There is a fully faithful $\otimes$-functor
\begin{align*}
 ^\sim: \mathrm{Mod}^{F,\varphi_q}_L & \longrightarrow \mathrm{Mod}^{F,\varphi}_{L_0 \otimes_{\mathbb{Q}_p} L}  \\
                             D & \longmapsto  \tilde D := L_0 \otimes_{\mathbb{Q}_p} D \ ,
\end{align*}
whose essential image consists of all {\em analytic} objects, i.e., those for which the filtration on the non-identity components is trivial. A quasi-inverse functor from the essential image is given by sending $\mathfrak{D}$ to the base extension $L \otimes_{L_0 \otimes_{\mathbb{Q}_p}L} \mathfrak{D}$ for the multiplication map $L_0 \otimes_{\mathbb{Q}_p} L \rightarrow L$.
\end{lemma}

Lemma \ref{filteredcategories} implies  that
\begin{equation}\label{f:Dcris}
 {D_{cris,L}(V)}^\sim  \cong D_{cris}(V)  \qquad\text{for any $V$ in $\Rep_L^{cris,an}(G_L)$}.
\end{equation}

We denote by $\mathfrak{M}^{et}(\mathbf{A}_L)$ the category of \'{e}tale $(\varphi_q,\Gamma_L)$-modules over $\mathbf{A}_L$  (cf.\ \cite[Def.\ 3.7]{SV15})  and by $\mathfrak{M}^{et}_{f}(\mathbf{A}_L)$  the full subcategory consisting of those objects, which are finitely generated free as $\mathbf{A}_L$-module. For $M$ in $\mathfrak{M}^{et}_{f}(\mathbf{A}_L)$, resp.\ for $T$ in $\Rep_{o_L,f}(G_L)$, we put $V(M):=(\mathbf{A} \otimes_{\mathbf{A}_L} M)^{\varphi_q \otimes \varphi_M =1}$, resp.\  $D_{LT}(T) := (\mathbf{A} \otimes_{o_L} T)^{\ker(\chi_{LT})}$.

Having defined all of the relevant categories (and most of the functors) we now contemplate the following diagram of functors:
\begin{equation*}
  \xymatrix{
     & & & \mathfrak{M}^{et}_{f}(\mathbf{A}_L) \ar@<-1.5ex>[d]_{V}   \\
    \mathrm{Mod}_{\mathbf{A}_L^+}^{\varphi_L,\Gamma_L,an} \ar@/_13ex/[dd]_{\mathcal{O} \otimes_{\mathbf{A}_L^+} -} \ar[d] \ar[rr]^{\simeq} \ar[rrru]^{\mathbf{A}_L \otimes_{\mathbf{A}_L^+} -} &  & \Rep_{o_L,f}^{cris,an}(G_L) \ar[dd]^{L \otimes_{o_L} -} \ar[r]^{\subseteq} & \Rep_{o_L,f}(G_L) \ar@<-1ex>[u]_{D_{LT}}^{\simeq} \\
    \mathrm{Mod}_{\mathbf{A}_L^+}^{\varphi_L,\Gamma_L,an} \otimes_{\mathbb{Z}_p} \mathbb{Q}_p \ar[d]_{\simeq}  &  & &  \\
    \mathrm{Mod}_{\mathcal{O}}^{\varphi_L,\Gamma_L,0} \ar[d]_{\subseteq} \ar@<1ex>[r]^{D}_{\simeq} & \mathrm{Mod}^{F,\varphi_q,wa}_L \ar[d]_{\subseteq} \ar@<1ex>[l]^{\mathcal{M}} \ar@<1ex>[r]^-{V_L}_-{\simeq} & \Rep_L^{cris,an}(G_L) \ar@<1ex>[l]^-{D_{cris,L}} & \\
    \mathrm{Mod}_{\mathcal{O}}^{\varphi_L,\Gamma_L,an} \ar@<1ex>[r]^{D}_{\simeq} & \mathrm{Mod}^{F,\varphi_q}_L \ar@<1ex>[l]^{\mathcal{M}} &  &  }
\end{equation*}
The arrows without decoration are the obvious natural ones. The following pairs of functors are quasi-inverse $\otimes$-equivalences of $\otimes$-categories:
\begin{itemize}
  \item[--] $(D_{LT},V)$ by \cite[Thm.\ 1.6]{KR};
  \item[--] $(D_{cris,L},V_L)$ by the crystalline comparison theorem (\cite[Rem.\ 3.6.7]{F1}) and Lemma \ref{filteredcategories};
  \item[--] $(D,\mathcal{M})$ by \cite[Prop.\ 2.2.6]{KR} (or \cite[Thm.\ 3.4.16]{BSX}) and \cite[Cor.\ 2.4.4]{KR}, to which we also refer for the definition of the functor $\mathcal{M}$.
\end{itemize}
In particular, all functors in the above diagram are $\otimes$-functors.
The second arrow in the left column, resp.\ the left arrow in the upper horizontal row, is an equivalence of categories by \cite[Cor.\ 2.4.2]{KR}, resp.\ by \cite[Cor.\ 3.3.8]{KR}. The lower square and the upper triangle are commutative for trivial reasons.

We list a few additional properties of these functors.

\begin{remark}\phantomsection\label{properties}
\begin{itemize}
  \item[i.] For any $M$ in $\mathfrak{M}^{et}_{f}(\mathbf{A}_L)$ the inclusion $V(M) \subseteq \mathbf{A}\otimes_{\mathbf{A}_L }M$ extends to an isomorphism
\begin{equation}
\label{f:compiso} \mathbf{A}\otimes_{o_L}V(M) \xrightarrow{\cong} \mathbf{A}\otimes_{\mathbf{A}_L }M \ ,
\end{equation}
      which is compatible with the $\varphi_q$- and $\Gamma_L$-actions on both sides.
  \item[ii.] The functors $D_{LT}$, $V$, and $V(\mathbf{A}_L \otimes_{\mathbf{A}_L^+} -)$ respect exact sequences (of abelian groups).
  \item[iii.] (\cite[Prop.\ 3.4.14]{BSX}) For any $\mathcal{M}$ in $\mathrm{Mod}_{\mathcal{O}}^{\varphi_L,\Gamma_L,an}$ the projection map $\mathcal{M}[\frac{\omega_{LT}}{t_{LT}}] \longrightarrow D(\mathcal{M})$ restricts to an isomorphism $\mathcal{M}[\frac{\omega_{LT}}{t_{LT}}]^{\Gamma_L} \xrightarrow{\;\cong\;} D(\mathcal{M})$ such that the diagram
\begin{equation*}
  \xymatrix{
    \mathcal{M}[\frac{\omega_{LT}}{t_{LT}}]^{\Gamma_L}  \ar[d]_{\varphi_{\mathcal{M}}} \ar[r]^-{\cong} & D(\mathcal{M}) \ar[d]^{\varphi_{D(\mathcal{M})}} \\
    \mathcal{M}[\frac{\omega_{LT}}{t_{LT}}]^{\Gamma_L}  \ar[r]^-{\cong} & D(\mathcal{M})   }
\end{equation*}
      is commutative; moreover, $\mathcal{M}[\frac{\omega_{LT}}{t_{LT}}] = \mathcal{O}[\frac{\omega_{LT}}{t_{LT}}] \otimes_L \mathcal{M}[\frac{\omega_{LT}}{t_{LT}}]^{\Gamma_L} \cong \mathcal{O}[\frac{\omega_{LT}}{t_{LT}}] \otimes_L D(\mathcal{M})$.
\end{itemize}
\end{remark}

Now we recall that $A_{cris}$ is the $p$-adic completion of a divided power envelope of $W(\widetilde{\mathbf{E}}^+)$ and let $A_{cris,L}:=A_{cris}\otimes_{L_0} L$. The inclusion $W(\widetilde{\mathbf{E}}^+)\subseteq A_{cris}$ induces an embedding
$\mathbf{A}^+_L \subseteq W(\widetilde{\mathbf{E}}^+)_L\subseteq A_{cris,L}.$

 We observe that $t_{LT}=\log_{LT}(\omega_{LT})$ belongs to $B_{cris,L}^\times$. Indeed, by
 \cite[\S III.2]{CoIw} we know that $\varphi_p(B_{max})\subseteq B_{cris}\subseteq  B_{max},$ whence we obtain \[\varphi_q(B_{max}\otimes_{L_0}L)\subseteq B_{cris,L} \subseteq  B_{max}\otimes_{L_0}L,\] where the definition of  $B_{max}$ can be found in (loc.\ cit.).  By
\cite[Prop.\ 9.10, Lem.\ 9.17,\S 9.7]{Co4} $t_{LT}$ and $\omega_{LT}$ are invertible in $ B_{max,L}\subseteq B_{max}\otimes_{L_0}L$ (This reference assumes that the power series $[\pi_L](Z)$ is a polynomial. But, by some additional convergence considerations, the results can be seen to hold in general (cf.\ \cite[\S2.1]{GAL} for more details)). Hence, by the above inclusions and using that $\varphi_q(t_{LT})=\pi_L t_{LT},$ we see that  $t_{LT}$ is a unit $B_{cris,L}.$ In particular, we have an inclusion $A_{cris,L}[\tfrac{1}{\pi_L},\tfrac{1}{t_{LT}}]\subseteq B_{cris,L} $. Moreover,  since $\varphi_q(\omega_{LT})=Q \omega_{LT}$  is invertible in $\varphi_q(B_{max}\otimes_{L_0}L)$, the elements  $\omega_{LT}$ and $Q$ are  units in $B_{cris,L}$ as well. In particular, we have an inclusion
\begin{eqnarray}\label{f:ABcris}
  \mathbf{A}^+[\tfrac{1}{\omega_{LT}}] &\subseteq B_{cris,L}.
\end{eqnarray}

  Next we shall recall in Lemma \ref{lem:OBcris} below that the above inclusion $\mathbf{A}^+_L \subseteq A_{cris,L}$ extends to a (continuous) ring homomorphism
\begin{align}\label{f:OBcris}
 \mathcal{O} &\rightarrow A_{cris,L}[\tfrac{1}{\pi_L}] \subseteq B_{cris,L}.
\end{align}

For $ \alpha\in \widetilde{\mathbf{E}}^+\cong  \projlim_{n} o_{\mathbb{C}_p}$ we denote by $\alpha^{(0)}$ as usual its zero-component.

\begin{lemma}
 The following diagram of $o_{L_0}$-modules is commutative
 \begin{equation}
 \label{f:diagtheta}
 \xymatrix{
   0 \ar[r]^{ } & J \ar[d]_{ } \ar[r]^{ } & W(\widetilde{\mathbf{E}}^+) \ar[d]_{u} \ar[r]^{\Theta} & o_{\mathbb{C}_p} \ar@{=}[d]_{} \ar[r]^{ } & 0  \\
  0 \ar[r]^{} & \ker(\Theta_L) \ar[r]^{} & W(\widetilde{\mathbf{E}}^+)_L \ar[r]^{\Theta_L} & o_{\mathbb{C}_p}\ar[r]^{} & 0 , }
 \end{equation}
  where $J:=\ker(\Theta)$,  $\Theta(\sum_{n\geq 0} [\alpha_n]p^n))=\sum_{n\geq 0} \alpha_n^{(0)}p^n$ and similarly  $\Theta_L(\sum_{n\geq 0} [\alpha_n]\pi_L^n))=\sum_{n\geq 0} \alpha_n^{(0)}\pi_L^n$, while $u$ denotes the canonical map as defined in \cite[Lem.\ 1.2.3]{FF}, it sends Teichm\"{u}ller lifts $[\alpha]$ with respect to $W(\widetilde{\mathbf{E}}^+)  $ to the Teichm\"{u}ller lift $[\alpha]$  with respect to $W(\widetilde{\mathbf{E}}^+)_L.  $
\end{lemma}

\begin{proof}
First of all we recall from \cite[Lem.\ 1.6.1]{GAL} that $\Theta$ and $\Theta_L$ are continuous and show that also $u$ is continuous, each time with respect to the weak topology, of which a fundamental system of open neighbourhoods consists of \[U_{\mathfrak{a},m}:=\{(b_0,b_1,\ldots)\in W(\widetilde{\mathbf{E}}^+)| b_0,\dots, b_{m-1}\in \mathfrak{a}\}=\sum_{i=0}^{m-1} V_p^i([\mathfrak{a}])+p^m W(\widetilde{\mathbf{E}}^+) \] and similarly $U_{\mathfrak{a},m}^L:=\{(b_0,b_1,\ldots)\in W(\widetilde{\mathbf{E}}^+)_L| b_0,\dots, b_{m-1}\in \mathfrak{a}\}$ for open ideals $\mathfrak{a}$ of $\widetilde{\mathbf{E}}^+$ and  $m\geq 0$; see \S 1.5 in (loc.\ cit.). By $o_{L_0}$-linearity, we see that $u(p^m W(\widetilde{\mathbf{E}}^+))\subseteq p^mW(\widetilde{\mathbf{E}}^+)_L$. Using the relation
\[u(V_px)=\frac{p}{\pi_L} V_{\pi_L}( u(F^{f-1}x)  )\] from \cite[Lem.\ 1.2.3]{FF}\footnote{Note the typos in (loc.\ cit.) where $u(V_px)=\frac{p}{\pi_L} V_{\pi_L}( F^{f-1}u(x)  ) $ is stated. Moreover, one has the relation $u(F^fx)=Fu(x).$}, where $V_?$ denotes the Verschiebung, one easily concludes that \[u(V_p^i([b]))=(\frac{p}{\pi_L})^iV_{\pi_L}^i([b^{p^{i(f-1)}}]),\] whence $u(U_{\mathfrak{a},m})\subseteq U_{\mathfrak{a},m}^L$ and continuity of $u$ follows.

Since the commutativity is clear on Teichm\"{u}ller lifts and on $p$ by $o_{L_0}$-linearity, which generate a dense ideal, the result follows by continuity.
\end{proof}

The following lemma generalizes parts from \cite[Prop.\ 1.5.2.]{PR}.

\begin{lemma}\label{lem:OBcris}
Sending $f=\sum_{n\geq 0} a_nZ^n$ to $f(\omega_{LT})$ induces a continuous map
\[\mathcal{O}\to A_{cris,L}[\frac{1}{\pi_L}],\] where the source carries the Fr\'{e}chet-topology while the target is a topological $o_{L_0}$-module, of which the topology is uniquely determined by requiring that $A_{cris,L}$ is open, i.e., the system $p^mA_{cris,L}$ with $m\geq 0$ forms a basis of open neighbourhoods of $0.$\Footnote{Apparently this coincides with the inductive limit topology (in the category of topological spaces or groups or $o_{L_0}$-modules).  }
\end{lemma}

\begin{proof}
First of all, the relation $J^p\subseteq p A_{cris}$ from \cite[\S 1.4.1, bottom of p.\ 96]{PR} (note that $J^p\subseteq W_p(R)$ regarding the notation in (loc.\ cit.) for the last object)  implies easily by flat base change
\begin{equation}
\label{f:PR}
J_L^p\subseteq pA_{cris,L}
\end{equation}
with $J_L:=J\otimes_{o_{L_0}}o_L.$ By \cite[Lem.\ 2.1.12]{GAL} we know that $\omega_{LT}$ belongs to $\ker(\Theta_L).$ Now we claim that there exists a natural number $r'$ such that $\omega_{LT}^{r'}$ lies in $W_1:=J_L+pW(\widetilde{\mathbf{E}}^+ )_L$, whence for $r:=pr'$ we have $\omega_{LT}^r\in W_p$ with $W_m:=W_1^m$ for all $m\geq 0.$ To this aim note that diagram \eqref{f:diagtheta} induces the following commutative diagram with exact lines
\begin{equation}\scriptsize
 \label{f:diagtheta+p}
 \xymatrix{
   0 \ar[r]^{ } & W_1 \ar[d]_{ } \ar[r]^{ } & W(\widetilde{\mathbf{E}}^+) \otimes_{o_{L_0}}o_L \ar[d]_{\cong} \ar[r]^(0.35){\Theta} & (o_{\mathbb{C}_p} \otimes_{o_{L_0}}o_L)/p(o_{\mathbb{C}_p} \otimes_{o_{L_0}}o_L) \ar[d]_{\mu} \ar[r]^{ } & 0  \\
  0 \ar[r]^{} & \ker(\Theta_L)+pW(\widetilde{\mathbf{E}}^+)_L \ar[r]^{} & W(\widetilde{\mathbf{E}}^+)_L \ar[r]^{\Theta_L} & o_{\mathbb{C}_p}/po_{\mathbb{C}_p}\ar[r]^{} & 0 , }
 \end{equation}
 where the map $\mu$ is induced by sending $a\otimes b$ to $ab$ and a reference for the middle vertical isomorphism is \cite[Prop.\ 1.1.26]{GAL}.
By the snake lemma the cokernel of the left vertical map is isomorphic to
\begin{align*}
  \ker(\mu) & \subseteq \ker\left( (o_{\mathbb{C}_p} \otimes_{o_{L_0}}o_L)/p(o_{\mathbb{C}_p} \otimes_{o_{L_0}}o_L)   \to \bar{k}\right) \\
   & = \ker\left( o_{\mathbb{C}_p}/p o_{\mathbb{C}_p}\otimes_{k}o_L/po_L   \to  o_{\mathbb{C}_p}/\mathfrak{m}_{\mathbb{C}_p}\otimes_{k}o_L/\pi_Lo_L \right)  \\
   & = \mathfrak{m}_{\mathbb{C}_p}\otimes_k o_L/po_L + o_{\mathbb{C}_p}/p o_{\mathbb{C}_p}\otimes_k \pi_L o_L/po_L
\end{align*}
and thus consists of nilpotent elements whence the claim follows. Here $\mathfrak{m}_{\mathbb{C}_p}$ denotes the maximal ideal of $o_{\mathbb{C}_p}.$

Now let $f=\sum_{n\geq 0} a_nZ^n$ satisfy that $|a_n|\rho^n$ tends to zero for all $\rho<1.$ Writing $n=q_nr+r_n$ with $0\leq r_n<r,$ we have
\[a_n\omega_{LT}^n=a_n \omega_{LT}^{r_n}(\omega_{LT}^r)^{q_n}\in a_n W_{pq_n}\subseteq a_n p^{q_n}A_{cris,L},\]
where the last inclusion follows from \eqref{f:PR}.
But $|a_np^{q_n}|\leq |a_n|_pp^{1-\frac{n}{r}}$ tends to $0$ for $n\to\infty.$ Thus the series $\sum_{n\geq 0} a_n \omega_{LT}^n$ converges in $A_{cris,L}[\frac{1}{\pi_L}].$

Moreover, since one has $sup |a_n p^{-1+\frac{n}{r}}|\leq p\|f\|_\rho$ for the usual norms $\|\cdot \|_\rho$ if $1>\rho>p^{-\frac{1}{r}}$, we obtain for any $m$ that
\[\{f|\; \|f\|_\rho<p^{-m-1}\}\subseteq \{f\in \mathcal{O}| f(\omega_{LT})\in p^mA_{cris,L}\},\] whence the latter set, which is the preimage of $p^mA_{cris,L},$ is open. This implies continuity.
\end{proof}

\begin{lemma}\label{square}
The big square in the middle is a commutative square of $\otimes$-functors (up to a natural isomorphism of $\otimes$-functors).
\end{lemma}
\begin{proof}
We have to establish a natural isomorphism
\begin{equation}\label{f:2}
  L \otimes_{o_L} V(\mathbf{A}_L \otimes_{\mathbf{A}_L^+} N) \cong V_L(D(\mathcal{O} \otimes_{\mathbf{A}_L^+} N))  \qquad\text{for any $N$ in $\mathrm{Mod}_{\mathbf{A}_L^+}^{\varphi_L,\Gamma_L,an}$}.
\end{equation}
In fact, we shall prove the dual statement, i.e., using  \eqref{f:1},  that
\begin{equation}
\label{f:3} (L \otimes_{o_L} V(\mathbf{A}_L \otimes_{\mathbf{A}_L^+} N))^* \cong V_L(D( N))^*,
\end{equation}
where $^*$ indicates the $L$-dual. From the canonical isomorphisms
\begin{align*}
 \mathrm{Hom}_{\mathbf{A}_L ,\varphi_q}(M,\mathbf{A} )&\cong\mathrm{Hom}_{\mathbf{A} ,\varphi_q}(\mathbf{A}\otimes_{\mathbf{A}_L}M,\mathbf{A} )\\
 &\cong \mathrm{Hom}_{\mathbf{A} ,\varphi_q}(\mathbf{A}\otimes_{o_L}V(M),\mathbf{A} )\\
 &\cong \mathrm{Hom}_{o_L}( V(M),\mathbf{A}^{\varphi_q=1} )\\
 &\cong \mathrm{Hom}_{o_L}( V(M),o_L ),
\end{align*}
where we used \eqref{f:compiso} for the second isomorphism and write $M$ for $\mathbf{A}_L\otimes_{\mathbf{A}_L^+}N$, we conclude that the left hand side of \eqref{f:3} is canonically isomorphic to $\mathrm{Hom}_{\mathbf{A}_L ,\varphi_q}(\mathbf{A}_L\otimes_{\mathbf{A}_L^+} N,\mathbf{A} )\otimes_{o_L}L. $
Let $\mathbf{A}^+ := \mathbf{A} \cap W(\widetilde{\mathbf{E}}^+)_L$. 
On the one hand, by \cite[Lem.\ 3.2.1]{KR}, base extension induces an isomorphism
\begin{equation*}
     \mathrm{Hom}_{\mathbf{A}_L^+ ,\varphi_q}(N,\mathbf{A}^+[\tfrac{1}{\omega_{LT}}]) \xrightarrow{\cong} \mathrm{Hom}_{\mathbf{A}_L ,\varphi_q}(\mathbf{A}_L\otimes_{\mathbf{A}_L^+} N,\mathbf{A} ) \ .
\end{equation*}
On the other hand, in \cite[Prop.\ 3.2.3]{KR} they construct a natural isomorphism
\begin{equation}\label{f:V2}
  \mathrm{Hom}_{\mathbf{A}_L^+ ,\varphi_q}(N,\mathbf{A}^+[\tfrac{1}{\omega_{LT}}]) \otimes_{o_L} L  \xrightarrow{\cong} \mathrm{Hom}_{L,\varphi_q,\mathrm{Fil}}((N/\omega_{LT}N)[\tfrac{1}{p}],B_{cris,L}) \ .
\end{equation}
\footnote{$\otimes_{o_L}L $ is missing in the reference!} Therefore, the left hand side of \eqref{f:3} becomes naturally isomorphic to
\begin{equation}\label{f:V3}
 \Hom_{L,\varphi_q,\mathrm{Fil}}(D(N),B_{cris,L})\cong V_L(D(N)^*),
\end{equation}
where the last isomorphism is straightforward. Thus the proof of \eqref{f:3} is reduced to the canonical identity
\begin{equation}\label{f:V4}
 V_L(D(N)^*)\cong V_L(D(N))^*.
\end{equation}
This can be proved in the same way as in \cite[Rem.\ 3.4.5 (iii), Rem.\ 3.6.7]{F1}: Since $V_L$ is a rigid $\otimes$-functor, it preserves inner Hom-objects, in particular duals.

In order to see that \eqref{f:2} is compatible with tensor products note that base change, taking $L$-duals or applying comparison isomorphisms are $\otimes$-compatible. Thus the claim is reduced to the tensor compatiblity of the isomorphism \eqref{f:V2} the construction of  which we therefore recall from \cite{KR}. It  is induced by a natural map
\begin{equation*}
  \mathrm{Hom}_{\mathbf{A}_L^+}(N,\mathbf{A}^+[\tfrac{1}{\omega_{LT}}]) \otimes_{o_L} L \longrightarrow \mathrm{Hom}_{L}((N/\omega_{LT}N)[\tfrac{1}{p}],B_{cris,L})
\end{equation*}
which comes about as follows. Let $f \in \mathrm{Hom}_{\mathbf{A}_L^+}(N,\mathbf{A}^+[\tfrac{1}{\omega_{LT}}])$.
%
%
%
%
 By composing $f$ with the inclusion \eqref{f:ABcris} we obtain $f_1 : N \rightarrow B_{cris,L}$.
By base extension to $\mathcal{O}$ via \eqref{f:OBcris} and then localization in $Q$ the map $f_1$ gives rise to a map $f_2 : (\mathcal{O} \otimes_{\mathbf{A}_L^+} N)[\tfrac{1}{Q}] \rightarrow B_{cris,L}$. This one we precompose with the isomorphism $1 \otimes \varphi_N$ to obtain
\begin{equation*}
  f_3 : (\mathcal{O} \otimes_{\mathbf{A}_L^+,\varphi_L} N)[\tfrac{1}{Q}] \xrightarrow{\cong} (\mathcal{O} \otimes_{\mathbf{A}_L^+} N)[\tfrac{1}{Q}] \rightarrow B_{cris,L} \ .
\end{equation*}
Now we observe the inclusions
\begin{multline*}
  (\mathcal{O} \otimes_{\mathbf{A}_L^+,\varphi_L} N)[\tfrac{1}{Q}] \subseteq (\mathcal{O} \otimes_{\mathbf{A}_L^+,\varphi_L} N)[\tfrac{\omega_{LT}}{t_{LT}}] \supseteq (\mathcal{O} \otimes_{\mathbf{A}_L^+,\varphi_L} N)[\varphi_L(\tfrac{\omega_{LT}}{t_{LT}})] \\
  = \mathcal{O} \otimes_{\mathcal{O},\varphi_L} \big( (\mathcal{O} \otimes_{\mathbf{A}_L^+} N)[\tfrac{\omega_{LT}}{t_{LT}}] \big) \ .
\end{multline*}
They only differ by elements which are invertible in $B_{cris,L}$. Therefore giving the map $f_3$ is equivalent to giving a map $f_4 : \mathcal{O} \otimes_{\mathcal{O},\varphi_L} \big( (\mathcal{O} \otimes_{\mathbf{A}_L^+} N)[\tfrac{\omega_{LT}}{t_{LT}}] \big) \rightarrow B_{cris,L}$. Finally we use Remark \ref{properties}.iii which gives the map
\begin{equation*}
  \xi : (N/\omega_{LT}N)[\tfrac{1}{p}] \xleftarrow[\pr]{\cong} \big( (\mathcal{O} \otimes_{\mathbf{A}_L^+} N)[\tfrac{\omega_{LT}}{t_{LT}}] \big)^{\Gamma_L} \xrightarrow{\subseteq} (\mathcal{O} \otimes_{\mathbf{A}_L^+} N)[\tfrac{\omega_{LT}}{t_{LT}}] \ .
\end{equation*}
By precomposing $f_4$ with $1 \otimes \xi$ we at last arrive at a map $f_5 : (N/\omega_{LT}N)[\tfrac{1}{p}] \rightarrow B_{cris,L}$.
From this description the compatibility with tensor products is easily checked.
\end{proof}

Suppose that $N$ is in $\mathrm{Mod}_{\mathbf{A}_L^+}^{\varphi_L,\Gamma_L,an}$ and put $T := V(\mathbf{A}_L \otimes_{\mathbf{A}_L^+} N)$ in $\Rep_{o_L,f}^{cris,an}(G_L)$. Then, by Remark \ref{properties}.iii and Lemma \ref{square}, we have a natural isomorphism of  $\otimes$-functors
\begin{equation}\label{f:compB}
  \mathrm{comp} : \mathcal{O}[\tfrac{\omega_{LT}}{t_{LT}}] \otimes_{\mathbf{A}_L^+} N \xrightarrow{\;\cong\;}  \mathcal{O}[\tfrac{\omega_{LT}}{t_{LT}}] \otimes_L D_{cris,L}(L \otimes_{o_L} T)
\end{equation}
which is compatible with the diagonal $\varphi$'s on both sides.

In the proof of \cite[Cor.\ 3.3.8]{KR} it is shown that, for any $T$ in $\Rep_{o_L,f}^{cris,an}(G_L)$, there exists an $\mathbf{A}_L^+$-submodule  $\mathfrak{M} \subseteq D_{LT}(T)$ which
\begin{enumerate}
 \item[(N1)] lies in $\mathrm{Mod}_{\mathbf{A}_L^+}^{\varphi_L,\Gamma_L,an}$ with $\varphi_{\mathfrak{M}}$ and the $\Gamma_L$-action on $\mathfrak{M}$ being induced by the $(\varphi_q,\Gamma_L)$-structure of $D_{LT}(T),$ and
 \item[(N2)] satisfies  $\mathbf{A}_L\otimes_{\mathbf{A}_L^+} \mathfrak{M} = D_{LT}(T)$.\Footnote{ Let $M$ be a (free) $\mathbf{A}_L$-module together with a (free) $\mathbf{A}_L^+$-submodule $\mathfrak{M}$    such that $\mathbf{A}_L\otimes_{\mathbf{A}_L^+} \mathfrak{M}\subseteq M$. Then, if $\mathfrak{M}$ is $p$-saturated in $M$, then $\mathbf{A}_L\otimes_{\mathbf{A}_L^+} \mathfrak{M}$ is so, too. Let $m\in M$ and  assume that $\pi_Lm = \sum_i f_i\otimes \mathfrak{m}_i$ belongs to $\mathbf{A}_L\otimes_{\mathbf{A}_L^+}\mathfrak{M} $. Choose $j\geq 0$ such that $X^jf_i=\pi_Lg_i +f_i'$ for some  $g_i\in \mathbf{A}_L$ and     $f_i'\in\mathbf{A}_L^+$ . Then $\pi_L X^jm=\pi_L m_0 +1\otimes \mathfrak{m}$ for some $m_0\in \mathbf{A}_L\otimes_{\mathbf{A}_L^+}\mathfrak{M}$ and $\mathfrak{m}\in \mathfrak{M}.$ It follows that $X^jm$ belongs to $\mathbf{A}_L\otimes_{\mathbf{A}_L^+}\mathfrak{M}  $ whence also $m$ as $X$ is invertible in $\mathbf{A}_L.$        }
\end{enumerate}
Note that property (N2) implies that $\mathfrak{M}  $ is \textit{$p$-saturated in $ D_{LT}(T)$}, i.e., $\mathfrak{M} [\frac{1}{p}] \cap D_{LT}(T) = \mathfrak{M}$, since $\mathbf{A}_L^+$ is obviously $p$-saturated in $\mathbf{A}_L^.$

We once and for all pick such an $N(T) := \mathfrak{M}$. This defines a functor
\begin{equation*}
  N : \Rep_{o_L,f}^{cris,an}(G_L) \longrightarrow \mathrm{Mod}_{\mathbf{A}_L^+}^{\varphi_L,\Gamma_L,an}
\end{equation*}
which is quasi-inverse to the upper left horizontal arrow in the above big diagram.  Note that $N$ is in a unique way a $\otimes$-functor by \cite[I.4.4.2.1]{Sa}.

\begin{remark}\label{positive}
For $T$ in $\Rep_{o_L,f}^{cris,an}(G_L)$ and $N := N(T)$ in $\mathrm{Mod}_{\mathbf{A}_L^+}^{\varphi_L,\Gamma_L,an}$ we have:
\begin{enumerate}
\item If $L \otimes_{o_L} T$ is a positive\footnote{i.e., the Hodge-Tate weights are non-positive, i.e., $gr^j D_{dR}(T)\neq 0$ implies that $j\geq 0.$ } analytic crystalline representation, then $N$ is stable under $\varphi_N$;
\item If the Hodge-Tate weights of $L \otimes_{o_L} T$ are all $\geq 0$, then we have $N \subseteq \mathbf{A}_L^+ \cdot \varphi_N(N)$, where the latter means the $\mathbf{A}_L^+ $-span generated by $ \varphi_N(N)$.
\end{enumerate}
\end{remark}
\begin{proof}  The corresponding assertions for $\mathcal{M} := \mathcal{O} \otimes_{\mathbf{A}_L^+}N$ are contained in \cite[Cor.\ 3.4.9]{BSX}. Let $n_1, \ldots, n_d$ be an $\mathbf{A}_L^+$-basis of $N$.

For (i) we have to show that $\varphi_N(n_j) \in N$ for any $1 \leq j \leq d$. Writing $\varphi_N(n_j) = \sum_{i=1}^d f_{ij} n_i$ we know from the definition of the category $\mathrm{Mod}_{\mathbf{A}_L^+}^{\varphi_L,\Gamma_L,an}$ that $f_{ij} \in \mathbf{A}_L^+[\frac{1}{Q}]$ and from the above observation that $f_{ij} \in \mathcal{O}$. This reduces us to showing that $o_L[[Z]][\frac{1}{Q}] \cap \mathcal{O} \subseteq o_L[[Z]]$. Suppose therefore that $Q^r h = f$ for some $r \geq 1$, $h \in \mathcal{O}$, and $f \in o_L[[Z]]$. The finitely many zeros of $Q \in o_L[[Z]]$, which are the nonzero $\pi_L$-torsion points of the Lubin-Tate formal group, all lie in the open unit disk. By Weierstrass preparation it follows that $Q$ must divide $f$ already in $o_L[[Z]]$. Hence $h \in o_L[[Z]]$.

For (ii) we have to show that $n_j = \sum_{i=1}^d f_{ij} \varphi_N(n_i)$, for any $1 \leq j \leq d$, with $f_{ij} \in \mathbf{A}_L^+$. For the same reasons as in the proof of (1) we have $n_j = \sum_{i=1}^d f_{ij}' \varphi_N(n_i) = \sum_{i=1}^d f_{ij}'' \varphi_N(n_i)$ with $f_{ij}' \in \mathbf{A}_L^+[\frac{1}{Q}]$ and $f_{ij}'' \in \mathcal{O}$. Then $\sum_{i=1}^d (f_{ij}' - f_{ij}'') \varphi_N(n_i) = 0$. But, again by the  definition of the category $\mathrm{Mod}_{\mathbf{A}_L^+}^{\varphi_L,\Gamma_L,an}$, the $\varphi_N(n_i)$ are linearly independent over $\mathbf{A}_L^+[\frac{1}{Q}]$ and hence over $\mathcal{O}[\frac{1}{Q}]$. It follows that $f_{ij}' = f_{ij}'' \in \mathbf{A}_L^+$.
\end{proof}

First we further investigate any $T$ in $\Rep_{o_L,f}^{cris,an}(G_L)$ whose {\bf Hodge-Tate weights are all $\leq 0$}, i.e., which is positive. For this purpose we need the ring $\mathbf{A}^+ = \mathbf{A} \cap W(\widetilde{\mathbf{E}}^+)_L$. One has the following general fact.

\begin{lemma}\label{fin-gen}
Let $F$ be any nonarchimedean valued field which contains $o_L/\pi_L o_L$, and let $o_F$ denote its ring of integers; we have:
\begin{itemize}
  \item[i.] Let $\alpha \in W(F)_L$ be any element; if the $W(o_F)_L$-submodule of $W(F)_L$ generated by $\{\varphi_q^i(\alpha)\}_{i \geq 0}$ is finitely generated then $\alpha \in W(o_F)_L$.
  \item[ii.] Let $X$ be a finitely generated free $o_L$-module, and let $M$ be a finitely generated $W(o_F)_L$-submodule of $W(F)_L \otimes_{o_L} X$; if $M$ is $\varphi_q \otimes \id$-invariant then $M \subseteq W(o_F)_L \otimes_{o_L} X$.
\end{itemize}
\end{lemma}
\begin{proof}
i. This is a simple explicit calculation as given, for example, in the proof of \cite[Lem.\ III.5]{Co}. ii. This is a straightforward consequence of i.
\end{proof}

\begin{proposition}\label{finite-height}
For positive $T$ in $\Rep_{o_L,f}^{cris,an}(G_L)$ we have
\begin{equation*}
  N(T) \subseteq D_{LT}^+(T) := (\mathbf{A}^+ \otimes_{o_L} T)^{\ker(\chi_{LT})} \ ,
\end{equation*}
and $N(T)$ is $p$-saturated in $D_{LT}^+(T)$.
\end{proposition}
\begin{proof}
By Remark \ref{positive}.i the $\mathbf{A}_L^+$-submodule $N(T)$ of $W(\widetilde{\mathbf{E}})_L \otimes_{o_L} T$ is $\varphi_q \otimes \id$-invariant (and finitely generated). Hence we may apply Lemma \ref{fin-gen}.ii to $M := W(\widetilde{\mathbf{E}}^+)_L \cdot N(T)$ and obtain that $N(T) \subseteq (W(\widetilde{\mathbf{E}}^+)_L \otimes_{o_L} T) \cap (\mathbf{A} \otimes_{o_L} T)^{\ker(\chi_{LT})} = D_{LT}^+(T)$. Since $N(T)$ is even $p$-saturated in $D_{LT}(T),$ the same holds with respect to the smaller $D_{LT}^+(T)$.
\end{proof}

\begin{corollary}\label{same-rank}
For positive $T$ in $\Rep_{o_L,f}^{cris,an}(G_L)$ the $\mathbf{A}_L^+$-module $D_{LT}^+(T)$ is free of the same rank as $N(T)$.
\end{corollary}
\begin{proof}
By the argument in the proof of \cite[Lem.\ III.3]{Co} the $\mathbf{A}_L^+$-module $D_{LT}^+(T)$ always is free of a rank less or equal to the rank of $N(T)$. The equality of the ranks in the positive case then is a consequence of Prop.\ \ref{finite-height}.
\end{proof}

Next we relate $N(T)$ to the construction in \cite[Prop.\ II.1.1]{Be}.

\begin{proposition}\label{Berger-Wach-module}
Suppose that $T$ in $\Rep_{o_L,f}^{cris,an}(G_L)$ is positive. For $N := N(T)$ we then have:
\begin{itemize}
  \item[i.] $N$ is the unique $\mathbf{A}_L^+$-submodule of $D_{LT}(T)$ which satisfies (N1) and (N2).
  \item[ii.] $N$ is also the unique $\mathbf{A}_L^+$-submodule of $D_{LT}^+(T)$ which satisfies:
\begin{itemize}
  \item[(a)] $N$ is free of rank equal to the rank of $D_{LT}^+(T)$;
  \item[(b)] $N$ is $\Gamma_L$-invariant;
  \item[(c)] the induced $\Gamma_L$-action on $N/\omega_{LT} N$ is trivial;
  \item[(d)] $\omega_{LT}^r D_{LT}^+(T) \subseteq N$ for some $r \geq 0$.
\end{itemize}
\end{itemize}
\end{proposition}
\begin{proof}
Let $P = P(\mathbf{A}_L^+)$ denote the set of height one prime ideals of $\mathbf{A}_L^+$. It contains the prime ideal $\mathfrak{p}_0 := (\omega_{LT})$.

\textit{Step 1:} We show the existence of a unique $\mathbf{A}_L^+$-submodule $N'$ of $D_{LT}^+(T)$ which satisfies (a) -- (d), and we show that this $N'$ is $\varphi_q$-invariant.

\textit{Existence:} We begin by observing that the $\mathbf{A}_L^+$-submodule $N := N(T)$ of $D_{LT}^+(T)$ has the properties (a), (b), and (c), but possibly not (d). In particular, the quotient $D_{LT}^+(T)/N$ is an $\mathbf{A}_L^+$-torsion module. Hence the localizations $N_{\mathfrak{p}} = D_{LT}^+(T)_{\mathfrak{p}}$ coincide for all but finitely many $\mathfrak{p} \in P$. By \cite[VII.4.3 Thm.\ 3]{B-CA} there exists a unique intermediate $\mathbf{A}_L^+$-module $N \subseteq N' \subseteq D_{LT}^+(T)$ which is finitely generated and reflexive and such that $N'_{\mathfrak{p}_0} = N_{\mathfrak{p}_0}$ and $N'_{\mathfrak{p}} = D_{LT}^+(T)_{\mathfrak{p}}$ for any $\mathfrak{p} \in P \setminus \{\mathfrak{p}_0\}$. Since $\mathbf{A}_L^+$ is a two dimensional regular local ring the finitely generated reflexive module $N'$ is actually free, and then, of course, must have the same rank as $N$ and $D_{LT}^+(T)$. We also have $N' = \bigcap_{\mathfrak{p}} N'_{\mathfrak{p}} = N_{\mathfrak{p}_0} \cap \bigcap_{\mathfrak{p} \neq \mathfrak{p}_0} D_{LT}^+(T)_{\mathfrak{p}}$. Since $\mathfrak{p}_0$ is preserved by $\varphi_{D_{LT}^+(T)}$ and $\Gamma_L$ it follows that $N'$ is $\varphi_{D_{LT}^+(T)}$- and $\Gamma_L$-invariant. Next the identities
\begin{equation*}
  L \otimes_{o_L} N/\omega_{LT} N = N_{\mathfrak{p}_0}/ \omega_{LT} N_{\mathfrak{p}_0} = N'_{\mathfrak{p}_0}/ \omega_{LT} N'_{\mathfrak{p}_0} = L \otimes_{o_L} N'/ \omega_{LT} N' \supseteq N'/ \omega_{LT} N'
\end{equation*}
show that the induced $\Gamma_L$-action on $N'/ \omega_{LT} N'$ is trivial. By using \cite[VII.4.4 Thm.\ 5]{B-CA} we obtain, for some $m_1, \ldots, m_d \geq 0$, a homomorphism of $\mathbf{A}_L^+$-modules $D_{LT}^+(T)/N' \longrightarrow \oplus_{i=1}^d \mathbf{A}_L^+/ \mathfrak{p}_0^{m_i} \mathbf{A}_L^+$ whose kernel is finite. Any finite $\mathbf{A}_L^+$-module is annihilated by a power of the maximal ideal in $\mathbf{A}_L^+$. We see that $D_{LT}^+(T)/N'$ is annihilated by a power of $\mathfrak{p}_0$, which proves (d).

\textit{Uniqueness:} Observing that $\gamma(\omega_{LT}) = [\chi_{LT}(\gamma)](\omega_{LT})$ for any $\gamma \in \Gamma_L$ (\cite[Lem.\ 2.1.15]{GAL}) this is exactly the same computation as in the uniqueness part of the proof of \cite[Prop.\ II.1.1]{Be}.

\textit{Step 2:} We show that $N'$ is $p$-saturated in $D_{LT}^+(T)$. By construction we have $(N')_{(\pi_L)} = D_{LT}^+(T)_{(\pi_L)}$. This implies that the $p$-torsion in the quotient $D_{LT}^+(T)/N'$ is finite. On the other hand, both modules, $N'$ and $D_{LT}^+(T)$, are free of the same rank. Hence the finitely generated $\mathbf{A}_L^+$-module $D_{LT}^+(T)/N'$ has projective dimension $\leq 1$ and therefore has no nonzero finite submodule (cf. \cite[Prop.\ 5.5.3(iv)]{NSW}).

\textit{Step 3:} We show that $N' = N$. Since both, $N$ and $N'$, are $p$-saturated in $D_{LT}^+(T)$ it suffices to show that the free $\mathbf{B}_L^+$-modules $N(V) := N[\frac{1}{p}]$ and $N'(V) := N'[\frac{1}{p}]$ over the principal ideal domain $\mathbf{B}_L^+ := \mathbf{A}_L^+[\frac{1}{p}]$ coincide. As they are both $\Gamma_L$-invariant, so is the annihilator ideal $I := ann_{\mathbf{B}_L^+}(N'(V)/N(V))$. Hence, by a standard argument as in \cite[Lem.\ 1.3.2]{Be}, the ideal $I$ is generated by an element $f$ of the form $\omega_{LT}^{\alpha_0} \prod_{n\geq 1}^s \varphi_L^{n-1}(Q)^{\alpha_n} $ with certain $\alpha_n \geq 0$, $0\leq n \leq s$, for some (minimal) $s\geq 0$. Since $N(V)_{(\omega_{LT} )} = N'(V)_{(\omega_{LT} )}$ by the construction of $N'$, it follows that $\alpha_0 = 0$. Assuming that $M := N'(V)/N(V) \neq 0$ we conclude that $s \geq 1$ (with $\alpha_s\geq 1$), i.e., that, with   $\mathfrak{p}_n := (\varphi_L^{n-1}(Q))$, we have $M_{\mathfrak{p}_s} \neq 0$ while $M_{\mathfrak{p}_{s+1}} = 0$. We claim that $(\varphi_L^*M)_{\mathfrak{p}_{s+1}} \neq 0. $ First note that we have have an exact sequence
\begin{equation*}
  0 \to (\mathbf{B}_L^+)^d \xrightarrow{A} (\mathbf{B}_L^+)^d \to M\to 0 \ ,
\end{equation*}
with $f$ dividing $\det(A) \in \mathbf{B}_L^+\setminus(\mathbf{B}_L^+)_{\mathfrak{p}_s}^\times$, which induces an exact sequence
\begin{equation*}
  0 \to (\mathbf{B}_L^+)^d \xrightarrow{\varphi_L(A)} (\mathbf{B}_L^+)^d \to \varphi_L^*M \to 0 \ .
\end{equation*}
Since $\varphi_L(f)= \prod_{n\geq 2}^{s+1} \varphi_L^{n-1}(Q)^{\alpha_{n-1}}$ divides $\det(\varphi_L(A))$ we conclude that $\det(\varphi_L(A))$ belongs to $\mathfrak{p}_{s+1} $ which implies the claim.

Now consider the following diagram with exact rows
\begin{equation*}
  \xymatrix{
 0  \ar[r]^{ }&  {(\varphi_L^*N(V))[\frac{1}{Q}]} \ar@{^(->}[d]_{ } \ar[r]^-{\cong } & N(V)[\frac{1}{Q}] \ar@{^(->}[d]_{ } \ar[r]^{ } & 0 \ar[d]_{ } \ar[r]^{} & 0     \\
  0 \ar[r]^{ } &  {(\varphi_L^*N'(V))[\frac{1}{Q}]} \ar[r]^{ } &  N'(V)[\frac{1}{Q}] \ar[r]^{ } & C \ar[r]^{} & 0.   }
\end{equation*}
The upper isomorphism comes from the definition of the category $\mathrm{Mod}_{\mathbf{A}_L^+}^{\varphi_L,\Gamma_L,an}$ in which $N$ lies. The map $(\varphi_L^*N'(V))[\frac{1}{Q}] \rightarrow N'(V)[\frac{1}{Q}]$ is injective since $\varphi_L^*N' \rightarrow N'$ is the restriction of the isomorphism $\varphi_L^* D_{LT}(T) \xrightarrow{\cong} D_{LT}(T)$. By the snake lemma and as $Q \notin \mathfrak{p}_{s+1}$ we obtain an injection
\begin{equation*}
  0 \neq (\varphi_L^*M)_{\mathfrak{p}_{s+1}} \hookrightarrow M_{\mathfrak{p}_{s+1}}=0 \ ,
\end{equation*}
which is a contradiction. Thus $M=0$ as had to be shown.
\end{proof}

\begin{remark}\phantomsection\label{twist}
\begin{enumerate}
  \item  $N(o_L(\chi_{LT}^{-1})) = \omega_{LT}\mathbf{A}_L^+ \otimes_{o_L} o_L\eta^{\otimes -1}$ and $N(o_L)= \mathbf{A}_L^+$.
  \item Let $o_L(\chi)=o_Lt_0$ with $\chi:G_L\to o_L^\times$ unramified. Then there exists an $a\in W(\bar{k}_L)_L^\times$ with $\sigma a=\chi^{-1}(\sigma) a$ for all $\sigma\in G_L$ by Remark \ref{rem:detcomp} $\phantom{m}$ \footnote{Since $\pi_L$ has trivial $G_L$-action the period $a$ there can be normalized such it becomes a unit in $ W(\bar{k}_L)_L$. };  in particular,
\begin{equation*}
  N(o_L(\chi)) = D^+_{LT}(o_L(\chi)) = \mathbf{A}_L^+n_0   \quad\text{for $n_0=a\otimes t_0$},
\end{equation*}
       where $\Gamma_L$ fixes $n_0$ and $\varphi_{N(o_L(\chi))}(n_0) = c n_0$ with $c := \frac{\varphi_L(a)}{a} \in o_L^\times$.
\end{enumerate}

\end{remark}
\begin{proof}
Each case belongs  to a positive representation $T$: in all cases the right hand side of the equality satisfies the properties   characterizing $N(T)$ in Prop.\ \ref{Berger-Wach-module}.ii (cf.\ \cite[Lem.\ 2.1.15]{GAL}).
\end{proof}

\begin{lemma}\label{twist-invariance}
For any $T\in \Rep_{o_L,f}^{cris,an}(G_L)$ we have:
\begin{itemize}
  \item[i.]  $N(T)$ is the unique $\mathbf{A}_L^+$-submodule of $D_{LT}(T)$ which satisfies (N1) and (N2);
  \item[ii.] $N(T(\chi_{LT}^{-r})) \cong \omega_{LT}^r N(T) \otimes_{o_L} o_L\eta ^{\otimes -r}$.
\end{itemize}
\end{lemma}
\begin{proof}
First we choose $r \geq 0$ such that $T(\chi_{LT}^{-r})$ is positive. Sending $N$ to $\omega_{LT}^{r} N(T) \otimes_{o_L} o_L\eta ^{\otimes -r} \subseteq D_{LT}(T)\otimes_{o_L} o_L\eta ^{\otimes -r}$ viewed in $D_{LT}(T)\otimes_{o_L} o_L\eta ^{\otimes -r} \cong D_{LT}(T(\chi_{LT}^{-r}))$ sets up a bijection between the $\mathbf{A}_L^+$-submodules of $D_{LT}(T)$ and $D_{LT}(T(\chi_{LT}^{-r}))$, respectively. One checks that $N$ satisfies (N1) and (N2) if and only if its image does. Hence i. and ii. (for such $r$) are a consequence of Prop.\ \ref{Berger-Wach-module}.i. That ii. holds in general follows from the obvious transitivity property of the above bijections.
\end{proof}

\begin{proposition}\label{elementarydiv}
Let $T$ be in $\Rep_{o_L,f}^{cris,an}(G_L)$ of $o_L$-rank $d$ and such that $V = L\otimes_{o_L} T$ is positive with Hodge-Tate weights $-r=-r_d \leq \cdots \leq - r_1 \leq 0$. Taking \eqref{f:compB} as an identification we then have
\begin{equation}\label{f:compBintegral}
    (\tfrac{t_{LT}}{\omega_{LT}})^{r}\mathcal{O}  \otimes_{\mathbf{A}_L^+} N(T) \subseteq\mathcal{O} \otimes_L D_{cris,L}(L \otimes_{o_L} T) \subseteq  \mathcal{O}  \otimes_{\mathbf{A}_L^+} N(T)
\end{equation}
with elementary divisors
\begin{equation*}
  [  \mathcal{O}  \otimes_{\mathbf{A}_L^+} N(T):\mathcal{O} \otimes_L D_{cris,L}(L \otimes_{o_L} T)]=[ (\tfrac{t_{LT}}{\omega_{LT}})^{r_1}:\cdots : (\tfrac{t_{LT}}{\omega_{LT}})^{r_d}].
\end{equation*}
\end{proposition}

\begin{proof}
We abbreviate $D := D_{cris,L}(V)$. By the definition of the functor $\mathcal{M}$ in \cite{KR} we have
\begin{equation}\label{f:defcalM}
  \mathcal{O}  \otimes_L D \subseteq \mathcal{M}(D) \subseteq  (\tfrac{t_{LT}}{\omega_{LT}})^{-r} \mathcal{O}  \otimes_L D \ .
\end{equation}
On the other hand, the commutativity of the big diagram before Remark \ref{properties} says that $\mathcal{M}(D) \cong \mathcal{O}  \otimes_{\mathbf{A}_L^+} N(T)$. This implies the inclusions \eqref{f:compBintegral}.

Concerning the second part of the assertion we first of all note that, although $\mathcal{O}$ is only a Bezout domain, it does satisfy the elementary divisor theorem (\cite[proof of Prop.\ 4.4]{ST1}). We may equivalently determine the elementary divisors of the $\mathcal{O}$-module $\mathcal{M}(D) / (\mathcal{O}  \otimes_L D)$. The countable set $\mathbb{S}$ of zeros of the function $\frac{t_{LT}}{\omega_{LT}} \in \mathcal{O}$ coincides with the set of nonzero torsion points of our Lubin-Tate formal group, each occurring with multiplicity one. The first part of the assertion implies that the module $\mathcal{O}$-module $\mathcal{M}(D) / (\mathcal{O}  \otimes_L D)$ is supported on $\mathbb{S}$. Let $\mathcal{M}_z(D)$, resp.\ $\mathcal{O}_z$, denote the stalk in $z \in \mathbb{S}$ of the coherent sheaf on $\mathbf{B}$ defined by $\mathcal{M}(D)$, resp.\ $\mathcal{O}$. The argument in the proof of \cite[Prop.\ 1.1.10]{BSX} then shows that we have
\begin{equation*}
  \mathcal{M}(D) / (\mathcal{O}  \otimes_L D) = \prod_{z \in \mathbb{S}} \mathcal{M}_z(D) / (\mathcal{O}_z  \otimes_L D) \ .
\end{equation*}
The ring $\mathcal{O}_z$ is a discrete valuation ring with maximal ideal $\mathfrak{m}_z$ generated by $\frac{t_{LT}}{\omega_{LT}}$. We consider on its field of fractions $\mathrm{Fr}(\mathcal{O}_z)$ the $\mathfrak{m}_z$-adic filtration and then on $\mathrm{Fr}(\mathcal{O}_z) \otimes_L D$ the tensor product filtration. By \cite[Lem.\ 1.2.1(2)]{Kis} (or \cite[Lem.\ 3.4.4]{BSX}) we have
\begin{equation*}
  \mathcal{M}_z(D) \cong \Fil^0(\mathrm{Fr}(\mathcal{O}_z) \otimes_L D)  \qquad\text{for any $z \in \mathbb{S}$},
\end{equation*}
and this isomorphism preserves $\mathcal{O}_z \otimes_L D$. At this point we let $0 \leq s_1 < \ldots < s_m < r$ denote the jumps of the filtration $\Fil^\bullet D$, i.e., the $r_j$ but without repetition. We write
\begin{equation*}
  D = D_1 \oplus \ldots \oplus D_m  \quad\text{such that}\quad  \Fil^{s_i} D = D_i \oplus \ldots \oplus D_m \ .
\end{equation*}

For the following computation let, for notational simplicity, $R$ denote any $L$-algebra which is a discrete valuation ring with maximal ideal $\mathfrak{m}$. We compute
\begin{align*}
  \Fil^0(\mathrm{Fr}(R) \otimes_L D) & = \sum_{j \in \mathbb{Z}} \mathfrak{m}^{-j} \otimes_L \Fil^j D = \sum_{j = 0}^r \mathfrak{m}^{-j} \otimes_L \Fil^j D  \\
   & = \sum_{i = 1}^m \mathfrak{m}^{-s_i} \otimes_L \Fil^{s_i} D = \sum_{i = 1}^m \sum_{j=i}^m  \mathfrak{m}^{-s_i} \otimes_L D_j  \\
   & = \sum_{j = 1}^m (\sum_{i=1}^j  \mathfrak{m}^{-s_i}) \otimes_L D_j =  \sum_{j = 1}^m \mathfrak{m}^{-s_j} \otimes_L D_j \ .
\end{align*}
Hence we obtain
\begin{equation*}
  \Fil^0(\mathrm{Fr}(R) \otimes_L D) / ( R \otimes_L D) = \oplus_{j=1}^m \mathfrak{m}^{-s_j}/R \otimes_L D_j \cong \oplus_{j=1}^m R / \mathfrak{m}^{s_j} \otimes_L D_j \ .
\end{equation*}

By combining all of the above we finally arrive at
\begin{align*}
  \mathcal{M}(D) / (\mathcal{O}  \otimes_L D) & = \prod_{z \in \mathbb{S}} \mathcal{M}_z(D) / (\mathcal{O}_z  \otimes_L D) \cong \prod_{z \in \mathbb{S}} \Fil^0(\mathrm{Fr}(\mathcal{O}_z) \otimes_L D) / (\mathcal{O}_z \otimes_L D)   \\
   & \cong \prod_{z \in \mathbb{S}} ( \oplus_{j=1}^m \mathcal{O}_z / \mathfrak{m}_z^{s_j} \otimes_L D_j )   =  \oplus_{j=1}^m ( \prod_{z \in \mathbb{S}} (\mathcal{O}_z / (\tfrac{t_{LT}}{\omega_{LT}})^{s_j}\mathcal{O}_z \otimes_L D_j ))   \\
   &  =  \oplus_{j=1}^m ( \prod_{z \in \mathbb{S}} \mathcal{O}_z / (\tfrac{t_{LT}}{\omega_{LT}})^{s_j}\mathcal{O}_z) \otimes_L D_j   =  \oplus_{j=1}^m  \mathcal{O} / (\tfrac{t_{LT}}{\omega_{LT}})^{s_j}\mathcal{O} \otimes_L D_j \ .
\end{align*}
\end{proof}

For a first application of this result we recall the comparison isomorphism
\begin{equation}\label{f:comp-iso}
\xymatrix{
  N(V) \ar[d]_{\subseteq} &  &  \\
  N(V)[\tfrac{1}{Q}] \ar[r]^-{\subseteq} & \mathcal{O}[\tfrac{\omega_{LT}}{t_{LT}}] \otimes_{\mathbf{A}_L^+} N(V) \ar[r]^-{\mathrm{comp}}_-{\cong}& \mathcal{O}[\tfrac{\omega_{LT}}{t_{LT}}] \otimes_L D_{cris,L}(V)   }
\end{equation}
for any $T$ in $\Rep_{o_L,f}^{cris,an}(G_L)$, $V := L \otimes_{o_L} T$, and $N(V):=N(T)[\frac{1}{\pi_L}]$. The left horizontal inclusion comes from the fact that $\frac{t_{LT}}{\omega_{LT}}$ is a multiple of $Q$ in $\mathcal{O}$. In particular, we have the commutative diagram
\begin{equation*}
\xymatrix{
  & N(V) \ar[dl]_{\varphi_{N(V)}} \ar[d]^{\varphi_{N(V)}} \ar[r]^-{\mathrm{comp}}  &  \mathcal{O}[\tfrac{\omega_{LT}}{t_{LT}}] \otimes_L D_{cris,L}(V) \ar[d]^{\varphi_L \otimes \varphi_{cris}}  \\
  N^{(\varphi)}(V) \ar[r]^{\subseteq} & N(V)[\tfrac{1}{Q}]  \ar[r]^-{\mathrm{comp}}  & \mathcal{O}[\tfrac{\omega_{LT}}{t_{LT}}] \otimes_L D_{cris,L}(V)   }
\end{equation*}
where $\varphi_{cris}$ denotes the $q$-Frobenius on $D_{cris,L}(V)$ and where
\begin{equation*}
  N^{(\varphi)}(V) := \text{the $\mathbf{A}_L^+$-submodule of $N(V)[\tfrac{1}{Q}]$ generated by the image of $N(V)$ under $\varphi_{N(V)}$}.
\end{equation*}
We note that, since $Q$ is invertible in $\mathbf{A}_L$, $N^{(\varphi)}(V)$ can also be viewed as the $\mathbf{A}_L^+$-submodule of $D_{LT}(V) = \mathbf{A}_L \otimes_{\mathbf{A}_L^+} N(V)$ generated by the image of $N(V)$ under $\varphi_{D_{LT}(V)}$. from this one easily deduces (use the projection formula for the $\psi$-operator) that the map $\psi_{D_{LT}(V)}$ on $D_{LT}(V)$ restricts to an operator
\begin{equation*}
  \psi_{N(V)} : N^{(\varphi)}(V) \longrightarrow N(V) \ .
\end{equation*}

\begin{corollary}
Assume that the Hodge-Tate weights of $V$ are all in $[0,r]$. Then we have
\begin{align}
  & \mathrm{comp}(N(V)) \subseteq \mathcal{O} \otimes_L D_{cris,L}(V),\ \mathrm{comp}(N^{(\varphi)}(V)) \subseteq \mathcal{O} \otimes_L D_{cris,L}(V), \ \text{and}  \\
 & \mathrm{comp}(N^{(\varphi)}(V)^{\psi_{N(V)}=0}) \subseteq \mathcal{O}^{\psi_L=0} \otimes_L D_{cris,L}(V) \ .   \label{f:phistarcomp}
\end{align}
\end{corollary}

\begin{proof}
Apply Prop.\  \ref{elementarydiv} to $T(\chi_{LT}^{-r})$, then divide the resulting (left) inclusion in \eqref{f:compBintegral} by $t_{LT}^r$ and tensor with $o_L(\chi^r_{LT})$. This gives the first inclusion by Lemma \ref{twist-invariance} upon noting that $t_{LT}^r D_{cris,L}(L \otimes_{o_L} T) \otimes_L L \eta^{\otimes -r} = D_{cris,L}(L \otimes_{o_L} T(\chi_{LT}^{-r}) )$. The second inclusion easily derives from the first by using that the map $\mathrm{comp}$ is compatible with the $\varphi$'s.

For the third inclusion we consider any element $x = \sum_i f_i \varphi_{N(V)}(x_i) \in N^{(\varphi)}(V)$, with $f_i \in \mathbf{A}_L^+$ and $x_i \in N(V)$, such that $\psi_{N(V)}(x) = \sum_i \psi_L(f_i) x_i = 0$. We choose an $L$-basis $e_1, \ldots, e_m$ of $D_{cris,L}(V)$ and write $\mathrm{comp}(x_i) = \sum_j f_{ij} \otimes e_j$ with $f_{ij} \in \mathcal{O}$. Then
\begin{equation*}
  0 = \mathrm{comp}(\psi_{N(V)}(x)) = \sum_i \psi_L(f_i) \mathrm{comp}(x_i) = \sum_i \sum_j \psi_L(f_i) f_{ij} \otimes e_j
\end{equation*}
and it follows that
\begin{equation*}
  \psi_L(\sum_i f_i \varphi_L(f_{ij})) = \sum_i \psi_L(f_i) f_{ij} = 0 \ ,
\end{equation*}
i.e., that $\sum_i f_i \varphi_L(f_{ij}) \in \mathcal{O}^{\psi_L=0}$. On the other hand we compute
\begin{align*}
  \mathrm{comp}(x) & = \sum_i f_i \varphi_{N(V)}(x_i) = \sum_i f_i (\varphi_L \otimes \varphi_{cris})(\mathrm{comp}(x_i))    \\
   & = \sum_i \sum_j f_i (\varphi_L(f_{ij}) \otimes \varphi_{cris}(e_j))   \\
   & = \sum_j \big( \sum_i f_i \varphi_L(f_{ij}) \big) \otimes \varphi_{cris}(e_j) \ .
\end{align*}
\end{proof}

\begin{corollary}\label{cor:Dcris}
In the situation of Prop.\  \ref{elementarydiv} we have
\begin{equation*}
  D_{cris,L}(V) \cong \left( \mathcal{O}\otimes_{\mathbf{A}_L^+} N(T) \right)^{\Gamma_L}.
\end{equation*}
\end{corollary}

\begin{proof}
We set $\mathcal{M} := \mathcal{O}\otimes_{\mathbf{A}_L^+} N(T)$ and identify $D(\mathcal{M})$ and $ D_{cris,L}(V)$ based on Lemma \ref{square} and using \eqref{f:compBintegral}. { The proof of \cite[Prop.\ (2.2.6)]{KR} combined with Remark \ref{properties} iii. implies the commutativity of the following diagram
\begin{equation*}
\xymatrix{
 D(\mathcal{M})\ar@{^(->}[r]_(0.4){incl.}\ar@{^(->}[rd]_{ incl.} & \mathcal{O}[\frac{\omega_{LT}}{t_{LT}}] \otimes_L D(\mathcal{M})  \ar[r]^(0.6){\xi}_(0.6){\cong} & \mathcal{M}[\frac{\omega_{LT}}{t_{LT}}] \\
 & \ar@{^(->}[u]_{ incl.} \mathcal{M}(D(\mathcal{M})) \ar[r]^{\cong} & \mathcal{M}, \ar@{^(->}[u]_{ incl.}   }
\end{equation*}
in which the right vertical map is the canonical inclusion while the left vertical map stems from the definition of the functor $\mathcal{M}$ as in \eqref{f:defcalM}  (which also implies the commutativity of the left triangle). Taking $\Gamma_L$-invariants and using the fact that the upper line induces the isomorphism $D(\mathcal{M}) \cong \mathcal{M}[\frac{\omega_{LT}}{t_{LT}}]^{\Gamma_L}$ in Remark \ref{properties} (iii)  the result follows.  }
\end{proof}

\begin{corollary}\label{cor:Q}
In the situation of Prop.\  \ref{elementarydiv} we have $Q^{r}N(V) \subseteq N^{(\varphi)}(V)$.
\end{corollary}

\begin{proof}
In the present situation $\varphi_{N(V)} : N(V) \rightarrow N(V)$ is an semilinear endomorphism of $N(V)$ by Remark \ref{positive}(i). Then $\varphi_{\mathcal{O} \otimes_{\mathbf{A}_L^+} N(V)} = \varphi_L \otimes \varphi_{N(V)} : \mathcal{O} \otimes_{\mathbf{A}_L^+} N(V) \rightarrow \mathcal{O} \otimes_{\mathbf{A}_L^+} N(V)$ is an endomorphism as well. The corresponding linearized maps are
\begin{align*}
  \varphi_{N(V)}^{lin} : \mathbf{A}_L^+ \otimes_{\mathbf{A}_L^+,\varphi_L} N(V) & \xrightarrow{\;\cong\;} N^{(\varphi)}(V) \subseteq  N(V) \\
  f \otimes x & \longmapsto f \varphi_{N(V)} (x)
\end{align*}
and
\begin{multline*}
  \varphi_{\mathcal{O} \otimes_{\mathbf{A}_L^+} N(V)}^{lin} = \id_{\mathcal{O}} \otimes \varphi_{N(V)}^{lin} : \mathcal{O} \otimes_{\mathbf{A}_L^+,\varphi_L} N(V) = \mathcal{O} \otimes_{\mathbf{A}_L^+} \big( \mathbf{A}_L^+ \otimes_{\mathbf{A}_L^+,\varphi_L} N(V) \big)  \\
   \longrightarrow \mathcal{O} \otimes_{\mathbf{A}_L^+} N(V) \ .
\end{multline*}
Since $\mathcal{O}$ is flat over $\mathbf{A}_L^+[\frac{1}{\pi_L}]$ it follows that
\begin{equation*}
  \mathcal{O} \otimes_{\mathbf{A}_L^+} N(V) / \im(\varphi_{\mathcal{O} \otimes_{\mathbf{A}_L^+} N(V)}^{lin}) = \mathcal{O} \otimes_{\mathbf{A}_L^+} (N(V) / N^{(\varphi)}(V)) \ .
\end{equation*}
But $\mathcal{O}$ is even faithfully flat over $\mathbf{A}_L^+[\frac{1}{\pi_L}]$. Hence the natural map
\begin{equation*}
  N(V) / N^{(\varphi)}(V) \longrightarrow \mathcal{O} \otimes_{\mathbf{A}_L^+} N(V) / \im(\varphi_{\mathcal{O} \otimes_{\mathbf{A}_L^+} N(V)}^{lin})
\end{equation*}
is injective. This reduces us to proving that
\begin{equation*}
  Q^r (\mathcal{O} \otimes_{\mathbf{A}_L^+} N(V)) \subseteq \im(\varphi_{\mathcal{O} \otimes_{\mathbf{A}_L^+} N(V)}^{lin}) \ .
\end{equation*}
As for any object in the category $\Mod_{\mathcal{O}}^{\varphi_L,\gamma_L,an}$, we do have
\begin{equation*}
  Q^h (\mathcal{O} \otimes_{\mathbf{A}_L^+} N(V)) \subseteq \im(\varphi_{\mathcal{O} \otimes_{\mathbf{A}_L^+} N(V)}^{lin}) \ .
\end{equation*}
for some sufficiently big integer $h$. On the other hand, \eqref{f:compBintegral} says that
\begin{equation*}
    (\tfrac{t_{LT}}{\omega_{LT}})^{r} \mathrm{comp} \big( \mathcal{O}  \otimes_{\mathbf{A}_L^+} N(V) \big) \subseteq \mathcal{O} \otimes_L D_{cris,L}(V) \subseteq  \mathrm{comp} \big(\mathcal{O}  \otimes_{\mathbf{A}_L^+} N(V) \big) \ .
\end{equation*}
Since $\varphi_{cris}$ is bijective we can sharpen the right hand inclusion to
\begin{equation*}
  \mathcal{O} \otimes_L D_{cris,L}(V) \subseteq  \mathrm{comp} \big( \im(\varphi_{\mathcal{O} \otimes_{\mathbf{A}_L^+} N(V)}^{lin}) \big) \ .
\end{equation*}
It follows that $(\tfrac{t_{LT}}{\omega_{LT}})^{r} (\mathcal{O}  \otimes_{\mathbf{A}_L^+} N(V)) \subseteq \im(\varphi_{\mathcal{O} \otimes_{\mathbf{A}_L^+} N(V)}^{lin})$. Since the greatest common divisor of $Q^h$ and $(\frac{t_{LT}}{\omega_{LT}})^r $ is $Q^{\min(h,r)}$ we finally obtain that $Q^r (\mathcal{O} \otimes_{\mathbf{A}_L^+} N(V)) \subseteq \im(\varphi_{\mathcal{O} \otimes_{\mathbf{A}_L^+} N(V)}^{lin})$.
\end{proof}

\begin{corollary}\label{cor:det}
In the situation of Prop.\  \ref{elementarydiv} we have, with regard to an $\mathbf{A}_L^+$-basis of $N := N(T)$ and with $s := \sum_{i=1}^d r_i$, that
\begin{equation*}
  \det (\varphi_N : N(T) \to N(T))= \det (\varphi_{N(V)} : N(V)\to N(V))= Q^{s}
\end{equation*}
up to an element in $o_L^\times \cdot (\varphi_L - 1)\big( (\mathbf{A}_L^+)^\times \big)$.
\end{corollary}

\begin{proof}
Note first that $N$ is $\varphi_N$-stable by Remark \ref{positive}(i). Moreover, the determinant of $\varphi_N$ acting on $N(V)$ equals the determinant of $\varphi_L \otimes \varphi_N$ acting on $\mathcal{O}\otimes_{\mathbf{A}_L^+} N(T)$, since we can take for both an $\mathbf{A}_L^+$-basis of $N(T)$. Since $\varphi_L(\frac{t_{LT}}{\omega_{LT}}) = \frac{\pi_L}{Q}\frac{t_{LT}}{\omega_{LT}}$, by proposition \ref{elementarydiv} the latter determinant equals $(\frac{\pi_L}{Q})^{-s} $ multiplied by the determinant of $\varphi_L \otimes \mathrm{Frob}$ acting on $\mathcal{O} \otimes_{L} D_{cris,L}(V)$. The latter is equal to the determinant of $\mathrm{Frob}$ on $D_{cris,L}(V)$, which is $\pi_L^s$ up to a unit in $o_L$ since the filtered Frobenius module $D_{cris,L}(V)$ is weakly admissible. This shows the claim up to an element in $o_L^\times\cdot(\varphi_L - 1)( \mathcal{O}^\times)$. But $\mathcal{O}^\times = \pi_L^{\mathbb{Z}} \times (\mathbf{A}_L^+)^\times$ by \cite[(4.8)]{Laz}. Hence $(\varphi_L - 1)( \mathcal{O}^\times) = (\varphi_L - 1)\big( (\mathbf{A}_L^+)^\times \big)$.
\end{proof}

\subsection{The determinant of the crystalline comparison isomorphism}

Let $T$ be any object in $\Rep_{o_L,f}^{cris,an}(G_L)$ of $o_L$-rank $d$ and such that $V=L\otimes_{o_L} T$ has Hodge-Tate weights $-r=-r_d\leq \cdots \leq - r_1$; we set $s:=\sum_{i=1}^d r_i$,  $N:=N(T)$ and $\mathcal{M}=\mathcal{O}\otimes_{ }N$. Consider the integral lattice
\begin{equation*}
  \mathcal{D} := \mathcal{D}(T) \subseteq D_{cris,L}(V)
\end{equation*}
which is defined as the image of $N/\omega_{LT}N \subseteq D(N)$ under the natural isomorphisms $D(N) \cong D(\mathcal{M}) \cong D_{cris,L}(V)$ arising from Lemma \ref{square} and \eqref{f:1}. Then with $N(-)$ also $\mathcal{D}(-)$ is a $\otimes$-functor. The aim of this subsection is to prove the following result.

\begin{proposition}\label{prop:detcomp}
With regard to bases of $T$ and $\mathcal{D}$ the determinant of the crystalline comparison isomorphism
\[B_{cris,L} \otimes_L V\cong B_{cris,L} \otimes_L D_{cris,L}(V)\]
belongs to $t_{LT}^s W(\bar{k}_L)_L^\times.$
\end{proposition}
We write $\bigwedge V$ for the highest exterior power of $V$ over $L.$

\begin{remark}
If $V$ is $L$-analyic (Hodge-Tate, crystalline), then so is $\bigwedge V.$
\end{remark}

Since $D_{cris,L}$ is a tensor functor, we are mainly reduced to consider characters $\rho:G_L\to L^\times,$ for which we denote by $V_\rho$ its representation space.

\begin{remark}\phantomsection\label{rem:rho}
\begin{enumerate}
\item If $V_\rho$ is Hodge-Tate, then $\rho$ coincides on an open subgroup of the inertia group $I_L$ of $G_L$ with
\begin{equation*}
  \prod_{\sigma\in \Sigma_L}\sigma^{-1}\circ \chi_{\sigma L,LT}^{n_\sigma},
\end{equation*}
      for some integers $n_\sigma,$ where $\Sigma_L$ denotes the set of embeddings of $L$ into $\bar{L}$
      and $\chi_{\sigma L,LT}$ is the Lubin-Tate character for $\sigma L$ and $\sigma(\pi_L).$
\item  If, in addition, $V_\rho$ is $L$-analytic, then $\rho$ coincides on an open subgroup of the inertia group $I_L$ with $\chi_{LT}^n$ for some integer $n.$
\end{enumerate}
\end{remark}

\begin{proof}
This follows from \cite[III.A4 Prop.\ 4 as well as III.A5 Thm.\  2 and its corollary]{Se0}.
\end{proof}

\begin{remark}\phantomsection\label{rem:detcomp} Let $\rho$  be a crystalline (hence Hodge-Tate) and $L$-analytic character. We then have:
\begin{enumerate}
\item If  $\rho$ factorizes through $G(L'/L)$ for some discretely valued Galois extension $L'$ of $L$, then the determinant of the crystalline comparison isomorphism for $V_\rho$ belongs to \linebreak $(W(\bar{k}_L)_L[\frac{1}{p}])^\times$ (with respect to arbitrary bases of $V$ and $ D_{cris,L}(V)$.)
\item If $\rho$  has Hodge-Tate   weight $-s$, then the determinant of the crystalline  comparison isomorphism for $V_\rho$ lies in $t_{LT}^s(W(\bar{k}_L)_L[\frac{1}{p}])^\times.$
\item   $\rho$ is of the form $\chi_{LT}^n\chi^{un}$ with an integer $n$ and an unramified character $\chi^{un}.$\footnote{ Also the converse statement is true: Indeed, any unramified character is locally algebraic (by definition, see \cite{Se0}), whence HT. By \cite[A3 Prop 3, A5 Thm 2]{Se0} the Hodge Tate weights are all zero, whence $\rho$ is $L$-analytic. It is crystalline as it is admissible, i.e. there is a period in $\mathbb{C}_p^\times$, which in this case needs to lie in the fixed field under inertia, which is contained in $B_{cris}$.}
\end{enumerate}
\end{remark}

\begin{proof} We shall write $K_0$ for the maximal absolutely unramified subextension of $K,$ any algebraic extension of $\Qp.$
Taking $G_{L'}$-invariants of the comparison isomorphism shows that the latter is already defined over
\begin{equation*}
  B_{cris,L}^{G_{L'}} = \left( L\otimes_{ L_{0}}B_{cris}\right)^{G_{L'}} = L\otimes_{L_0}\left(B_{cris}\right)^{G_{L'}} = L\otimes_{L_0}\widehat{L'_0}\subseteq W(\bar{k}_L)_L[\tfrac{1}{p}],
\end{equation*}
whence (i). Using Remark \ref{rem:rho} (ii) and applying (i) to $\rho\chi_{LT}^{-n}$ gives (ii). By the same argument it suffices to prove (iii) in the case of Hodge-Tate  weight $0.$ Then its period lies in the completion of the maximal unramified extension of $L$ by (i), whence the claim that $\rho$ is unramified follows, as the inertia subgroup of $G_L$ must act trivially.
\end{proof}

By Prop.\  \ref{finite-height} we have
\[N(T)\subseteq D^+_{LT}(T)\subseteq \mathbf{A}^+ \otimes_{o_L} T\]
if $T$ is positive. Using (N2) and the isomorphism
\[\mathbf{A}\otimes_{\mathbf{A}_L}D_{LT}(T)\cong\mathbf{A} \otimes_{o_L} T\] we obtain a canonical injection
\begin{equation}\label{f:NinjT}
\mathbf{A}^+ \otimes_{\mathbf{A}_L^{+}}N(T) \hookrightarrow \mathbf{A}^+ \otimes_{o_L} T.
\end{equation}

\begin{proposition}\label{prop:NinjT}
If $T$ is positive, then the determinant of \eqref{f:NinjT} with respect to bases of $N(T)$ and $T$ is contained in $\omega_{LT}^s(\mathbf{A}_L^+)^\times\cdot W(\bar{k}_L)_L^\times$.
\end{proposition}

\begin{proof}
Let $M\in M_d(\mathbf{A}^+)$ be the matrix of a basis of $N(T)$ with respect to a basis of $T$ and $P\in M_d(\mathbf{A}_L^+)$ the matrix of $\varphi_L$ with respect to the same basis of $N(T)$. Then we have $\varphi_L(M) = MP$.
By Corollary \ref{cor:det} we have $\det(P) = Q^s \varphi_L(f) f^{-1}u$ for some $f \in (\mathbf{A}_L^+)^\times$ and $u\in o_L^\times.$ But $Q = \varphi_L(\omega_{LT}) \omega_{LT}^{-1}$. We deduce that
\begin{equation*}
  \varphi_L(\det(M)) = \varphi_L(\omega_{LT}^s f) (\omega_{LT}^s f)^{-1}u \det(M) \ \text{, i.e., that}\  (\omega_{LT}^s fa)^{-1} \det(M) \in \mathbf{A}^{\varphi_L = 1} = o_L \ .
\end{equation*}
with $a\in W(\bar{k}_L)_L^\times$ such that $\varphi_L(a)/a=u.$
It follows that $\det(M) \in \omega_{LT}^s o_L (\mathbf{A}_L^+)^\times\cdot W(\bar{k}_L)_L^\times $. But we also have $\det(M) \in \mathbf{A}^\times$. Hence we finally obtain $\det(M) \in \omega_{LT}^s o_L (\mathbf{A}_L^+)^\times \cdot W(\bar{k}_L)_L^\times\cap \mathbf{A}^\times = \omega_{LT}^s(\mathbf{A}_L^+)^\times\cdot W(\bar{k}_L)_L^\times$.

\footnote{The argument also shows that $\omega_{LT}^{s}\mathbf{A}^+ \otimes_{o_L} T\subseteq\mathbf{A}^+ \otimes_{\mathbf{A}_L}N(T)\subseteq \mathbf{A}^+ \otimes_{o_L} T$. But the left inclusion does even hold with $r$ instead of $s$ (cf.\ \cite{Ste} following \cite{Be} in the cyclotomic case). }
\end{proof}

\begin{remark}\label{rem:unram}
For $T=o_L(\chi)$ with unramified $\chi$ as in Remark \ref{twist} the map \eqref{f:NinjT} maps the basis $n_0$ to $a\otimes t_0.$
\end{remark}

\begin{lemma}\label{lem:DinjN}
If $T$ is positive, then we have:
\begin{itemize}
  \item[i.] $\mathcal{O} \otimes_{o_L} \mathcal{D}(T) = \mathcal{O} \otimes_L D_{cris,L}(V) \subseteq \mathrm{comp}(\mathcal{O}\otimes_{\mathbf{A}_L^+} N(T))$;
  \item[ii.] the determinant of the inclusion in i. with respect to bases of $\mathcal{D}(T)$ and $N(T)$ belongs to $(\frac{t_{LT}}{\omega_{LT}})^s(\mathbf{A}_L^+)^\times$.
  \item[iii.] for $T=o_L(\chi)$ with unramified $\chi$ as in Remark \ref{twist}: $\mathrm{comp}(n_0)=\varphi_L(a)\otimes t_0=ca\otimes t_0\in D_{cris,L}(V)$ with $c=\frac{\varphi_L(a)}{a}\in o_L^\times;$ in particular, the element $a\otimes t_0$ is a basis of $ \mathcal{D}(T).$
\end{itemize}
\end{lemma}

\begin{proof}
By construction the comparison isomorphism \eqref{f:compB} is of the form
\begin{equation*}
  \mathrm{comp} = \id_{\mathcal{O}[\frac{\omega_{LT}}{t_{LT}}]} \otimes_L\, \mathrm{comp}_0
\end{equation*}
with
\begin{equation*}
  \mathrm{comp}_0 : \big( \mathcal{O} \otimes_{\mathbf{A}_L^+} N[\tfrac{\omega_{LT}}{t_{LT}}] \big)^{\Gamma_L} \xrightarrow[\pr]{\;\cong\;} N/\omega_{LT} N [\tfrac{1}{p}]=D(N) \xrightarrow{\;\cong\;} D_{cris,L}(V)
\end{equation*}
the right hand arrow being the natural isomorphism from Lemma \ref{square}. For positive $T$ we know in addition from the proof of Lemma \ref{cor:Dcris} that $\big( \mathcal{O} \otimes_{\mathbf{A}_L^+} N \big)^{\Gamma_L} = \big( \mathcal{O} \otimes_{\mathbf{A}_L^+} N[\tfrac{\omega_{LT}}{t_{LT}}] \big)^{\Gamma_L}$. We deduce that
\begin{equation*}
  \mathrm{comp}(\mathcal{O}\otimes_{\mathbf{A}_L^+} N) \supseteq \mathcal{O} \otimes_L \mathrm{comp}_0(\big( \mathcal{O} \otimes_{\mathbf{A}_L^+} N \big)^{\Gamma_L}) = \mathcal{O} \otimes_L D_{cris,L}(V) \ .
\end{equation*}
By Prop.\  \ref{elementarydiv} we know that the determinant in ii. is of the form $(\frac{t_{LT}}{\omega_{LT}})^s f(\omega_{LT})$ with $f(\omega_{LT}) \in \mathcal{O}^\times$. On the other hand, if we base change the inclusion in i. to $L = \mathcal{O}/\omega_{LT} \mathcal{O}$ then we obtain the base change from $o_L$ to $L$ of the isomorphism $\mathcal{D} \cong N/\omega_{LT} N$. By our choice of bases the determinant of the latter lies in $o_L^\times$. Since evaluation in zero maps $(\frac{t_{LT}}{\omega_{LT}})^s f(\omega_{LT})$ to $f(0)$ it follows that $f(0)$ belongs to $o_L^\times$ and hence (\cite[(4.8)]{Laz}) that $f(\omega_{LT})$ belongs to $(\mathbf{A}_L^+)^\times$.

Now we prove iii.: By the above description of $\mathrm{comp}_0$ we have to show that the image $\bar{n}_0\in D(N(T))$ of $n_0$ is mapped to $ca\otimes t_0$ under the natural isomorphism from Lemma \ref{square}. Since under the crystalline comparison isomorphisms these elements are sent to $a\otimes (a^{-1}\otimes \bar{n}_0)\in B_{cris,L}\otimes_L V_L(D(N))$ and $ca\otimes t_0\in B_{cris,L}\otimes_{o_L} T,$ respectively, it suffices to show that the map \eqref{f:2} sends $a^{-1}\otimes n_0\in L\otimes_{o_L} V(\mathbf{A}_L\otimes_{\mathbf{A}_L^+}N)$ (which corresponds to $t_0$ under the canonical isomorphism $T\cong V(\mathbf{A}_L\otimes_{\mathbf{A}_L^+}N)$) to $(ca)^{-1}\otimes \bar{n}_0\in V_L(D(N)).$ Dualizing, this is equivalent to the claim that the map \eqref{f:3} sends the dual basis $\delta_{a^{-1}\otimes n_0}\in (L\otimes_{o_L}V(M))^*$ of $a^{-1}\otimes n_0 $ to $\delta_{(ca)^{-1}\otimes \bar{n}_0}\in V_L(D(N))^*.$ Note that the isomorphism
\begin{equation*}
(L\otimes_{o_L}V(M))^*\cong L\otimes_{o_L} \mathrm{Hom}_{\mathbf{A}_L ,\varphi_q}(\mathbf{A}_L\otimes_{\mathbf{A}_L^+} N,\mathbf{A} )\cong  L\otimes_{o_L} \mathrm{Hom}_{\mathbf{A}_L^+ ,\varphi_q}(N,\mathbf{A}^+[\tfrac{1}{\omega_{LT}}])
\end{equation*}
sends $ \delta_{a^{-1}\otimes n_0}$ to $a\delta_{n_0}. $ Thus it suffices to show that the map \eqref{f:V2} sends $a\delta_{n_0} $  to  $ca\delta_{\bar{n}_0}$ in $\mathrm{Hom}_{L,\varphi_q,\mathrm{Fil}}((N/\omega_{LT}N)[\tfrac{1}{p}],B_{cris,L})$,
since the latter corresponds under \eqref{f:V3} to $ca\otimes\delta_{\bar{n}_0}\in V_L(D(N)^*) $ which in turn corresponds to $ \delta_{(ca)^{-1}\otimes \bar{n}_0}$ under \eqref{f:V4}.

If $f= a\delta_{n_0},$ which is the map which sends $n_0$ to $a$, then - in the notation of the proof of Lemma \ref{square} - $f_1$ and $f_2$ share this property, while $f_3$ (and hence $f_4$) sends $c^{-1}{n_0}$
to $a$, because $\varphi_N(c^{-1}n_0)= c^{-1}\varphi_L( a)a^{-1} n_0=n_0.$ Then $f_5$ sends $c^{-1}\bar{n}_0$ to $a,$  because $\xi(c^{-1}\bar{n}_0)=c^{-1}n_0.$ Altogether this means, that $a\delta_{n_0}$ is mapped to $\varphi_L(a)\delta_{n_0}=ca\delta_{\bar{n}_0}$ as claimed.
\end{proof}

\begin{proof}[Proof of Prop.\  \ref{prop:detcomp}]

 The functor $D_{cris,L}(-)$ on crystalline Galois representations is a $\otimes$-functor and commutes with exterior powers, and the crystalline comparison isomorphism is compatible with tensor products and exterior powers. The analogous facts hold for the functor $N(-)$ and hence for the functor $\mathcal{D}(-)$ (by base change). The case of the functor $N(-)$ reduces, by using the properties (N1) and (N2) in Lemma \ref{twist-invariance} (i), to the case of the functor $D_{LT}(-).$ Here the properties can easily be seen by the comparison isomorphism \eqref{f:compiso}.

 Upon replacing $T$ by its highest exterior power\Footnote{Siehe Bourbaki, Algebra I, Chap. III, §7, f\"{u}r die folgenden Fakten: \"{A}u{\ss}ere Potenz ist definiert als Quotient der entsprechenden  Tensorpotenz, die \"{a}u{\ss}ere Potenz eines freien Moduls ist wieder frei, Bildung der \"{a}u{\ss}eren Potenz kommutiert mit Basiswechsel, die Vertauschbarkeit von $D_{cris,L}$ mit der \"{a}u{\ss}eren Potenz folgt aus Vergleichsisomorphismus, ebenso f\"{u}r $D_{LT}.$}
we may and do assume that the $o_L$-module $T$ has rank  $1.$ In addition by twisting $T$ if necessary with a power of $\chi_{LT}$ we may and do assume that $T$ is positive with $s=0,$ i.e., unramified by \ref{rem:detcomp}. In this case it is clear that - using the notation of Lemma \ref{lem:DinjN} iii.  - the crystalline comparison isomorphism sends $t_0$ to $a\otimes t_0.$ Since the latter is also a basis of $\mathcal{D}(T)$ by the same Lemma, the proposition follows.\Footnote{But apparently it is not true in general that the determinant of the crystalline comparison isomorphism is the product of the determinants in Lemma \ref{lem:DinjN} and Prop.\  \ref{prop:NinjT} (with compatible basis)! In the unramified rank one case they differ by $c$ or its inverse!}

\end{proof}

\subsection{Non-negative Hodge-Tate weights}

Now assume that for $T$ in $\Rep_{o_L,f}^{cris,an}(G_L)$ the {\bf Hodge-Tate weights are all $\geq 0$} and set $N := N(T)$.  By \cite[Remark 3.2.i.-ii.]{SV15} the map $\psi_L$ preserves $\mathbf{A}_L^+$. It follows that $\psi_{D_{LT}(T)}$ maps $\mathbf{A}_L^+ \cdot \varphi_N(N)$ - and hence $N$ by Remark \ref{positive}.i(2)  - into $N$. The following lemmata generalize those of  \cite[Appendix A]{B}.

\begin{lemma}\label{psi-omega}
For $m\geq 1,$ there exists $Q_m\in o_L[[Z]]$ such that
\[
\psi_L(\frac{1}{\omega_{LT}^m})=\frac{\pi_L^{m-1}+\omega_{LT}Q_m(\omega_{LT})}{\omega_{LT}^m}.
\]
\end{lemma}

\begin{proof}
According to the paragraph after Remark 2.1 in \cite{SV15} combined with Remark 3.2 ii. in (loc.\ cit.) we have that
\begin{equation*}
   h(\omega_{LT}) := \omega_{LT}^{m}\psi_L(\frac{1}{\omega_{LT}^m}) = \psi_L(\frac{[\pi_L]^m}{\omega_{LT}^m})\in \mathbf{A}_L^+.
\end{equation*}
Obviously there exists $Q_m\in o_L[[Z]]$ such that
\[
h(\omega_{LT})-h(0)=\omega_{LT}Q_m(\omega_{LT}).
\]
Thus the claim follows from
  \begin{align*}
 h(0) & = \varphi_L(h(\omega_{LT}))_{|\omega_{LT}=0}= \varphi_L\circ \psi_L(\frac{[\pi_L]^m}{\omega_{LT}^m})_{|\omega_{LT}=0}= \pi_L^{-1} \sum_{a\in LT_1}\left(\frac{[\pi_L]^m(a+_{LT}\omega_{LT})}{(a+_{LT}\omega_{LT})^m}\right)_{|\omega_{LT}=0}  \\
      & = \pi_L^{-1} \sum_{a\in LT_1}\left(\frac{[\pi_L]( \omega_{LT})}{a+_{LT}\omega_{LT}}\right)^m_{|\omega_{LT}=0}=\pi_L^{m-1},
\end{align*}
because $\left(\frac{[\pi_L]( \omega_{LT})}{a+_{LT}\omega_{LT}}\right)_{|\omega_{LT}=0} = \pi_L$ for $a=0$ and $= 0$ otherwise.
\end{proof}

\begin{lemma}\label{psi-shrinking}
We have
\[
\psi_{D_{LT}(T)}(\pi_LD_{LT}(T)+\omega_{LT}^{-1}N(T))\subseteq \pi_L D_{LT}(T)+\omega_{LT}^{-1}N(T)
\]
and, for $k\geq 1,$
\[
\psi_{D_{LT}(T)}(\pi_L D_{LT}(T)+\omega_{LT}^{-(k+1)}N(T))\subseteq \pi_L D_{LT}(T)+\omega_{LT}^{-k}N(T).
\]
\end{lemma}

\begin{proof}
By Remark \ref{positive} (2) we can write any $x\in N(T)$ in the form $x=\sum a_i\varphi_N(x_i)$ with $a_i\in \mathbf{A}_L^+$ and $x_i\in N(T).$ Therefore $\psi_{D_{LT}(T)}(\omega_{LT}^{-(k+1)}x)=\sum \psi_L(\omega_{LT}^{-(k+1)}a_i)x_i$ by the projection formula. Since $\psi_L$ preserves $\mathbf{A}_L^+$ and is $o_L$-linear we conclude by Lemma \ref{psi-omega} that $\psi_L(\omega_{LT}^{-(k+1)}a_i)$ belongs to $\pi_L \mathbf{A}_L+\omega_{LT}^{-k}\mathbf{A}_L^+$, whenever $k\geq 1,$ from which the second claim follows as $\psi_{D_{LT}(T)}(\pi_L D_{LT}(T))\subseteq \pi_LD_{LT}(T)$ by $o_L$-linearity of $\psi_{D_{LT}(T)}.$ For $k=0$ finally, $\psi_L(\omega_{LT}^{-1}a_i)$ belongs to $ \omega_{LT}^{-1}\mathbf{A}_L^+$, from which the first claim follows.
\end{proof}

\begin{lemma}\label{psi-id-shrinking}
If $k\geq 1$ and $x\in D_{LT}(T)$ satisfies $\psi_{D_{LT}(T)}(x)-x\in \pi_L D_{LT}(T)+\omega_{LT}^{-k}N(T),$ then $x$ belongs to $\pi_L D_{LT}(T)+\omega_{LT}^{-k}N(T).$
\end{lemma}

\begin{proof}
Since $D_{LT}(T)/\pi_LD_{LT}(T)$ is a finitely generated (free) $k_L((\omega_{LT}))$-module there exists an integer $m\geq 0$ such that $x\in \pi_LD_{LT}(T)+\omega_{LT}^{-m}N(T);$ let $l$ denote the smallest among them.
Assume that $l>k.$ Then Lemma \ref{psi-shrinking} shows that
\[
\psi_{D_{LT}(T)}(x)\in \pi_LD_{LT}(T)+\omega_{LT}^{-(l-1)}N(T).
\]
Hence $\psi_{D_{LT}(T)}(x)-x$ would belong to  $\pi_LD_{LT}(T)+\omega_{LT}^{-l}N(T)$ but not to $(\pi_LD_{LT}(T)+\omega_{LT}^{-(l-1)}N(T)),$ a contradiction to our assumption. It follows that $l\leq k$, and we are done.
\end{proof}

\begin{lemma}\label{psi-Wach}
It holds $D_{LT}(T)^{\psi_{D_{LT}(T)}=1}\subseteq \omega_{LT}^{-1}N(T),$ i.e., \[D_{LT}(T)^{\psi_{D_{LT}(T)}=1}=\left(\omega_{LT}^{-1}N(T)\right)^{\psi_{D_{LT}(T)}=1}.\]
\end{lemma}
\begin{proof}
By induction on $k\geq 1$ we will show that $ D_{LT}(T)^{\psi_{D_{LT}(T)}=1}\subseteq \pi_L^kD_{LT}(T)+ \omega_{LT}^{-1}N(T),$ i.e., writing $x=\pi_L^ky_k+n_k\in D_{LT}(T)^{\psi_{D_{LT}(T)}=1} $ the sequence $n_k$ will $\pi_L$-adically converge in $\omega_{LT}^{-1}N(T)$ with limit $x.$

In order to show the claim assume $x\in D_{LT}(T)^{\psi_{D_{LT}(T)}=1}$. As in the previous proof there exists some minimal  integer $m\geq 0$ such that $x\in \pi_LD_{LT}(T)+\omega_{LT}^{-m}N(T).$  Then $m=1$ and we are done since otherwise  Lemma \ref{psi-id-shrinking} implies that $m$ can be   decreased by $1$.  This proves the claim for $k=1.$

By our induction hypothesis we can write $x\in D_{LT}(T)^{\psi_{D_{LT}(T)}=1}$ as $x=\pi_L^ky+n$ with $y\in D_{LT}(T)$ and $n \in \omega_{LT}^{-1}N(T).$ The equation $\psi_{D_{LT}(T)}(x)=x$ implies that $\psi_{D_{LT}(T)}(n)-n=\pi_L^k(\psi_{D_{LT}(T)}(y)-y).$  In the proof of Lemma \ref{psi-shrinking} we have seen that $\psi_{D_{LT}(T)}(n) - n \in \omega_{LT}^{-1}N(T)$. Note that $\pi_L^kD_{LT}(T)\cap \omega_{LT}^{-1}N(T)=\pi_L^k\omega_{LT}^{-1}N(T)$ because $\mathbf{A}_L/\omega_{LT}^{-1}\mathbf{A}_L^+$ has no $\pi_L$-torsion. Therefore $\psi_{D_{LT}(T)}(y)-y \in \omega_{LT}^{-1}N(T)$, whence $y$, by Lemma \ref{psi-id-shrinking}, belongs to $\pi_L D_{LT}(T)+\omega_{LT}^{-1}N(T)$ so that we can write $x=\pi_L^k(\pi_Ly'+n')+n=\pi_L^{k+1}y'+(\pi_L^kn'+n)$ as desired.
\end{proof}

Set $V := T\otimes_{o_L}L$.

\begin{lemma}\label{trivialquotient}
If $D_{cris,L}(V)^{\varphi_q=1}\neq 0,$ then $V$ has the trivial representation $L$ as quotient, i.e., the co-invariants $V_{G_L}$ are non-trivial.
\end{lemma}

\begin{proof} Let $W=V^*$ be the $L$-dual of $V.$ Then, by \cite[(51)]{SV15} we have \Footnote{ Achtung: dort arbeiten wir mit $B_{max,L}\subsetneqq B_{max}\otimes L$ die Aussage hier sollte daher aber erst recht gelten, wir m\"{u}ssen hier nur sorgf\"{a}ltig sein! Au{\ss}erdem: Hier ist $D_{cris,L}(V) = (B_{cris} \otimes_{L_0} V)^{G_L}$, aber in \cite{SV15} war $D_{cris,L}(V) = (B_{max,\mathbb{Q}_p} \otimes_{\mathbb{Q}_p} V)^{G_L}$.}
\begin{equation*}
  (V_{G_L})^* \cong H^0(L,W) \cong D_{cris,L}(W)^{\varphi_q=1}\cap (B_{dR}^+ \otimes_{L} W)^{G_L} = D_{cris,L}(W)^{\varphi_q=1} \neq 0,
\end{equation*}
because $(B_{dR}^+ \otimes_{L} W)^{G_L}=(B_{dR} \otimes_{L} W)^{G_L}\supseteqq D_{cris,L}(W)$ since the Hodge-Tate weights of $W$ are $\leq 0.$
\end{proof}

\begin{lemma}\label{notrivialquot}
If $V$ does not have any quotient isomorphic to the trivial representation $L$, then $D_{LT}(T)^{\psi_{D_{LT}(T)}=1} \subseteq N(T)$, i.e.,
\begin{equation*}
  D_{LT}(T)^{\psi_{D_{LT}(T)}=1} = N(T)^{\psi_{D_{LT}(T)}=1} \ .
\end{equation*}
\end{lemma}
\begin{proof}
Because of Lemma \ref{psi-Wach} it suffices to show that $(\omega_{LT}^{-1} N(T))^{\psi_{D_{LT}(T)}=1} \subseteq N(T)$. Let $e_1,\ldots ,e_d$ be a basis of $N := N(T)$ over $\mathbf{A}_L^+.$ Then, by Remark \ref{positive} (ii) there exist $\beta_{ij} = \sum_{\ell \geq 0} \beta_{ij,\ell}\omega_{LT}^\ell \in \mathbf{A}_L^+$ such that $e_i=\sum_{j=1}^d\beta_{ij}\varphi_N(e_j).$ Now assume that $\omega_{LT}^{-1}n=\sum_{i=1}^d \alpha_{i}e_i=\sum_{i,j} \alpha_{i}\beta_{ij}\varphi_N(e_j)$ belongs to $(\omega_{LT}^{-1}N)^{\psi_{D_{LT}(T)}=1}$ with $\alpha_i = \sum_{\ell \geq -1} \alpha_{i,\ell}\omega_{LT}^\ell \in \omega_{LT}^{-1}\mathbf{A}_L^+$. By the projection formula this implies, for $1 \leq j\leq d$,
\begin{equation*}
  \alpha_j=  \psi_L(\sum_{i=1}^d \alpha_i \beta_{ij})\equiv \omega_{LT}^{-1}\sum_{i=1}^d \alpha_{i,-1} \beta_{ij,0}  \mod  \mathbf{A}_L^+
\end{equation*}
because $\psi_L(\omega_{LT}^{-1})\equiv\omega_{LT}^{-1} \mod  \mathbf{A}_L^+$ by Lemma \ref{psi-omega}, whence
\begin{align*}
  \varphi_L (\omega_{LT}) \varphi_L(\alpha_j)&\equiv\sum_{i=1}^d \alpha_{i,-1} \beta_{ij,0}  \mod  \omega_{LT}\mathbf{A}_L^+.
\end{align*}
It follows from the definition of $\beta_{ij}$ that
\begin{align*}
    \varphi_N(n)&= \sum_j \varphi_L (\omega_{LT})\varphi_L(\alpha_j)\varphi_N(e_j)\equiv\sum_{j,i} \alpha_{i,-1} \beta_{ij,0}\varphi_N(e_j)\equiv \sum_{i} \alpha_{i,-1} e_i\equiv n \mod \omega_{LT}N,
\end{align*}
i.e., that $D_{cris,L}(V)\cong N/\omega_{LT}N[\frac{1}{p}]$ (by \eqref{f:1} and Lemma \ref{square}) contains an eigenvector for $\varphi_q$ with eigenvalue $1,$ if $\omega_{LT}^{-1}n$ does not belong to $N.$ Now the result follows from Lemma \ref{trivialquotient}.
\end{proof}

\newpage

\section{\texorpdfstring{$(\varphi_L,\Gamma_L)$}{(phi,Gamma)}-modules over the Robba ring }

\subsection{Robba rings of character varieties}\label{sec:basics}

Throughout our coefficient field $K$ is a complete intermediate extension $L \subseteq K \subseteq \mathbb{C}_p$.
For any reduced affinoid variety $\mathfrak{Y}$ over $\mathbb{Q}_p$ of $L$ we let $|\ |_{\mathfrak{Y}}$ denote the supremum norm on the affinoid algebra $\mathcal{O}_K(\mathfrak{Y})$ of $K$-valued holomorphic functions on $\mathfrak{Y}$. It is submultiplicative and defines the intrinsic Banach topology of this algebra.

\subsubsection{The additive character variety and its Robba ring}\label{sec:add-Robba}

 Let $\mathbf{B}_1$ denote the rigid $\mathbb{Q}_p$-analytic open disk of radius one around the point $1 \in \mathbb{Q}_p$. The rigid analytic group variety
\begin{equation*}
    \mathfrak{X}_0 := \mathbf{B}_1 \otimes_{\mathbb{Z}_p} \Hom_{\mathbb{Z}_p}(o_L,\mathbb{Z}_p)
\end{equation*}
over $\mathbb{Q}_p$ (which noncanonically is a $d$-dimensional open unit polydisk) parametrizes the locally $\mathbb{Q}_p$-analytic characters of the additive group $o_L$: the point $z \otimes \beta$ is sent to the character $\chi_{z \otimes \beta}(a) := z^{\beta(a)}$. It is shown in \cite[\S2]{ST2} that the rigid analytic group variety $\mathfrak{X}$ over $L$, which parametrizes the locally $L$-analytic characters of $o_L$, is the common zero set in $\mathfrak{X}_{0/L}$ of the functions
\begin{equation*}
    \sum_{j=1}^d z_j \otimes \beta_j \longmapsto \sum_{j=1}^d (\beta_j(t_i) - t_i \cdot \beta_j(1)) \cdot \log(z_j)
\end{equation*}
for $1 \leq i \leq d$; here $t_1, \ldots, t_d$ is a $\mathbb{Z}_p$-basis of $o_L$ and $\beta_1,\ldots,\beta_d$ is the corresponding dual basis. It is one dimensional, smooth, and connected. As a closed analytic subvariety of the Stein space $\mathfrak{X}_0$ the rigid variety $\mathfrak{X}$ is Stein as well.

For any $a \in o_L$ the map $b \longmapsto ab$ on $o_L$ is locally $L$-analytic. This induces an action of the multiplicative monoid $o_L \setminus \{0\}$ first on the $\mathbb{Z}_p$-module $\Hom_{\mathbb{Z}_p}(o_L,\mathbb{Z}_p)$ and then on the varieties $\mathfrak{X}_0$ and $\mathfrak{X}$. The latter actions further induce actions on the rings of $K$-valued holomorphic functions $\mathcal{O}_K(\mathfrak{X}_0) \twoheadrightarrow \mathcal{O}_K(\mathfrak{X})$, which we will denote by $(a,f) \mapsto a_*(f)$.

We also have induced translation actions of $o_L \setminus \{0\}$ on the vectors spaces $C^{an}_{\mathbb{Q}_p}(o_L,K)$, resp.\ $C^{an}(o_L,K)$, of $K$-valued locally $\mathbb{Q}_p$-analytic, resp. $L$-analytic, functions on $o_L$ and then by duality on the spaces $D_{\mathbb{Q}_p}(o_L,K) \twoheadrightarrow D(o_L,K)$ of locally $\mathbb{Q}_p$-analytic and locally $L$-analytic distributions on $o_L$, respectively; they will be denoted by $(a,\lambda) \mapsto a_*(\lambda)$. By \cite[Thm.\ 2.3]{ST2} we have the Fourier isomorphism
\begin{align}\label{f:Fourier}
  D(o_L,K) & \xrightarrow{\;\cong\;} \mathcal{O}_K(\mathfrak{X}) \\
  \lambda & \longmapsto F_\lambda(\chi) = \lambda(\chi) \ .\notag
\end{align}
One easily checks that this isomorphism is $o_L \setminus \{0\}$-equivariant. In the following we will denote the endomorphism $(\pi_L)_*$ in all situations also by $\varphi_L$. The Fourier isomorphism maps the Dirac distribution $\delta_a$, for any $a \in o_L$, to the evaluation function $\ev_a(\chi) := \chi(a)$.\\

\noindent
{\bf{The $\psi$-operator and the Mellin transform}} \\

\begin{lemma}\label{phi-free}
The endomorphism $\varphi_L$ makes $\mathcal{O}_K(\mathfrak{X})$ into a free module over itself of rank equal to the cardinality of $o_L/\pi_L o_L$; a basis is given by the functions $\ev_a$ for $a$ running over a fixed system of representatives for the cosets in $o_L/\pi_L o_L$.
\end{lemma}
\begin{proof}
This is most easily seen by using the Fourier isomorphism which reduces the claim to the corresponding statement about the distribution algebra $D(o_L,K)$. But here the ring homomorphism $\varphi_L$ visibly induces an isomorphism between $D(o_L,K)$ and the subalgebra $D(\pi_L o_L,K)$ of $D(o_L,K)$. Let $R \subseteq o_L$ denote a set of representatives for the cosets in $o_L/\pi_L o_L$. Then the Dirac distributions $\{\delta_a\}_{a \in R}$ form a basis of $D(o_L,K)$ as a $D(\pi_L o_L,K)$-module.
\end{proof}

\begin{lemma}\label{module-struc}
The $o_L^\times$-action on $D(o_L,K) \cong \mathcal{O}_K(\mathfrak{X})$ extends naturally to a (jointly) continuous $D(o_L^\times,K)$-module structure.
\end{lemma}
\begin{proof}
In a first step we consider the case $K = L$, so that $K$ is spherically complete. By \cite[Cor.\ 3.4]{ST1} it suffices to show that $C^{an}(G,K)$ as an $o_L^\times$-representation is locally analytic. This means we have to establish that, for any $f \in C^{an}(G,K)$, the orbit map $a \longmapsto a^*(f)$ on $o_L^\times$ is locally analytic. But this map is the image of the locally analytic function $(a,g) \longmapsto f(ag)$ under the isomorphism $C^{an}(o_L^\times \times G,K) = C^{an}(o_L^\times,C^{an}(G,K))$ in \cite[Lem\. A.1]{ST3}.

Now let $K$ be general. All tensor products in the following are understood to be formed with the projective tensor product topology. By the universal property of the latter the jointly continuous bilinear map $D(o_L^\times,L) \times \mathcal{O}_L(\mathfrak{X}) \rightarrow \mathcal{O}_L(\mathfrak{X})$ extends uniquely to a continuous linear map $D(o_L^\times,L) \widehat{\otimes}_L \mathcal{O}_L(\mathfrak{X}) \rightarrow \mathcal{O}_L(\mathfrak{X})$. This further extends to the right hand map in the sequence of continuous $K$-linear maps
\begin{equation*}
  \big(K \widehat{\otimes}_L D(o_L^\times,L) \big) \widehat{\otimes}_K \big(K \widehat{\otimes}_L \mathcal{O}_L(\mathfrak{X}) \big) \rightarrow
  K \widehat{\otimes}_L \big( D(o_L^\times,L) \widehat{\otimes}_L \mathcal{O}_L(\mathfrak{X}) \big) \rightarrow  K \widehat{\otimes}_L \mathcal{O}_L(\mathfrak{X}) \ .
\end{equation*}
The left hand map is the obvious canonical one. We refer to \cite[\S 10.6]{PGS} for the basics on scalar extensions of locally convex vector spaces. The same reasoning as in the proof of \cite[Prop.\ 2.5.ii]{BSX} shows that $K \widehat{\otimes}_L \mathcal{O}_L(\mathfrak{X}) = \mathcal{O}_K(\mathfrak{X})$. It remains to check that $K \widehat{\otimes}_L D(o_L^\times,L) = D(o_L^\times,K)$ holds true as well. For any open subgroup $U \subseteq o_L^\times$ we have $D(o_L^\times,-) = \oplus_{a \in o_L^\times/U} \delta_a D(U,-)$. Hence it suffices to check that $K \widehat{\otimes}_L D(U,L) = D(U,K)$ for one appropriate $U$. But $o_L^\times$ contains such a subgroup $U$ which is isomorphic to the additive group $o_L$ so that $D(U,-) \cong D(o_L,-) \cong \mathcal{O}_-(\mathfrak{X})$. In this case we had established our claim already.
\end{proof}

The operator $\varphi_L$ has a distinguished $K$-linear continuous left inverse ${\psi_L^D}$ which is defined to be the dual of the map
\begin{align*}
  C^{an}(o_L,K) & \longrightarrow C^{an}(o_L,K)  \\
  f & \longmapsto (\pi_L)_!(f)(a) :=
  \begin{cases}
  f(\pi_L^{-1}a) & \text{if $a \in \pi_L o_L$}, \\
  0 & \text{otherwise},
  \end{cases}
\end{align*}
and then, via the Fourier transform, induces  an operator $\psi_L^\mathfrak{X}$ on $\mathcal{O}_K(\mathfrak{X})$. One checks that for Dirac distributions we have
\begin{equation}\label{f:psi-Dirac}
  {\psi_L^D}(\delta_a) =
    \begin{cases}
  \delta_{\pi_L^{-1}a} & \text{if $a \in \pi_L o_L$}, \\
  0 & \text{otherwise}.
  \end{cases}
\end{equation}
Together with Lemma \ref{phi-free} this implies the following.

\begin{lemma}\label{psi-ker}
If $R_0 \subseteq o_L$ is a set of representatives for the nonzero cosets in $o_L/\pi_L o_L$ then
\begin{equation*}
  \ker({\psi_L^\mathfrak{X}}) = \oplus_{a \in R_0} \ev_a \cdot \varphi_L(\mathcal{O}_K(\mathfrak{X})) \ .
\end{equation*}
\end{lemma}

We also recall the resulting projection formula
\begin{equation*}
  {\psi_L^\mathfrak{X}}(\varphi_L(F_1) F_2) = F_1 {\psi_L^\mathfrak{X}}(F_2) \qquad\text{for any $F_1, F_2 \in \mathcal{O}_K(\mathfrak{X})$}.
\end{equation*}

Sometimes it will be useful to view $\psi_L^\mathfrak{X}$ as a normalized trace operator. Since $\mathcal{O}_K(\mathfrak{X})$ is a free module over $\varphi_L(\mathcal{O}_K(\mathfrak{X}))$ of rank $q$ we have the corresponding trace map
\begin{equation*}
  trace_{\mathcal{O}_K(\mathfrak{X}) / \varphi_L(\mathcal{O}_K(\mathfrak{X}))} : \mathcal{O}_K(\mathfrak{X}) \longrightarrow \varphi_L(\mathcal{O}_K(\mathfrak{X})) \ .
\end{equation*}

\begin{remark}\label{rem:trace}
   $\psi_L^\mathfrak{X} = \frac{1}{q} \varphi_L^{-1} \circ trace_{\mathcal{O}_K(\mathfrak{X}) / \varphi_L(\mathcal{O}_K(\mathfrak{X}))}$.
\end{remark}
\begin{proof}
Since the functions $\ev_a$ generate a dense subspace in $\mathcal{O}_K(\mathfrak{X})$ (\cite[Lem.\ 3.1]{ST1} the proof of which remains valid for general $K$ by \cite[Cor.\ 4.2.6 and Thm.\ 11.3.5]{PGS}) it suffices, by the continuity of all operators involved, to check the asserted equality on the functions $\ev_a$. As before we choose a set of representatives $R \subseteq o_L$ for the cosets $o_L/\pi_Lo_L$, so that the functions $\ev_c$, for $c \in R$, form a basis of $\mathcal{O}_K(\mathfrak{X})$ over $\varphi_L(\mathcal{O}_K(\mathfrak{X}))$. \textit{Case 1}: Let $a \in o_L^\times$. Then $\psi_L^\mathfrak{X}(\ev_a) = 0$ by \eqref{f:psi-Dirac}. On the other hand $\ev_a \cdot \ev_c = \ev_{a+c}  \in \ev_{c'} \cdot \varphi_L(\mathcal{O}_K(\mathfrak{X}))$ for some $c \neq c' \in R$. Hence the matrix of multiplication by $\ev_a$ w.r.t.\ to our choice of basis has only zero entries on the diagonal. This means that $trace_{\mathcal{O}_K(\mathfrak{X}) / \varphi_L(\mathcal{O}_K(\mathfrak{X}))}(\ev_a) = 0$. \textit{Case 2}: Let $a \in \pi_L o_L$. Then  $\psi_L^\mathfrak{X}(\ev_a) = \ev_{\pi_L^{-1} a}$. On the other hand the matrix of multiplication by $\ev_a$ now is the diagonal matrix with constant entry $\ev_a = \varphi_L(\ev_{\pi_L^{-1} a})$. We see that $\frac{1}{q} \varphi_L^{-1} (trace_{\mathcal{O}_K(\mathfrak{X}) / \varphi_L(\mathcal{O}_K(\mathfrak{X}))}(\ev_a)) = \frac{1}{q} \varphi_L^{-1}(q \varphi_L(\ev_{\pi_L^{-1} a})) = \ev_{\pi_L^{-1} a}$.
\end{proof}

In order to establish a formula for the composition $\varphi_L \circ {\psi_L^\mathfrak{X}}$ we let $\mathfrak{X}[\pi_L] := \ker (\mathfrak{X} \xrightarrow{\pi_L^*} \mathfrak{X})$. Then $\mathfrak{X}[\pi_L](\mathbb{C}_p)$ is the character group of the finite group $o_L/\pi_L o_L$. The points in $\mathfrak{X}[\pi_L](\mathbb{C}_p)$ are defined over some finite extension $K_1/K$. For any $\zeta \in \mathfrak{X}[\pi_L](\mathbb{C}_p)$ we have the continuous translation operator
\begin{align*}
  \mathcal{O}_{K_1}(\mathfrak{X}) & \longrightarrow \mathcal{O}_{K_1}(\mathfrak{X}) \\
  F & \longmapsto ({_\zeta}F)(\chi) := F(\chi\zeta) \ .
\end{align*}

\begin{proposition}\phantomsection\label{phi-psi}
\begin{itemize}
  \item[i.] For any $F \in \mathcal{O}_{K_1}(\mathfrak{X})$ we have
  \begin{equation*}
    [o_L : \pi_L o_L] \cdot \varphi_L \circ {\psi_L^\mathfrak{X}} (F) = \sum_{\zeta \in \mathfrak{X}[\pi_L](\mathbb{C}_p)} {_\zeta}F \ .
  \end{equation*}
  \item[ii.] $\varphi_L (\mathcal{O}_K(\mathfrak{X})) = \{ F \in \mathcal{O}_K(\mathfrak{X}) : {_\zeta}F = F \ \text{for any $\zeta \in \mathfrak{X}[\pi_L](\mathbb{C}_p)$}\}$.
\end{itemize}
\end{proposition}
\begin{proof}
i. Again it suffices to consider any $F = \ev_a$. We compute
\begin{align*}
  (\sum_\zeta {_\zeta}\ev_a)(\chi) & = \sum_\zeta \ev_a(\chi\zeta) = \chi(a) \sum_\zeta \zeta(a) \\
   & = \begin{cases}
           [o_L : \pi_L o_L] \cdot \chi(a) & \text{if $a \in \pi_L o_L$}, \\
            0  & \text{otherwise}
       \end{cases}   \\
   & = \begin{cases}
           [o_L : \pi_L o_L] \cdot \ev_a(\chi) & \text{if $a \in \pi_L o_L$}, \\
            0  & \text{otherwise}.
       \end{cases}
\end{align*}
On the other hand
\begin{equation*}
  \varphi_L ({\psi_L^\mathfrak{X}} (\ev_a)) = \varphi_L \Big(
  \begin{cases}
     \ev_{\pi_L^{-1}a} & \text{if $a \in \pi_L o_L$}, \\
     0   & \text{otherwise}
  \end{cases}
  \Big) =
  \begin{cases}
     \ev_a & \text{if $a \in \pi_L o_L$}, \\
     0   & \text{otherwise}.
  \end{cases}
\end{equation*}
ii. If ${_\zeta}F = F$ for any $\zeta \in \mathfrak{X}[\pi_L](\mathbb{C}_p)$ then $\varphi_L({\psi_L^\mathfrak{X}}(F)) = F$ by i. On the other hand
\begin{equation*}
  ({_\zeta}\varphi_L(F))(\chi) = \varphi_L(F)(\chi\zeta)) = F(\pi_L^*(\chi) \pi_L^*(\zeta)) = F(\pi_L^*(\chi)) = \varphi_L(F)(\chi) \ .
\end{equation*}
\end{proof}

We have observed in the above proof that the functions $\ev_a$, for $a \in o_L$, generate a dense subspace of $\mathcal{O}_K(\mathfrak{X})$. Considering the topological decomposition
\begin{align}\label{f:decomp}
  \mathcal{O}_K(\mathfrak{X}) & = \varphi_L(\mathcal{O}_K(\mathfrak{X})) \oplus \mathcal{O}_K(\mathfrak{X})^{{\psi_L^\mathfrak{X}} = 0} \\
  F & = \varphi_L({\psi_L^\mathfrak{X}}(F)) + (F - \varphi_L({\psi_L^\mathfrak{X}}(F)))   \nonumber
\end{align}
we see, using \eqref{f:psi-Dirac}, that the $\ev_a$ for $a \in \pi_L o_L$, resp.\ the $\ev_u$ for $u \in o_L^\times$, generate a dense subspace of $\varphi_L(\mathcal{O}_K(\mathfrak{X}))$, resp.\ of $\mathcal{O}_K(\mathfrak{X})^{{\psi_L^\mathfrak{X}} = 0}$. In view of Lemma \ref{module-struc} the obvious formula $u_*(\ev_a) = \ev_{ua}$ together with the fact, that the Dirac distributions $\delta_u$, for $u \in o_L^\times$, generate a dense subspace of $D(o_L^\times,K)$, then imply that the decomposition \eqref{f:decomp} is $D(o_L^\times,K)$-invariant.

\begin{lemma}\label{Mellin}
(Mellin transform) The natural inclusion $D(o_L^\times,K)\hookrightarrow D(o_L,K)$ combined with the Fourier isomorphism induces the map
\begin{align*}
  \mathfrak{M} : D(o_L^\times,K) & \xrightarrow{\; \cong \;} D(o_L,K)^{\psi_L^D=0}\cong \mathcal{O}_K(\mathfrak{X})^{{\psi_L^\mathfrak{X}} = 0} \\
  \lambda & \longmapsto \phantom{mmmmll} \lambda(\delta_1)\corresponds\lambda(\ev_1)
\end{align*}
which is a topological isomorphism of $D(o_L^\times,K)$-modules.
\end{lemma}
\begin{proof}
The disjoint decomposition into open sets $o_L = \pi_L o_L\, \cup \, o_L^\times$ induces the linear topological decomposition $D(o_L,K) = \varphi_L(D(o_L,K)) \oplus D(o_L^\times,K)$. The assertion follows by comparing this with the decomposition \eqref{f:decomp}.\footnote{The map $D(o_L^\times,K)\to D(G,K)$ sending $\lambda$ to $\lambda(\delta_1)$ is the inclusion map since $\delta_u(\delta_1) = \delta_u$.}
\end{proof}

\noindent{\bf The Robba ring}\\

We recall a few facts from \cite{BSX} about the analytic structure of the character variety $\mathfrak{X}$. As a general convention all \textbf{radii $r$} which will occur throughout the paper are \textbf{assumed to lie in  $(0,1) \cap p^\mathbb{Q}$}. Let $\mathbf{B}_1(r)$, resp.\ $\mathbf{B}(r)$, denote the $\mathbb{Q}_p$-affinoid disk of radius $r$ around $1$, resp.\ around $0$, and let $\mathbf{B}_1^-(r)$ be the open disk of radius $r$ around $1$. We put
\begin{equation*}
  \mathfrak{X}_0(r) := \mathbf{B}_1(r) \otimes_{\mathbb{Z}_p} \Hom_{\mathbb{Z}_p}(o_L,\mathbb{Z}_p) \quad\text{and}\quad    \mathfrak{X}(r) := \mathfrak{X} \cap  \mathfrak{X}_0(r)_{/L} \ .
\end{equation*}
These are affinoid subgroups of $\mathfrak{X}_0$ and $\mathfrak{X}$, respectively, which are respected by the action of the monoid $o_L \setminus \{0\}$. Since $\mathfrak{X}(r) \hookrightarrow \mathfrak{X}_0(r)_{/L}$ is a closed immersion of affinoid varieties the restriction map between the affinoid algebras $\mathcal{O}_K(\mathfrak{X}_0(r)) \twoheadrightarrow \mathcal{O}_K(\mathfrak{X}(r))$ is a strict surjection of Banach algebras. The families $\{\mathfrak{X}_0(r)\}_r$, resp.\ $\{\mathfrak{X}(r)\}_r$, form an increasing admissible covering of $\mathfrak{X}_0$, resp.\ $\mathfrak{X}$, which exhibits the latter as a quasi-Stein space. Hence $\mathcal{O}_K(\mathfrak{X}_0(r))$, resp.\ $\mathcal{O}_K(\mathfrak{X}(r))$, is the completion of $\mathcal{O}_K(\mathfrak{X}_0)$, resp.\ $\mathcal{O}_K(\mathfrak{X})$, in the supremum norm $|\ |_{\mathfrak{X}_0(r)}$, resp.\ $|\ |_{\mathfrak{X}(r)}$.

The structure of the affinoid variety $\mathfrak{X}(r_0)$ is rather simple for any radius $r_0 < p^{-\frac{d}{p-1}}$. Then (\cite[Lem.\ 1.16]{BSX}) the map
\begin{align}\label{f:small-disk}
    \mathbf{B}(r_0)_{/L} & \xrightarrow{\;\cong\;} \mathfrak{X}(r_0) \\
    y & \longmapsto \chi_y(a) := \exp(ay)   \nonumber
\end{align}
is an isomorphism of $L$-affinoid groups. Taking, somewhat unconventionally, $\exp - 1$ as coordinate function on $\mathbf{B}(r_0)$ we may view $\mathcal{O}_K(\mathbf{B}(r_0))$ as the Banach algebra of all power series $f = \sum_{i \geq 0} c_i (\exp - 1)^i$ such that $c_i \in K$ and $\lim_{i \rightarrow \infty} |c_i|r_0^i = 0$; the norm is $|f|_{\mathbf{B}(r_0)} := \max_i |c_i|r_0^i$. Since $\exp - 1$ corresponds under the above isomorphism to the function $\ev_1 - 1$ on $\mathfrak{X}(r_0)$ we deduce that
\begin{equation}\label{f:r0-expansion}
  \mathcal{O}_K(\mathfrak{X}(r_0)) = \{f = \sum_{i \geq 0} c_i (\ev_1 - 1)^i : c_i \in K \ \text{and}\ \lim_{i \rightarrow \infty} |c_i|r_0^i = 0\}
\end{equation}
is a Banach algebra with the supremum norm $|f|_{\mathfrak{X}(r_0)} = \max_i |c_i|r_0^i$.

Next we need to explain the admissible open subdomains $\mathfrak{X}_I$ of $\mathfrak{X}$, where the $I \subseteq (0,1)$ are certain intervals (cf.\ \cite[\S2.1]{BSX}). First of all we have the admissible open subdomains
\begin{equation*}
  \mathfrak{X}_{(r,1)} := \mathfrak{X} \setminus \mathfrak{X}(r) \ .
\end{equation*}
To introduce the relevant affinoid subdomains we also need the open disk $\mathbf{B}_1^-(r)$ of radius $r$ around $1$. This allows us to first define the admissible open subdomains $\mathfrak{X}_0^-(r) := (\mathbf{B}_1^-(r) \otimes_{\mathbb{Z}_p} \Hom_{\mathbb{Z}_p}(o_L,\mathbb{Z}_p))_{/L}$ and $\mathfrak{X}^-(r) := \mathfrak{X} \cap \mathfrak{X}_0^-(r)$ of $\mathfrak{X}_0$ and $\mathfrak{X}$, respectively. For $r \leq s$ we then have the admissible open subdomains
\begin{equation*}
  \mathfrak{X}_0[r,s] := \mathfrak{X}_0(s) \setminus \mathfrak{X}_0^-(r) \subseteq \mathfrak{X}_0 \quad\text{and}\quad  \mathfrak{X}_{[r,s]} := \mathfrak{X}(s) \setminus \mathfrak{X}^-(r) = \mathfrak{X} \cap \mathfrak{X}_0[r,s] \subseteq \mathfrak{X} \ .
\end{equation*}
We recall that the $\mathfrak{X}_{[r,s]}$ are actually affinoid varieties. There are the obvious inclusions $\mathfrak{X}_{[r,s]} \subseteq \mathfrak{X}(s)$ and $\mathfrak{X}_{[r,s]} \subseteq \mathfrak{X}_{(r',1)}$ provided $r' < r$. Moreover, $\mathfrak{X}_{(r',1)}$ is the increasing admissible union of the $\mathfrak{X}_{[r,s]}$ for $r' < r \leq s < 1$. Hence
\begin{equation*}
  \mathcal{O}_K(\mathfrak{X}_{(r',1)}) = \varprojlim_{r' < r \leq s < 1} \mathcal{O}_K(\mathfrak{X}_{[r,s]}) \ ,
\end{equation*}
which exhibits the Fr\'echet algebra structure of the left side.

We point out that these subdomains $\mathfrak{X}_I$ all are invariant under $o_L^\times$. Their behaviour with respect to $\pi_L^*$ is more complicated. We recall from \cite[Lem.\ 2.11]{BSX}  that, for any radius $p^{-\frac{dp}{p-1}} \leq r < 1$ we have
\begin{equation}\label{f:pi-1}
  (\pi_L^*)^{-1} (\mathfrak{X}_{(r,1)}) \subseteq \mathfrak{X}_{(r^{1/p},1)} \ .
\end{equation}

It is technically necessary in the following to sometimes only work with a smaller set of radii. We put
\begin{equation*}
  S_0 := [p^{-\frac{d}{e} - \frac{d}{p-1}},p^{-\frac{d}{p-1}}) \cap p^{\mathbb{Q}} \subseteq [p^{-\frac{dp}{p-1}},p^{-\frac{d}{p-1}}) \ ,
\end{equation*}
$S_n := S_0^{\frac{1}{p^n}}$ for $n \geq 1$, and $S_\infty := \bigcup_{n \geq 1} S_n$. Note that the sets $S_n$ are pairwise disjoint. The point is (\cite[Prop.\ 1.20]{BSX}) that for $s \in S_\infty$ we know that $\mathfrak{X}(s)$ becomes isomorphic to a closed disk over $\mathbb{C}_p$. Let $s_n$ for $n \geq 0$, denote the left boundary point of the set $S_n$. Then we have the following result (\cite[Prop.\ 2.1]{BSX}).

\begin{proposition}\label{quasi-Stein}
  For any $n \geq 0$ the rigid variety $\mathfrak{X}_{(s_n,1)}$ is quasi-Stein with respect to the admissible covering $\{\mathfrak{X}_{[r,s]}\}$ where $s_n < r \leq s < 1$, $r \in S_n$, and $s \in \bigcup_{m \geq n} S_m$. In particular, the affinoid algebra $\mathcal{O}_K(\mathfrak{X}_{[r,s]})$ is the completion of $\mathcal{O}_K(\mathfrak{X}_{(s_n,1)})$ with respect to the supremum norm $|\ |_{\mathfrak{X}_{[r,s]}}$.
\end{proposition}

Obviously, with $\mathfrak{X}$ each $\mathfrak{X}_{(s_n,1)}$ is one dimensional and smooth. But, in order to be able to apply later on Serre duality to the spaces $\mathfrak{X}_{(s_n,1)}$, we need to show that they are actually Stein spaces. This means that we have to check that the admissible covering in Prop.\ \ref{quasi-Stein} has the property that $\mathfrak{X}_{[r',s']}$ is relatively compact in $\mathfrak{X}_{[r,s]}$ over $L$ (\cite[\S 9.6.2]{BGR}) for any $r < r' \leq s' < s$. We simply write $U \Subset X$ for an affinoid subdomain $U$ being relatively compact over $L$ in an $L$-affinoid variety $X$.

\begin{lemma}\label{lem:rel-compact}
  Let $U \subseteq X \subseteq X'$ be affinoid subdomains of the affinoid variety $X'$; we then have:
\begin{itemize}
  \item[i.] If $U \Subset X$ then $U \Subset X'$;
  \item[ii.] suppose that $U = U_1 \cup \ldots \cup U_m$ is an affinoid covering; if $U_i \Subset X$ for any $1 \leq i \leq m$ then $U \Subset X$.
\end{itemize}
\end{lemma}
\begin{proof}
Let $A \rightarrow B$ be the homomorphism of affinoid algebras which induces the inclusion $U = \Sp(B) \subseteq X = \Sp(A)$. It is not difficult to see that the property $U \Subset X$ is equivalent to the homomorphism $A \rightarrow B$ being inner w.r.t. $L$ in the sense of \cite[Def.\ 2.5.1]{Ber}. Therefore i., resp.\ ii., is a special case of Cor.\ 2.5.5, resp.\ Lemma 2.5.10, in \cite{Ber}.
\end{proof}

\begin{proposition}\label{Stein}
   $\mathfrak{X}_{(s_n,1)}$, for any $n \geq 0$, is a Stein space.
\end{proposition}
\begin{proof}
   Since, by Prop.\ \ref{quasi-Stein}, we already know that the $\mathfrak{X}_{(s_n,1)}$ are quasi-Stein. Hence it remains to show that $\mathfrak{X}_{[r',s']} \Subset \mathfrak{X}_{[r,s]}$ for any $r < r' \leq s' < s$. Looking first at $\mathfrak{X}_0$, let $\mathfrak{B}_1[r,s] \subseteq \mathfrak{B}_1$ denote the affinoid annulus of inner radius $r$ and outer radius $s$. Fixing coordinate functions $z_1, \ldots, z_d$ on $\mathfrak{X}_0$ we have the admissible open covering
\begin{equation*}
  \mathfrak{X}_0[r,s] = \bigcup_{i = 1}^d \mathfrak{X}_0^{(i)}[r,s]  \quad\text{with}\quad  \mathfrak{X}_0^{(i)}[r,s] := \{x \in \mathfrak{X}_0(s) : |z_i(x)| \geq r \}.
\end{equation*}
The affinoid varieties of this covering have the direct product structure
\begin{equation*}
  \mathfrak{X}_0^{(i)}[r,s] = \mathfrak{B}_1(s) \times \ldots \mathfrak{B}_1(s) \times \mathfrak{B}_1[r,s] \times \mathfrak{B}_1(s) \times \ldots \mathfrak{B}_1(s)
\end{equation*}
with the annulus being the \textit{i}th factor. It immediately follows that $\mathfrak{X}_0^{(i)}[r',s'] \Subset \mathfrak{X}_0^{(i)}[r,s]$ (\cite[Lem.\ 9.6.2.1]{BGR}). Since relative compactness is preserved by passing to closed subvarieties we deduce that $\mathfrak{X} \cap \mathfrak{X}_0^{(i)}[r',s'] \Subset \mathfrak{X} \cap \mathfrak{X}_0^{(i)}[r,s]$ for any $1 \leq i \leq d$. Applying now Lemma \ref{lem:rel-compact} we conclude first that $\mathfrak{X} \cap \mathfrak{X}_0^{(i)}[r',s'] \Subset \mathfrak{X}_0[r,s]$ and then that $\mathfrak{X}_{[r',s']} \Subset \mathfrak{X}_{[r,s]}$.
\end{proof}

We finally recall that the Robba ring of $\mathfrak{X}$ over $K$ is defined as the locally convex inductive limit $\varinjlim_{\mathfrak{Y}} \mathcal{O}_K(\mathfrak{X} \setminus \mathfrak{Y})$ where $\mathfrak{Y}$ runs over all affinoid subdomains of $\mathfrak{X}$. Since any such $\mathfrak{Y}$ is contained in some $\mathfrak{X}(r)$ we have
\begin{equation*}
  \mathcal{R}_K(\mathfrak{X}) = \varinjlim_{n \geq 0} \mathcal{O}_K(\mathfrak{X}_{(s_n,1)}) \ ,
\end{equation*}
and we view $\mathcal{R}_K(\mathfrak{X})$ as the locally convex inductive limit of the Fr\'echet algebras $\mathcal{O}_K(\mathfrak{X}_{(s_n,1)})$. By \cite{BSX} Prop.\ 1.20 the system $\mathfrak{X}_{(s_n,1)/\mathbb{C}_p}$ is isomorphic to a decreasing system of one dimensional annuli. This implies:
\begin{itemize}
   \item[--] $\mathcal{R}_K(\mathfrak{X})$ is the increasing union of the rings $\mathcal{O}_K(\mathfrak{X}_{(s_n,1)})$ and contains $\mathcal{O}_K(\mathfrak{X})$;\Footnote{
In functional analysis the Hahn-Banach theorem is crucial for many applications. Therefore many references assume that the base field $K$ is spherically complete. Unfortunately, in this article we frequently have to deal with bigger extensions of $\mathbb{Q}_p$ not satisfying this property. From \cite{PGS} we recall that a locally convex vector space $V$ is said to be of {\it countable type}, if for every continuous
seminorm $p$ on $V$ its completion $V_p$ at $p$ has a dense subspace of countable algebraic dimension. They are  stable under forming subspaces, linear images,
projective limits, and countable inductive limits, cf.\    theorem 4.2.13 in (loc.\ cit.). By  corollary 4.2.6 in (loc.\ cit.) for such vector spaces the Hahn-Banach theorem holds, too. By \cite[Prop.\ 5.4.3]{Th} the Robba ring over any intermediate field $\mathbb{Q}_p\subseteq K\subseteq \mathbb{C}_p$ is of countable type as $K$-vector space.
}
   \item[--] each $\mathcal{O}_K(\mathfrak{X}_{(s_n,1)})$ as well as $\mathcal{R}_K(\mathfrak{X})$ are integral domains.
\end{itemize}
The action of the monoid $o_L \setminus \{0\}$ on $\mathcal{O}_K(\mathfrak{X})$ extends naturally to a continuous action on $\mathcal{R}_K(\mathfrak{X})$ (\cite[Lem.\ 2.12]{BSX}). In fact, this action extends further uniquely to a separately continuous action of $D(o_L^\times,K)$-action on $\mathcal{R}_K(\mathfrak{X})$. This is a special case of the later Prop.\ \ref{prop:distributionaction} which implies that we will have such an action on any $L$-analytic $(\varphi_L,\Gamma_L)$-module over $\mathcal{R}_K(\mathfrak{X})$. Via the isomorphism $\chi_{LT} : \Gamma_L \xrightarrow{\cong} o_L^\times$ we later on will view this as a $D(\Gamma_L,K)$-action.

In order to extend the $\psi$-operator to the Robba ring we need the following fact.

\begin{lemma}\label{lem:piL-finite}
   The morphism $\pi_L^* : \mathfrak{X} \rightarrow \mathfrak{X}$ is finite, faithfully flat, and etale.
\end{lemma}
\begin{proof}
The character variety $\mathfrak{X}'$ of the subgroup $\pi_L o_L \subseteq o_L$ is isomorphic to $\mathfrak{X}$ via
\begin{align*}
  \mathfrak{X} & \xrightarrow{\;\cong\;} \mathfrak{X}' \\
  \chi & \longmapsto \chi'(\pi_L a) := \chi(a) \ .
\end{align*}
We have the commutative diagram
\begin{equation*}
  \xymatrix{
  \mathfrak{X} \ar[d]_{\pi_L^*} \ar[dr]^{\chi \mapsto \chi | \pi_L o_L}        \\
  \mathfrak{X} \ar[r]_-{\chi \mapsto \chi'}^-{\cong}  & \mathfrak{X}'   \ .           }
\end{equation*}
The oblique arrow is finite and faithfully flat by the proof of \cite[Prop.\ 6.4.5]{Eme}. For its etaleness it remains to check that all its fibers are unramified. This can be done after base change to $\mathbb{C}_p$. Then, since this arrow is a homomorophism of rigid groups, all fibers are isomorphic. But the fiber in the trivial character of $\pi_L o_L$ is isomorphic to $\Sp(\mathbb{C}_p[o_L/\pi_L o_L]) \cong \Sp(\mathbb{C}_p^q)$. It follows that $\pi_L^*$ has these properties as well.
\end{proof}

Since the subsequent reasoning will be needed again in the next section in an analogous situation we proceed in an axiomatic way. \textbf{Suppose} that:
\begin{itemize}
  \item[-] $\rho : \mathfrak{Y} \rightarrow \mathfrak{Z}$ is a finite and faithfully flat morphism of quasi-Stein spaces over $K$. In particular, the induced map $\rho^* : \mathcal{O}_K(\mathfrak{Z}) \rightarrow \mathcal{O}_K(\mathfrak{Y})$ is injective. Moreover, the finiteness of $\rho$ implies that the preimage under $\rho$ of any affinoid subdomain in $\mathfrak{Z}$ is an affinoid subdomain in $\mathfrak{Y}$ (\cite[Prop.\ 9.4.4.1(i)]{BGR}) and hence that $\rho^*$ is continuous.
  \item[-] $\mathcal{O}_K(\mathfrak{Y})$ is finitely generated free as a $\rho^*(\mathcal{O}_K(\mathfrak{Z}))$-module. Fix a corresponding basis $f_1, \ldots, f_h \in \mathcal{O}_K(\mathfrak{Y})$.
\end{itemize}

\begin{proposition}\label{crossed-product}
For any admissible open subset $\mathfrak{U} \subseteq \mathfrak{Z}$ we have
\begin{equation*}
  \mathcal{O}_K(\rho^{-1}(\mathfrak{U})) = \mathcal{O}_K(\mathfrak{Y}) \otimes_{\mathcal{O}_K(\mathfrak{Z})} \mathcal{O}_K(\mathfrak{U})
\end{equation*}
is free with basis $f_1, \ldots, f_h$ over $\mathcal{O}_K(\mathfrak{U})$.
\end{proposition}
\begin{proof}
Since $\rho$ is finite, $\rho_* \mathcal{O}_{\mathfrak{Y}}$ is a coherent $\mathcal{O}_{\mathfrak{Z}}$-module by \cite[Prop.\ 9.4.4.1(ii)]{BGR}. Gruson's theorem (cf.\ \cite[Prop.\ 1.13]{BSX}) then tells us that $\rho_* \mathcal{O}_{\mathfrak{Y}}$ is, in fact, a free $\mathcal{O}_{\mathfrak{Z}}$-module with basis $f_1, \ldots, f_h$.
\end{proof}

We observe that the definition of the Robba ring $\mathcal{R}_K(\mathfrak{X})$ above was completely formal and works precisely the same way for any quasi-Stein space. Hence we have available the Robba rings $\mathcal{R}_K(\mathfrak{Y})$ and $\mathcal{R}_K(\mathfrak{Z})$. Since the morphism $\rho : \mathfrak{Y} \rightarrow \mathfrak{Z}$ is finite the preimage under $\rho$  of any affinoid subdomain in $\mathfrak{Z}$ is an affinoid subdomain in $\mathfrak{Y}$ (\cite[Prop.\ 9.4.4.1(i)]{BGR}).
We note again that the preimage under $\rho$  of any affinoid subdomain in $\mathfrak{Z}$ is an affinoid subdomain in $\mathfrak{Y}$.  The injective map $\rho^* : \mathcal{O}_K(\mathfrak{Z}) \subseteq \mathcal{O}_K(\mathfrak{Y})$ therefore extends to a natural homomorphism of rings
\begin{equation}\label{f:rho-Robba}
  \rho^* : \mathcal{R}_K(\mathfrak{Z}) \longrightarrow \mathcal{R}_K(\mathfrak{Y}) \ .
\end{equation}

\begin{remark}\label{inj}
The homomorphism \eqref{f:rho-Robba} is injective.
\end{remark}
\begin{proof}
We fix an admissible covering $\mathfrak{Z} = \bigcup_{j \geq 1} \mathfrak{U}_j$ by an increasing sequence of affinoid subdomains $\mathfrak{U}_j \subseteq \mathfrak{Z}$. As $\rho$ is a finite map, $\mathfrak{Y} = \bigcup_{j \geq 1} \rho^{-1}(\mathfrak{U}_j)$ again is an admissible covering by affinoid  subdomains. It follows that $\mathcal{R}_K(\mathfrak{Y}) = \varinjlim_j \mathcal{O}_K(\mathfrak{Y} \setminus \rho^{-1}(\mathfrak{U}_j))$, and therefore it suffices to show the injectivity of the maps $\rho^* : \mathcal{O}_K(\mathfrak{Z} \setminus \mathfrak{U}_j) \rightarrow \mathcal{O}_K(\mathfrak{Y} \setminus \rho^{-1}(\mathfrak{U}_j))$. But this is clear since the map $\rho : \mathfrak{Y} \setminus \rho^{-1}(\mathfrak{U}_j) \rightarrow \mathfrak{Z} \setminus \mathfrak{U}_j$ is faithfully flat.
\end{proof}

\begin{corollary}\label{crossed-product-R}
$\mathcal{R}_K(\mathfrak{Y}) = \mathcal{O}_K(\mathfrak{Y}) \otimes_{\mathcal{O}_K(\mathfrak{Z})}  \mathcal{R}_K(\mathfrak{Z})$ is free over $\rho^*(\mathcal{R}_K(\mathfrak{Z}))$ with the basis $f_1, \ldots, f_h$. In fact, the map
\begin{align*}
  \mathcal{R}_K(\mathfrak{Z})^h & \xrightarrow{\;\cong\;} \mathcal{R}_K(\mathfrak{Y}) \\
                      (z_1,\ldots,z_h) & \longmapsto \sum_{i=1}^h \rho^*(z_i)f_i
\end{align*}
is a homeomorphism.
\end{corollary}
\begin{proof}
By passing to locally convex limits this follows from Prop.\ \ref{crossed-product} which says that the map
\begin{align*}
  \mathcal{O}_K(\mathfrak{U})^h & \xrightarrow{\;\cong\;} \mathcal{O}_K(\rho^{-1}(\mathfrak{U})) \\
                      (z_1,\ldots,z_h) & \longmapsto \sum_{i=1}^h \rho^*(z_i)f_i
\end{align*}
is a continuous bijection between Fr\'echet spaces and hence a homeomorphism by the open mapping theorem.
\end{proof}

By the Lemmas \ref{phi-free} and \ref{lem:piL-finite} the above applies to the morphism $\pi_L^* : \mathfrak{X} \rightarrow \mathfrak{X}$ and we obtain the following result.

\begin{proposition}\label{phi-free-Robba}
   Let $R \subseteq o_L$ be a set of representatives for the cosets in $o_L/\pi_L o_L$. Then the Robba ring $\mathcal{R}_K(\mathfrak{X})$ is a free module over $\varphi_L(\mathcal{R}_K(\mathfrak{X}))$ with basis $\{\ev_a\}_{a \in R}$.
\end{proposition}

In particular we have the trace map
\begin{equation*}
  trace_{\mathcal{R}_K(\mathfrak{X}) / \varphi_L(\mathcal{R}_K(\mathfrak{X}))} : \mathcal{R}_K(\mathfrak{X}) \longrightarrow \varphi_L(\mathcal{R}_K(\mathfrak{X}))
\end{equation*}
and therefore may introduce the operator
\begin{equation*}
  \psi_L^\mathfrak{X} := \frac{1}{q} \varphi_L^{-1} \circ trace_{\mathcal{R}_K(\mathfrak{X}) / \varphi_L(\mathcal{R}_K(\mathfrak{X}))} : \mathcal{R}_K(\mathfrak{X}) \longrightarrow \mathcal{R}_K(\mathfrak{X}) \ .
\end{equation*}
Because of Remark \ref{rem:trace} it extends the operator $\psi_L^\mathfrak{X}$ on $\mathcal{O}_K(\mathfrak{X})$, which justifies denoting it by the same symbol. By construction it is a left inverse of $\varphi_L$ and satisfies the projection formula.
 Furthermore, as a consequence of Cor.\ \ref{crossed-product-R}, $\psi_L^\mathfrak{X}$ is continuous.

\subsubsection{The multiplicative character variety and its Robba ring}\label{sec:mult-Robba}

In this section we consider the multiplicative group $o_L^\times$ as a locally $L$-analytic group. We introduce the open subgroups $U_n := 1 + \pi_L^n o_L$ for $n \geq 1$. Correspondingly we have the inclusion of distribution algebras $D(U_{n+1},K) \subseteq D(U_n,K) \subseteq D(o_L^\times,K)$. There is an integer $n_0 \geq 1$ such that, for any $n \geq n_0$, the logarithm series induces an isomorphism of locally $L$-analytic groups $\log : U_n \xrightarrow{\cong} \pi_L^n o_L$. We then introduce the isomorphisms $\ell_n := \pi_L^{-n} \log : U_n \xrightarrow{\cong} o_L$ together with the algebra isomorphisms
\begin{equation}\label{f:ell}
  \ell_{n*} : D(U_n,K) \xrightarrow{\cong} D(o_L,K) \cong \mathcal{O}_K(\mathfrak{X})
\end{equation}
which they induce.

As for $o_L$ in the previous section we have rigid analytic varieties (over $L$) of locally $L$-analytic characters $\mathfrak{X}^\times$ for $o_L^\times$ and $\mathfrak{X}^\times_n$ for $U_n$ as well (cf.\ \cite[Thm.\ 2.3, Lemma 2.4, Cor.\ 3.7]{ST2}  and \cite[Prop.s 6.4.5 and 6.4.6]{Eme}):
\begin{itemize}
  \item[--] $\ell_n^* : \mathfrak{X} \xrightarrow{\cong} \mathfrak{X}^\times_n$ is, for $n \geq n_0$, an isomorphism of group varieties. \Footnote{Note also that the inclusion $o_L^\times \hookrightarrow L^\times$ is a distinguished point in $\mathfrak{X}^\times(L)$.}
  \item[--] The restriction map $\rho_n : \mathfrak{X}^\times  \longrightarrow \mathfrak{X}^\times_n$ is a finite faithfully flat covering (\cite{Eme} proof of Prop.\ 6.4.5).
\Footnote{The following fact is a consequence. For the convenience of the reader we give an elementary argument.

\begin{remark}
The restriction map $\mathfrak{X}^\times(\mathbb{C}_p) \longrightarrow \mathfrak{X}^\times_n(\mathbb{C}_p)$ is surjective.
\end{remark}
\begin{proof}
Since $\mathbb{C}_p^\times$ is divisible and hence injective as an abelian group, any abstract homomorphism $U_n \rightarrow \mathbb{C}_p^\times$ extends to a homomorphism $o_L^\times \rightarrow \mathbb{C}_p^\times$. But the continuity as well as the local $L$-analyticity of a homomorphism $\chi : o_L^\times \rightarrow \mathbb{C}_p^\times$ can be tested on its restriction $\chi | U_n$. Observe for this that the diagram
\begin{equation*}
  \xymatrix{
    U_n \ar[d]_{u \cdot} \ar[r]^{\chi} & \mathbb{C}_p^\times \ar[d]^{\cdot \chi(u)} \\
    u U_n \ar[r]^{\chi} & \mathbb{C}_p^\times   }
\end{equation*}
is commutative for any $u \in o_L^\times$.
\end{proof}
}
  \item[--] $\mathfrak{X}^\times$ and $\mathfrak{X}^\times_n$ are one dimensional Stein spaces. (As group varieties they are separated and equidimensional.)
  \item[--] For $n \geq n_0$ the variety $\mathfrak{X}^\times_n$ is smooth and $\mathcal{O}_L(\mathfrak{X}^\times_n)$ is an integral domain.
  \item[--] The Fourier transforms
\begin{equation*}
  D(o_L^\times,K) \xrightarrow{\cong} \mathcal{O}_K(\mathfrak{X}^\times)  \qquad\text{and}\qquad   D(U_n,K) \xrightarrow{\cong} \mathcal{O}_K(\mathfrak{X}^\times_n)
\end{equation*}
      sending a distribution $\mu$ to the function $F_\mu(\chi) := \mu(\chi)$ are isomorphisms of Fr\'echet algebras.
\end{itemize}

As a consequence of the properties of the morphism $\rho := \rho_n: \mathfrak{X}^\times  \rightarrow \mathfrak{X}^\times_n$ the homomorphism $\rho^* : \mathcal{O}_K(\mathfrak{X}^\times_n) \rightarrow \mathcal{O}_K(\mathfrak{X}^\times)$ is injective and extends to an injective homomorphism $\rho^* : \mathcal{R}_K(\mathfrak{X}^\times_n) \rightarrow \mathcal{R}_K(\mathfrak{X}^\times)$ (cf.\ Remark \ref{inj}).

\begin{lemma}\phantomsection\label{lem:crossed-product-R}
\begin{itemize}
  \item[i.] $\mathcal{O}_K(\mathfrak{X}^\times) = \mathbb{Z}[o_L^\times] \otimes_{\mathbb{Z}[U_n]}  \mathcal{O}_K(\mathfrak{X}^\times_n)$.
  \item[ii.] $\mathcal{R}_K(\mathfrak{X}^\times) = \mathcal{O}_K(\mathfrak{X}^\times) \otimes_{\mathcal{O}_K(\mathfrak{X}_n^\times)} \mathcal{R}_K(\mathfrak{X}_n^\times) = \mathbb{Z}[o_L^\times] \otimes_{\mathbb{Z}[U_n]}  \mathcal{R}_K(\mathfrak{X}^\times_n)$.
\end{itemize}
\end{lemma}
\begin{proof}
i. Let $u_1, \ldots, u_h \in o_L^\times$ be a set of representatives for the cosets of $U_n$ in $o_L^\times$. We then have the decomposition into open subsets $o_L^\times = u_1 U_n \cup \ldots \cup u_h U_n$. It follows that
\begin{equation*}
  D(o_L^\times,K) = \delta_{u_1} D(U_n,K) \oplus \ldots \oplus \delta_{u_h} D(U_n,K) = \mathbb{Z}[o_L^\times] \otimes_{\mathbb{Z}[U_n]} D(U_n,K)
\end{equation*}
is, in particular, a free $D(U_n,K)$-module of rank $h = [o_L^\times :U_n]$. Using the Fourier isomorphism we obtain that $\mathcal{O}_K(\mathfrak{X}^\times)$ is a free $\mathcal{O}_K(\mathfrak{X}^\times_n)$-module over the basis $\ev_{u_1}, \ldots, \ev_{u_h}$.

ii. Because of i. the assumptions before Prop.\ \ref{crossed-product} are satisfied and the present assertion is a special case of Cor.\ \ref{crossed-product-R}.
\end{proof}

\begin{lemma}\label{etale}
  The morphism $\rho$ is etale.
\end{lemma}
\begin{proof}
This is the same argument as in the proof of Lemma \ref{lem:piL-finite}.
\end{proof}

\begin{corollary}\label{cor:etale1}
  $\mathfrak{X}^\times$ is smooth.
\end{corollary}
\begin{proof}
This follows from the lemma since $\mathfrak{X}_n^\times$ is smooth for $n \geq n_0$.
\end{proof}

\begin{remark}\label{rem:crossed-product}
If $n \geq m$ then all the above assertions hold analogously for the finite morphism $\rho_{m,n} : \mathfrak{X}^\times_m \longrightarrow \mathfrak{X}^\times_n$. In particular, all $\mathfrak{X}^\times_n$ are smooth.
\end{remark}

Suppose that $n \geq n_0$. Then, due to the isomorphism $\ell_n^* : \mathfrak{X} \xrightarrow{\cong} \mathfrak{X}^\times_n$, everything which was defined for and recalled about $\mathfrak{X}$ in section \ref{sec:add-Robba} holds correspondingly for $\mathfrak{X}^\times_n$. In particular we have the admissible open subdomains $\mathfrak{X}^\times_n(r)$, $\mathfrak{X}^\times_{n,(r,1)}$, and $\mathfrak{X}^\times_{n,[r,s]}$. For $n \geq m \geq n_0$ we have the commutative diagram
\begin{equation}\label{diag:m,n}
  \xymatrix{
    \mathfrak{X} \ar[d]_{(\pi_L^*)^{n-m}} \ar[r]^-{\ell_m^*}_-{\cong} & \mathfrak{X}^\times_m \ar[d]^{\rho_{m,n}} \\
    \mathfrak{X} \ar[r]^{\ell_n^*}_-{\cong} & \mathfrak{X}^\times_n.   }
\end{equation}

\begin{lemma}\label{p-map}
Let $n \geq n_0$ and $m \geq 0$; for any $p^{-\frac{dp}{p-1}} \leq r < 1$ the map $\rho_{n,n+me}^* : \mathcal{O}_K(\mathfrak{X}^\times_{n+me}) \rightarrow \mathcal{O}_K(\mathfrak{X}^\times_n)$ extends to an isometric homomorphism of Banach algebras
\begin{equation*}
  (\mathcal{O}_K(\mathfrak{X}^\times_{n+me}(r)),|\ |_{\mathfrak{X}^\times_{n+me}(r)}) \longrightarrow        (\mathcal{O}_K(\mathfrak{X}^\times_n(r^{\frac{1}{p^m}})),|\ |_{\mathfrak{X}^\times_n(r^{\frac{1}{p^m}})}) \ .
\end{equation*}
\end{lemma}
\begin{proof}
By the above commutative diagram \eqref{diag:m,n} our assertion amounts to the statement that the map $(\pi_L^{me})^* : \mathfrak{X} \rightarrow \mathfrak{X}$ restricts to a surjection $\mathfrak{X}(r^{\frac{1}{p^m}}) \rightarrow \mathfrak{X}(r)$. In \cite[Lem.\ 3.3]{ST2} this is shown to be the case for the map $(p^m)^*$. But $p^m$ and $\pi_L^{me}$ differ by a unit $u \in o_L^\times$, and $u^*$ preserves $\mathfrak{X}(r)$.
\end{proof}


\subsubsection{Twisting}\label{subsec:twisting}

Consider any locally $L$-analytic group $G$ and fix a locally $L$-analytic character $\chi : G \rightarrow L^\times$. Then multiplication by $\chi$ is a $K$-linear topological isomorphism $C^{an}(G,K) \xrightarrow[\cong]{\chi \cdot} C^{an}(G,K)$. We denote the dual isomorphism by
\begin{equation*}
  Tw_\chi^D : D(G,K) \xrightarrow{\;\cong\;} D(G,K) \ ,
\end{equation*}
i.e., $Tw_\chi^D(\mu) = \mu(\chi -)$, and call it the twist by $\chi$. For Dirac distributions we obtain $Tw_\chi^D(\delta_g) = \chi(g) \delta_g$.

Suppose now that $G$ is one of the groups $o_L$ or $U_n \subseteq   o_L^\times$ of the previous subsections, and let $\mathfrak{X}_G$ denote its character variety. Then $\chi$ is an $L$-valued point $z_\chi \in \mathfrak{X}_G(L)$. Using the product structure of the variety $\mathfrak{X}_G$ we similarly have the twist operator
\begin{equation*}
  Tw_z^{\mathfrak{X}_G} : \mathcal{O}_K(\mathfrak{X}_G) \xrightarrow{\;\cong\;} \mathcal{O}_K(\mathfrak{X}_G) \ , \ f \longmapsto f(z -) \ .
\end{equation*}
As any rigid automorphism multiplication by a rational point respects the system of affinoid subdomains and hence the system of their complements. Hence $Tw_z^{\mathfrak{X}_G}$ extends straightforwardly to an automorphism $Tw_z^{\mathfrak{X}_G} : \mathcal{R}_K(\mathfrak{X}_G) \xrightarrow{\cong} \mathcal{R}_K(\mathfrak{X}_G)$. The following properties are straightforward to check:
\begin{itemize}
  \item[1.] Under the Fourier isomorphism $Tw_\chi^D$ and $Tw_{z_\chi}^{\mathfrak{X}_G}$ correspond to each other.
  \item[2.] $Tw_{z_1}^{\mathfrak{X}_G} \circ Tw_{z_2}^{\mathfrak{X}_G} = Tw_{z_1 \cdot z_2}^{\mathfrak{X}_G}$.
  \item[3.] If $\alpha : G_1 \xrightarrow{\cong} G_2$ is an isomorphism between two of our groups then, for any $z \in \mathfrak{X}_{G_2}(L)$, the twist operators $Tw_{\alpha^*(z)}^{\mathfrak{X}_{G_1}}$ and $Tw_z^{\mathfrak{X}_{G_2}}$ correspond to each other under the isomorphism $\alpha_* : \mathcal{R}_K(\mathfrak{X}_{G_1}) \xrightarrow{\cong} \mathcal{R}_K(\mathfrak{X}_{G_2})$.
\end{itemize}

\subsubsection{The LT-isomorphism, part 1}\label{subsec:LT}

We write $\mathbf{B}$ for the  open unit ball over $L$. The Lubin-Tate formal $o_L$-module gives $\mathbf{B}$ an $o_L$-action via $(a,z)\mapsto [a](z).$ If $\mathcal{O}_K(\mathbf{B})$ is the ring of power series in $Z$ with coefficients in $K$ which converge on $\mathbf{B}(\mathbb{C}_p)$, then the above $o_L$-action on $\mathbf{B}$ induces an action of the monoid $o_L\setminus \{0\}$ on $\mathcal{O}_K(\mathbf{B})$ by $(a,F)\mapsto F\circ [a].$ Similarly  as before we let $\varphi_L$ denote the action of $\pi_L$. Next we consider the continuous operator
\begin{align*}
  tr : \mathcal{O}_K(\mathbf{B}) & \longrightarrow \mathcal{O}_K(\mathbf{B})   \\
    f(z) & \longmapsto \sum_{y \in \ker([\pi_L])} f(y +_{LT} z) \ .
\end{align*}
Coleman has shown (cf.\ \cite[\S 2]{SV15}) that $tr(Z^i) \in \im(\varphi_L)$ for any $i \geq 0$. Hence, since $\varphi_L$ is a homeomorphism onto its image, we have $\im(tr) \subseteq \im(\varphi_L)$ and hence, since $\varphi_L$ is injective, we may introduce the $K$-linear operator
\begin{equation*}
  \psi_L : \mathcal{O}_K(\mathbf{B}) \longrightarrow \mathcal{O}_K(\mathbf{B})    \qquad\text{such that $\pi_L^{-1} tr = \varphi_L \circ \psi_L$}.
\end{equation*}
One easily checks that $\psi_L$ is equivariant for the $o_L^\times$-action and satisfies the projection formula $\psi_L(f_1 \varphi_L(f_2)) = \psi_L(f_1) f_2$ as well as $\psi_L \circ \varphi_L = \frac{q}{\pi_L}$.

Furthermore, we fix a generator  $\eta'$ of $\TLT'$ as $o_L$-module and denote by $\Omega=\Omega_{\eta'}$ the corresponding period.  In the following we \textbf{assume that $\Omega$ belongs to $K$.} From \cite[Thm.\ 3.6]{ST2} we recall the LT-isomorphism
\begin{align}\label{f:LTiso}
 \kappa^*:\mathcal{O}_K(\mathfrak{X})&\xrightarrow{\cong}\mathcal{O}_K(\mathbf{B})\\ F &\mapsto [z\mapsto F(\kappa_z)],\notag
\end{align}
where $\kappa_z(a)=1+F_{\eta'}([a](z))$ with $ 1+F_{\eta'}(Z):=\exp \left(\Omega_{}\log_{LT}(Z)\right).$ It is an isomorphism of topological rings which is equivariant with respect to the action by the monoid $o_L\setminus \{0\}$ (as a consequence of \cite[Prop.\ 3.1]{ST2}). Moreover, Lemma \ref{module-struc} implies that the $o_L^\times$-action on $\mathcal{O}_K(\mathbf{B})$ extends to a jointly continuous $D(o_L^\times,K)$-module structure (by descent even for general $K$) and that the LT-isomorphism is an isomorphism of $D(o_L^\times,K)$-modules.

 By the construction of the LT-isomorphism we have
\begin{equation*}
  \kappa^*(\ev_a) = \exp(a \Omega \log_{LT}(Z)) \in o_{\mathbb{C}_p}[[Z]]  \qquad\text{for any $a \in o_L$}.
\end{equation*}
Hence Lemma \ref{psi-ker} implies that
\begin{equation*}
  \kappa^*(\ker(\psi_L^{\mathfrak{X}})) = \sum_{a \in R_0} \exp(a \Omega \log_{LT}(Z)) \varphi_L(\mathcal{O}_K(\mathbf{B}))
\end{equation*}
where $R_0 \subseteq o_L$ denotes a set of representatives for the nonzero cosets in $o_L/\pi_L o_L$. Using that $\log_{LT} (Z_1 +_{LT} Z_2) = \log_{LT}(Z_1) + \log_{LT}(Z_2)$ we compute
\begin{align*}
  tr (\exp(a \Omega \log_{LT}(Z)) & = \sum_{y \in \ker([\pi_L])} \exp(a \Omega \log_{LT}(y +_{LT}Z))   \\
   & = \big( \sum_{y \in \ker([\pi_L])} \exp(a \Omega \log_{LT}(y)) \big) \exp(a \Omega \log_{LT}(Z))  \\
   & = \big( \sum_{y \in \ker([\pi_L])} \kappa_y(a) \big) \exp(a \Omega \log_{LT}(Z)) \ .
\end{align*}
But the $\kappa_y$ for $y \in \ker([\pi_L])$ are precisely the characters of the finite abelian group $o_L/\pi_L oL$. Hence $\sum_{y \in \ker([\pi_L])} \kappa_y(a) = 0$ for $a \in R_0$.
It follows that $\kappa^*(\ker(\psi_L^{\mathfrak{X}})) = \ker(\psi_L)$. We conclude that under the LT-isomorphism $\psi_L$ corresponds to $\frac{q}{\pi_L} \psi_L^{\mathfrak{X}}$ using the fact that we also have a decomposition
\begin{equation}\label{f:Odecomp}
  \mathcal{O}_K(\mathbf{B}) = \sum_{a \in o_L/\pi_L} \exp(a \Omega \log_{LT}(Z)) \varphi_L(\mathcal{O}_K(\mathbf{B})).
\end{equation}

In the following we denote by
\[\mathfrak{M}_{LT} : D(\Gamma_L,K) \xrightarrow{\cong} \mathcal{O}_{K}(\mathbf{B})^{{\psi_L}=0}\] the composite
\begin{align*}
  D(\Gamma_L,K)\cong D(o_L^\times,K)\cong\mathcal{O}_{K}(\mathfrak{X})^{{\psi_L^\mathfrak{X}} = 0}\cong \mathcal{O}_{K}(\mathbf{B})^{{\psi_L}=0}
\end{align*}
where the first isomorphism is induced by the character $\chi_{LT}:\Gamma_L\xrightarrow{\cong} o_L^\times,$ the second one is the Mellin transform $\mathfrak{M}$ from Lemma \ref{Mellin} while the third one is the LT-isomorphism. By inserting the definitions we obtain the explicit formula
\begin{equation*}
  \mathfrak{M}_{LT} (\lambda)(z) = \lambda(\kappa_z \circ \chi_{LT}) \ .
\end{equation*}
The construction of the above map $\mathfrak{M}_{LT}$ is related to the pairing
\begin{align*}
 \{\;,\;\} : \mathcal{O}_{K}(\mathbf{B}) \times C^{an}(o_L ,K) & \to K \\
 (F,f) & \mapsto \mu(f), \quad \text{ where $\mu \in D(o_L,K)$ is such that $\mu(\kappa_z)=F(z),$ }
\end{align*}
in \cite[lem.\ 4.6]{ST2} by the following commutative diagram:
\begin{equation*}
  \xymatrix{
    D(\Gamma_L,K) \ar[d]_{\mathfrak{M}_{LT}} & \times & C^{an}(\Gamma_L,K) \ar[d]^{(\chi_{LT})_*} \ar[rr]^-{(\lambda,f) \mapsto \lambda(f)} && K \ar@{=}[dd]^{} \\
    \mathcal{O}_{K}(\mathbf{B})^{\psi_L = 0} \ar[d]_{\subseteq}  & \times & C^{an}(o_L^\times,K) \ar[d]^{extension\ by\ 0}  &  \\
    \mathcal{O}_{K}(\mathbf{B}) & \times & C^{an}(o_L ,K) \ar[rr]^-{\{\;,\;\}} && K   }
\end{equation*}

\begin{remark}\label{rem:inc}
   ($\Omega \in K$) For any $F \in \mathcal{O}_{K}(\mathbf{B})^{\psi_L = 0}$ and any $f \in C^{an}(o_L,K)$ such that $f|o_L^\times = 0$ we have $\{F,f\} = 0$.
\end{remark}
\begin{proof}
We have seen above that under the LT-isomorphism $\psi_L$ corresponds, up to a nonzero constant, to $\psi_L^{\mathfrak{X}}$ and hence further under the Fourier isomorphism to $\psi_L^D$. It therefore suffices to show that for any $\mu \in D(o_L,K)^{\psi_L^D = 0}$ we have $\mu(f) = 0$. For this we define $\tilde{f} := f(\pi_L -) \in C^{an}(o_L,K)$ and note that $(\pi_L)_!(\tilde{f}) = f$. By the definition of $\psi_L^D$ we therefore obtain, under our assumption on $\mu$, that $\mu(f) = \mu(f) - \psi_L^D(\mu)(\tilde{f}) = \mu(f - (\pi_L)_!(\tilde{f})) = \mu(0) = 0$.
\end{proof}

\begin{lemma}\label{phi}
($\Omega \in K$) For any $F \in \mathcal{O}_{K}(\mathbf{B})^{\psi_L=1}$ and $n \geq 0$ we have
\begin{equation*}
  \mathfrak{M}_{LT}^{-1}((1-\frac{\pi_L}{q}\varphi_L)F)(\chi_{LT}^n) = \Omega^{-n}(1-\frac{\pi_L^{n+1}}{q}) (\partial_{\mathrm{inv}} ^n F)_{|Z=0} .
\end{equation*}
\end{lemma}
\begin{proof}
Note that $ (1-\frac{\pi_L}{q}\varphi_L)F$ belongs to $\mathcal{O}_{K}(\mathbf{B})^{\psi_L=0}$. Let $\mathrm{inc}_! \in C^{an}(o_L,K)$ denote the extension by zero of the inclusion $o_L^\times \subseteq o_L$, and let $\id : o_L \rightarrow K$ be the identity function. Using the above commutative diagram the assertion reduces to the equality
\begin{equation*}
  \{ (1-\frac{\pi_L}{q}\varphi_L)F, \mathrm{inc}_!^n \} = \Omega^{-n}(1-\frac{\pi_L^{n+1}}{q}) (\partial_{\mathrm{inv}} ^n F)_{|Z=0} .
\end{equation*}
By Remark \ref{rem:inc} we may replace on the left hand side the function $\mathrm{inc}_!^n$ by the function $\id^n$. Next we observe that $x \, \id^n(x) = \id^{n+1}(x)$. Hence, by \cite[Lem.\ 4.6(8)]{ST2}, i.e., $\{F,xf(x)\}=\{\Omega^{-1}\partial_{\mathrm{inv}}  F,f\}$ and induction, the left hand side is equal to
\begin{align*}
  \{ (1-\frac{\pi_L}{q}\varphi_L)F, \id^n \} & = \{ \Omega^{-n}\partial_{\mathrm{inv}}^n ((1-\frac{\pi_L}{q}\varphi_L)F), \id^0 \}    \\
      & = \{ \Omega^{-n} (1-\frac{\pi_L^{n+1}}{q}\varphi_L) (\partial_{\mathrm{inv}}^n F), \id^0 \}      \quad\text{since $\partial_{\mathrm{inv}} \varphi_L = \pi_L\varphi_L \partial_{\mathrm{inv}}$ by \eqref{f:dlog}}      \\
      & = \Omega^{-n} (1-\frac{\pi_L^{n+1}}{q}\varphi_L) (\partial_{\mathrm{inv}}^n F)_{|Z=0}   \quad\text{since $\id^0$ is the trivial character of $o_L$}  \\
      & = \Omega^{-n} (1-\frac{\pi_L^{n+1}}{q}) (\partial_{\mathrm{inv}}^n F)_{|Z=0}   \quad\text{since $[\pi_L](0)=0$}  .
\end{align*}
\end{proof}

In the course of the previous proof we have seen that, for $F$ in $\mathcal{O}_{K}(\mathbf{B})^{\psi_L=0}$ and $n \geq 0$,
\begin{equation}\label{f:evaluation}
   \mathfrak{M}_{LT}^{-1}(F)(\chi_{LT}^n) = \Omega^{-n} (\partial_{\mathrm{inv}}^n F)_{|Z=0} \ .
\end{equation}

\begin{lemma}\label{log}
($\Omega \in K$) For any $F \in \mathcal{O}_{K}(\mathbf{B})^{\psi_L=0}$ and $n \geq 1$ we have
\begin{equation*}
  \mathfrak{M}_{LT}^{-1}(\log_{LT}\cdot F)(\chi_{LT}^n) = n \Omega^{-1} \mathfrak{M}_{LT}^{-1}(F)(\chi_{LT}^{n-1}) \ .
\end{equation*}
\end{lemma}
\begin{proof}
First, using \eqref{f:dlog}, observe that
\begin{equation*}
  \psi_L (\log_{LT}\cdot F) = \psi_L (\pi_L^{-1} \varphi_L(\log_{LT}) \cdot F) = \pi_L^{-1} \varphi_L(\log_{LT}) \psi_L(F) = 0 \ .
\end{equation*}
Secondly note that $\partial_{\mathrm{inv}} \log_{LT}=1$, i.e., $\partial_{\mathrm{inv}}^i \log_{LT} = 0$ for $i \geq 2$; also $\log_{LT}(0)=0$. Using \eqref{f:evaluation} twice we have
\begin{align*}
\mathfrak{M}_{LT}^{-1}(\log_{LT}F)(\chi_{LT}^n) & = \Omega^{-n} (\partial_{\mathrm{inv}}^n (\log_{LT}F))_{|Z=0}    \\
                   & = \Omega^{-n} \left( \sum_{i+j=n} \left(
                                 \begin{array}{c}
                                   n \\
                                   i \\
                                 \end{array}
                               \right)
   (\partial_{\mathrm{inv}}^i \log_{LT}) (\partial_{\mathrm{inv}}^j F) \right)_{|Z=0}    \\
   &= \Omega^{-n} n (\partial_{\mathrm{inv}}^{n-1} F)_{|Z=0}    \\
   &= n \Omega^{-1} \mathfrak{M}^{-1}(F)(\chi_{LT}^{n-1}) \ .
\end{align*}
\end{proof}

For the rest of this section we \textbf{assume} not only that $K$ contains $\Omega$ but also that the action of  $G_L$ on $\mathbb{C}_p$ leaves $K$ invariant.

The LT-isomorphism is a topological ring isomorphism
\begin{align*}
K\widehat{\otimes}_L\mathcal{O}_L(\mathfrak{X})=\mathcal{O}_K(\mathfrak{X})\cong \mathcal{O}_K(\mathbf{B})=K\widehat{\otimes}_L\mathcal{O}_L(\mathbf{B})
\end{align*}
(see \cite[Prop.\ 2.1.5 ii.]{BSX} for the outer identities).

On both sides we have the obvious coefficientwise $G_L$-action induced by the Galois-action on the tensor-factor $K$.
We use the following notation:

\begin{itemize}
\item $\sigma\in G_L$ acting {\em coefficientwise} on $\mathcal{O}_K(\mathbf{B})$ is denoted by: $F\mapsto {^\sigma F};$ the corresponding fixed ring is $\mathcal{O}_K(\mathbf{B})^{G_L}=\mathcal{O}_L(\mathbf{B}).$
\item The coefficientwise action on $\mathcal{O}_K(\mathfrak{X})$ transfers to the {\em twisted action } on $\mathcal{O}_K(\mathbf{B})$ by \cite[before Cor.\ 3.8]{ST2} given as $F\mapsto {    ^{\sigma*}F}:={^\sigma F}\circ [\tau(\sigma^{-1})];$ the corresponding fixed ring is $\mathcal{O}_K(\mathbf{B})^{G_L,*}=\mathcal{O}_L(\mathfrak{X})=D(o_L,L).$
\end{itemize}

\begin{remark}\label{comactions}
 Note that the   $o_L \setminus \{0\}$-action and hence the $D(o_L^\times,L)$-module structure commute with both $G_L$-actions.   Moreover, $\psi_L$ commutes with the $G_L$-actions as well.
\end{remark}

Recall that using the notation from \cite[Lem. 4.6, 1./2.]{ST2} the function $1+F_{a\eta'}(Z)=\exp \left(a\Omega_{\eta'}\log_{LT}(Z)\right)$ corresponds to the Dirac distribution $\delta_a$ of $a\in o_L$ under the Fourier isomorphism.

\begin{lemma}\label{lem:actionperiod} Let $\sigma$ be in $G_L, $  $t'\in \TLT'$ and $a\in o_L.$ Then
\begin{enumerate}
\item $\sigma(\Omega_{t'})=\Omega_{\tau(\sigma)t'}=\Omega_{t'}\tau(\sigma)$ and
\item ${^\sigma F_{a\eta'}}=F_{a\eta'}\circ[\tau(\sigma)]=F_{a\tau(\sigma)\eta'}$.
\end{enumerate}
\end{lemma}

\begin{proof}
(i) The Galois equivariance of the pairing $(\; ,\; ):\TLT'\otimes_{o_L}\mathbb{C}_p\to\mathbb{C}_p , $  from (loc.~cit.\ before Prop.~ 3.1) with $(t',x)=\Omega_{t'}x$ implies that \[\sigma(\Omega_{t'})=\Omega_{\sigma(t')}=\Omega_{\tau(\sigma)t'}\] while the $o_L$-invariance of that pairing implies that the latter expression equals $\Omega_{t'}\tau(\sigma)$.

(ii) This is immediate from (i) and the definition of $F_a$ taking equation \eqref{f:dlog} into account.
\end{proof}

\begin{proposition}\phantomsection\label{prop:twistinvariance}
\begin{enumerate}
\item The      isomorphism (LT together with Fourier) $ \textfrak{L}:   D(o_L,K)\cong \mathcal{O}_K(\mathfrak{X})\cong \mathcal{O}_K(\mathbf{B})$ restricts to an isomorphism
\begin{equation*}
  D(o_L,K)^{G_L,\tilde{*}}=\mathcal{O}_K(\mathfrak{X})^{G_L}\cong \mathcal{O}_K(\mathbf{B})^{G_L}= \mathcal{O}_L(\mathbf{B})
\end{equation*}
of $D(o_L^\times,L)$-modules.
\item The Mellin transform restricts to an isomorphism of $D(o_L^\times,L)$-modules
\begin{equation*}
  D(o_L^\times,K)^{G_L,*}=\mathcal{O}_K(\mathfrak{X})^{G_L,\psi_L=0}\cong \mathcal{O}_L(\mathbf{B})^{\psi_L=0}.
\end{equation*}
\end{enumerate}
Here the $G_L$-action on the distribution rings on the left hand sides is induced from the coefficientwise action on $\mathcal{O}_K(\mathbf{B})$ and $\mathcal{O}_K(\mathbf{B})^{\psi_L=0}$ via  the LT-isomorphism and Mellin transform, respectively.\Footnote{It is not true that the $*$-action of $G_L$ on $ D(o_L^\times,K)$ coincides on some $D(U_n,K)$ with the action, which arises on the latter via transport of structures along $\ell_n$ from the $\tilde{*}$-action of $G_L$ on  $D(o_L,K)$ as in (i)!}
\end{proposition}

\begin{proof} (i) and (ii) follow from passing to the fixed vectors with respect to the coefficientwise $G_L$-action and Remark \ref{comactions}.
\end{proof}



In order to express the $D(o_L^\times,L)$-module $D(o_L^\times,K)^{G_L,*} $ in the above proposition more explicitly we describe the previous two actions on $\mathcal{O}_K(\mathbf{B})$ now on $D(o_L,K)$:
\begin{itemize}
\item The coefficientwise $G_L$-action on $D(o_L,K)=K\widehat{\otimes}_L D(o_L,L),$ which corresponds to the twisted action on $\mathcal{O}_K(\mathbf{B}),$ will be written as $\lambda\mapsto{^\sigma\lambda}.$
\item The $G_L$-action given by $\lambda\mapsto \tau(\sigma)_*({^\sigma}\lambda)$ corresponds to the coefficientwise action on $\mathcal{O}_K(\mathbf{B}).$
\end{itemize}
Note that for  $\lambda\in D(o_L^\times,K)$ we have
    $\tau(\sigma)_*(\lambda)=\delta_{\tau(\sigma)}\lambda$, where the right hand side refers to the product of $\lambda$ and the Dirac distribution $ \delta_{\tau(\sigma)}$ in the ring
$ D(o_L^\times,K).$  Then we conclude that
\[D(o_L^\times,K)^{G_L,*}=\{\lambda\in D(o_L^\times,K)| \; {^\sigma\lambda}=\delta_{\tau(\sigma)^{-1}}\lambda \mbox{ for all }\sigma\in G_L\}.  \]

\newpage

\subsection{Consequences of Serre duality}\label{sec:duality}

Recall that in any quasi-separated rigid analytic variety the complement of any affinoid subdomain is admissible open (\cite[\S3 Prop.\ 3(ii)]{Sch}). This applies in particular to quasi-Stein spaces since they are separated by definition. For a rigid analytic variety $\mathfrak{Y}$ we will denote by $\Aff(\mathfrak{Y})$ the set of all affinoid subdomains of $\mathfrak{Y}$.

We have seen that $\mathfrak{X}$, $\mathfrak{X}^\times$, and $\mathfrak{X}_n^\times$ for $n \geq 1$ all are 1-(equi)dimensional smooth Stein spaces.

\subsubsection{Cohomology with compact support}\label{subsec:compact-supp}

We slightly rephrase the definition of cohomology with compact support given in \cite[\S1]{vdP} in the case of a Stein space $\mathfrak{Y}$ over $L$. For any abelian sheaf $\mathcal{F}$ on $\mathfrak{Y}$ and any $\mathfrak{U} \in \Aff(\mathfrak{Y})$ we put
\begin{equation*}
  H^0_\mathfrak{U}(\mathfrak{Y}, \mathcal{F}) := \ker (\mathcal{F}(\mathfrak{Y}) \longrightarrow \mathcal{F}(\mathfrak{Y} \setminus \mathfrak{U})) \ .
\end{equation*}
This is a left exact functor in $\mathcal{F}$, and we denote by $H^*_\mathfrak{U}(\mathfrak{Y}, \mathcal{F})$ its right derived functors. Since quasi-Stein spaces have no coherent cohomology the relative cohomology sequence (\cite[Lem.\ 1.3]{vdP}) gives rise, for a coherent sheaf $\mathcal{F}$, to the exact sequence
\begin{equation}\label{f:relative}
  0 \rightarrow  H^0_\mathfrak{U}(\mathfrak{Y}, \mathcal{F}) \rightarrow \mathcal{F}(\mathfrak{Y}) \rightarrow \mathcal{F}(\mathfrak{Y} \setminus \mathfrak{U}) \rightarrow H^1_\mathfrak{U}(\mathfrak{Y}, \mathcal{F}) \rightarrow 0 \ .
\end{equation}
We then define the cohomology with compact support as
\begin{equation*}
  H^*_c(\mathfrak{Y}, \mathcal{F}) := \varinjlim_{\mathfrak{U} \in \Aff(\mathfrak{Y})} H^*_\mathfrak{U}(\mathfrak{Y}, \mathcal{F}) \ .
\end{equation*}
Again, if $\mathcal{F}$ is a coherent sheaf we obtain the exact sequence
\begin{equation*}
  0 \rightarrow  H^0_c(\mathfrak{Y}, \mathcal{F}) \rightarrow \mathcal{F}(\mathfrak{Y}) \rightarrow \varinjlim_{\mathfrak{U} \in \Aff(\mathfrak{Y})} \mathcal{F}(\mathfrak{Y} \setminus \mathfrak{U}) \rightarrow H^1_c(\mathfrak{Y}, \mathcal{F}) \rightarrow 0 \ .
\end{equation*}

Suppose in the following that $j : \mathfrak{Y}_0 \rightarrow \mathfrak{Y}$ is an open immersion of Stein spaces (over $L$) which, for simplicity we view as an inclusion.

\begin{lemma}\label{covering}
  For any $\mathfrak{U} \in \Aff(\mathfrak{Y}_0)$ the covering $\mathfrak{Y} = (\mathfrak{Y} \setminus \mathfrak{U}) \cup \mathfrak{Y}_0$ is admissible.
\end{lemma}
\begin{proof}
This follows from \cite[Lem.\ 1.1]{vdP}.
\end{proof}

\begin{lemma}\label{c-functoriality}
  For any $\mathfrak{U} \in \Aff(\mathfrak{Y}_0)$ and any sheaf $\mathcal{F}$ on $\mathfrak{Y}$ the natural map
\begin{equation*}
  H^*_\mathfrak{U}(\mathfrak{Y}, \mathcal{F}) \xrightarrow[\cong]{\;\mathrm{res}\;} H^*_\mathfrak{U}(\mathfrak{Y}_0, \mathcal{F})
\end{equation*}
is bijective.
\end{lemma}
\begin{proof}
Recall that for an injective sheaf on $\mathfrak{Y}$ its restriction to $\mathfrak{Y}_0$ is injective as well. Hence, by using an injective resolution, it suffices to proof the assertion for $* = 0$. \textit{Injectivity}: Let $f \in H^0_\mathfrak{U}(\mathfrak{Y}, \mathcal{F})$ such that $f | \mathfrak{Y}_0 = 0$. Since $f | \mathfrak{Y} \setminus \mathfrak{U} = 0$ as well it follows from Lemma \ref{covering} that $f = 0$. \textit{Surjectivity}: Let $g \in H^0_\mathfrak{U}(\mathfrak{Y}_0, \mathcal{F})$ so that $g | \mathfrak{Y}_0 \setminus \mathfrak{U} = 0$. Using Lemma \ref{covering} again we may define a preimage $f$ of $g$ by $f | \mathfrak{Y} \setminus \mathfrak{U} := 0$ and $f | \mathfrak{Y}_0 := g$.
\end{proof}

By passing to inductive limits we obtain the composed map
\begin{multline*}
  j_! : H^*_c(\mathfrak{Y}_0, \mathcal{F}) = \varinjlim_{\mathfrak{U} \in \Aff(\mathfrak{Y}_0)} H^*_\mathfrak{U}(\mathfrak{Y}_0, \mathcal{F}) \xrightarrow[\cong]{\mathrm{res}^{-1}} \varinjlim_{\mathfrak{U} \in \Aff(\mathfrak{Y}_0)} H^*_\mathfrak{U}(\mathfrak{Y}, \mathcal{F})    \\
  \rightarrow \varinjlim_{\mathfrak{U} \in \Aff(\mathfrak{Y})} H^*_\mathfrak{U}(\mathfrak{Y}, \mathcal{F}) =    H^*_c(\mathfrak{Y}, \mathcal{F}) \ .
\end{multline*}
For later purposes we have to analyze the following situation. Let
\begin{equation*}
  \mathfrak{U}_1 \subseteq  \mathfrak{U}_2 \subseteq \ldots \subseteq  \mathfrak{U}_n \subseteq \ldots \subseteq \mathfrak{Y} = \bigcup_n  \mathfrak{U}_n
\end{equation*}
be a Stein covering. We \textbf{assume} that each admissible open subset $\mathfrak{Y}_n := \mathfrak{Y} \setminus \mathfrak{U}_n$ also is a Stein space. Since $\ldots \subseteq \mathfrak{Y}_n \subseteq \ldots \subseteq \mathfrak{Y}_1 \subseteq \mathfrak{Y}$ we then have the projective system
\begin{equation*}
  \ldots \rightarrow H^*_c(\mathfrak{Y}_n,\mathcal{F}) \rightarrow \ldots \rightarrow H^*_c(\mathfrak{Y}_1,\mathcal{F}) \rightarrow H^*_c(\mathfrak{Y},\mathcal{F}) \ .
\end{equation*}
By Lemma \ref{c-functoriality} we may rewrite it as the projective system
\begin{equation}\label{f:projsystem}
   \ldots \rightarrow \varinjlim_{\mathfrak{U} \in \Aff(\mathfrak{Y}_n)} H^*_\mathfrak{U}(\mathfrak{Y}, \mathcal{F}) \rightarrow \ldots \rightarrow \varinjlim_{\mathfrak{U} \in \Aff(\mathfrak{Y}_1)} H^*_\mathfrak{U}(\mathfrak{Y}, \mathcal{F}) \rightarrow \varinjlim_{\mathfrak{U} \in \Aff(\mathfrak{Y})} H^*_\mathfrak{U}(\mathfrak{Y}, \mathcal{F}) \ .
\end{equation}

\begin{lemma}\label{projlim}
   In the above situation we assume in addition that $\mathcal{F}$ is a coherent sheaf and that the restriction maps $\mathcal{F}(\mathfrak{Y}) \rightarrow \mathcal{F}(\mathfrak{Y} \setminus \mathfrak{U})$ for any $\mathfrak{U} \in \Aff(\mathfrak{Y})$ are injective. We then have
\begin{equation}\label{f:projlim}
  \varprojlim_n H^1_c(\mathfrak{Y}_n,\mathcal{F}) = \big( \varprojlim_n \varinjlim_{\mathfrak{U} \in \Aff(\mathfrak{Y}_n)} \mathcal{F}(\mathfrak{Y} \setminus \mathfrak{U}) \big) /\mathcal{F}(\mathfrak{Y}) \ .
\end{equation}
\end{lemma}
\begin{proof}
This is immediate from \eqref{f:projsystem} and the relative cohomology sequence.
\end{proof}

For coherent sheaves $\mathcal{F}$ all the above cohomology vector spaces carry natural locally convex topologies, which we briefly recall. The global sections $\mathcal{F}(\mathfrak{Y})$ and $\mathcal{F}(\mathfrak{Y} \setminus \mathfrak{U})$, for $\mathfrak{U} \in \Aff(\mathfrak{Y})$, are Fr\'echet spaces. Using the relative cohomology sequence \eqref{f:relative} we equip $H^1_{\mathfrak{U}}(\mathfrak{Y},\mathcal{F})$ with the quotient topology from $\mathcal{F}(\mathfrak{Y} \setminus \mathfrak{U})$ (which might be non-Hausdorff) and then $H^1_c(\mathfrak{Y},\mathcal{F})$ with the locally convex inductive limit topology (w.r.t.\ varying $\mathfrak{U}$).

\begin{remark}\label{rem:BS}
  Let
\begin{equation*}
  \xymatrix{
    M^0 \ar[d] \ar[r]^{\alpha} & M^1 \ar[d] \\
    N^0 \ar[r]^{\beta} & N^1   }
\end{equation*}
  be a commutative diagram of Fr\'echet spaces such that the induced map $\coker(\alpha) \rightarrow \coker(\beta)$ is bijective; then the latter map is a topological isomorphism for the quotient topologies.
\end{remark}
\begin{proof}
A more general statement can be found in \cite[Chap.\ VII Lem.\ 1.32]{BS}.
\end{proof}

Using this Remark we see that the bijection $H^1_\mathfrak{U}(\mathfrak{Y}, \mathcal{F}) \xrightarrow[\cong]{\mathrm{res}} H^1_\mathfrak{U}(\mathfrak{Y}_0, \mathcal{F})$ from Lemma \ref{c-functoriality} is a topological isomorphism. It follows that, in degree one at least, the above map $j_!$ is as well as the transition maps in the projective system \eqref{f:projsystem} are continuous.

\begin{lemma}\label{projlim-topol}
  The isomorphism \eqref{f:projlim} is topological.
\end{lemma}
\begin{proof}
The assertion has to be understood, of course, in such a way that forming projective limits, inductive limits, and quotient spaces on both sides of \eqref{f:projlim} is meant in the topological sense. First of all one checks that forming the projective limit on the right hand side commutes with passing to the quotient space by $\mathcal{F}(\mathfrak{Y})$ (compare (ii) in the proof of \cite[Thm.\ 4.3]{Pro} for a more general statement). Secondly, as a special case of \cite[II.28 Cor.\ 2]{B-TVS}, forming inductive limits commutes with passing to quotient spaces. This reduces us to $H^1_{\mathfrak{U}}(\mathfrak{Y} \setminus \mathfrak{U},\mathcal{F}) = \mathcal{F}(\mathfrak{Y} \setminus \mathfrak{U}) /\mathcal{F}(\mathfrak{Y})$ being a topological isomorphism, but which holds by definition.
\end{proof}

In the following we compute \eqref{f:projlim} further in two concrete cases.\\

\noindent
{\bf The open unit disk}\\

Let $\mathbf{B} = \mathbf{B}_{[0,1)}$ denote the open unit disk over $L$. We recall our convention that all radii are assumed to lie in $(0,1) \cap p^{\mathbb{Q}}$. For any radii $r \leq s$ we introduce the affinoid disk $\mathbf{B}_{[0,r]}$ as well as the open disk $\mathbf{B}_{[0,r)}$ of radius $r$ around $0$ and the affinoid annulus $\mathbf{B}_{[r,s]} := \{ r \leq |x| \leq s \}$. We put $\mathbf{B}_{(r,1)} := \mathbf{B} \setminus \mathbf{B}_{[0,r]}$, which are Stein spaces. By the identity theorem for Laurent series the assumptions of Lemma \ref{projlim} are satisfied for the structure sheaf $\mathcal{O} = \mathcal{O}_{\mathbf{B}}$ of $\mathbf{B}$. We first fix a radius $r$ and compute the cohomology $H^1_c(\mathbf{B}_{(r,1)},\mathcal{O})$. By Lemma \ref{c-functoriality} and the relative cohomology sequence  we have
\begin{equation*}
  H^1_c(\mathbf{B}_{(r,1)},\mathcal{O}) = \varinjlim_{\mathfrak{U} \in \Aff(\mathbf{B}_{(r,1)})} H^*_\mathfrak{U}(\mathbf{B}, \mathcal{O}) = \big( \varinjlim_{\mathfrak{U} \in \Aff(\mathbf{B}_{(r,1)})} \mathcal{O}(\mathbf{B} \setminus \mathfrak{U}) \big) /\mathcal{O}(\mathbf{B})  \ .
\end{equation*}
Of course, it suffices to take the inductive limit over a cofinal sequence of larger and larger affinoid annuli in $\mathbf{B}_{(r,1)}$. For this we choose two sequences of radii $r < \ldots < r_m < \ldots < r_1\bigcup < s_1 < \ldots < s_m < \ldots < 1$ with $(r_m)_m$ and $(s_m)_m$ converging to $r$ and $1$, respectively. Each space $\mathbf{B} \setminus \mathbf{B}_{[r_m,s_m]} = \mathbf{B}_{(s_m,1)} \dot{\cup} \mathbf{B}_{[0,r_m)}$ has two connected components. We see that
\begin{equation*}
  H^1_c(\mathbf{B}_{(r,1)},\mathcal{O}) = \big( \varinjlim_{m \rightarrow \infty} \mathcal{O}(\mathbf{B}_{(s_m,1)}) \oplus \varinjlim_{m \rightarrow \infty} \mathcal{O}(\mathbf{B}_{[0,r_m)}) \big) / \mathcal{O}(\mathbf{B}) \ .
\end{equation*}
As explained in Lemma \ref{projlim-topol} this is a topological equality. We observe that $\mathcal{R}_L = \mathcal{R}_L(\mathbf{B}) = \varinjlim_{m \rightarrow \infty} \mathcal{O}(\mathbf{B}_{(s_m,1)})$ is the usual Robba ring (over $L$) whereas $\mathcal{O}^\dagger(\mathbf{B}_{[0,r]}) = \varinjlim_{m \rightarrow \infty} \mathcal{O}(\mathbf{B}_{[0,r_m)})$ is the ring of overconvergent analytic functions on $\mathbf{B}_{[0,r]}$. Hence
\begin{equation*}
  H^1_c(\mathbf{B}_{(r,1)},\mathcal{O}) = \big( \mathcal{R}_L \oplus \mathcal{O}^\dagger(\mathbf{B}_{[0,r]}) \big) / \mathcal{O}(\mathbf{B}) \ .
\end{equation*}
Passing now to the projective limit w.r.t.\ $r \rightarrow 1$ of the continuous restriction maps $\mathcal{O}(\mathbf{B}) \rightarrow \mathcal{O}^\dagger(\mathbf{B}_{[0,r]}) \rightarrow \mathcal{O}(\mathbf{B}_{[0,r]})$ we observe that $\varprojlim_{r \rightarrow 1} \mathcal{O}^\dagger(\mathbf{B}_{[0,r]}) = \mathcal{O}(\mathbf{B})$ holds true topologically. We finally deduce that
\begin{equation}\label{f:projlim=Robba}
  \varprojlim_{r \rightarrow 1} H^1_c(\mathbf{B}_{(r,1)},\mathcal{O})=\varprojlim_{r \rightarrow 1}\big( \varinjlim_{\mathfrak{U} \in \Aff(\mathbf{B}_{(r,1)})} \mathcal{O}(\mathbf{B} \setminus \mathfrak{U}) \big) /\mathcal{O}(\mathbf{B}) = \big( \mathcal{R}_L \oplus \mathcal{O}(\mathbf{B}) \big) / \mathcal{O}(\mathbf{B}) \cong \mathcal{R}_L
\end{equation}
as locally convex vector spaces.\\

\noindent
{\bf The character variety $\mathfrak{X}$}\\

Since $\mathfrak{X}_{/\mathbb{C}_p} \cong \mathbf{B}_{/\mathbb{C}_p}$ by \cite{ST2} the injectivity of the restriction maps $\mathcal{O}(\mathfrak{X}) \rightarrow \mathcal{O}(\mathfrak{X} \setminus \mathfrak{U})$ for any $\mathfrak{U} \in \Aff(\mathfrak{X})$ follows from the corresponding fact for $\mathbf{B}$, which we saw already. According to Prop.\ \ref{Stein} the admissible open subdomains $\mathfrak{X}_{(s_n,1)}$ of $\mathfrak{X}$ are Stein spaces. In order to compute their cohomology with compact support in the structure sheaf $\mathcal{O} = \mathcal{O}_{\mathfrak{X}}$ we fix an $n \geq 0$. We choose a sequence of radii $r_{n+1} > r_{n+2} > \ldots > s_n$ in $S_n$ converging to $s_n$. Furthermore we observe that the increasing sequence $(s_m)_{m > n}$ in $S_\infty$ converges to $1$. By Prop.\ \ref{quasi-Stein} we then have the Stein covering
\begin{equation*}
  \mathfrak{X}_{(s_n,1)} = \bigcup_{m > n} \mathfrak{X}_{[r_m,s_m]} \ .
\end{equation*}
Hence, by Lemma \ref{c-functoriality}, the relative cohomology sequence, and the explanation in Lemma \ref{projlim-topol}, we have the topological equality
\begin{equation*}
  H^1_c(\mathfrak{X}_{(s_n,1)},\mathcal{O}) = \big( \varinjlim_{m \rightarrow \infty} \mathcal{O}(\mathfrak{X} \setminus \mathfrak{X}_{[r_m,s_m]}) \big) / \mathcal{O}(\mathfrak{X}) \ .
\end{equation*}
The obvious set theoretic decomposition
\begin{equation*}
  \mathfrak{X} \setminus \mathfrak{X}_{[r_m,s_m]} = \mathfrak{X} \setminus (\mathfrak{X}(s_m) \setminus \mathfrak{X}^-(r_m)) = (\mathfrak{X} \setminus \mathfrak{X}(s_m)) \, \dot\cup \, \mathfrak{X}^-(r_m) = \mathfrak{X}_{(s_m,1)} \, \dot\cup \, \mathfrak{X}^-(r_m)
\end{equation*}
is in fact the decomposition of the space $\mathfrak{X} \setminus \mathfrak{X}_{[r_m,s_m]}$ into its connected components. This can be checked after base change to $\mathbb{C}_p$ where, by \cite[Prop.\ 1.20 and proof of Prop.\ 2.1]{BSX}, the setting becomes isomorphic to the setting for the open unit disk which we discussed in the previous section. Entirely similarly as in the previous subsection it follows now that
\begin{equation*}
  H^1_c(\mathfrak{X}_{(s_n,1)},\mathcal{O}) = \big( \mathcal{R}_L(\mathfrak{X}) \oplus \mathcal{O}^\dagger(\mathfrak{X}(s_n)) \big) / \mathcal{O}(\mathfrak{X})
\end{equation*}
where $\mathcal{O}^\dagger(\mathfrak{X}(s_n)) := \varinjlim_{m \rightarrow \infty} \mathfrak{X}^-(r_m)$, and then
\begin{equation}\label{f:projlim=RobbaX}
  \varprojlim_{n \rightarrow \infty} H^1_c(\mathfrak{X}_{(s_n,1)},\mathcal{O}) =\varprojlim_{n \rightarrow \infty}\big( \varinjlim_{m \rightarrow \infty} \mathcal{O}(\mathfrak{X} \setminus \mathfrak{X}_{[r_m,s_m]}) \big) / \mathcal{O}(\mathfrak{X}) = \big( \mathcal{R}_L(\mathfrak{X}) \oplus \mathcal{O}(\mathfrak{X}) \big) / \mathcal{O}(\mathfrak{X}) \cong \mathcal{R}_L(\mathfrak{X})
\end{equation}
as locally convex vector spaces.

\subsubsection{Serre duality for Stein spaces}\label{subsec:duality}

In the following the continuous dual of a locally convex vector space is always equipped with the strong topology.

The Serre duality for smooth Stein spaces is established in \cite{Chi} and \cite{Bey}. Let $\mathfrak{Y}$ be a 1-(equi)dimensional smooth Stein space over $L$.

\begin{theorem}\label{Serre}
   For any coherent sheaf $\mathcal{F}$ on $\mathfrak{Y}$ we have:
\begin{itemize}
  \item[i.] $H^1_c(\mathfrak{Y}, \mathcal{F})$ is a complete reflexive Hausdorff space.
  \item[ii.] $\Hom_\mathfrak{Y}(\mathcal{F}, \Omega^1_\mathfrak{Y}) = H^0(\mathfrak{Y}, \underline{\Hom}_\mathfrak{Y}(\mathcal{F}, \Omega^1_\mathfrak{Y}))$, being the global sections of another coherent sheaf, is a reflexive Fr\'echet space strictly of countable type (\cite[Def.\ 4.2.3]{PGS}).
  \item[iii.] There is a canonical trace map
\begin{equation*}
  t_\mathfrak{Y} : H^1_c(\mathfrak{Y}, \Omega^1_\mathfrak{Y}) \longrightarrow L
\end{equation*}
such that the Yoneda pairing
\begin{equation*}
  H^1_c(\mathfrak{Y}, \mathcal{F}) \times \Hom_\mathfrak{Y}(\mathcal{F}, \Omega^1_\mathfrak{Y}) \longrightarrow H^1_c(\mathfrak{Y}, \Omega^1_\mathfrak{Y})
\end{equation*}
composed with the trace map induces isomorphisms of topological vector spaces
\begin{equation*}
  \Hom_L^{cont}(H^1_c(\mathfrak{Y}, \mathcal{F}),L) \xrightarrow{\cong} \Hom_\mathfrak{Y}(\mathcal{F}, \Omega^1_\mathfrak{Y}) \ and \ \Hom_L^{cont}(\Hom_\mathfrak{Y}(\mathcal{F}, \Omega^1_\mathfrak{Y}),L) \xrightarrow{\cong} H^1_c(\mathfrak{Y}, \mathcal{F})
\end{equation*}
which are natural in $\mathcal{F}$.
\end{itemize}
\end{theorem}
\begin{proof}
See \cite[Thm.\ 7.1]{Bey}  \footnote{  This references depends on the results in the article \cite{aBey}, which unfortunately  contains the following gaps. Firstly, in the proof of Lemma 4.2.2.  he quotes a result of Bosch concerning the connectedness of formal fibers without verifying the required assumptions. This is repaired by \cite[Thm.\ 22/Cor.\ 23]{Mal}. Secondly,   Beyer claims implicitly and without proof,  that {\it special affinoid wide-open spaces} are {\it affinoid wide-open spaces}
in the sense of  Definition 4.1.1 and Remark 4.1.2 in (loc.\ cit.). This crucial ingredient has now been shown  explicitly in \cite[\S 2.5]{Mal}. } and \cite[Thm.\ 4.21]{Chi} (as well as \cite[Prop.\ 3.6]{vdP}).
\end{proof}

If we specialize the above assertion to the case of the structure sheaf $\mathcal{F} = \mathcal{O}_\mathfrak{Y}$ then we have $\Hom_\mathfrak{Y}(\mathcal{O}_\mathfrak{Y}, \Omega^1_\mathfrak{Y}) = \Omega^1_\mathfrak{Y}(\mathfrak{Y})$ for trivial reasons. On the other hand the relative cohomology sequence implies that
\begin{equation}\label{f:Robba}
  H^1_c(\mathfrak{Y}, \mathcal{O}_\mathfrak{Y}) = \mathcal{R}_L(\mathfrak{Y}) / \mathcal{O}_\mathfrak{Y}(\mathfrak{Y}) \ .
\end{equation}
Hence we have the following consequence of Serre duality.

\begin{corollary}\label{Serre2}
  Serre duality gives rise to an isomorphism of topological vector spaces
\begin{equation*}
  \Hom_L^{cont}(\mathcal{R}_L(\mathfrak{Y}) / \mathcal{O}_\mathfrak{Y}(\mathfrak{Y}),L) \xrightarrow{\;\cong\;} \Omega^1_\mathfrak{Y}(\mathfrak{Y}) \ .
\end{equation*}
\end{corollary}

\begin{lemma}\label{Serre-funct}
  Let $\alpha : \mathfrak{Y} \longrightarrow \mathfrak{Y}'$ be a finite etale morphism of 1-dimensional smooth Stein spaces over $L$. We then have, for any coherent sheaf $\mathcal{F}$ on $\mathfrak{Y}'$, the commutative diagram of Serre duality pairings
\begin{equation*}
  \xymatrix{
    H^1_c(\mathfrak{Y}, \alpha^*\mathcal{F})  & \times & \Hom_\mathfrak{Y}(\alpha^*\mathcal{F}, \Omega^1_\mathfrak{Y}) \ar[d] \ar[r]^{} & H^1_c(\mathfrak{Y}, \Omega^1_\mathfrak{Y})  \ar[d] \ar[r]^-{t_\mathfrak{Y}} & L \ar@{=}[d] \\
    H^1_c(\mathfrak{Y}', \alpha_*\alpha^*\mathcal{F}) \ar[u]^{\cong} & \times & \Hom_{\mathfrak{Y}'}(\alpha_*\alpha^*\mathcal{F}, \Omega^1_{\mathfrak{Y}'}) \ar[d] \ar[r]^{} & H^1_c(\mathfrak{Y}', \Omega^1_{\mathfrak{Y}'}) \ar@{=}[d] \ar[r]^-{t_{\mathfrak{Y}'}} & L \ar@{=}[d] \\
    H^1_c(\mathfrak{Y}', \mathcal{F}) \ar[u] & \times & \Hom_{\mathfrak{Y}'}(\mathcal{F}, \Omega^1_{\mathfrak{Y}'})  \ar[r]^{} & H^1_c(\mathfrak{Y}', \Omega^1_{\mathfrak{Y}'}) \ar[r]^-{t_{\mathfrak{Y}'}} & L .  }
\end{equation*}
\end{lemma}
\begin{proof}
The vertical arrows in the lower part of the diagram are induced by the adjunction homomorphism $\mathcal{F} \rightarrow \alpha_*\alpha^* \mathcal{F}$. It is commutative by the naturality of Serre duality in the coherent sheaf.

For the upper part we consider more generally a coherent sheaf $\mathcal{G}$ on $\mathfrak{Y}$ and check the commutativity of the diagram
\begin{equation*}
  \xymatrix{
    H^1_c(\mathfrak{Y}, \mathcal{G}) & \times & \Hom_\mathfrak{Y}(\mathcal{G}, \Omega^1_\mathfrak{Y}) \ar[d]^{f \mapsto \alpha_*(f)} \ar[r]^{} & H^1_c(\mathfrak{Y}, \Omega^1_\mathfrak{Y})  \ar[r]^-{t_\mathfrak{Y}} & L \ar@{=}[dd] \\
    H^1_c(\mathfrak{Y}', \alpha_*\mathcal{G}) \ar[u]^{\cong} & \times & \Hom_{\mathfrak{Y}'}(\alpha_*\mathcal{G}, \alpha_*\Omega^1_\mathfrak{Y})  \ar[d]  \ar[r]^{} & H^1_c(\mathfrak{Y}', \alpha_*\Omega^1_\mathfrak{Y})  \ar[d]  \ar[u]^{\cong} &  \\
     H^1_c(\mathfrak{Y}', \alpha_*\mathcal{G})  \ar@{=}[u] & \times & \Hom_{\mathfrak{Y}'}(\alpha_*\mathcal{G}, \Omega^1_{\mathfrak{Y}'}) \ar[r]^{} & H^1_c(\mathfrak{Y}', \Omega^1_{\mathfrak{Y}'}) \ar[r]^-{t_{\mathfrak{Y}'}} & L .    }
\end{equation*}
Here the second and third lower vertical arrows are induced by the relative trace map $t_\alpha : \alpha_*\Omega^1_\mathfrak{Y} \rightarrow \Omega^1_{\mathfrak{Y}'}$ (see below). The commutativity of the Yoneda pairings (before applying the horizontal trace maps) is a trivial consequence of functoriality properties. That the first and third upper vertical arrows are isomorphisms follows from the fact that for a finite morphism the functor $\alpha_*$ is exact on quasi-coherent sheaves. This reduces us to showing that the diagram
\begin{equation}\label{f:diagramMilan}
  \xymatrix@R=0.5cm{
                &         H^1_c(\mathfrak{Y}, \Omega^1_\mathfrak{Y}) \ar[dr]^-{t_\mathfrak{Y}}   &  \\
  H^1_c(\mathfrak{Y}', \alpha_*\Omega^1_\mathfrak{Y}) \ar[ur]^{\cong} \ar[dr]_-{H^1_c(\mathfrak{Y}', t_\alpha )}   && L              \\
                & H^1_c(\mathfrak{Y}', \Omega^1_{\mathfrak{Y}'}) \ar[ur]_-{t_{\mathfrak{Y}'}}   &                      }
\end{equation}
is commutative.

For the convenience of the reader we briefly explain the definition of the relative trace map $t_\alpha$. But first we need to recall that any coherent $\alpha_* \mathcal{O}_\mathfrak{Y}$-module $\mathcal{M}$ can  naturally be viewed (\cite[Prop.\ I.9.2.5]{EGA}) as a coherent $\mathcal{O}_\mathfrak{Y}$-module $\widetilde{\mathcal{M}}$ such that $\alpha_* \widetilde{\mathcal{M}} = \mathcal{M}$ (for any open affinoid subdomain $\mathfrak{V} \subseteq \mathfrak{Y}'$ one has $\widetilde{\mathcal{M}}(\alpha^{-1}(\mathfrak{V})) = \mathcal{M}(\mathfrak{V})$). Since $\alpha$ is etale we have (\cite[Thm.\ 25.1]{Mat}) that
\begin{equation*}
  (\alpha_* \mathcal{O}_\mathfrak{Y} \otimes_{\mathcal{O}_{\mathfrak{Y}'}} \Omega^1_{\mathfrak{Y}'})^\sim \xrightarrow{\cong} \Omega^1_\mathfrak{Y} \ .
\end{equation*}
Since $\alpha$ is finite flat the natural map
\begin{equation*}
  \underline{\Hom}_{\mathfrak{Y}'} (\alpha_* \mathcal{O}_\mathfrak{Y}, \mathcal{O}_{\mathfrak{Y}'}) \otimes_{\mathcal{O}_{\mathfrak{Y}'}} \Omega^1_{\mathfrak{Y}'}  \xrightarrow{\; \cong\;}
  \underline{\Hom}_{\mathfrak{Y}'} (\alpha_* \mathcal{O}_\mathfrak{Y}, \Omega^1_{\mathfrak{Y}'})
\end{equation*}
is an isomorphism. Finally, since $\alpha$ is finite etale the usual trace pairing is nondegenerate and induces an isomorphism \footnote{Compare https://stacks.math.columbia.edu/tag/0BVH}
\begin{equation*}
   \underline{\Hom}_{\mathfrak{Y}'} (\alpha_* \mathcal{O}_\mathfrak{Y}, \mathcal{O}_{\mathfrak{Y}'}) \xrightarrow{\cong} \alpha_* \mathcal{O}_\mathfrak{Y} \ .
\end{equation*}
The relative trace map is now defined to be the composite map
\begin{multline*}
  \alpha_* \Omega^1_\mathfrak{Y} \cong \alpha_* (\alpha_* \mathcal{O}_\mathfrak{Y} \otimes_{\mathcal{O}_{\mathfrak{Y}'}} \Omega^1_{\mathfrak{Y}'})^\sim \cong \alpha_* ( \underline{\Hom}_{\mathfrak{Y}'} (\alpha_* \mathcal{O}_\mathfrak{Y}, \mathcal{O}_{\mathfrak{Y}'}) \otimes_{\mathcal{O}_{\mathfrak{Y}'}} \Omega^1_{\mathfrak{Y}'})^\sim \cong \\
   \alpha_* \underline{\Hom}_{\mathfrak{Y}'} (\alpha_* \mathcal{O}_\mathfrak{Y}, \Omega^1_{\mathfrak{Y}'})^\sim =
  \underline{\Hom}_{\mathfrak{Y}'} (\alpha_* \mathcal{O}_\mathfrak{Y}, \Omega^1_{\mathfrak{Y}'}) \xrightarrow{\; f \mapsto f(1) \;} \Omega^1_{\mathfrak{Y}'} \ .
\end{multline*}
The commutativity of \eqref{f:diagramMilan}   is shown in detail in \cite{Mal} and should also be consequence of \cite[Prop.\ 6.5.1 (2)]{AL} upon showing that their general construction boils down to the above description of the relative trace map.
\end{proof}

We make the last lemma more explicit for the structure sheaf. Let $\rho : \mathfrak{Y} \rightarrow \mathfrak{Z}$ be a finite, faithfully flat, and etale morphism of $1$-dimensional smooth Stein spaces over $L$ such that $\mathcal{O}_{\mathfrak{Y}}(\mathfrak{Y})$ is finitely generated free as an $\mathcal{O}_{\mathfrak{Z}}(\mathfrak{Z})$-module. Fix a basis $f_1, \ldots, f_h \in \mathcal{O}_{\mathfrak{Y}}(\mathfrak{Y})$. Going through the proof of Lemma \ref{Serre-funct} one checks that on global sections the relative trace map is given by
\begin{align}\label{f:reltrace}
  t_\rho : \Omega^1_{\mathfrak{Y}}(\mathfrak{Y}) = \mathcal{O}_{\mathfrak{Y}}(\mathfrak{Y}) \otimes_{\mathcal{O}_{\mathfrak{Z}}(\mathfrak{Z})} \Omega^1_{\mathfrak{Z}}(\mathfrak{Z}) & \longrightarrow \Omega^1_{\mathfrak{Z}}(\mathfrak{Z}) \\
  \omega = \sum_{i=1}^h f_i \otimes \omega_i & \longmapsto \sum_{i=1}^h trace_{\mathcal{O}_{\mathfrak{Y}}(\mathfrak{Y})/\mathcal{O}_{\mathfrak{Z}}(\mathfrak{Z})} (f_i) \omega_i \ .   \nonumber
\end{align}
Hence we have the commutative diagram of duality pairings
\begin{equation}\label{f:Serre2-funct}
  \xymatrix{
    \Hom_L^{cont}(\mathcal{R}_L(\mathfrak{Y}) / \mathcal{O}_{\mathfrak{Y}}(\mathfrak{Y}),L) \ar[d]^{\Hom(\rho^*,L)} \ar[r]^-{\cong} & \Omega^1_{\mathfrak{Y}}(\mathfrak{Y}) \ar[d]^{t_\rho} \\
     \Hom_L^{cont}(\mathcal{R}_L(\mathfrak{Z}) / \mathcal{O}_{\mathfrak{Z}}(\mathfrak{Z}),L)  \ar[r]^-{\cong} & \Omega^1_{\mathfrak{Z}}(\mathfrak{Z})  . }
\end{equation}
It remains to explicitly compute the relative trace map in the cases of interest to us. But first we observe that, by the explanation at the end of section 2.3 in \cite{BSX}, the sheaf of differentials $\Omega^1_{\mathfrak{Y}}$ on a smooth $1$-dimensional Stein group variety $\mathfrak{Y}$ is a free $\mathcal{O}_{\mathfrak{Y}}$-module. Furthermore, if $\mathfrak{Y}$ is one of our character varieties, say of the group $G$, then by the construction before Def. \ 1.27 in \cite{BSX} we have the embedding
\begin{align*}
  L = \Lie(G) & \longrightarrow \mathcal{O}_{\mathfrak{Y}}(\mathfrak{Y}) \\
  \mathfrak{x} & \longmapsto [\chi \mapsto d\chi(\mathfrak{x})]
\end{align*}
and the function $\log_{\mathfrak{Y}}$ defined as the image of $1 \in L = \Lie(G)$.

\begin{remark}\label{rem:logXlogLT}
The function $\log_{\mathfrak{X}}$ corresponds under the LT-isomorphism $\kappa$ to the function $\Omega \log_{LT}$ by the commutative diagram after \cite[Lemma 3.4]{ST2}. In particular, $d\log_{\mathfrak{X}}$ corresponds to $\Omega d\log_{LT}$ and $\partial_{inv}^\mathfrak{X}$ to $\frac{1}{\Omega}\partial_{inv}$, where $df=\partial_{inv}^\mathfrak{X}d\log_{\mathfrak{X}}$ defines the invariant derivation on $\mathcal{O}_K(\mathfrak{X})$ similarly as for $\partial_{inv}$ in \eqref{f:inv}.
\end{remark}

\begin{proposition}\phantomsection\label{explicit-reltrace}
\begin{enumerate}
  \item[1.] For $\pi_L^* : \mathfrak{X} \rightarrow \mathfrak{X}$ we have $\Omega^1_{\mathfrak{X}}(\mathfrak{X}) = \mathcal{O}_{\mathfrak{X}}(\mathfrak{X}) d\log_{\mathfrak{X}}$ and
\begin{equation*}
  t_{\pi_L^*}(f d\log_{\mathfrak{X}}) = \frac{q}{\pi_L} \psi_L^{\mathfrak{X}}(f) d\log_{\mathfrak{X}} \ .
\end{equation*}
  \item[2.] For $n \geq n_0$ and $\ell_n^* : \mathfrak{X} \xrightarrow{\cong} \mathfrak{X}^\times_n$ we have $\Omega^1_{\mathfrak{X}^\times_n}(\mathfrak{X}^\times_n) = \mathcal{O}_{\mathfrak{X}^\times_n}(\mathfrak{X}^\times_n) d\log_{\mathfrak{X}^\times_n}$ and
\begin{equation*}
  t_{\ell_n^*}(f d\log_{\mathfrak{X}}) = \pi_L^n (\ell_n^*)_*(f) d\log_{\mathfrak{X}^\times_n} \ .
\end{equation*}
  \item[3.] For $n \geq m \geq 1$ and $\rho_{m,n} : \mathfrak{X}^\times_m \rightarrow \mathfrak{X}^\times_n$ we have $\Omega^1_{\mathfrak{X}^\times_m}(\mathfrak{X}^\times_m) = \mathcal{O}_{\mathfrak{X}^\times_m}(\mathfrak{X}^\times_m) d\log_{\mathfrak{X}^\times_m}$ and
\begin{equation*}
  t_{\rho_{m,n}}(f d\log_{\mathfrak{X}^\times_m}) = q^{n-m} f_1 d\log_{\mathfrak{X}^\times_n} \quad\text{if $f = \sum_{i=1}^h \ev_{u_i} \rho_{m,n}^*(f_i)$} \ .
\end{equation*}
  where $u_1 = 1, u_2, \ldots , u_h \in U_m$ are representatives for the cosets of $U_n$ in $U_m$ (with $h := q^{n-m}$).
  \item[4.] For $n \geq 1$ and $\rho_n : \mathfrak{X}^\times \rightarrow \mathfrak{X}^\times_n$ we have $\Omega^1_{\mathfrak{X}^\times}(\mathfrak{X}^\times) = \mathcal{O}_{\mathfrak{X}^\times}(\mathfrak{X}^\times) d\log_{\mathfrak{X}^\times}$ and
\begin{equation*}
  t_{\rho_n}(f d\log_{\mathfrak{X}^\times}) = (q-1)q^{n-1} f_1 d\log_{\mathfrak{X}^\times_n} \quad\text{if $f = \sum_{i=1}^h \ev_{u_i} \rho_n^*(f_i)$} \ .
\end{equation*}
  where $u_1 = 1, u_2, \ldots , u_h \in U_n$ are representatives for the cosets of $U_n$ in $o_L^\times$ (with $h := (q-1)q^{n-1}$).
  \item[5.] For the multiplication $\mu_\chi : \mathfrak{X}^\times \xrightarrow{\cong} \mathfrak{X}^\times$ by a fixed point $\chi \in \mathfrak{X}^\times(L)$ we have
\begin{equation*}
  t_{\mu_\chi}(f d\log_{\mathfrak{X}^\times}) = \mu_{\chi*}(f) d\log_{\mathfrak{X}^\times} = \mu_{\chi^{-1}}^*(f) d\log_{\mathfrak{X}^\times} \ .
\end{equation*}
\end{enumerate}
\end{proposition}
\begin{proof}
All subsequent computations start, of course, from the formula \eqref{f:reltrace} for the relative trace map.

1. The assumptions are satisfied by Lemmas \ref{phi-free} and \ref{lem:piL-finite}. As explained at the end of section 2.3 in \cite{BSX}, the sheaf of differentials $\Omega^1_{\mathfrak{X}}$ on $\mathfrak{X}$ is a free $\mathcal{O}_{\mathfrak{X}}$-module of rank one with basis the global differential $d\log_{\mathfrak{X}}$. By \cite[Lem.\ 1.28.ii]{BSX} we have $\varphi_L(\log_{\mathfrak{X}}) = \pi_L \log_{\mathfrak{X}}$. The formula for $t_{\pi_L^*}$ now follows from Remark \ref{rem:trace}.

2. The assumptions are trivially satisfied. The map $d\ell_n : L = \Lie(U_n) \rightarrow L = \Lie(o_L)$ is multiplication by $\pi_L^{-n}$. It follows that the isomorphism $(\ell_n^*)^* : \mathcal{O}_{\mathfrak{X}^\times_n}(\mathfrak{X}^\times_n) \rightarrow \mathcal{O}_{\mathfrak{X}}(\mathfrak{X})$ sends $\log_{\mathfrak{X}^\times_n}$ to $\pi_L^{-n} \log_{\mathfrak{X}}$. This implies the assertions.

3. The assumptions are satisfied by Remark \ref{rem:crossed-product}. The inclusion $U_n \hookrightarrow U_m$ of an open subgroup induces the identity map on the Lie algebras. It follows that $\rho_{m,n}^*(\log_{\mathfrak{X}^\times_n}) = \log_{\mathfrak{X}^\times_m}$. Since $\rho_{m,n}$ is etale we first may apply this with some $n \geq n_0$ and, using 2., deduce that $\Omega^1_{\mathfrak{X}^\times_m}(\mathfrak{X}^\times_m) = \mathcal{O}_{\mathfrak{X}^\times_m}(\mathfrak{X}^\times_m) d\log_{\mathfrak{X}^\times_m}$ for any $m \geq 1$. The formula for $t_{\rho_{m,n}}$ follows by the same argument as in the proof of Remark \ref{rem:trace}.

4. The argument is the same as the one for 3.

5. The assumptions are trivially satisfied. Using that $d\chi(1) = \frac{d}{dt}\chi(\exp(t))_{|t=0}$ we check that $\log_{\mathfrak{X}^\times}(\chi_1\chi_2) = \log_{\mathfrak{X}^\times}(\chi_1) + \log_{\mathfrak{X}^\times}(\chi_2)$ holds true. It follows that $\mu_\chi^* (d\log_{\mathfrak{X}^\times}) = d\log_{\mathfrak{X}^\times}$ and hence the formula for $t_{\mu_\chi}$.
\end{proof}

We briefly remark on the case where our Stein space is the open unit disk $\mathbf{B}$ around zero. Then $\mathcal{R}_L(\mathbf{B})$ is the usual Robba ring of all Laurent series $f(Z) = \sum_{i \in \mathbb{Z}} c_i Z^i$ with coefficients $c_i \in L$ which converge in some annulus near $1$. Analogously to \eqref{f:Robba} we have
\begin{equation*}
  H^1_c(\mathbf{B},\Omega^1_\mathbf{B}) = \mathcal{R}_L(\mathbf{B}) dZ / \mathcal{O}_\mathbf{B}(\mathbf{B})dZ
\end{equation*}
and the trace map sends $\sum_{i \in \mathbb{Z}} c_i Z^i dZ$ to its residue which is the coefficient $c_{-1}$ (\cite[\S3.1]{Bey}).

\subsubsection{Duality for boundary sections}\label{subsec:boundary}

First we recall another functoriality property of Serre duality.

\begin{proposition}\label{open-im}
  Let $j : \mathfrak{Y}_0 \rightarrow \mathfrak{Y}$ be an open immersion of 1-(equi)dimensional smooth Stein spaces over $L$, and let $\mathcal{F}$ be a coherent sheaf on $\mathfrak{Y}$. Then the diagram
\begin{equation*}
  \xymatrix{
    H^1_c(\mathfrak{Y}, \mathcal{F}) & \times & \Hom_\mathfrak{Y}(\mathcal{F}, \Omega^1_\mathfrak{Y}) \ar[d]^{\mathrm{res}} \ar[r]^{} & H^1_c(\mathfrak{Y}, \Omega^1_\mathfrak{Y})  \ar[r]^-{t_\mathfrak{Y}} & L \ar@{=}[d] \\
    H^1_c(\mathfrak{Y}_0, \mathcal{F}) \ar[u]^{j_!} & \times & \Hom_{\mathfrak{Y}_0}(\mathcal{F}, \Omega^1_{\mathfrak{Y}_0})    \ar[r]^{} & H^1_c(\mathfrak{Y}_0, \Omega^1_{\mathfrak{Y}_0})   \ar[u]^{j_!}  \ar[r]^-{t_{\mathfrak{Y}_0}}   &  L
  }
\end{equation*}
is commutative.
\end{proposition}
\begin{proof}
The commutativity of the Yoneda pairing (before applying trace maps) is immediate from the functoriality of the cohomology with compact support in the coefficient sheaf. The assertion that $t_\mathfrak{Y} \circ j_! = t_{\mathfrak{Y}_0}$ holds true is shown in \cite[Thm.\ 3.7]{vdP}.
\end{proof}

In order to combine the above functoriality property with Lemma \ref{projlim} in the case of the structure sheaf $\mathcal{F} = \mathcal{O}_\mathfrak{Y}$ we first recall the setting of that lemma.
\begin{itemize}
  \item[1.] $\mathfrak{Y} = \bigcup_n \mathfrak{U}_n$ is a Stein covering of the Stein space $\mathfrak{Y}$ such that the $\mathfrak{Y}_n = \mathfrak{Y} \setminus \mathfrak{U}_n$ are Stein spaces as well. In particular, $\mathcal{R}_L(\mathfrak{Y}) = \varinjlim_n \mathcal{O}_{\mathfrak{Y}}(\mathfrak{Y}_n)$ with the locally convex inductive limit topology.
  \item[2.] The restriction maps $\mathcal{O}_\mathfrak{Y}(\mathfrak{Y}) \rightarrow \mathcal{O}_\mathfrak{Y}(\mathfrak{Y} \setminus \mathfrak{U})$ are injective for any $\mathfrak{U} \in \Aff(\mathfrak{Y})$.
  \item[3.] The inductive system of Fr\'echet spaces $\Omega^1_{\mathfrak{Y}}(\mathfrak{Y}_1) \rightarrow \ldots \rightarrow \Omega^1_{\mathfrak{Y}}(\mathfrak{Y}_n) \rightarrow \Omega^1_{\mathfrak{Y}}(\mathfrak{Y}_{n+1}) \rightarrow \ldots$ is regular (\cite[Def.\ 11.1.3(ii)]{PGS}). By \cite[Thm.\ 11.2.4(ii)]{PGS} the locally convex inductive limit
\begin{equation*}
  \Omega^1_{\mathcal{R}_L(\mathfrak{Y})} := \varinjlim_n \Omega^1_{\mathfrak{Y}}(\mathfrak{Y}_n) = \varinjlim_{\mathfrak{U} \in \Aff(\mathfrak{Y})} \Omega^1_{\mathfrak{Y}}(\mathfrak{Y} \setminus \mathfrak{U})
\end{equation*}
       is a locally convex Hausdorff space.
\end{itemize}

\begin{proposition}\label{Robba-duality}
  In the above setting 1.-3. we have a natural topological isomorphism
\begin{equation*}
  \Hom_L^{cont}(\Omega^1_{\mathcal{R}_L(\mathfrak{Y})},L) \cong \big( \varprojlim_n \varinjlim_{\mathfrak{U} \in \Aff(\mathfrak{Y}_n)} \mathcal{O}_{\mathfrak{Y}}(\mathfrak{Y} \setminus \mathfrak{U}) \big) /\mathcal{O}_{\mathfrak{Y}}(\mathfrak{Y})
\end{equation*}
(where the left hand side is equipped with the strong topology of bounded convergence).
\end{proposition}
\begin{proof}
The asserted isomorphism is the composite of the isomorphisms
\begin{align*}
  \Hom_L^{cont}(\Omega^1_{\mathcal{R}_L(\mathfrak{Y})},L) & = \Hom_L^{cont}(\varinjlim_n \Omega^1_{\mathfrak{Y}}(\mathfrak{Y}_n),L) \xrightarrow{\cong} \varprojlim_n \Hom_L^{cont}(\Omega^1_{\mathfrak{Y}}(\mathfrak{Y}_n),L) \\
     & = \varprojlim_n H^1_c(\mathfrak{Y}_n,\mathcal{O}_{\mathfrak{Y}})   \\
     & = \big( \varprojlim_n \varinjlim_{\mathfrak{U} \in \Aff(\mathfrak{Y}_n)} \mathcal{O}_{\mathfrak{Y}}(\mathfrak{Y} \setminus \mathfrak{U}) \big) /\mathcal{O}_{\mathfrak{Y}}(\mathfrak{Y}) \ .
\end{align*}
The isomorphism in the first line comes from \cite[Thm.\ 11.1.13]{PGS}. The equality in the second, resp.\ third, line is a consequence of Thm.\ \ref{Serre}  and Prop.\ \ref{open-im}, resp.\ Lemmata \ref{projlim} and \ref{projlim-topol}.
\end{proof}

We now evaluate this latter result in the same concrete cases as in section \ref{subsec:compact-supp}.\\

\noindent
{\bf The open unit disk}\\

First of all, the sheaf of differentials $\Omega^1_{\mathbf{B}}$ on the open unit disk $\mathbf{B}$ is a free $\mathcal{O}_{\mathbf{B}}$-module of rank one. Hence, by choosing, for example the global differential $dZ$ for a coordinate function $Z$, as a basis we obtain a topological isomorphism $\mathcal{R}_L(\mathbf{B}) \cong \Omega^1_{\mathcal{R}_L(\mathbf{B})}$ as $\mathcal{R}_L(\mathbf{B})$-modules. The regularity assumption in 3. above therefore is reduced to the corresponding property for $\mathcal{R}_L(\mathbf{B})$, which is established in the proof of \cite[Prop.\ 2.6.i]{BSX}. Hence Prop.\ \ref{Robba-duality} is available. By combining its assertion with \eqref{f:projlim=Robba} we obtain a natural topological isomorphism
\begin{equation}\label{f:Robba-selfdual}
  \Hom_L^{cont}(\mathcal{R}_L(\mathbf{B}),L) \cong \Hom_L^{cont}(\Omega^1_{\mathcal{R}_L(\mathbf{B})},L) \cong \mathcal{R}_L(\mathbf{B}) \ .
\end{equation}
This shows that $\mathcal{R}_L(\mathbf{B})$ is topologically selfdual. By going through the definitions and using the explicit description of the trace map in this case as the usual residue map (end of section \ref{subsec:duality}) one checks that this selfduality comes from the pairing
\begin{align*}
  \mathcal{R}_L(\mathbf{B}) \times \mathcal{R}_L(\mathbf{B}) & \longrightarrow L \\
        (f_1(Z),f_2(Z)) & \longmapsto \text{residue of $f_1(Z) f_2(Z) dZ$}.
\end{align*}
This latter form of the result was known (\cite{CR}, \cite{MS}) before Serre duality in rigid analysis was established. In this paper it is more natural to use the selfduality given by the pairing
\begin{align*}
  < \ ,\ >_{\mathbf{B}} \, : \mathcal{R}_L(\mathbf{B}) \times \mathcal{R}_L(\mathbf{B}) & \longrightarrow L \\
        (f_1,f_2) & \longmapsto \text{residue of $f_1 f_2 d\log_{LT}$}.
\end{align*}
We will denote by $\res_\mathbf{B} : \Omega^1_{\mathcal{R}_L(\mathbf{B})} \rightarrow L$ the linear form which corresponds to $1 \in \mathcal{R}_L(\mathbf{B})$ under the second isomorphism in \eqref{f:Robba-selfdual}.
\\

\noindent
{\bf The character variety $\mathfrak{X}$}\\

We recall that the sheaf of differentials $\Omega^1_{\mathfrak{X}}$ on $\mathfrak{X}$ is a free $\mathcal{O}_{\mathfrak{X}}$-module of rank one with basis the global differential $d\log_{\mathfrak{X}}$. Hence again we have a topological isomorphism $\mathcal{R}_L(\mathfrak{X}) \cong \Omega^1_{\mathcal{R}_L(\mathfrak{X})}$ as $\mathcal{R}_L(\mathfrak{X})$-modules. The regularity assumption in 3. above therefore holds by \cite[Prop.\ 2.6.i]{BSX}. Hence Prop.\ \ref{Robba-duality} is available. By combining its assertion with \eqref{f:projlim=RobbaX} we obtain a natural topological isomorphism
\begin{equation}\label{f:RobbaX-selfdual}
  \Hom_L^{cont}(\mathcal{R}_L(\mathfrak{X}),L) \cong \Hom_L^{cont}(\Omega^1_{\mathcal{R}_L(\mathfrak{X})},L) \cong \mathcal{R}_L(\mathfrak{X}) \ .
\end{equation}
This shows that $\mathcal{R}_L(\mathfrak{X})$ is topologically selfdual. Let $\res_{\mathfrak{X}} : \Omega^1_{\mathcal{R}_L(\mathfrak{X})} \rightarrow L$ be the linear form which corresponds to $1 \in \mathcal{R}_L(\mathfrak{X})$ under the above isomorphism. Then, as a consequence of the naturality of the Yoneda pairing, this selfduality comes from the pairing
\begin{align}\label{f:pairing-RobbaX}
  < \ ,\ >_{\mathfrak{X}} \, : \mathcal{R}_L(\mathfrak{X}) \times \mathcal{R}_L(\mathfrak{X}) & \longrightarrow L \\ \notag
        (f_1,f_2) & \longmapsto \res_{\mathfrak{X}} (f_1 f_2 d\log_{\mathfrak{X}}) \ .
\end{align}
\footnote{WARNING: If $L = \mathbb{Q}_p$ then $\mathfrak{X} = \mathbf{B}_1 \cong \mathbf{B}$ with the latter isomorphism given by $z \mapsto z - 1$; but the selfdualities \eqref{f:Robba-selfdual} and \eqref{f:RobbaX-selfdual} do not correspond to each other since in this case $d\log_{\mathfrak{X}}$ corresponds to $d\log(1+Z) = \frac{1}{1+Z} dZ$.}

Next we consider $\mathfrak{X}^\times_n$ for some $n \geq n_0$, where $n_0 \geq 1$ is the integer from section \ref{sec:mult-Robba}. We then have the isomorphism of Stein group varieties $\ell_n^* : \mathfrak{X} \xrightarrow{\cong} \mathfrak{X}^\times_n$. Hence all we have established for $\mathfrak{X}$ holds true correspondingly for $\mathfrak{X}^\times_n$. In particular, we have a natural topological isomorphism
\begin{equation}\label{f:RobbaX-selfdual-mult}
  \Hom_L^{cont}(\mathcal{R}_L(\mathfrak{X}^\times_n),L) \cong \Hom_L^{cont}(\Omega^1_{\mathcal{R}_L(\mathfrak{X}^\times_n)},L) \cong \mathcal{R}_L(\mathfrak{X}^\times_n) \ .
\end{equation}
Let $\res_{\mathfrak{X}^\times_n} : \Omega^1_{\mathcal{R}_L(\mathfrak{X}^\times_n)} \rightarrow L$ be the linear form which corresponds to $1 \in \mathcal{R}_L(\mathfrak{X}^\times_n)$ under this isomorphism. We obtain that $\mathcal{R}_L(\mathfrak{X}^\times_n)$ is topologically selfdual w.r.t.\ the pairing
\begin{align*}
  < \ ,\ >_{\mathfrak{X}^\times_n} \, : \mathcal{R}_L(\mathfrak{X}^\times_n) \times \mathcal{R}_L(\mathfrak{X}^\times_n) & \longrightarrow L \\
        (f_1,f_2) & \longmapsto \res_{\mathfrak{X}_n^\times} (f_1 f_2 d\log_{\mathfrak{X}^\times_n}) \ .
\end{align*}

It follows from \cite[Thm.\ 3.7]{vdP}  that the diagram
\begin{equation*}
  \xymatrix@R=0.5cm{
  \Omega^1_{\mathcal{R}_L(\mathfrak{X})} \ar[dd]_{(\ell_n^*)_*}^{\cong} \ar[dr]^{\res_{\mathfrak{X}}}             \\
                & L          \\
  \Omega^1_{\mathcal{R}_L(\mathfrak{X}^\times_n)}    \ar[ur]_{\res_{\mathfrak{X}^\times_n}}             }
\end{equation*}
is commutative. But in the proof of Prop.\ \ref{explicit-reltrace}.2 we have seen that $(\ell_n^*)_*(\log_{\mathfrak{X}}) = \pi_L^n \log_{\mathfrak{X}^\times_n}$. Therefore the diagram of pairings
\begin{equation}\label{f:XmultXn}
  \xymatrix{
    \mathcal{R}_L(\mathfrak{X})  & \times & \mathcal{R}_L(\mathfrak{X}) \ar[d]_{\pi_L^n(\ell_n^*)_*}^{\cong} \ar[r]^-{< \ ,\ >_{\mathfrak{X}}} & L \ar@{=}[d] \\
    \mathcal{R}_L(\mathfrak{X}^\times_n) \ar[u]^{(\ell_n^*)^*}_{\cong} & \times & \mathcal{R}_L(\mathfrak{X}^\times_n)  \ar[r]^-{< \ ,\ >_{\mathfrak{X}^\times_n}}   &  L
  }
\end{equation}
is commutative. Alternatively we could have used the following observation.

\begin{remark}\label{reltrace-open}
   Let $\rho : \mathfrak{Y} \rightarrow \mathfrak{Z}$ be one of the morphisms in Prop.\ \ref{explicit-reltrace}. Then Prop.\ \ref{crossed-product} applies, and it follows that, for any admissible open subset $\mathfrak{U} \subseteq \mathfrak{Z}$ which is Stein, the relative trace map $t_{\rho|\rho^{-1}(\mathfrak{U})}$ is given by the same formula as for $t_\rho$.
\end{remark}

In the case of the morphisms $\pi_L^*$ and $\rho_{m,n}$ for $n \geq m \geq n_0$ this immediately leads to the following equalities of pairings.

\begin{lemma}\phantomsection\label{RobbaX-pairing-phi-psi}
\begin{itemize}
  \item[i.] $< \varphi_L(f_1),f_2 >_{\mathfrak{X}} \; = \, \frac{q}{\pi_L} < f_1, \psi_L^{\mathfrak{X}}(f_2) >_{\mathfrak{X}}$ for any $f_1, f_2 \in \mathcal{R}_L(\mathfrak{X})$.
  \item [ii.] Let $n \geq m \geq n_0$ and let $u_1 = 1, u_2, \ldots , u_h \in U_m$ be representatives for the cosets of $U_n$ in $U_m$ (with $h := q^{n-m}$); for any $f' \in \mathcal{R}_L(\mathfrak{X}^\times_n)$ and any $f = \sum_{i=1}^h \ev_{u_i} \rho_{m,n}^*(f_i) \in \mathcal{R}_L(\mathfrak{X})^\times_m$ (cf.\ Remark \ref{rem:crossed-product}) we have
\begin{equation*}
    < \rho_{m,n}^*(f'),f >_{\mathfrak{X}^\times_m} \, = \, q^{n-m} < f', f_1 >_{\mathfrak{X}^\times_n} \ .
\end{equation*}
\end{itemize}
\end{lemma}

\noindent
{\bf The multiplicative character variety $\mathfrak{X}^\times$}\\

We fix an $n \geq n_0$ for the moment as well as representatives $u_1 = 1, u_2, \ldots , u_h \in o_L^\times$ for the cosets of $U_n$ in $o_L^\times$ (with $h := (q-1)q^{n-1}$). Recalling from Lemma \ref{lem:crossed-product-R}.ii that
\begin{equation*}
  \mathcal{R}_L(\mathfrak{X}^\times) = \mathbb{Z}[o_L^\times] \otimes_{\mathbb{Z}[U_n]}  \rho_n^*(\mathcal{R}_L(\mathfrak{X}^\times_n))
\end{equation*}
we may write any $f \in \mathcal{R}_L(\mathfrak{X}^\times)$ as $f = \sum_{i=1}^h \ev_{u_i} \rho_n^*(f_i)$ with uniquely determined $f_i \in \mathcal{R}_L(\mathfrak{X}^\times_n)$. We now define
\begin{equation}\label{f:res-mult-charvar}
  \res_{\mathfrak{X}^\times} : \Omega^1_{\mathcal{R}_L(\mathfrak{X}^\times)} \longrightarrow L  \quad\text{by $\res_{\mathfrak{X}^\times}(f d\log_{\mathfrak{X}^\times}) := (q-1)q^{n-1} \res_{\mathfrak{X}^\times_n}(f_1 d\log_{\mathfrak{X}^\times_n})$}
\end{equation}
and then the pairing
\begin{align}\label{f:pairing-mult-charvar}
  < \ ,\ >_{\mathfrak{X}^\times} :  \mathcal{R}_L(\mathfrak{X}^\times) \times \mathcal{R}_L(\mathfrak{X}^\times) & \longrightarrow L  \\\notag
     (f_1,f_2) & \longmapsto  \res_{\mathfrak{X}^\times}(f_1 f_2) \ .
\end{align}
These definitions are obviously independent of the choice of the representatives $u_i$. Moreover, due to Lemma \ref{RobbaX-pairing-phi-psi}.ii they are independent of the choice of $n$ as well and we have
 \begin{equation}\label{f:RobbaX-pairing-Ln}
    < \rho_{n}^*(f'),f >_{\mathfrak{X}^\times} \, = \,(q-1) q^{n-1} < f', f_1 >_{\mathfrak{X}^\times_n} \
\end{equation}
 for any $f' \in \mathcal{R}_L(\mathfrak{X}^\times_n)$ and any $f = \sum_{i=1}^{(q-1) q^{n-1}} \ev_{u_i} \rho_{n}^*(f_i) \in \mathcal{R}_L(\mathfrak{X})^\times$ where $u_i   \in o_L^\times$ runs through representatives for the cosets of $U_n$ in $o_L^\times.$ The topological selfduality of $\mathcal{R}_L(\mathfrak{X}^\times_n)$ easily implies that this pairing makes $\mathcal{R}_L(\mathfrak{X}^\times)$ topological selfdual.

\begin{lemma}\phantomsection\label{RobbaXtimes-pairing-twist}
The twist morphism $\mu_\chi : \mathfrak{X}^\times \rightarrow \mathfrak{X}^\times$, for any $\chi \in \mathfrak{X}^\times (L)$, satisfies
\begin{equation*}
  < \mu_{\chi}^*(f_1), \mu_{\chi}^*(f_2) >_{\mathfrak{X}^\times} \; = \; < f_1, f_2 >_{\mathfrak{X}^\times}  \qquad\text{ for any $f_1, f_2 \in \mathcal{R}_L(\mathfrak{X}^\times)$}.
\end{equation*}
\end{lemma}
\begin{proof}
The assertion immediately reduces to checking the equality $\res_{\mathfrak{X}^\times} \circ \mu_{\chi}^* = \res_{\mathfrak{X}^\times}$. Obviously there are twist morphisms on $\mathfrak{X}^\times_n$ as well. One easily checks that
\begin{equation*}
  \mu_\chi^* \circ \rho_n^* = \rho_n^* \circ \mu_{\chi | U_n}^*
\end{equation*}
and that
\begin{equation*}
  \mu_\chi^*(\ev_u) = \chi(u) \ev_u  \qquad\text{for any $u \in o_L^\times$}.
\end{equation*}
Using Lemma \ref{lem:crossed-product-R}.ii we write an $f \in \mathcal{R}_L(\mathfrak{X}^\times)$ as $f = \sum_{i=1}^h \ev_{u_i} \rho_n^*(f_i)$ and compute
\begin{equation*}
  \mu_\chi^*(f) = \sum_{i=1}^h \mu_\chi^*(\ev_{u_i}) \mu_\chi^*(\rho_n^*(f)) = \sum_{i=1}^h \ev_{u_i} \rho_n^*(\chi(u_i) \mu_{\chi | U_n}^*(f_i)) \ .
\end{equation*}
This shows that $\mu_\chi^*(f)_1 = \mu_{\chi | U_n}^*(f_1)$. This further reduces us to showing that $\res_{\mathfrak{X}^\times_n} \circ \mu_{\chi | U_n}^* = \res_{\mathfrak{X}^\times_n}$. But this follows from \cite[Thm.\ 3.7]{vdP} or, alternatively, from a version of Prop.\ \ref{explicit-reltrace}.5 for $\mathcal{R}_L(\mathfrak{X}^\times_n)$.
\end{proof}

Of course, everything in this entire section \ref{sec:duality} remains valid over any complete extension field $K$ of $L$ contained in $\mathbb{C}_p$.
Moreover, our constructions above are compatible under (complete) base change: Let $\mathfrak{Y}_K$ denote the base change of  $\mathfrak{Y}$ over $L$ to $K$ (and similarly for affinoids). Then we obtain a commutative diagram
\begin{equation}\label{f:basechangeSerre}
   \xymatrix{
     \mathcal{R}_L(\mathfrak{Y}) \ar[d]_{ }   & \times   & \Omega^1_{\mathcal{R}_L(\mathfrak{Y})}  \ar[d]_{ } \ar[r]^{ } & L \ar@{^(->}[d]^{ } \\
     \mathcal{R}_K(\mathfrak{Y}_K) & \times   &  \Omega^1_{\mathcal{R}_K(\mathfrak{Y}_K)}  \ar[r]^{ } & K.  }
\end{equation}
Indeed, it is shown in \cite{Mal} that Serre-duality is compatible with base change  in the sense that there is the following commutative diagram    for any $n$, in which the horizontal lines are the Serre-dualities over $L$ and $K,$ respectively:
\begin{equation*}
  \xymatrix{
    H^1_c(\mathfrak{Y}_n, \mathcal{O}_\mathfrak{Y}) \ar[d]_{ }   & \times    & \Omega^1_\mathfrak{Y}(\mathfrak{Y}_n) \ar[d]_{ } \ar[r]^{ } & L \ar@{^(->}[d]^{ } \\
    H^1_c(\mathfrak{Y}_{n,K}, \mathcal{O}_{\mathfrak{Y}_K})   & \times & \Omega^1_{\mathfrak{Y}_K}(\mathfrak{Y}_{n,K})  \ar[r]^{ } & K   }
\end{equation*}
Hence, taking limits as in the proof of Proposition \ref{Robba-duality} the claim follows upon observing that also the relative cohomology sequence \eqref{f:relative} is compatible with base change. By inserting $1\in\mathcal{R}_L(\mathfrak{Y})\subseteq \mathcal{R}_K(\mathfrak{Y}_K) $ into the pairings of \eqref{f:basechangeSerre} we see that in any example discussed above the residue maps $\mathrm{res}_\mathfrak{Y}$ are compatible under base change as well as the pairings $< \ ,\ >_{\mathfrak{Y}}$. Since  $\mathcal{R}_K(\mathfrak{Y}_K)\cong K\hat{\otimes}_{L,\iota}\mathcal{R}_L(\mathfrak{Y})$ (and hence  $\Omega^1_{\mathcal{R}_K(\mathfrak{Y}_K)}\cong K\hat{\otimes}_{L,\iota}\Omega^1_{\mathcal{R}_L(\mathfrak{Y})}$) by \cite[Cor.\ 2.8]{BSX} with respect to the (completed) inductive tensor product, we see that
\begin{equation}\label{f:basechangeres}
 \mathrm{res}_{\mathfrak{Y}_K}= K\hat{\otimes}_{L,\iota}\mathrm{res}_\mathfrak{Y}
\end{equation}
and that we have the commutative diagram
\begin{equation*}
  \xymatrix{
    K\hat{\otimes}_{L,\iota} \Hom_L^{cont}(\mathcal{R}_L(\mathfrak{Y}),L)  \ar[d]_{\cong}^{\id_K\hat{\otimes}_{L,\iota}\text{duality}_{\mathfrak{Y}_L}} \ar[r]^{} &  \Hom_K^{cont}(K\hat{\otimes}_{L,\iota} \mathcal{R}_L(\mathfrak{Y}),K)  \ar[r]^-{\cong} & \Hom_K^{cont}(\mathcal{R}_K(\mathfrak{Y}_K),K) \ar[d]_{\cong}^{\text{duality}_{\mathfrak{Y}_K}} \\
    K\hat{\otimes}_{L,\iota}\mathcal{R}_L(\mathfrak{Y}) \ar[rr]^{\cong} & & \mathcal{R}_K(\mathfrak{Y}_K) .}
\end{equation*}
\newpage

\subsection{\texorpdfstring{$(\varphi_L,\Gamma_L)$}{(phi,Gamma)}-modules
}\label{sec:modulesRobba}

As before we let $L \subseteq K \subseteq \mathbb{C}_p$ be a complete intermediate field, and we denote by $o_K$ its ring of integers.

\subsubsection{The usual  Robba ring}\label{sec:usualRobba}

In sections \ref{subsec:duality} and \ref{subsec:boundary} we already had introduced the usual Robba ring $\mathcal{R} = \mathcal{R}_K = \mathcal{R}_K(\mathbf{B})$ of the Stein space $\mathbf{B}_{/K}$ in connection with Serre duality. We briefly review its construction in more detail. The ring of $K$-valued global holomorphic functions $\mathcal{O}_K(\mathbf{B})$\footnote{In the notation from \cite[\S 1.2]{Co2} this is the ring $\mathcal{R}^+. $ }  on $\mathbf{B}$ is the Fr\'echet algebra of all power series in the variable $Z$ with coefficients in $K$ which converge on the open unit disk $\mathbf{B}(\mathbb{C}_p)$. The Fr\'echet topology on ${\mathcal{O}_K(\mathbf{B})}$ is given by the family of norms
\begin{equation*}
  |\sum_{i \geq 0} c_i Z^i|_r := \max_i |c_i| r^i    \qquad\text{for}\ 0 < r < 1 \ .
\end{equation*}
In the commutative integral domain ${\mathcal{O}_K(\mathbf{B})}$ we have the multiplicative subset $Z^{\mathbb{N}} = \{Z^j : j \in \mathbb{N}\}$, so that we may form the corresponding localization ${\mathcal{O}_K(\mathbf{B})}_{Z^{\mathbb{N}}}$. Each norm $|\ |_r$ extends to this localization ${\mathcal{O}_K(\mathbf{B})}_{Z^{\mathbb{N}}}$ by setting $|\sum_{i \gg -\infty} c_iZ^i|_r := \max_i |c_i| r^i$.

The Robba ring $\mathcal{R} \supseteq {\mathcal{O}_K(\mathbf{B})}$ is constructed as follows. For any $s > 0$, resp.\ any $0 < r \leq s$, in $p^{\mathbb{Q}}$ let $\mathbf{B}_{[0,s]}$, resp.\ $\mathbf{B}_{[r,s]}$, denote the affinoid disk around $0$ of radius $s$, resp.\ the affinoid annulus of inner radius $r$ and outer radius $s$, over $K$. For $I=[0,s]$ or  $ [r,s]$ we denote by
\[
  \cR^I := \cR_K^I(\mathbf{B}) := \mathcal{O}_K(\mathbf{B}_{I})
\]
the affinoid $K$-algebra of $\mathbf{B}_{I}$. The Fr\'echet algebra $\mathcal{R}^{[r,1)} := \varprojlim_{r < s < 1} \cR^{[r,s]}$ is the algebra of (infinite) Laurent series in the variable $Z$ with coefficients in $K$ which converge on the half-open annulus $\mathbf{B}_{[r,1)} := \bigcup_{r < s < 1} \mathbf{B}_{[r,s]}$. The Banach algebra $\mathcal{R}^{[0,s]}$ is the completion of ${\mathcal{O}_K(\mathbf{B})}$ with respect to the norm $|\ |_s$. The Banach algebra $\mathcal{R}^{[r,s]}$ is the completion of ${\mathcal{O}_K(\mathbf{B})}_{Z^{\mathbb{N}}}$ with respect to the norm $|\ |_{r,s} := \max(|\ |_r,|\ |_s)$. It follows that the Fr\'echet algebra $\mathcal{R}^{[r,1)}$ is the completion of ${\mathcal{O}_K(\mathbf{B})}_{Z^{\mathbb{N}}}$ in the locally convex topology defined by the family of norms $(|\ |_{r,s})_{r < s < 1}$.  Finally, the Robba ring is $\mathcal{R} = \bigcup_{0 < r < 1} \mathcal{R}^{[r,1)}$.\Footnote{Let $K$ be a complete field which contains $\Omega$ and $L$. We recall the definition of the Robba ring  and  various related rings in the {\bf additive} notation: We set $r_n:=\frac{1}{e(q-1)q^{n-1}}=v_p(t_{0,n})$ for a generator $t_0=(t_{0,n})$ of the Tate module of LT.

Using the notation from \cite[\S 1.2]{Co2} we denote by $\cR^{[r,s]}=\cE^{[r,s]}\subseteq K[[Z,Z^{-1}]]$ for $0<r\leq s$ the analytic functions on the annulus $C_{[r,s]}=\{\frac{1}{p^s}\leq |z|_p\leq \frac{1}{p^r}\}.$ It is a principal ideal domain and a Banach algebra with valuation $v_{[r,s]}$ defined by
\[v_{[r,s]}=\inf_{r \leq v_p(z)\leq s} v_p(f(z)) =\min\{\inf_{k\in\z}(v_p(a_k)+rk),\inf_{k\in\z}(v_p(a_k)+sk)\},\]
if $f=\sum_{k\in\z} a_kZ^k.$ We set $|f |_{[r,s]}:=p^{-v_{[r,s]}(f)}.$

 We also allow $r=\infty$ with $s<\infty$ obtaining the full disks.
Then, for $0\leq r<s,$ the ring $\cE^{]r,s]}$,  defined as the $\varprojlim_{r<t\leq s} \cE^{[t,s]},$ describes the analytic functions on $C_{]r,s]}=\{\frac{1}{p^s}\leq |z|_p<\frac{1}{p^r}\}.$ We set ${\cE}^{s}=\cE^{{]0,s]}}$ and define the Robba ring $\cR$ as $\varinjlim_{s>0}\cE^s.$
Finally we set $\cR^+:=\cR\cap K[[Z]]$}

\begin{remark}\label{rem:denseB}
  ${\mathcal{O}_K(\mathbf{B})}_{Z^{\mathbb{N}}}$ is dense in $\mathcal{R}_K(\mathbf{B})$.
\end{remark}

Let $p^{-\frac{d}{(q-1)e}} < r \leq s < 1$. Then we have a surjective map
\begin{align}
 \mathbf{B}_{[r,s]}&\to \mathbf{B}_{[r^q,s^q]}\\
 z&\mapsto [\pi_L](z)\notag
\end{align}
according to \cite[proof of Lem.\ 2.6]{FX}\footnote{The proof there is only written for the special Lubin-Tate group, but generalizes easily to the general case by using the fact that $[\pi_L]= X^q+\pi_LXf(X)$ with $f(X)\in o_L[[X]]^\times.$} It induces a ring homomorphism
\begin{align}\label{f:isometry}
\varphi_L^{[r^q,s^q]}: \cR^{[r^q,s^q]}\to\cR^{[r,s]}
\end{align}
which is isometric with respect to the supremum norms, i.e., $|\varphi_L^{[r^q,s^q]}(f)|_{[r,s]}=|f|_{[r^q,s^q]}$ for any $f\in \cR^{[r^q,s^q]}.$ In particular, by taking first inverse and then direct limits we obtain a continuous ring homomorphism $\varphi_L:\cR\to \cR.$ We shall often omit the interval in $\varphi_L^{[r,s]}$ and just write $\varphi_L$.

Similarly, we obtain a continuous $\Gamma_L$-action on $\cR$: According to (loc.\ cit.) we have a bijective map
\begin{align}
 \mathbf{B}_{[r,s]}&\to \mathbf{B}_{[r,s]}\\
 z&\mapsto [\chi_{LT}(\gamma)](z)\notag
\end{align}
for any $\gamma\in \Gamma_L,$ whence we obtain an isometric isomorphism
\begin{align}
\gamma: \cR^{[r,s]}\to\cR^{[r,s]}
\end{align}
 with respect to the supremum norms, i.e., $|\gamma(f)|_{[r,s]}=|f|_{[r,s]}$ for any $f\in \cR^{[r,s]}.$

 Finally, we extend the operator $\psi_L$ to $\cR:$ For $y \in \ker([\pi_L])$ we have the isomorphism
 \begin{align}
 \mathbf{B}_{[r,s]}&\to \mathbf{B}_{[r,s]}\\
 z&\mapsto z+_{LT} y\notag
 \end{align}
of affinoid varieties, because $|z+_{LT}y|=|z+y|=|z|$. The latter equality comes from $|z| \geq r > p^{-\frac{d}{(q-1)e}} = q^{-\frac{1}{q-1}} = |y|$ for $y \neq 0$. Setting $tr(f):=\sum_{y\in \ker([\pi_L])} f(z+_{LT} y)$ we obtain a norm decreasing linear map $tr: \cR^{[r,s]}\to \cR^{[r,s]}$.
We claim that the image of $tr$ is contained in the (closed) image of the isometry $\varphi_L^{[r^q,s^q]},$ whence there is a norm decreasing map
\[\psi_{Col}:\cR^{[r,s]}\to \cR^{[r^q,s^q]},\] such that $\varphi_L\circ \psi_{Col} =tr.$ Indeed, by continuity it suffices to show that  $tr(Z^i)$ belongs to the image for any $i\in \z.$ For $i\geq 0,$ Coleman has shown that $tr(Z^i)=\varphi_L(\psi_{Col}(Z^i))$ with $\psi_{Col}(Z^i)\in o_L[[Z]]\subseteq \cR^{[r^q,s^q]}$, see  \cite[\S 2]{SV15}. For $i<0$, we calculate
\begin{align*}
\varphi_L( Z^i\psi_{Col}([\pi_L](Z)^{-i}Z^i))&=\varphi_L(Z^i)\left(\sum_{y\in \ker([\pi_L])}[\pi_L](Z)^{-i}Z^i\right)(Z+_{LT}y)\\
&=\varphi_L(Z^i)\left(\sum_{y\in \ker([\pi_L])}[\pi_L]((Z+_{LT}y))^{-i}(Z+_{LT}y)^i\right) \\
&=\varphi_L(Z^i)\left(\sum_{y\in \ker([\pi_L])}\varphi_L(Z)^{-i}(Z+_{LT}y)^i\right) \\
&=\sum_{y\in \ker([\pi_L])} (Z+_{LT}y)^i=tr(Z^i),
\end{align*}
whence the claim follows. We put $\psi_L^{[r,s]}:=\frac{1}{\pi_L}\psi_{Col}:\cR^{[r,s]}\to \cR^{[r^q,s^q]}$ which induces the   continuous operator $\psi_L:\cR\to \cR$ by taking first inverse limits and then direct limits. By definition of $tr$ the operators  $\psi_L^{[r,s]}$ and hence $\psi_L$ satisfy the projection formula.
We shall often omit the interval in $\psi_L^{[r,s]}$ and just write $\psi_L$.

As in section \ref{subsec:LT} we fix a generator $\eta'$ of the dual Tate module $T_\pi'$  and denote by $\Omega$ the corresponding period. For the rest of this subsection we \textbf{assume} in addition that $K$ contains $\Omega$. Following Colmez in the notation we introduce the power series $\eta(a,Z) := \exp(a\Omega\log_{LT}(Z)) \in o_K[[Z]]$ for $a\in o_L$. As noted in section \ref{subsec:LT} the power series $\eta(a,Z)$ is nothing else than the image under the LT-isomorphism $\kappa^*$ of the holomorphic function $\ev_a \in \mathcal{O}_K(\mathfrak{X})$. Generalizing the equality \eqref{f:Odecomp} we have the following decompositions of Banach spaces
\begin{equation}
 \label{f:cRdecompI}\cR^{[r,s]} =\bigoplus_{a\in o_L/\pi_L^n} \varphi_L^n(\cR^{[r^q,s^q]})\eta(a,Z)
\end{equation}
and hence
\begin{equation}
 \label{f:cRdecomp}\cR=\bigoplus_{a\in o_L/\pi_L^n} \varphi_L^n(\cR)\eta(a,Z)
\end{equation}
of   LF-spaces using the formula
\begin{equation}
\label{f:projectorformula}
r={ (\frac{\pi_L}{q})^n }\sum_a  \varphi_L^n\psi_L^n\left(\eta(-a,Z)r\right)\eta(a,Z).
\end{equation}
This can easily be reduced by induction on $n$ to the case $n=1.$ Using the definition of $tr$ and the orthogonality relations for the characters $\kappa_y$ for $y\in \ker([\pi_L]),$  the formula follows and, moreover, defines a continuous inverse to the continuous map
\begin{align*}
\z[o_L]\otimes_{\mathbb{Z}[\pi_L o_L]} \cR^{[r^q,s^q]} & \xrightarrow{\;\cong\;} \cR^{[r,s]}\\
a \otimes f & \longmapsto a \varphi_L(f) .
\end{align*}
Inductively, we obtain canonical isomorphisms
\begin{align}\label{f:cRtensorI}
\z[o_L]\otimes_{\mathbb{Z}[\pi_L^n o_L]} \cR^{[r^{q^n},s^{q^n}]} & \xrightarrow{\:\cong\:} \cR^{[r,s]}\\
a \otimes f & \longmapsto a\varphi_L^n(f).\notag
\end{align}

Moreover, immediately from the definitions we have
\begin{align}
\label{f:etaphi}\varphi_L(\eta(a,Z))&=\eta(\pi_L a,Z),\\
\label{f:etasigma}\sigma(\eta(a,Z))&=\eta(\chi_{LT}(\sigma)a,Z) \ \text{for $\sigma \in \Gamma_L$,} \\
\label{f:etapsi}\psi_L(\eta(a,Z)) & = \frac{q}{\pi_L} \eta(\frac{a}{\pi_L},Z) \ \text{for $a\in \pi_Lo_L$ and $= 0$ otherwise.}
\end{align}

\begin{remark}\label{rem:traceB}
  We have $\psi_L = \frac{1}{\pi_L} \varphi_L^{-1} \circ trace_{\cR/\varphi_L(\cR)}$.
\end{remark}
\begin{proof}
Both maps $tr$ and $trace_{\cR/\varphi_L(\cR)}$ are easily seen to be $\varphi_L(\cR)$-linear and to be multiplication by $q$ on $\varphi_L(\cR)$. Hence, by \eqref{f:cRdecomp}, it suffices to compare their values on the elements $\eta(a,Z)$. By \eqref{f:etapsi} we have
\begin{equation*}
  tr(\eta(a,Z)) =
  \begin{cases}
   q \varphi_L(\eta(\pi_L^{-1} a,Z)) & \text{if $a \in \pi_L o_L$},  \\
   0 &  \text{otherwise} .
  \end{cases}
\end{equation*}
On the other hand a computation as in the proof of Remark \ref{rem:trace} shows that $trace_{\cR/\varphi_L(\cR)}$ has exactly the same values. But $\psi_L = \frac{1}{\pi_L} \varphi_L^{-1} \circ tr$.
\end{proof}

For uniformity of notation we put $\psi_L^\mathbf{B} := \frac{\pi_L}{q} \psi_L = \frac{1}{q} \varphi_L^{-1} \circ trace_{\cR/\varphi_L(\cR)}$.

\subsubsection{The LT-isomorpism, part 2}\label{subsec:LT-Robba}

We \textbf{assume} throughout this subsection that $\Omega$ is contained in $K$. The map $\kappa : \mathbf{B} \xrightarrow{\cong} \mathfrak{X}$ being an isomorphism of rigid varieties it preserves the systems of affinoid subdomains on both sides. Hence the LT-isomorphism \eqref{f:LTiso} extends to a topological isomorphism
\begin{equation}\label{f:LTisoR}
  \kappa^* : \cR_K(\mathfrak{X}) \xrightarrow{\;\cong\;} \cR_K(\mathbf{B}) \ .
\end{equation}
In order to have a uniform notation, we usually write from now on
\begin{equation*}
  \cR_K^I(\mathfrak{X}) := \mathcal{O}_K(\kappa(\mathbf{B}_I))
\end{equation*}
for any closed interval $I \subseteq (0,1)$ so that we have the isomorphism of Banach algebras
\begin{equation}\label{f:LTisoRI}
  \kappa^* : \cR_K^I(\mathfrak{X}) \xrightarrow{\;\cong\;} \cR^I_K(\mathbf{B}) \ .
\end{equation}
We warn the reader that only for specific closed intervals $I$ there is another closed interval $I'$ given by a complicated but explicit rule such that $\kappa(\mathbf{B}_I) = \mathfrak{X}_{I'}$. The precise statement can be worked out from \cite[Prop.\ 1.20]{BSX}.

In the following we list a few compatibilities under this extended LT-isomorphism.

First of all under this isomorphism the $\Gamma_L \cong o_L^\times$-action and the maps $\varphi_L$ on both sides correspond to each other (cf.\ \cite[\S2.2]{BSX}). Then it follows from Remarks \ref{rem:trace} \ref{rem:traceB} that the operators $\psi_L^\mathfrak{X}$ (defined at the end of section \ref{sec:add-Robba}) and $\psi_L^\mathbf{B}$ (defined in previous section) also correspond under $\kappa^*$.

Secondly, as a consequence of \cite[Thm.\ 3.7]{vdP} we have the commutative diagram
\begin{equation}\label{f:vdP}
  \xymatrix@R=0.5cm{
  \Omega^1_{\mathcal{R}_K(\mathfrak{X})} \ar[dd]_{\kappa^*}^{\cong} \ar[dr]^{\res_{\mathfrak{X}}}             \\
                & K \ .          \\
  \Omega^1_{\mathcal{R}_K(\mathbf{B})}    \ar[ur]_{\res_{\mathbf{B}}}             }
\end{equation}
This combined with Remark \ref{rem:logXlogLT}  implies the explicit formula
\begin{equation}\label{f:res-XB}
  \res_{\mathfrak{X}}(f d\log_\mathfrak{X}) = \Omega  \res_\mathbf{B} (\kappa^*(f) g_{LT} dZ) \ .
\end{equation}

\subsubsection{\texorpdfstring{$\varphi_L$}{phi}-modules}

Let $\mathfrak{Y}$ be either $\mathfrak{X}$ or $\mathbf{B}$ and $\cR := \cR_{K}(\mathfrak{Y})$. Henceforth we will use the operator $\psi_L:=\frac{q}{\pi_L}\psi_L^\mathfrak{Y}$ on $\cR$. We also put
\begin{equation*}
  q_\mathfrak{Y} :=
  \begin{cases}
  p & \text{if $\mathfrak{Y} = \mathfrak{X}$},   \\
  q & \text{if $\mathfrak{Y} = \mathbf{B}$}.
  \end{cases}
\end{equation*}

\begin{definition}\label{def:phimodule}
  A $\varphi_L$-module $M$ over $\cR$ is a finitely generated free $\cR$-module $M$ equipped with a semilinear endomorphism $\varphi_M$ such that the $\cR$-linear map
\begin{align*}
  \varphi_M^{lin} :\cR \otimes_{\cR,\varphi_L} M & \xrightarrow{\;\cong\;} M\\
                                                      f\otimes m & \longmapsto f\varphi_M(m)
\end{align*}
is bijective.
\end{definition}

Technically important is the following fact, which for $\mathfrak{X}$ is part of the proof of \cite[Prop.\ 2.24]{BSX}. The proof for $\mathbf{B}$ is entirely analogous. It allows to extend the above maps and decompositions from the previous sections to $\varphi_L$-modules. For $r > 0$ we introduce the intervals
\begin{equation*}
  I(r,\mathfrak{Y}) :=
  \begin{cases}
  (r,1) & \text{if $\mathfrak{Y} = \mathfrak{X}$},  \\
  [r,1) & \text{if $\mathfrak{Y} = \mathbf{B}$}.
  \end{cases}
\end{equation*}

\begin{proposition}\label{descentphi}
  Let $M$ be a $\varphi_L$-module $M$ over $\mathcal{R}$. There exists a radius
\begin{equation*}
  r_0 \geq
  \begin{cases}
  p^{-\frac{dp}{p-1}} & \text{if $\mathfrak{Y} = \mathfrak{X}$},   \\
  p^{-\frac{dq}{(q-1)e}} & \text{if $\mathfrak{Y} = \mathbf{B}$}
  \end{cases}
\end{equation*}
  and a finitely generated free $\cO_K(\mathfrak{Y}_{I(r_0,\mathfrak{Y})})$-module $M_0$ equipped with
  a semilinear continuous homomorphism
\begin{equation*}
  \varphi_{M_0} : M_0 \longrightarrow \cO_K(\mathfrak{Y}_{I(r_0,\mathfrak{Y})^{1/q_\mathfrak{Y}}}) \otimes_{\cO_K(\mathfrak{Y}_{I(r_0,\mathfrak{Y})})} M_0
\end{equation*}
such that the induced $\cO_K(\mathfrak{Y}_{I(r_0,\mathfrak{Y})^{1/q_\mathfrak{Y}}})$-linear map
\begin{equation*}
  \varphi_{M_0}^{lin} : \cO_K(\mathfrak{Y}_{I(r_0,\mathfrak{Y})^{1/q_\mathfrak{Y}}}) \otimes_{\cO_K(\mathfrak{Y}_{I(r_0,\mathfrak{Y})}), \varphi_L} M_0 \xrightarrow{\;\cong\;} \cO_K(\mathfrak{Y}_{I(r_0,\mathfrak{Y})^{1/q_\mathfrak{Y}}}) \otimes_{\cO_K(\mathfrak{Y}_{I(r_0,\mathfrak{Y})})} M_0
\end{equation*}
is an isomorphism and such that
\begin{equation*}
  \mathcal{R} \otimes_{\cO_K(\mathfrak{Y}_{I(r_0,\mathfrak{Y})})} M_0 = M
\end{equation*}
with  $\varphi_L \otimes \varphi_{M_0}$ and $\varphi_M$ corresponding to each other.
\end{proposition}

The continuity condition for the $\varphi_{M_0}$, of course, refers to the product topology on $M_0 \cong (\cO_K(\mathfrak{Y}_{I(r_0,\mathfrak{Y})}))^d$.

In the following we fix a $\varphi_L$-module $M$ over $\mathcal{R}$ and a pair $(r_0,M_0)$ as in Prop.\ \ref{descentphi}. For any $r_0 \leq r' < 1$ and any closed interval $I = [r,s] \subseteq I(r',\mathfrak{Y})$ we then have the finitely generated free modules
\begin{equation*}
  M^{I(r',\mathfrak{Y})} := \cO_K(\mathfrak{Y}_{I(r',\mathfrak{Y})}) \otimes_{\cO_K(\mathfrak{Y}_{I(r_0,\mathfrak{Y})})} M_0 \quad \text{over}\ \cR^{[r,1)}
\end{equation*}
and
\begin{equation*}
   M^I := \cO_K(\mathfrak{Y}_I) \otimes_{\cO_K(\mathfrak{Y}_{I(r',\mathfrak{Y})})} M^{I(r',\mathfrak{Y})} \quad \text{over}\ \cO_K(\mathfrak{Y}_I) \ .
\end{equation*}
They satisfy
\begin{equation}\label{f:identities}
  M^{I(r',\mathfrak{Y})} = \varprojlim_{s > r} M^I \qquad\text{and}\qquad  M = \varinjlim_{r'} M^{I(r',\mathfrak{Y})} \ .
\end{equation}

We equip $M^I$ with the Banach norm $|-|_{M^I}$ given by the maximum norm with respect to any fixed basis (the induced topology  does not depend on the choice of basis) which is submultiplicative with respect to scalar multiplication and the norm $|-|_I$ on $\cO_K(\mathfrak{Y}_I)$.

Furthermore, base change with $\cO_K(\mathfrak{Y}_{I^{1/q_\mathfrak{Y}}})$ over $\cO_K(\mathfrak{Y}_{I(r_0,\mathfrak{Y})^{1/q_\mathfrak{Y}}})$ induces isomorphisms of Banach spaces
\begin{align*}
\varphi_{lin}^I = \cO_K(\mathfrak{Y}_{I^{1/q_\mathfrak{Y}}}) \otimes_{\cO_K(\mathfrak{Y}_{I(r_0,\mathfrak{Y})^{1/q_\mathfrak{Y}}})} \varphi_{M_0}^{lin} :  \cO_K(\mathfrak{Y}_{I^{1/q_\mathfrak{Y}}}) \otimes_{\cO_K(\mathfrak{Y}_I),\varphi_L} M^I  \xrightarrow{\;\cong\;} M^{I ^{1/q_\mathfrak{Y}}}
\end{align*}
and hence injective, continuous maps
\begin{align*}
\varphi^I:   M^I \xrightarrow{} M^{I^{1/q_\mathfrak{Y}}}
\end{align*}
by restriction.

Assuming that $I^{q_\mathfrak{Y}} \subseteq I(r',\mathfrak{Y})$ we define the additive, $K$-linear, continuous map $\psi^I:   M^I \rightarrow M^{I^{q_\mathfrak{Y}}}$ as the composite
\begin{align*}
   \psi^I : M^I \xrightarrow{(\varphi_{lin}^{I^{q_\mathfrak{Y}}})^{-1}} \cO_K(\mathfrak{Y}_I) \otimes_{\cO_K(\mathfrak{Y}_{I^{q_\mathfrak{Y}}}),\varphi_L} M^{I^{q_\mathfrak{Y}}} \to M^{I^{q_\mathfrak{Y}}},
\end{align*}
where the last map sends $f\otimes m$ to $\psi^I(f)m.$  By construction, it satisfies the projection formulas
\begin{equation}\label{f:projectionrs}
  \psi^I(\varphi^{I^{q_\mathfrak{Y}}}(f)m) = f \psi^I(m) \qquad \text{and}\qquad  \psi^I(g\varphi^{I^{q_\mathfrak{Y}}}(m')) = \psi^I(g) m' \ ,
\end{equation}
for any $f \in \cO_K(\mathfrak{Y}_{I^{q_\mathfrak{Y}}})$, $g \in \cO_K(\mathfrak{Y}_I)$ and $m \in  M^I $,  $m' \in  M^{I^{q_\mathfrak{Y}}}$ as well as the formula
\begin{equation*}
  \psi^I \circ \varphi^{I^{q_\mathfrak{Y}}} = \frac{q}{\pi_L} \cdot \id_{M^{I^{q_\mathfrak{Y}}}} \ .
\end{equation*}

Using Prop.\ \ref{phi-free-Robba} in case $\mathfrak{Y} = \mathfrak{X}$, resp.\ the decomposition \eqref{f:cRdecompI} in case $\mathfrak{Y} = \mathbf{B}$ (under the assumption that $\Omega$ is contained in $K$), combined with (iterates of) $\varphi_{lin}^I $ gives rise to decompositions
\begin{equation}\label{f:decompMI}
   M^{I^{\frac{1}{q_\mathfrak{Y}^n}}} =
   \begin{cases}
   \bigoplus_{a\in(o_L/\pi_L^n)} \  \ev_a \varphi^n_L(M^I)      & \text{if $\mathfrak{Y} = \mathfrak{X}$}, \\
   \bigoplus_{a\in(o_L/\pi_L^n)} \  \eta(a,Z)\varphi^n_L(M^I)  & \text{if $\mathfrak{Y} = \mathbf{B}$ and $\Omega \in K$}
   \end{cases}
\end{equation}
of Banach spaces and
\begin{equation}\label{f:decompM}
   M  =
   \begin{cases}
   \bigoplus_{a\in(o_L/\pi_L^n)} \ \ev_a \varphi^n_L(M)  & \text{if $\mathfrak{Y} = \mathfrak{X}$},    \\
   \bigoplus_{a\in(o_L/\pi_L^n)} \ \eta(a,Z)\varphi^n_L(M)  & \text{if $\mathfrak{Y} = \mathbf{B}$ and $\Omega \in K$}
   \end{cases}
   \end{equation}
of LF-spaces, again given by the formula
\begin{equation}\label{f:projectorformulaM}
   m =
   \begin{cases}
   (\frac{\pi_L}{q})^n \sum_a  \varphi_M\psi_M\left(\ev_{-a} m\right)\ev_a       &  \text{if $\mathfrak{Y} = \mathfrak{X}$},   \\
   (\frac{\pi_L}{q})^n \sum_a  \varphi_M\psi_M\left(\eta(-a,Z)m\right)\eta(a,Z)  &  \text{if $\mathfrak{Y} = \mathbf{B}$ and $\Omega \in K$}.
   \end{cases}
\end{equation}

\subsubsection{The Robba ring of a group}\label{sec:groupRobba}

Recall that $L_n = L(\ker([\pi_L^n]))$. We set
\[
\Gamma_n := G(L_\infty/L_n)=\ker\left(\Gamma_L\xrightarrow{\chi_{LT}} o_L^\times \to (o_L/\pi_L^n)^\times\right).
\]
Also recall from section \ref{sec:mult-Robba} the notation $U_n := 1 + \pi_L^n o_L$ for $n \geq 1$ and the isomorphisms $\log : U_n \xrightarrow{\cong} \pi_L^n o_L$ and $\ell_n = \pi_L^{-n} \log : U_n \xrightarrow{\cong} o_L$ for $n \geq n_0$, where $n_0 \geq 1$\Footnote{ $n_0>\frac{e}{p-1}.$} is  minimal among $n$ such that $\log:1+\pi_L^no_L\to \pi_L^no_L$  and $\exp:\pi_L^no_L\to 1+\pi_L^no_L$ are mutually inverse isomorphisms.

Obviously $\chi_{LT}$ restricts to isomorphisms $\Gamma_n \cong U_n$ for any $n\geq 1$. Consider the composed maps
\[
\hat{\ell} := \log \circ \chi_{LT} : \Gamma_{L}\to L   \qquad\text{and} \qquad \hat{\ell}_n := \ell_n \circ \chi_{LT} : \Gamma_n \xrightarrow{\cong} o_L \quad\text{for $n \geq n_0$}.
\]
The latter isomorphisms induce isomorphisms of Fr\'echet algebras $D(\Gamma_n,K) \xrightarrow{\;\cong\;} D(o_L,K)$.

Because of the isomorphisms $\Gamma_L \cong o_L^\times$ and $\Gamma_n \cong U_n$ the formalism of character varieties and corresponding Robba rings applies to the groups $\Gamma_L$ and $\Gamma_n$ as well giving us the corresponding character varieties $\mathfrak{X}_{\Gamma_L}$ and $\mathfrak{X}_{\Gamma_n}$, and the results of section \ref{sec:mult-Robba} transfer to this setting. To make a clear distinction we put $\cR_K(\Gamma_L) := \cR_K(\mathfrak{X}_{\Gamma_L})$ and $\cR_K(\Gamma_n) := \cR_K(\mathfrak{X}_{\Gamma_n})$ and call them the Robba rings of the groups $\Gamma_L$ and $\Gamma_n$. Clearly the Lubin-Tate character $\chi_{LT}$ induces topological ring isomorphisms
\begin{equation}\label{f:chi-iso}
  \cR_K(\Gamma_L) \xrightarrow{\cong} \mathcal{R}_K(\mathfrak{X}^\times) \qquad\text{and} \qquad  \mathcal{R}_K(\Gamma_n) \xrightarrow{\cong} \mathcal{R}_K(\mathfrak{X}_{_n}^\times) \ \text{for $n \geq 1$}.
\end{equation}
If $\Gamma$ denotes any of these groups then we will very often view, via the Fourier isomorphism, $K[\Gamma] \subseteq D(\Gamma,K)$ as subrings of $\cR_K(\Gamma)$. In particular we consider elements $\gamma \in \Gamma$ as elements of the Robba ring writing them in any of the forms $\gamma \corresponds \delta_{\gamma} \corresponds \ev_\gamma$.

Let $n \geq m \geq 1$. The inclusions $\iota_n : \Gamma_n \hookrightarrow \Gamma_L$ and $\iota_{n,m} : \Gamma_n \hookrightarrow \Gamma_m$ induce, by the transfer of the results in section \ref{sec:mult-Robba}, ring monomorphisms $\iota_{n*} : \cR_K(\Gamma_n) \hookrightarrow \cR_K(\Gamma_L)$ and $\iota_{n,m*} : \cR_K(\Gamma_n) \hookrightarrow \cR_K(\Gamma_m)$. More precisely we have (Lemma \ref{lem:crossed-product-R} and Remark \ref{rem:crossed-product}) topological ring isomorphisms
\begin{equation}\label{f:cross}
  \mathbb{Z}[\Gamma_L] \otimes_{\mathbb{Z}[\Gamma_n]} \cR_K(\Gamma_n) \xrightarrow{\cong} \cR_K(\Gamma_L)
\end{equation}
and
\begin{equation}\label{f:Rcrossproduct}
 \mathbb{Z}[\Gamma_m] \otimes_{\mathbb{Z}[\Gamma_n]} \cR_K(\Gamma_n) \xrightarrow{\cong} \cR_K(\Gamma_m) \ .
\end{equation}
Here the left hand sides are viewed as free $\cR_K(\Gamma_n)$-modules endowed with the product topology.

We also note that, for $ n \geq m \geq n_0$, the commutative diagram
\begin{equation*}
  \xymatrix{
   \Gamma_n \ar[d]_{\iota_{n,m}} \ar[r]^{\hat{\ell}_n} & o_L \ar[d]^{\pi_L^{n-m}} \\
   \Gamma_m \ar[r]^{\hat{\ell}_m} & o_L    }
\end{equation*}
induces the commutative diagrams
\begin{equation}\label{f:diagphi}
    \xymatrix{
     D(\Gamma_n,K) \ar[r]^{\hat{\ell}_{n*}} \ar@{^(->}[d]_{\iota_{n,m*}} & D(o_L,K) \ar[rr]^{\cong}_{Fourier} \ar[d]_{(\pi_L^{n-m})_*} && \mathcal{O}_K(\mathfrak{X}) \ar[d]_{\varphi_L^{n-m}}  \\
     D(\Gamma_m,K) \ar[r]^{\hat{\ell}_{m*}} & D(o_L,K) \ar[rr]^{\cong}_{Fourier}  &&  \mathcal{O}_K(\mathfrak{X})  }
\end{equation}
and
\begin{equation}\label{f:diagphiR}
   \xymatrix{
   \cR_K(\Gamma_n) \ar@{^(->}[d]_{\iota_{n,m*}} \ar[rr]^{\hat{\ell}_{n*}}_{\cong} && \cR_K(\mathfrak{X}) \ar@{^(->}[d]^{\varphi_L^{n-m}}    \\
   \cR_K(\Gamma_m) \ar[rr]^{\hat{\ell}_{m*}}_{\cong} && \cR_K(\mathfrak{X}) . }
\end{equation}

For the rest of this subsection we \textbf{assume} that $\Omega$ is contained in $K$. Let $n \geq n_0$.
We then have the isomorphisms of rigid varieties
\begin{equation*}
  \mathbf{B} \xrightarrow[\kappa]{\simeq} \mathfrak{X} \xrightarrow[\hat{\ell}_n^*]{\simeq} \mathfrak{X}_{\Gamma_n} \ .
\end{equation*}
For any closed interval $I \subseteq (0,1)$ we therefore have the affinoid subdomain $\hat{\ell}_n^* \circ \kappa (\mathbf{B}_I)$ in $\mathfrak{X}_{\Gamma_n}$ and we may introduce the Banach algebra $\cR^I_K(\Gamma_n) := \cO_K(\hat{\ell}_n^* \circ \kappa (\mathbf{B}_I))$.

By its very construction the diagram
\begin{equation}\label{f:diagphiRI}
   \xymatrix{
   \cR_K^{I^{q^{n-m}}}(\Gamma_n) \ar@{^(->}[d]_{\iota_{n,m*}} \ar[rr]^{\hat{\ell}_{n*}}_{\cong} && \cR_K^{I^{q^{n-m}}}(\mathfrak{X}) \ar@{^(->}[d]^{\varphi_L^{n-m}} \ar[r]^{\kappa^*}_{\cong} & \cR_K^{I^{q^{n-m}}}(\mathbf{B})  \ar@{^(->}[d]^{\varphi_L^{n-m}}   \\
   \cR_K^I(\Gamma_m) \ar[rr]^{\hat{\ell}_{m*}}_{\cong} && \cR_K^I(\mathfrak{X}) \ar[r]^{\kappa^*}_{\cong} & \cR_K^I(\mathbf{B}) , }
\end{equation}
for $n \geq m \geq n_0$, is commutative. Together with \eqref{f:cRtensorI} it implies the canonical isomorphism
\begin{equation}\label{f:RIcrossproduct}
  \mathbb{Z}[\Gamma_m] \otimes_{\mathbb{Z}[\Gamma_n]} \cR_K^{I^{q^{n-m}}}(\Gamma_n) \xrightarrow{\cong} \cR_K^I(\Gamma_m) .
\end{equation}

We will denote the composite of Fourier and LT-isomorphism by
\begin{equation*}
  \textfrak{L}: D(o_L,K) \xrightarrow{\;\cong\;} \mathcal{O}_K(\mathbb{\mathfrak{X}}) \xrightarrow{\;\cong\;}  \mathcal{O}_K(\mathbf{B}) \ .
\end{equation*}
Recall that $\mathcal{O}_K(\mathbf{B})$ is a space of certain power series in the variable $Z$. We put
\begin{equation*}
  X := \textfrak{L}^{-1}(Z) \in D(o_L,K)  \qquad\text{and}\qquad  Y_n := \hat{\ell}_{n*}^{-1}(X) \in D(\Gamma_n,K)  \ \text{for $n \geq n_0$}.
\end{equation*}
In this way we can express elements in these distribution algebras as power series in these variables. This will later on be an important technical tool for our proofs.

\Footnote{
\paragraph{The cyclotomic case}
We set $\Gamma_1=1+p\zp$ with topological generator $\gamma_1=1+p,$ we identify $\cR(\Gamma_1)$ with $\cR$ by $Z=\gamma_1-1$ and observe that then the differential $\frac{d\gamma_1}{\gamma_1}$ corresponds to $\frac{dZ}{Z+1}.$  Furthermore, corresponding to the twist with $\chi_{cyc}$ we have $Tw(Z)=\gamma_1Z +p$ (where $\gamma_1$ is now considered as coefficient!) and
\[Tw(Z)^{-1}=\sum_{j\leq -1} p^{-(j+1)}\gamma_1^{j} Z^j.\]
It follows that for $p^{-1}<r<s<1$ the element $Tw(Z)^{-1}$ belongs to $\cR^{[r,s]}$, whence we obtain a map
\[\cR^+_{Z^\mathbb{N}}\to\cR^{[r,s]},\]
which seems to be continuous \com{??} as   $|Z^{-1}|_{r,s}=|Tw(T)^{-1}|_{r,s}=r^{-1}$ and  $|Z|_{r,s}=|Tw(T)|_{r,s}=s$. Therefore we obtain a continuous map
\[Tw:\cR^{[r,s]}\to\cR^{[r,s]},\] which is a ring automorphism because $Tw^{-1}$ corresponding to the twist with $\chi_{cyc}^{-1}$ is its inverse. Taking projective limits with respect to $s$ and then inductive limits with respect to $r$ induces a ring automorphism $\cR(\Gamma_1)$, which is the desired twist operator. Setting $\omega_\ell:=\frac{d\gamma_1}{\ell(\gamma_1)\gamma_1}\in \Omega^1_{\cR(\Gamma_1)}$ we obtain a differential, which is invariant with respect to $Tw.$

\begin{lemma}For all $\omega\in \Omega^1_{\cR(\Gamma_1)}$ we have
\[Res_{\Gamma_1}(Tw(\omega))=Res_{\Gamma_1}(\omega).\]
\end{lemma}

\begin{proof} we calculate with respect to $\cR$ and $\omega= f dZ$, $f=\sum_{i\in\z} a_iZ^i\in \cR:$
\[Tw(\omega)=Tw(f)dTw(Z)=Tw(f) \gamma_1dZ\]
Since   the terms $Tw(Z)^i$ only contribute to the $Z^{-1}$ term for $i={-1}$,   the verification is reduced to the case $f=Z^{-1},$ in which case  one sees that
\[\frac{\gamma_1}{Tw(Z)}=\sum_{j\leq -1} p^{-(j+1)}\gamma_1^{j+1} Z^j \]
has residuum $1.$
\end{proof}

\begin{corollary} For all $\mu,\lambda\in \cR(\Gamma_1)$ we have
\[Res_{\Gamma_1}(\mu Tw(\lambda)\omega_\ell)=Res_{\Gamma_1}(Tw^{-1}(\mu) \lambda\omega_\ell).\]
\end{corollary}

\begin{proof}
Apply the lemma to $\omega:=Tw^{-1}(\mu) \lambda$ using that $Tw(\omega_\ell)=\omega_\ell.$
\end{proof}

Thus the pairing
\[<\;,\;>:\cR(\Gamma_1)\times \cR(\Gamma_1) \to K,\;\; (\mu,\lambda)\mapsto Res_{\Gamma_1}(\mu\lambda \omega_\ell)\]
satisfies  \eqref{cor:twist<}.
}

As an immediate consequence of Remark \ref{rem:denseB} we have the following.

\begin{remark}\phantomsection\label{rem:denseX}
\begin{itemize}
  \item[i.] $D(o_L,K)_{Z^\mathbb{N}}$ is dense in $\cR_K(\mathfrak{X})$.
  \item[ii.] $D(\Gamma_n,K)_{Y_n^\mathbb{N}}$ is dense in $\cR_K(\Gamma_n)$ for $n \geq n_0$.
\end{itemize}
\end{remark}
\Footnote{ Can we show this for arbitrary $K$? Runge!?}

\subsubsection{Locally \texorpdfstring{$\mathbb{Q}_p$}{Qp}-analytic versus locally \texorpdfstring{$L$}{L}-analytic distribution algebras.}


We fix a $\zp$-basis $h_1=1,\ldots, h_d$ of $o_L$ and set $b_i:=h_i-1$ and, for any multiindex $\mathbf{k}=(k_1,\ldots,k_d)\in \mathbb{N}_0^d,$ $\mathbf{b}^\alpha:=\prod_{i=1}^db_i^{\alpha_i}\in \zp[o_L].$ We write $D_{\qp}(G,K)$ for the algebra of $K$-valued locally $\qp$-analytic distributions on a $\qp$-Lie group $G$. Any $\lambda\in D_{\qp}(o_L,K)$ has a unique convergent expansion $\lambda=\sum_{\mathbf{k} \in \mathbb{N}_0^d} \alpha_{\mathbf{k}}\mathbf{b}^{\mathbf{k}}$ with $\alpha_{\mathbf{k}}\in K$ such that, for any $0<r<1,$ the set $\{\alpha_{\mathbf{k}}r^{{ \diamond}|\mathbf{k}|}_{\mathbf{k} \in \mathbb{N}_0^d}\}$ is bounded, { where $\diamond:=2$ if $p=2$ and $\diamond:=1$ otherwise}. The completion with respect to the norm
\[\|\lambda\|_{\mathbb{Q}_p,r}:=\sup_{\mathbf{k} \in \mathbb{N}_0^d} |\alpha_{\mathbf{k}}|r^{{\diamond}|\mathbf{k}|}\]
for $0<r<1$ is denoted by
\[
D_{\qp,r}(o_L,K)=\{\sum_{\mathbf{k} \in \mathbb{N}_0^d} \alpha_{\mathbf{k}}\mathbf{b}^{\mathbf{k}}|\alpha_{\mathbf{k}}\in K \mbox{ and } |\alpha_{\mathbf{k}}|r^{{\diamond}|\mathbf{k}|}\to 0 \mbox{ as } |\mathbf{k}|\to \infty \}.
\]
By \cite[Prop.\ 2.1]{Sc1} the group $o_L$ satisfies the hypothesis $(HYP)$ of \cite{ST} with $p$-valuation $\omega$ satisfying $\omega(h_i)=\diamond.$ Thus by \cite[Thm.\ 4.5]{ST}, restricting to the subfamily $q^{-e}<r<1$, $r\in p^\mathbb{Q},$ the norms $\|-\|_{\mathbb{Q}_p,r}$ are multiplicative.

If not otherwise specified, we denote by $V\otimes_K W$ the projective tensor product of locally convex $K$-vector spaces $V,W.$

\begin{lemma}\label{lem:completion}
Let
\[\xymatrix@C=0.5cm{
  0 \ar[r] & V \ar[rr]^{} && W \ar[rr]^{ } && X \ar[r] & 0 }\] be a strict exact sequence of locally convex topological $K$-vector spaces with $W$ metrizable and $X$ Hausdorff, then
\begin{enumerate}
\item the  sequence of the associated Hausdorff completed spaces
  \[\xymatrix@C=0.5cm{
  0 \ar[r] & \hat{V} \ar[rr]^{} && \hat{W} \ar[rr]^{ } && \hat{X}\ar[r] & 0 }\]
  is again strict exact,
\item for a complete valued field extension $F$ of $K$ the associated sequence of completed base extension
    \[\xymatrix@C=0.5cm{
  0 \ar[r] & F\widehat{\otimes}_K V \ar[rr]^{} && F\widehat{\otimes}_K W \ar[rr]^{ } && F\widehat{\otimes}_K X \ar[r] & 0 }\]
    is again strict exact.
\item If $W$ is a $K$-Banach space, $V$ a closed subspace with induced norm and $X=W/V$ endowed with the quotient norm, then in (ii) the quotient norm coincides with the tensor product norm on $F\hat{\otimes}_K X.$
\end{enumerate}
\end{lemma}
\begin{proof}
By \cite[I.17 \S2]{B-TVS} with $W$ also $V$, $X$ and all their completions are metrizable. Hence the first statement follows from \cite[IX.26 Prop.\ 5]{B-TG}. For the second statement we first obtain the exact sequence
\[\xymatrix@C=0.5cm{
  0 \ar[r] & F{\otimes}_K V \ar[rr]^{} && F{\otimes}_K W \ar[rr]^{ } && F{\otimes}_K X \ar[r] & 0 }\] of metrizable locally convex spaces (\cite[Thm.\ 10.3.13]{PGS}). The first non-trivial map is strict by Thm.\ 10.3.8 in (loc.\ cit.). Regarding the strictness of the second map one easily checks that   $F{\otimes}_K W /F{\otimes}_K V $ endowed with the quotient topology satisfies   the universal property of the projective tensor product $F{\otimes}_K X.$ Now apply (i). The third item is contained in \cite[\S 3, $n^\circ$ 2, Thm.\ 1]{G}, see also \cite[Thm.\ 4.28]{R}.
\end{proof}

The kernel of the surjection of Fr\'echet spaces $D_{\mathbb{Q}_p}(o_L,K) \twoheadrightarrow D(o_L,K)$ is generated as a closed ideal by $\mathfrak{a} := \ker(L \otimes_{\mathbb{Q}_p} \Lie_{\mathbb{Q}_p}(o_L) \xrightarrow{a \otimes \mathfrak{x} \mapsto a\mathfrak{x}} \Lie_L(o_L))$. For $K = L$ this is \cite[Lemma 5.1]{Sc1}. As seen in the proof of Lemma \ref{module-struc} we have $K \widehat{\otimes}_L D_{\mathbb{Q}_p}(o_L,L) = D_{\mathbb{Q}_p}(o_L,K)$ and $K \widehat{\otimes}_L D(o_L,L) = D(o_L,K)$. Hence the assertion for general $K$ follows from Lemma \ref{lem:completion}(ii). We write $D_{r}(o_L,K)$ for the completion of $D(o_L,K)$ with respect to the quotient norm $\|-\|_r $ of $\|-\|_{\mathbb{Q}_p,r}$. By the proof of \cite[Prop.\ 3.7]{ST} we have the exact sequence of $K$-Banach algebras
\begin{equation}\label{f:Drquotient}
  0 \longrightarrow \hat{\mathfrak{a}}_r \longrightarrow D_{\mathbb{Q}_p,r}(o_L,K) \longrightarrow D_r(o_L,K) \longrightarrow 0
\end{equation}
where $\hat{\mathfrak{a}}_r$ denotes the closed ideal generated by $\mathfrak{a}$. Moreover, the $K$-Banach algebras $D_{r}(o_L,K)$ realize a Fr\'{e}chet-Stein structure on $D(o_L,K)$. For convenience we set $\mathfrak{r}_0:= q^{-\diamond e}$ and $\mathfrak{r}_m:= q^{-\frac{\diamond e}{p^m}}$ for $m\geq 1$. We, of course, have
\[
D(o_L,K) = \varprojlim_m D_{\mathfrak{r}_m}(o_L,K) \ .
\]
Moreover according to \cite[Cor.\ 5.13]{Sc} one has $D_{\mathfrak{r}_m}(o_L,K) = \mathbb{Z}[o_L] \otimes_{\mathbb{Z}[p^m o_L]}   D_{\mathfrak{r}_0}( {p^m}o_L,K)$.

We have corresponding results and will be using analogous notation for groups isomorphic to $o_L$. This applies, in particular, to $\Gamma_n$ for any $n \geq n_0$.  Note that $\Gamma_n^{p^m}=\Gamma_{n+me}$.

\subsubsection{\texorpdfstring{$(\varphi,\Gamma)$}{(phi,Gamma)}-modules}\label{subsec:phiGamma}

We recall the definition of as well as a few known facts about $(\varphi_L,\Gamma_L)$-modules (cf.\ \cite{BSX}).
Let $\mathfrak{Y}$ be either $\mathfrak{X} $ or $\mathbf{B}$ and $\cR := \cR_{K}(\mathfrak{Y})$.
Any $(\varphi_L,\Gamma_L)$-module $M$ over $\cR$ is, by definition, in particular an $\mathcal{R}$-module with a semilinear action of the group $\Gamma_L$. Our aim in this section is to show that these two structures on $M$ give rise to a module structure on $M$ under the 'group' Robba ring $\mathcal{R}_K(\Gamma_L)$.

\begin{definition}
  A $(\varphi_L,\Gamma_L)$-module $M$ over $\cR$ is a $\varphi_L$-module $M$ (see Definition \ref{def:phimodule}) equipped with a semilinear continuous action of $\Gamma_L$ which commutes with the endomorphism $\varphi_M.$ We shall write $\mathcal{M}(\cR)$ for the category of $(\varphi_L,\Gamma_L)$-modules over $\cR.$
\end{definition}

The continuity condition for the $\Gamma_L$-action on $M$, of course, refers to the product topology on $M \cong \cR^d$.

%

According to \cite[Prop.\ 2.25]{BSX}  the $\Gamma_L$-action on a $(\varphi_L,\Gamma_L)$-module $M$ is differentiable so that the derived action of the Lie algebra $\Lie(o_L^\times)$ on $M$ is available.

\begin{definition}
  The $(\varphi_L,\Gamma_L)$-module $M$ over $\cR$ is called $L$-analytic, if the derived action $\Lie(\Gamma_L) \times M \rightarrow M$ is $L$-bilinear, i.e., if the induced action $\Lie(\Gamma_L)\to \End(M)$ of the Lie algebra $\Lie(\Gamma_L)$ of $\Gamma_L$ is $L$-linear (and not just $\Qp$-linear). We shall write $\mathcal{M}^{an}(\cR)$ for the category of $L$-analytic $(\varphi_L,\Gamma_L)$-modules over $\cR.$
\end{definition}


In \cite{BSX} a $(\varphi_L,\Gamma_L)$-module $M$ over $\cR$ is only required to be projective instead of free as in our definition. Since  throughout this paper we are exclusively interested in $L$-analytic modules, that makes no difference as by \cite[Thm.\ 3.17]{BSX}  any $L$-analytic $(\varphi_L,\Gamma_L)$-module $M$ is actually a free $\cR$-module.

We have the following variant of Prop.\ \ref{descentphi} (cf.\ \cite[Prop.\ 2.24]{BSX}).

\begin{proposition}\label{descentphiGamma}
  Let $M$ be a $(\varphi_L,\Gamma_L)$-module over $\mathcal{R}$. Then there exists a model $(M_0,r_0)$ as in Prop.\ \ref{descentphi} equipped with a semilinear continuous action of $\Gamma_L$   such that
\begin{equation*}
  \mathcal{R} \otimes_{\cR^{[r_0,1)}} M_0 = M
\end{equation*}
respects the $\Gamma_L$-actions (acting diagonally on the left hand side).
\end{proposition}

From now on in this subsection \textbf{we fix a $(\varphi,\Gamma)$-module $M$ over $\cR$ and a pair $(r_0,M_0)$ as in Prop.\ \ref{descentphiGamma}}. We then have available the objects introduced after Prop.\ \ref{descentphi}. But now the finitely generated free modules $M^{I(r',\mathfrak{Y})}$ and $M^I$ are each in addition equipped with a semilinear continuous $\Gamma_L$-action, compatible with the identities \eqref{f:identities}. Moreover, the $\Gamma_L$-actions commutes with the $\psi^I$-operators, and the decompositions \eqref{f:decompM} and \eqref{f:decompMI} are $\Gamma_L$-equivariant.

\textbf{Assume} henceforth in this subsection that $M$ is an $L$-{\it analytic} $(\varphi_L,\Gamma_L)$-module over $\cR.$

\begin{proposition}\label{prop:distributionaction}
   The $\Gamma_L$-action on $M$ extends uniquely to a separately continuous action of the locally $L$-analytic distribution algebra $D(\Gamma_L,K)$ of $\Gamma_L$ with coefficients in $K$. If $M\xrightarrow{f}N$ is a homomorphism of $L$-analytic $(\varphi_L,\Gamma_L)$-modules, then $f$ is $D(\Gamma_L,K)$-equivariant with regard to this action.
\end{proposition}
\begin{proof}
First of all we observe that the Dirac distributions generate a dense $L$-subspace in $D(\Gamma_L,L)$ by \cite[Lem.\ 3.1]{ST1}. Since $\Gamma_L \cong o_L^\times$ we have seen in the proof of Lemma \ref{module-struc} that $D(\Gamma_L,K) = K \widehat{\otimes}_L D(\Gamma_L,L)$. Hence the Dirac distributions also generate a dense $K$-linear subspace of $D(\Gamma_L,K)$. Therefore the extended action is unique provided it exists.

Our assertion is easily reduced to the analogous statement concerning the Banach spaces $M^I$ for a closed interval $I=[r,s]$. From \cite[Prop.\ 2.16 and Prop.\ 2.17]{BSX}  we know that the $\Gamma_L$-action on $M^I$ is locally $\mathbb{Q}_p$-analytic. But since we assume $M$ to be $L$-analytic it is actually locally $L$-analytic (cf.\ Addendum to Prop.\ 2.25 and the argument at the end of the proof of Prop.\ 2.17 in \cite{BSX} ).

For our purpose we show more generally the existence, for any $K$-Banach space $W$, of a continuous $K$-linear map
\begin{equation*}
  I : \cC^{an}(\Gamma_L,W) \to \cL_b(D(\Gamma_L,K),W)
\end{equation*}
satisfying $I(f)(\delta_g)=f(g)$. Note that this map, if it exists is unique by our initial observation. Recall (cf.\ \cite[\S12]{pLG}) that the locally convex vector space $\cC^{an}(\Gamma_L,W)$ is the locally convex inductive limit of finite products of Banach spaces of the form $B \widehat{\otimes}_K W$ with a Banach space $B$, and that its strong dual $D(\Gamma_L,K)$ is the corresponding projective limit of the finite sums of dual Banach spaces $B'$. We therefore may construct the map $I$ as the inductive limit of finite products
of maps of the form
\begin{align*}
  B \widehat{\otimes}_K W & \longrightarrow \cL_b(B',W) \\
  x \otimes y & \longmapsto [\ell \mapsto \ell(x)y] \ .
\end{align*}
Since $B$ as a Banach space is barrelled this map is easily seen to be continuous (cf.\ the argument in the proof of \cite[Lem.\ 9.9]{NFA}).

Now suppose that $W$ carries a locally $L$-analytic $\Gamma_L$-action (e.g., $W = M^I$). For $y \in W$ let $\rho_y(g) := gy$ denote the orbit map in $\cC^{an}(\Gamma_L,W)$. We then define
\begin{align*}
  D(\Gamma_L,K) \times W & \longrightarrow W \\
  (\mu,y) & \longmapsto I(\rho_y)(\mu) \ .
\end{align*}
Due to our initial observation the proof of \cite[Prop.\ 3.2]{ST1}, that the above is a separately continuous module structure, remains valid even so $K$ is not assumed to be spherically complete.

By \cite[Rem.\ 2.20]{BSX} the homomorphism $f$ is continuous and hence the $D(\Gamma_L,K)$-equi\-variance of $f$ follows from the $\Gamma_L$-invariance by the first paragraph of this proof.
\end{proof}

Recall that $M^{I}$, for each $I = [r,s]$ with $r \geq r_0$, bears a natural $\Gamma_L$-action. Now, for each $n\geq 1,$ we will define a different action of  $\Gamma_n$ on $M^{[r,s]} $, which is motivated by Lemma \ref{lem:operatorH} below and which is crucial for analysing the structure of $M^{\psi_M=0}$ in the next subsection. To this end consider for each $\gamma\in \Gamma_n$ the operator $H_n(\gamma)$ on $M^{[r,s]}$ defined by
\begin{equation*}
  H_n(\gamma)(m) :=
  \begin{cases}
  \ev_{\pi_L^{-n}(\chi_{LT}(\gamma)-1)} \gamma m    & \text{if $\mathfrak{Y} = \mathfrak{X}$},  \\
  \eta(\pi_L^{-n}(\chi_{LT}(\gamma)-1),Z) \gamma m  & \text{if $\mathfrak{Y} = \mathbf{B}$ and $\Omega \in K$}.
  \end{cases}
\end{equation*}
Note that, since $\Gamma_n$ acts on ${\mathcal{O}_K(\mathfrak{Y})}$ via $\chi_{LT}$ and the $o_L^\times$-action, we may form the skew group ring ${\mathcal{O}_K(\mathfrak{Y})}[\Gamma_n]$, which due to the semi-linear action of $\Gamma_L$ on $M$ maps into the $K$-Banach algebra $\mathcal{E}nd_K(M^I)$ of continuous $K$-linear endomorphisms of $M^I,$ endowed with the operator norm $\|\ \|_{M^I}$. Hence we obtain the ring homomorphism
\begin{align*}
  H_n : K[\Gamma_n] & \longrightarrow {\mathcal{O}_K(\mathfrak{Y})}[\Gamma_n] \longrightarrow \mathcal{E}nd_K(M^I) \\
  \gamma & \longmapsto
  \begin{cases}
  \ev_{\pi_L^{-n}(\chi_{LT}(\gamma)-1)} \gamma     & \text{if $\mathfrak{Y} = \mathfrak{X}$},  \\
  \eta(\pi_L^{-n}(\chi_{LT}(\gamma)-1),Z) \gamma   & \text{if $\mathfrak{Y} = \mathbf{B}$ and $\Omega \in K$}.
  \end{cases}
\end{align*}

The next lemma holds true in both cases. We spell it out only in the $\mathbf{B}$-case since we technically need it only there.

\begin{lemma}\label{lem:operatorH}
Suppose that $\Omega$ is contained in $K$, and let $n \geq m \geq 1$.
\begin{enumerate}
\item We have for all $\sigma\in\Gamma_n$ \[\sigma\big(\eta(1,Z)\varphi_L^n(y)\big)=\eta(1,Z)\varphi_L^n(H_n(\sigma)(y)),\] i.e., the isomorphisms
 \begin{align*}
 M & \xrightarrow{\cong} \eta(1,Z)\varphi_L^n(M),\\
   M^{[r,s]} & \xrightarrow{\cong} \eta(1,Z)\varphi_L^n(M^{[r,s]})\\
   y&\mapsto \eta(1,Z)\varphi_L^n(y)
 \end{align*}
are $\Gamma_n$-equivariant with respect to the natural action on the right hand side and the action via $H_n$ on the left hand side.
\item The map
\begin{align*}
\z[\Gamma_m]\otimes_{\z[\Gamma_n],H_n}M^{[r,s]}& \to M^{[r^{1/q^{n-m}},s^{1/q^{n-m}}]}\\
\gamma \otimes y&\mapsto \eta(\frac{\chi_{LT}(\gamma)-1}{\pi_L^m},Z)\varphi_M^{n-m}(\gamma y)
\end{align*}
is a homeomorphism of Banach-spaces, where the left hand side is viewed as the direct sum of Banach-spaces $ \bigoplus_{\gamma\in\Gamma_m/\Gamma_n} \gamma \otimes M^{[r,s]}.$   Moreover, the  map is $\Gamma_m$-equivariant with respect to the $H_m$-action on the right hand side.
\item If the homomorphism  $H_n: K[\Gamma_n] \to \mathcal{E}nd_K(M^I) $ extends to a continuous homomorphism $R_K^I(\Gamma_n)\to \mathcal{E}nd_K(M^I)$, then    $H_m: K[\Gamma_m] \to \mathcal{E}nd_K(M^{I^{1/q^{n-m}}})$ extends to a continuous homomorphism $R_K^{I^{1/q^{n-m}}}(\Gamma_m)\to \mathcal{E}nd_K(M^{I^{1/q^{n-m}}})$. If the first extension is unique, so is the second one.
\end{enumerate}
\end{lemma}
\begin{proof}
(i) Setting $b:=\frac{\chi_{LT}(\sigma)-1}{\pi_L^n}$ we calculate
\begin{align*}
\sigma\big(\eta(1,Z)\varphi_L^n(m)\big)&=\sigma(\eta(1,Z))\varphi_L^n(\sigma m) \\
&=\eta(1+\pi_L^nb,Z)\varphi_L^n(\sigma m) \\
&=\eta(1,Z)\eta(\pi^nb,Z)\varphi_L^n(\sigma m)\\
&=\eta(1,Z)\varphi_L^n\left(\eta(b,Z)\sigma m \right).
\end{align*}
where we used the multiplicativity of $\eta$ in the first variable in the third and \eqref{f:etaphi} in the last equality.

(ii) By a straight forward computation one first checks that the map is well defined. The bijectivity  follows from \eqref{f:decompMI} using the bijection $1+\pi_L^m o_L/1+\pi_L^n o_L\xrightarrow{\cong}   o_L/\pi_L^{n-m} o_L,$ $\gamma\mapsto \frac{\chi_{LT}(\gamma)-1}{\pi_L^m}$ and that $M^{[r,s]}={\gamma}M^{[r,s]}.$

(iii) Base change induces the $R_K^{I^{1/q^{n-m}}}(\Gamma_m)$-action on
\begin{align}\label{f:basechange}
R_K^{I^{1/q^{n-m}}}(\Gamma_m) \otimes_{R_K^I(\Gamma_n),H_n} M^I&\cong   \z[\Gamma_m]\otimes_{\z[\Gamma_n]} R_K^I(\Gamma_n)\otimes_{R_K^I(\Gamma_n),H_n} \notag M^I\\
&\cong \z[\Gamma_m]\otimes_{\z[\Gamma_n],H_n} \otimes M^I\\
&\cong M^{I^{1/q^{n-m}}},\notag
\end{align}
where we used \eqref{f:RIcrossproduct} and (ii). The continuity is easily checked by considering 'matrix entries' which are built by composites of the original continuous map by other continuous transformations. Here we use that the identifications \eqref{f:RIcrossproduct} and \eqref{f:basechange} are  homeomorphisms when we endow the left hand side with the maximum norm. Finally, the claim regarding uniqueness follows from \eqref{f:basechange} as the action of $\Gamma_m$ is already determined by the original $H_m.$
\end{proof}

For the rest of this subsection we \textbf{assume that $\Omega$ is contained in $K$} and we will work exclusively in the $\mathbf{B}$-case, i.e., $\cR = \cR_K(\mathbf{B})$ and $\cR^I = \cR_K^I(\mathbf{B})$. The consequences for the $\mathfrak{X}$-case will be given in section \ref{subsec:descent}.



 There is a natural ring homomorphism $\cR^I\to\mathcal{E}nd_K(M^I)$ by assigning to $f\in\cR^I$ the multiplication- with-$f$-operator, which we denote   by the same symbol $f$.   Part (iii) of the following remark means that this ring homomorphism has operator norm $1.$
\begin{remark}\label{rem:normI}
\begin{enumerate}\phantom{section}
\item We have $ \sup_{x\in o_L} |\eta(x,Z)-1|_I <1$  and $| \eta(x,Z)|_I=1$ for all $x\in o_L$.\footnote{$  |\eta(x,Z)-1|_I =|\eta(1,Z)-1|_I<1$ for all $x\in o_L^\times$ because   any $x\in o_L^\times$ induces an isomorphism $[x](-)$ of $B_I$.}
\item $|\eta(px,Z)-1|_I\leq \max\{|\eta(x,Z)-1|_I^p,\frac{1}{q^e}|\eta(x,Z)-1|_I\} (=|\eta(x,Z)-1|_I^p,\mbox{ if }|\eta(x,Z)-1|_I{\geq q^{-\frac{e}{p-1}})}.$
\item   $|f|_I=\|f\|_{M^I}$ for all $f\in \cR^I$.
\end{enumerate}
\end{remark}

\begin{proof}
  It is known (\cite{ST2}) that $\eta(x,Z)=\eta(1,[x](Z))$ belongs to $1+Zo_{\mathbb{C}_p}[[Z]]$, whence we have, for any $x\in o_L,$ that $|\eta(x,Z)-1|_I <1$  from the definition of $|-|_I$, and (i) follows from the fact that the map $o_L\to \mathbb{R},$ $x\mapsto |\eta(x,Z)-1|_I$ is continuous with compact source. Affirmation (ii) is a consequence of  the expansion
  \begin{align*}
  \eta(px,Z)-1&=(\eta(x,Z)-1+1)^p-1\\
  &=(\eta(x,Z)-1)^p+\sum_{k=1}^{p-1}\begin{pmatrix}
                                      p \\
                                      k \\
                                    \end{pmatrix}
   (\eta(x,Z)-1)^k
  \end{align*} and $|\begin{pmatrix}
                                      p \\
                                      k \\
                                    \end{pmatrix}|=q^{ -e}$ for $k=1,\ldots, p-1.$
  (iii) follows from the submultiplicativity of $|- |_I$ plus the fact that $1\in\cR^I,$ which implies the statement on $M\cong(\cR^I)^m.$
\end{proof}

The above Remark allows us to fix a natural number $m_0=m_0(r_0)$ such that for all $m\geq m_0$ we have that
\begin{align}\label{f:m0}
|\eta(x,Z)-1|_I&< \mathfrak{r}_m \mbox { for all } x\in o_L \mbox{ and } |\eta(x,Z)-1|_I\leq \mathfrak{r}_0 \mbox { for all } x\in p^mo_L,\\
r_0^{{1}/{q}}&< \mathfrak{r}_m,\label{f:m0r0}
\end{align}
for any of the intervals $I=[r_0,r_0], $ $[r_0,r_0^{{1}/{q}}]$ and $ [r_0^{{1}/{q}},r_0^{{1}/{q}}]$. In the following let $I$ always denote one of those intervals.

\begin{lemma}\label{lem:Ked}
Let $\epsilon>0$ arbitrary. Then there exists $n_1\gg0$ such that, for any $n\geq n_1,$ the operator norm $\|-\|_{M^I}$ on $M^I$ satisfies 
\begin{equation}\label{f:n1}
 \|\gamma-1\|_{M^I}\leq  \epsilon \mbox{ for all } \gamma\in \Gamma_n.
\end{equation}
\end{lemma}

\begin{proof} We first  prove the statement for the module $M=\cR.$
For the convenience of the reader we adopt the proof of  \cite[Lem.\ 5.2]{Ked}. First note that for any fixed $f\in R^I$ by continuity of the action of $\Gamma_L$ there exists an open normal subgroup $H$ of $\Gamma_L$ such that
\begin{equation}
\label{f:Ked}
|(\gamma-1)f|_I<\epsilon| f|_I
\end{equation}
holds for all $\gamma\in H$. So me may assume that the latter inequality holds for $Z$ and $Z^{-1}$ simultaneously. Using the twisted Leibniz rule
\[(\gamma-1)(gf)=(\gamma-1)(g)f+\gamma(g)(\gamma -1)(f)\] and induction
we get \eqref{f:Ked} for all powers $Z^\z$. Since the latter form an orthogonal basis, the claim follows using that $|\gamma(g)|_I=|g|_I$ for any $\gamma\in H, g\in\cR^I$. If $M\cong \bigoplus_{i=1}^d\cR e_i$ and $m=\sum f_ie_i,$ we may assume that
\begin{equation}
\label{f:KedM}
|(\gamma-1)e_i|_{M^I}<\epsilon| e_i|_{M^I}
\end{equation}
 holds for   $1\leq i\leq d,$ and  apply the same Leibniz rule to $f_ie_i$ instead of $gf,$ whence the result follows, noting that $| e_i|_{M^I}=1$ by the definition of the maximum norm and that $|\gamma(e_i)|_{M^I}=| e_i|_{M^I}=1$ for any $\gamma\in H$ and $1\leq i\leq d$ as a consequence of \eqref{f:KedM}.
\end{proof}

We fix $n_1=n_1(r_0)\geq n_0$ such that the Lemma holds for $\epsilon=\mathfrak{r}_0.$ Then, for any $n\geq n_1, m\geq m_0,$  the above $H_n$ extends to   continuous ring homomorphisms
\begin{align*}
\tilde{H}_{n}: D_{\qp,\mathfrak{r}_m}(\Gamma_n,K)&\to \mathcal{E}nd_K(M^I),\\ \sum_{\mathbf{k} \in \mathbb{N}_0^d} \alpha_{\mathbf{k}}\hat{\ell}_{n,*}^{-1}(\mathbf{b})^{\mathbf{k}}&\mapsto \sum_{k\geq 0} \alpha_{\mathbf{k}} \prod_{i=1}^dH_n(\hat{\ell}_{n,*}^{-1}(b_i))^{k_i},
\end{align*}
and
\begin{align*}
\tilde{\mathbb{H}}_{n}:=\tilde{H}_{n}\circ \hat{\ell}_{n,*}^{-1}: D_{\qp,\mathfrak{r}_m}(o_L,K)\xrightarrow{\hat{\ell}_{n,*}^{-1}}D_{\qp,\mathfrak{r}_m}(\Gamma_n,K)&\to \mathcal{E}nd_K(M^I),\\ \sum_{\mathbf{k} \in \mathbb{N}_0^d} \alpha_{\mathbf{k}} \mathbf{b}^{\mathbf{k}}&\mapsto \sum_{k\geq 0} \alpha_{\mathbf{k}} \prod_{i=1}^dH_n(\hat{\ell}_{n,*}^{-1}(b_i))^{k_i}.
\end{align*}
  Indeed, we have
  \[{H}_n(\hat{\ell}_{n,*}^{-1}(b_i))=  \eta(\frac{\hat{\ell}_n^{-1}(h_i)-1}{\pi_L^n},Z)-1 +\eta(\frac{\hat{\ell}_n^{-1}(h_i)-1}{\pi_L^n},Z)(\hat{\ell}_n^{-1}(h_i) -1 )\]
 and since
\begin{align*}
\| \eta(\frac{\hat{\ell}_n^{-1}(h_i)-1}{\pi_L^n}&,Z)-1 +\eta(\frac{\hat{\ell}_n^{-1}(h_i)-1}{\pi_L^n},Z)(\hat{\ell}_n^{-1}(h_i) -1 ) \|_{M^I}\leq \\ &\max \{\| \eta(\frac{\hat{\ell}_n^{-1}(h_i)-1}{\pi_L^n},Z)-1    \|_{M^I},\|  \eta(\frac{\hat{\ell}_n^{-1}(h_i)-1}{\pi_L^n},Z)(\hat{\ell}_n^{-1}(h_i) -1 )  \|_{M^I}\}\leq \mathfrak{r}_m
\end{align*}
  by \eqref{f:n1},\eqref{f:m0} and Remark \ref{rem:normI} (i), the above defining sum converges with respect to the operator norm. Moreover, we have
\begin{align}\label{f:operatorHntilde}
\|\tilde{\mathbb{H}}_{n}(\lambda)\|_{M^I}\leq \sup_{\mathbf{k}} |\alpha_{\mathbf{k}}|\mathfrak{r}_m^{|\mathbf{k}|}=\|\lambda\|_{\mathbb{Q}_p,\mathfrak{r}_m}
\end{align}
for all $\lambda\in D_{\qp,\mathfrak{r}_m}(o_L,K),$ i.e., the operator norm of $\tilde{\mathbb{H}}_{n}$ is less or equal to $1.$

Since $M$ is assumed to be $L$-analytic, $\tilde{H}_n$ factorises over the desired ring homomorphism
 \begin{align*}
{H}_{n}: \bigg(D(\Gamma_n,K)\subseteq\bigg )D_{\mathfrak{r}_m}(\Gamma_n,K)\to \mathcal{E}nd_K(M^I)
\end{align*}
and $\tilde{\mathbb{H}}_n$   over
\begin{align*}
{\mathbb{H}}_{n}: \bigg(D(o_L,K)\subseteq\bigg )D_{\mathfrak{r}_m}(o_L,K)\to \mathcal{E}nd_K(M^I)
\end{align*}
by \eqref{f:Drquotient}. As $D_{\mathfrak{r}_m}(o_L,K)$ carries the quotient norm of $D_{\qp,\mathfrak{r}_m}(o_L,K) $ we obtain from \eqref{f:operatorHntilde}
\begin{align}\label{f:operatorHn}
\| {\mathbb{H}}_{n}(\lambda)\|_{M^I}\leq \inf_{\tilde{\lambda},pr(\tilde{\lambda})={\lambda}}\|\lambda\|_{\mathbb{Q}_p,\mathfrak{r}_m}=\|\lambda\|_{\mathfrak{r}_m}
\end{align}
for all $\lambda\in D_{\mathfrak{r}_m}(o_L,K),$ i.e., the operator norm of ${\mathbb{H}}_{n}$ is again less or equal to $1.$ By a similar, but simpler reasoning one shows the following
\begin{lemma}\label{lem:LTFourier}
The  isomorphism (LT together with Fourier) $\textfrak{L}:D(o_L,K)\cong\mathcal{O}_K(\mathbf{B}), \delta_a\mapsto \eta(a,Z),$ induces, for all $ m\geq m_0,$     a commutative diagram of continuous maps
\[\xymatrix{
  D_{\qp,\mathfrak{r}_m}(o_L,K) \ar[d]_{pr} \ar[dr]^{ } &  \\
  D_{\mathfrak{r}_m}(o_L,K) \ar[r]^{\textfrak{L} } & \mathcal{R}^I }\]
with operator norms less or equal to $1.$ Moreover, the operator norm of the {\it scalar} action via $\textfrak{L}$
\begin{align}\label{f:scal}
scal:(D(o_L,K)\subseteq ) D_{\mathfrak{r}_m}(o_L,K)\xrightarrow{\textfrak{L}} \mathcal{R}^I \to \mathcal{E}nd_K(M^I)
\end{align}
is also  bounded by $1$, see Remark \ref{rem:normI} (iii).
\end{lemma}

\begin{remark}\label{rem:uniqueness-H_n}
The maps $\tilde{H}_{n}$ and $\tilde{\mathbb{H}}_{n},$ as well as $  {H}_{n}$ and $ {\mathbb{H}}_{n}$ are uniquely determined by their restriction to $K[\Gamma_n]$ and $K[o_L]$, respectively, because these group algebras are dense in the sources $D_{\qp,\mathfrak{r}_m}(\Gamma_n,K)$, $D_{\mathfrak{r}_m}(\Gamma_n,K)$ and $D_{\qp,\mathfrak{r}_m}(o_L,K) $, $D_{\mathfrak{r}_m}(o_L,K)$, respectively.
\end{remark}

Applying our convention before Remark \ref{rem:normI} we usually  shall abbreviate  $scal(\mu)$ by $\textfrak{L}(\mu)$ for $\mu\in D(o_L,K) $ below  when we refer to this scalar action on $M^I.$ For the proof of Thm.\  \ref{Rrs-extension} below it will be crucial to compare the two actions $scal$ and $\mathbb{H}_n$ of $D(o_L,K)$ on $M^I$.

Finally,   for $n\geq n_1,$  we obtain similar maps for the original (multiplicative) action of $\Gamma_n\subseteq \Gamma$ on $M^I:$
\begin{align}\label{f:originalaction}
  D_{\qp,\mathfrak{r}_m}(\Gamma_n,K)&\to \mathcal{E}nd_K(M^I),\\ \sum_{\mathbf{k} \in \mathbb{N}_0^d} \alpha_{\mathbf{k}}\hat{\ell}_{n,*}^{-1}(\mathbf{b})^{\mathbf{k}}&\mapsto \sum_{k\geq 0} \alpha_{\mathbf{k}} \prod_{i=1}^d \hat{\ell}_{n,*}^{-1}(b_i))^{k_i}
\end{align}
with operator norm bounded by $1$.

A special case of the following lemma was pointed out to us by Rustam Steingart.

\begin{lemma}\phantomsection\label{lem:rustam}
Let $m\in\mathbb{N}$ be arbitrary. Setting $u_n(a):=\frac{\exp(\pi_L^na)-1}{\pi^n_La}$ for $a\in o_L\setminus \{0\}$ and $u_n(0)=1$ there exist $n_2 = n_2(m)$\Footnote{$=m+2e$, $v_p(n!)\leq \frac{n}{p-1}$.} such that $u_n(a)\in 1+\pi_L^mo_L$ for all $a\in o_L$ and $n\geq n_2.$
\end{lemma}
\begin{proof}

This is easily checked using $v_p(n!)\leq \frac{n}{p-1}$.
\end{proof}

In order to distinguish Dirac distributions for elements $\gamma$ in the {\it multiplicative} group $\Gamma_n$ from those for elements $a$ in the {\it additive} group $o_L$ we often shall write $\delta^\times_\gamma$ in contrast to $\delta_a.$

\begin{lemma}\label{lem:estimateX}
Let  $0 < \epsilon <1$ be arbitrary and $\Delta=\sum_k c_k(\delta_{a_k}-1)\in D(o_L,K)$ a finite sum with $a_k\in o_L$. Then there exists  $n_3=n_3(\epsilon,\Delta,r_0)$ such that for all $n\geq n_3$ it holds
\[\|\textfrak{L}(\Delta)-\mathbb{H}_{n}(\Delta)\|_{M^I}<\epsilon.\]
\end{lemma}
\begin{proof}
 Put $\xi:= \sup_k|c_k|$ and choose $\epsilon'\leq \epsilon$ such that $\epsilon'\xi<\epsilon.$ Then choose $m_1\geq m_0$ such that
\[\|\delta_\gamma^\times-1\|_{M^I}<\epsilon' \mbox{ and }  \|\delta_\gamma^\times-1\|_{\cR^I}<\epsilon'\] for all $\gamma\in\Gamma_{m_1}$ (see Lemma \ref{lem:Ked}). Now according to Lemma \ref{lem:rustam} we choose   $n_3:=n_2(m_1)\geq m_1$. Observing that for $a\in o_L$
\[\mathbb{H}_n(\delta_a)=\textfrak{L}(\delta_{u_n(a)a})\circ \delta^\times_{\hat{\ell}_n^{-1}(a)}\]
we estimate, for  $ n\geq  n_3,$
{\footnotesize\begin{align*}
\|\textfrak{L}( \Delta)-\mathbb{H}_{n}(\Delta)\|_{M^I}&=\|\sum_{k } c_k \left\{\textfrak{L}(\delta_{a_k})-1  - \left( \textfrak{L}(\delta_{u_n(a_k)a_k})\circ\delta^\times_{\hat{\ell}_n^{-1}(a_k)} -1 \right)\right\}\|_{M^I}\\
&=\|\sum_{k } c_k \left\{\textfrak{L}(\delta_{a_k})-1  - \left( \textfrak{L}(\delta_{u_n(a_k)a_k})\circ(\delta^\times_{\hat{\ell}_n^{-1}(a_k)} -1)+\textfrak{L}(\delta_{u_n(a_k)a_k})-1 \right)\right\}\|_{M^I}\\
&=\|\sum_{k } c_k \left\{\textfrak{L}(\delta_{a_k})-1  - \left( \textfrak{L}(\delta_{u_n(a_k)a_k})\circ(\delta^\times_{\hat{\ell}_n^{-1}(a_k)} -1)+\delta^\times_{u_n(a_k)}(\textfrak{L}(\delta_{a_k})-1 )\right)\right\}\|_{M^I}\\
&=\|\sum_{k } c_k \left\{(1-\delta^\times_{u_n(a_k)})(\textfrak{L}(\delta_{a_k})-1)  - \left(\textfrak{L}( \delta_{u_n(a_k)a_k})\circ(\delta^\times_{\hat{\ell}_n^{-1}(a_k)} -1) \right)\right\}\|_{M^I}\\
&\leq\sup_{k}\left( | c_k |\max\left\{\|(1-\delta^\times_{u_n(a_k)})(\textfrak{L}(\delta_{a_k})-1) \|_{M^I},  \| \textfrak{L}(\delta_{u_n(a_k)a_k})\circ(\delta^\times_{\hat{\ell}_n^{-1}(a_k)} -1) \|_{M^I} \right\}\right)\\
&\leq\epsilon'\sup_{k} | c_k |=\epsilon'\xi<\epsilon,
\end{align*}}
where we used for the  last line the estimate
\begin{align*}
 \|(1-\delta^\times_{u_n(a_k)})(\textfrak{L}(\delta_{a_k})-1) \|_{M^I} & =|(1-\delta^\times_{u_n(a_k)})(\textfrak{L}(\delta_{a_k})-1) |_{I} \\
    & \leq \|(1-\delta^\times_{u_n(a_k)})\|_{\cR^I} |\textfrak{L}(\delta_{a_k})-1 |_{I}\\
   & \leq \epsilon'|\eta(a_k,Z)-1 |_{I}\leq \epsilon'
\end{align*} (by Remark \ref{rem:normI}(i)/(iii) and due to the choice of $m_1$ and $n_3$) for the first term as well as the estimate
\begin{align*}
\| \textfrak{L}(\delta_{u_n(a_k)a_k})\circ(\delta^\times_{\hat{\ell}_n^{-1}(a_k)} -1) \|_{M^I}&\leq | \eta(u_n(a_k)a_k,Z)   |_I\| \delta^\times_{\hat{\ell}_n^{-1}(a_k)} -1  \|_{M^I}\\
&=\| \delta^\times_{\hat{\ell}_n^{-1}(a_k)} -1  \|_{M^I}<\epsilon'
\end{align*}
for the second term (again by Remark \ref{rem:normI}(i)/(iii) and due to the choice of $m_1$ and $n_3\geq m_1$).
\end{proof}
\begin{lemma}\label{lem:norm}
Let   $0 < \epsilon < 1$ be arbitrary and $\mu \in D(o_L,K)$ be any element. Then there exists  $\Delta=\sum_k c_k(\delta_{a_k}-1)\in D(o_L,K)$ a finite sum with $a_k\in o_L$ such that $\|\mu-\Delta\|_{\mathfrak{r}_{m_0}}<\epsilon.$ Moreover, for  $n_3=n_3(\epsilon,\Delta,r_0)$ from the previous Lemma and all $n\geq n_3$  we have
\[\|\textfrak{L}(\mu)-\mathbb{H}_{{n}}(\mu)\|_{M^I}<\epsilon.\]
In particular, if $\textfrak{L}(\mu)$ is invertible in $\mathcal{E}nd_K(M^I)$  or equivalently invertible as an element of $\cR^I$,\footnote{$M^I$ being a free $\cR^I$-module on which $\textfrak{L}(\mu)$ acts via the diagonal matrix with all diagonal entries  equal to $\textfrak{L}(\mu)$. }
  then firstly there exists $n_4=n_4(\mu,r_0)$ such that $\|\textfrak{L}(\mu)-\mathbb{H}_{n}(\mu)\|_{M^I}<|\textfrak{L}(\mu)^{-1}|_{I}^{-1}$ and $\|\mathbb{H}_{n}(\mu)^{-1}-\textfrak{L}(\mu)^{-1}\|_{M^I}< |\textfrak{L}(\mu)^{-1}|_I$ for any $n\geq n_4$ and secondly   the operator $ \mathbb{H}_{n}(\mu)$ is invertible, too.
\end{lemma}
\begin{proof} The existence of such $\Delta$ is clear because such elements form a dense subset of $D(o_L,K)$ in the Fr\'{e}chet topology   (as noted at the beginning of the proof of Prop.\ \ref{prop:distributionaction}). Consider the estimation   for $n\geq {n_3}$
\begin{align*}
  \|\textfrak{L}(\mu)-\mathbb{H}_n(\mu)\|_{M^I}& \leq \max\left(  \|\textfrak{L}(\mu-\Delta)\|_{M^I}, \|\textfrak{L}(\Delta)-\mathbb{H}_n(\Delta)\|_{M^I}, \| \mathbb{H}_n(\mu-\Delta)\|_{M^I} \right) <\epsilon,
\end{align*}
where we use the estimate
\[\|\textfrak{L}(\mu-\Delta)\|_{M^I}=|\textfrak{L}(\mu-\Delta)|_{I}\leq \|\mu-\Delta\|_{\mathfrak{r}_{m_0} }<\epsilon\]
by \eqref{f:scal} for the first, Lemma \ref{lem:estimateX} for the second and \eqref{f:operatorHn} for the last term.

Now suppose that $\textfrak{L}(\mu)$ as an operator on $M^I$ is invertible. We choose $\epsilon\leq \|\textfrak{L}(\mu)^{-1}\|_{M^I}^{-1}$,    $\Delta$ accordingly and put $n_4=n_3(\epsilon,\Delta,r_0).$ Then, for $n\geq n_4,$ we have $\| 1- \textfrak{L}(\mu)^{-1}\mathbb{H}_n(\mu)\|_{M^I}=\| \textfrak{L}(\mu)^{-1}(\textfrak{L}(\mu ) -\mathbb{H}_n(\mu))\|_{M^I}< 1$, whence $\sum_{k\geq 0} (1-\textfrak{L}(\mu)^{-1}\mathbb{H}_n(\mu) )^k$ converges in $\mathcal{E}nd_K(M^I)$ and $\mathbb{H}_n(\mu)^{-1}:=\left(\sum_{k\geq 0} (1-\textfrak{L}(\mu)^{-1}\mathbb{H}_n(\mu) )^k\right)\textfrak{L}(\mu)^{-1}$ is the  inverse of $\mathbb{H}_n(\mu)= \mu(1-(1-\textfrak{L}(\mu)^{-1}\mathbb{H}_n(\mu))).$

  Furthermore,
 \begin{align*}
\|\mathbb{H}_n(\mu)^{-1}-\textfrak{L}(\mu)^{-1}\|_{M^I}&=\|\left(\sum_{k\geq 1} (1-\textfrak{L}(\mu)^{-1}\mathbb{H}_n(\mu) )^k\right)\textfrak{L}(\mu)^{-1}\|_{M^I}\\
&\leq\sup_{k\geq 1} \|1-\textfrak{L}(\mu)^{-1}\mathbb{H}_n(\mu)\|_{M^I}^k|\textfrak{L}(\mu)^{-1}|_I<|\textfrak{L}(\mu)^{-1}|_I.
 \end{align*}
\end{proof}

Note that the above lemma applies to the variable $X$ and from now on we set $ n_4:=n_4(X,r_0)$.
In view of Lemma \ref{lem:operatorH} (iii)  the following lemma is crucial for the main result  Thm.\  \ref{Rrs-extension} of this section.
%

\begin{lemma}\label{lem:technical}\phantom{section} For $n\geq n_4$
\begin{enumerate}
\item the map $\Theta_n:  D(\Gamma_n,K)\xrightarrow{H_n}\mathcal{E}nd_K(M^I)$ extends  uniquely to a continuous ring homomorphism
\begin{align*}
 \mathcal{R}_K^{I }({\Gamma_{n }}) \to \mathcal{E}nd_K(M^I).
\end{align*} If $M\xrightarrow{f}N$ is a homomorphism of $L$-analytic $(\varphi_L,\Gamma_L)$-modules, then $f^I:M^I\to N^I$ is $\mathcal{R}_K^{I }({\Gamma_{n }})$-equivariant with regard to this action.
\item $M^I$ is a free $\mathcal{R}_K^I({\Gamma_n})$-module of rank $\mathrm{rk}_{\mathcal{R}}M.$  Any basis as $\mathcal{R}^I$-module also is a basis as $\mathcal{R}_K^I({\Gamma_n})$-module .
\item  The natural maps
\begin{align*}
 \cR_K^{[r_0,r_0]}({\Gamma_n})\otimes_{\cR_K^{[r_0,r_0^{\frac{1}{q}}]}({\Gamma_n})} M^{[r_0,r_0^{\frac{1}{q}}]} & \xrightarrow{\cong} M^{[r_0,r_0]},      \\
  \cR_K^{[r_0^{\frac{1}{q}},r_0^{\frac{1}{q}}]}({\Gamma_n})\otimes_{\cR_K^{[r_0,r_0^{\frac{1}{q}}]}({\Gamma_n})} M^{[r_0,r_0^{\frac{1}{q}}]} & \xrightarrow{\cong} M^{[r_0^{\frac{1}{q}},r_0^{\frac{1}{q}}]}
\end{align*}
\end{enumerate}
are isomorphisms.
\end{lemma}

\begin{proof}
(i) Inductively, for $n\geq n_4,$ we obtain from   Lemma \ref{lem:norm} - by expressing $(\mathbb{H}_n(\mu)^\pm)^k-(\textfrak{L}(\mu)^\pm)^k$ as  $\sum_{l=1}^k \begin{pmatrix}
                k \\
                l \\
              \end{pmatrix}
(\mathbb{H}_n(\mu)^\pm-\textfrak{L}(\mu)^\pm)^l(\textfrak{L}(\mu)^\pm)^{k-l}$ -  that
\begin{align}\label{f:estimatek}
 \|\mathbb{H}_n(\mu)^k-\textfrak{L}(\mu)^k\|_{M^I}<\left\{
                            \begin{array}{ll}
                               |\textfrak{L}(\mu)|_{I}^k, & \hbox{for $k\geq 0;$} \\
                              \leq |\textfrak{L}(\mu)^{-1}|_{I}^{-k}\leq |\textfrak{L}(\mu)|_{I}^{k}\leq |\textfrak{L}(\mu)^{k}|_{I}, & \hbox{for $k<0.$}
                            \end{array}
                          \right.
\end{align}
for all $k\in \mathbb{Z}.$ It follows for $\mu=X$ that, if $\sum_{k\in \z} a_kZ^k\in \cR^{I}$ with $a_i\in K,$
 then $\sum_{k\in \z}  a_k\mathbb{H}_n(X)^k$ converges in $\mathcal{E}nd_K(M^I)$,  because
\begin{align*}
\|a_k \mathbb{H}_n(X)^k\|_{M^I}\leq \max \{\|a_k (\mathbb{H}_n(X)^k-\textfrak{L}(\mu)^k)\|_{M^I},\|a_k \textfrak{L}(\mu)^k\|_{M^I}\}\leq \left\{
                            \begin{array}{ll}
                              |a_k| |Z|_{I}^k, & \hbox{for $k\geq 0;$} \\
                              |a_k| |Z^{-1}|_{I}^{-k}, & \hbox{for $k<0.$}
                            \end{array}
                          \right.
\end{align*}
goes to zero for $k$ going to $\pm\infty.$ In other words, we have extended the continuous ring homomorphism   $\Theta_n$  to a continuous ring homomorphism
\begin{align*}
 \cR^I \to \mathcal{E}nd_K(M^I),\;\; Z\mapsto \mathbb{H}_n(X).
\end{align*}
As by definition $\kappa^* \circ \hat{\ell}_{n ,*}$ extends to a continuous ring  isomorphism $ \mathcal{R}^{I }_K({\Gamma_{n}}) \xrightarrow{\cong}  \cR_K^{I }(\mathbf{B}) = \cR^I$ we have constructed a continuous ring homomorphism
\begin{align*}
 \mathcal{R}^{I }_K({\Gamma_{n }}) \to \mathcal{E}nd_K(M^I)
\end{align*}
 as claimed.

The uniqueness is a consequence of the fact that $\mathcal{R}^{I }_K({\Gamma_{n }}) $ is the completion of the localization $D(\Gamma_n,K)_{Y_n^{\mathbb{N}}}   $ by a certain norm, for which the extended map is continuous.

Concerning functoriality observe that the maps $f$ and $f^I$ are automatically continuous by \cite[Rem.\ 2.20]{BSX} (with respect to the canonical topologies). Without loss of generality we may assume that the estimates of Lemma \ref{lem:norm} hold for $M$ and $N$ simultaneously. By the invariance under the distribution algebra and $\cR$-linearity of $f$, the map $f^I$ is compatible with respect to the operators $\mathbb{H}_n(X)^\pm$ of $M^I$ and $N^I$. By continuity this extends to arbitrary elements of $\mathcal{R}_K^{I }({\Gamma_{n }})$.

 (ii)   follows similarly as in \cite{KPX}: Recall  that $(e_k)$ denotes an $\cR^I$-basis of $M^I$ and consider the maps
\[\Phi:\bigoplus_{k=1}^m \cR_K^I(\mathbf{B}) \cong M^I, (f_k)\mapsto \sum_{k=1}^m f_ke_k,\]
\[\Phi':\bigoplus_{k=1}^m \cR_K^I({\Gamma_n})  \to M^I, (f_k)\mapsto \sum_{k=1}^m f_k(e_k),\]
and
\[\Upsilon:\bigoplus_{k=1}^m \cR_K^I(\mathbf{B}) \xrightarrow{\cong} \bigoplus_{k=1}^m \cR_K^I({\Gamma_n}),\]
which in each component is given by $(\kappa^* \circ \hat{\ell}_{n ,*})^{-1}$.
Then we have from \eqref{f:estimatek} that
\[|\Phi'\circ\Upsilon\circ \Phi^{-1}(m)-m|_I<|m|_I,\]
i.e.,
\[\|\Phi'\circ\Upsilon\circ \Phi^{-1}-\id\|_I<1,\]
whence with $\Phi$ and $\Upsilon$ also $\Phi'$ is an isomorphism because $\Phi'\circ\Upsilon\circ \Phi^{-1}$ is invertible by the usual argument using the geometric series.

(iii) The base change property   follows from the fact that $\Phi'$ is compatible with changing the interval.
\end{proof}

\begin{theorem}\label{Rrs-extension}
Suppose that $\Omega$ is contained in $K$.
\begin{enumerate}
\item Let $J$ be any of the intervals
\[
[r_0,r_0]^{1/q^n} \mbox{ or } [r_0,r_0^{1/q}]^{1/q^n} \mbox{ for } n\geq 0 .
\]
Then the $\Gamma_{n_4}$-action on $M^J$ via $H_{n_4}$ extends uniquely to a continuous $\cR_K^J(\Gamma_{n_4})$-module structure. Moreover, $M^J$ is a finitely generated free $\cR_K^J(\Gamma_{n_4})$-module; any $\cR^{[r_0,1)}$-basis of $M_0$ is also an $\cR_K^J(\Gamma_{n_4})$-basis of $M^J$. If $M\xrightarrow{f}N$ is a homomorphism of $L$-analytic $(\varphi_L,\Gamma_L)$-modules, then $f^J : M^J \to N^J$ is $\cR_K^J(\Gamma_{n_4})$-equivariant with regard to this action.
\item The $\Gamma_1$-action on $M$ via $H_1$ extends uniquely to a   separately  continuous $\cR_K(\Gamma_1)$-module structure. Moreover,  $M$ is a finitely generated free  $\cR_K(\Gamma_1)$-module; any $\cR^{[r_0,1)}$-basis of $M_0$ is also an $\cR_K(\Gamma_1)$-basis of  $M$.  If $M\xrightarrow{f}N$ is a homomorphism of $L$-analytic $(\varphi_L,\Gamma_L)$-modules, then $f$ is $\cR_K(\Gamma_1)$-equivariant with regard to this action.
\end{enumerate}
\end{theorem}

\begin{proof}
(i) From Lemma \ref{lem:technical} we obtain, for any $n \geq n_4$, the $H_n$-action of $\cR_K(\Gamma_n)$ on $M^I$ for the original three intervals $I$. Using Lemma \ref{lem:operatorH}(iii) we deduce the $H_{n_4}$-action of $\cR_K^{I^{1/q^{n-n_4}}}(\Gamma_{n_4})$-action on $M^{I^{1/q^{n-n_4}}}$. The asserted properties of these actions follow from the same lemmas.

(ii) By the uniqueness part in (i) we may glue the $\cR_K^J(\Gamma_{n_4})$-actions on the $M^J$ to a continuous $\cR_K^{[r_0,1)}(\Gamma_{n_4})$-action on $M^{[r_0,1)}$. By Remark \ref{rem:denseX}.ii it is uniquely determined by the $\Gamma_{n_4}$-action. Therefore we may vary $r_0$ now and obtain in the inductive limit a separately continuous $H_{n_4}$-action of $\cR_K(\Gamma_{n_4})$ on $M$. Using \eqref{f:Rcrossproduct} and Lemma \ref{lem:operatorH}(ii) we deduce the separately continuous $H_1$-action of $\cR_K(\Gamma_1)$ on $M$. Again by Remark \ref{rem:denseX} this action is uniquely determined by the $\Gamma_1$-action. The remaining assertions follow from the corresponding ones in (i).
\end{proof}


\subsubsection{The structure of \texorpdfstring{$ M^{\psi_M=0}$}{psi=0 on M}}

We still \textbf{assume that $\Omega$ is contained in $K$ and  let $M$ be an $L$-analytic $(\varphi_L,\Gamma)$-module over $\cR = \cR_K(\mathbf{B})$.} We want to show that $M^{\psi=0}$ carries a natural $\cR_K(\Gamma_L)$-action extending the action of $D(\Gamma_L,K)$.\\

From \eqref{f:decompM} and using formula \eqref{f:etapsi} and \eqref{f:etasigma} we have
\begin{align}\label{f:Mspiequal0}
M^{\psi_L=0}&=\bigoplus_{a\in(o_L/\pi_L)^\times} \eta(a,Z)\varphi_M(M) =\mathbb{Z}[\Gamma_L]\otimes_{\mathbb{Z}[\Gamma_{1}]}\left(\eta(1,Z)\varphi_M(M) \right).
 \end{align}

\begin{theorem}\label{thm:Mpsiequal0}
The $\Gamma_L$ action on $M$ extends   to a unique separately continuous $\cR_K(\Gamma_L)$-action on $M^{\psi_L=0}$ (with respect to the $LF$-topology on $\cR_K(\Gamma_L) $ and the subspace topology on $M^{\psi_L=0} $); moreover the latter is a free $\cR_K(\Gamma_L)$-module of rank $\mathrm{rk}_{\cR}M.$ If $e_1,\ldots, e_r$ is a basis of $M$ over $\cR,$ an $\cR_K(\Gamma_L)$-basis of $M^{\psi_L=0}$ is given by $\eta(1,Z)\varphi_M(e_1),\ldots, \eta(1,Z)\varphi_M(e_r).$
  If $M\xrightarrow{f}N$ is a homomorphism of $L$-analytic $(\varphi_L,\Gamma_L)$-modules, then $M\xrightarrow{f^{\psi_L=0}}N$ is $\cR_K(\Gamma_L)$-equivariant with regard to this action.
\end{theorem}



\begin{proof}
By Lemma \ref{lem:operatorH} (i) we transfer the $\cR_K(\Gamma_1)$-action on $M$ from Thm.\  \ref{Rrs-extension}(ii) to the space $\eta(1,Z)\varphi_M(M).$ Note that the resulting action is separately continuous for the subspace topology of $\eta(1,Z)\varphi_M(M),$ because the map
$\varphi_L : M \rightarrow M$ is a homeomorphism onto its image. The latter is a consequence of the existence of the continuous operator  $\psi_L$ and   the relation $\psi_L\circ \varphi_L=\frac{q}{\pi_L}\id_M.$ Finally, because of \eqref{f:cross} and \eqref{f:Mspiequal0} the $\cR_K(\Gamma_1)$-action extends to the asserted $\cR_K(\Gamma_L)$-action. Similarly as before, since $\Gamma_L$ spans a dense subspace of $D(\Gamma_L,K)$, the uniqueness of the action follows from Remark \ref{rem:denseX}.
\end{proof}

\subsubsection{Descent}\label{subsec:descent}

For the proof of Thm.\ \ref{thm:Mpsiequal0} we had to work over a field $K$ containing the period $\Omega$ since only then we were able to write elements in $\cR_K(\Gamma_L)$ or rather $\cR_K(\Gamma_n)$ as certain Laurent series in one variable $Y_n$ by means of the Lubin-Tate isomorphism $\cR_K(\mathfrak{X})\cong\cR_K(\mathbf{B})$,  which in general does not exist  over   $L.$ In this section we are going to explore to which extent the structure theorem over $K$ descends to $L.$ We shall consider two situations, i.e., we now start with an $L$-analytic $(\varphi_L,\Gamma)$-module $M$ over $\cR_L(\mathfrak{X})$ or $\cR_L(\mathbf{B})$, respectively. Thus, in what follows let $\mathfrak{Y}$ be either $\mathfrak{X}$ or $\mathbf{B}$ and $\cR = \cR_L(\mathfrak{Y})$. Then we consider the functor
\begin{align*}
 \mathfrak{M}^{an}(\cR_{L}(\mathfrak{Y}))&\to \mathfrak{M}^{an}(\cR_{K}(\mathfrak{Y})), \\
  M&\mapsto M_K:=\cR_{K}(\mathfrak{Y})\otimes_{\cR_{L}(\mathfrak{Y})}M\cong K\hat{\otimes}_{\iota,L} M,
\end{align*}
where the last isomorphism  and the well-definedness of the functor are established in \cite[Lem.\ 2.23]{BSX}. Moreover, there is a natural action of $G_L$ on both $\cR_{K}(\mathfrak{Y})\cong K\hat{\otimes}_{\iota,L}\cR_{L}(\mathfrak{Y})$ and $M_K$ via the first tensor factor (and the identity on the second). We have
\begin{align}\label{f:decentR}
\cR_{K}(\mathfrak{Y})^{G_L}\cong  \cR_{L}(\mathfrak{Y})
\end{align}
by \cite[Prop.\ 2.7 (iii)]{BSX},
whence also
\begin{align}\label{f:decentM_K}
(M_K)^{G_L}=M
\end{align}
because $M$ is finitely generated free over $\cR_{L}(\mathfrak{Y}) $ (by definition or \cite[Thm.\ 3.17]{BSX}) and hence $M_K$ has a $G_L$-invariant basis over $\cR_{K}(\mathfrak{Y}) $.

Since the $\varphi_L$-operator on $\cR_{K}(\mathfrak{Y})$ is induced from that on $\cR_{L}(\mathfrak{Y})$, it commutes with the action of $G_L.$ Similarly, one checks that this action commutes with the operator $\psi_L$ of $\cR_{K}(\mathfrak{Y})$. Indeed, by Lemma \ref{phi-free-Robba} there exists a $G_L$-invariant basis of $\cR_{K}(\mathfrak{Y})$ over $\varphi_L(\cR_{K}(\mathfrak{Y})),$ whence the trace commutes with the $G_L$-action. From this and the construction of the operator $\psi_M$, one derives easily that also $\psi_M$ commutes with the $G_L$-action. As a consequence we obtain natural isomorphisms
\begin{align}
\label{f:G_Linvariants} M^{\psi_M=0}\cong \left((M_K)^{G_L}\right)^{\psi_M=0}\cong\left((M_K)^{\psi_M=0}\right)^{G_L}.
\end{align}
Since the    Lubin-Tate isomorphism $\cR_K(\mathfrak{X})\cong\cR_K(\mathbf{B})$ respects the $(\varphi_L,\Gamma_L)$-module structure, Thm.\  \ref{thm:Mpsiequal0} applies for both choices of $\mathfrak{Y},$ i.e., we obtain a separately  continuous action
\begin{align}\label{f:actionRonpsiequal0}
  \cR_K(\Gamma_L)\times (M_K)^{\psi=0} \to (M_K)^{\psi=0}.
\end{align}
Moreover, if $M=\bigoplus_{i=1}^r \cR_{L}(\mathfrak{Y})e_i$, then the families $(\eta(1,Z)\varphi_M(e_i))$ and $(\mathrm{ev}_1\varphi_M(e_i))$ form basis of $(M_K)^{\psi=0}$ as
$ \cR_K(\Gamma_L)$-modules in case $\mathbf{B}$ and $\mathfrak{X},$ respectively. Therefore, we consider next a natural $G_L$-action on $ \cR_K(\Gamma_L)$ and show that \eqref{f:actionRonpsiequal0} is $G_L$-equivariant. To the first aim we use the canonical isomorphisms \eqref{f:cross} and \eqref{f:diagphiR}
\begin{equation*}
  \cR_K(\Gamma_L) \xleftarrow{\cong} \mathbb{Z}[\Gamma_L]\otimes_{\mathbb{Z}[\Gamma_n]} \cR_K(\Gamma_n) \xrightarrow[\hat{\ell}_{n*}]{\cong} \mathbb{Z}[\Gamma_L]\otimes_{\mathbb{Z}[\Gamma_n]} \cR_K(\mathfrak{X})
\end{equation*}
to extend the $G_L$-action from $\cR_K(\mathfrak{X})$ to $\cR_K(\Gamma_L);$ clearly, we obtain from \eqref{f:decentR} and the fact that the isomorphism $\cR_K(\Gamma_n)\xrightarrow[\cong]{\ell_{n*}}\cR_K(\mathfrak{X})$ is defined over $L$, that\Footnote{On $\cR_K(\Gamma_L)$ we only consider the natural, but not any twisted action by $G_L!$}
\begin{align}\label{f:decentRGamma}
 \cR_K(\Gamma_L)^{G_L} & \cong \cR_L(\Gamma_L).
\end{align}
To the second aim we proof the following
\begin{lemma}
The action \eqref{f:actionRonpsiequal0} is $G_L$-equivariant.
\end{lemma}

\begin{proof}
 We fix $\sigma\in G_L$ and define a second separately continuous action
 \begin{align*}
  \cR_K(\Gamma_L)\times (M_K)^{\psi=0} \to (M_K)^{\psi=0}
\end{align*}
by sending $(\lambda,x)$ to $\sigma^{-1}\left(\sigma(\lambda)(\sigma(x))  \right)$ (using that $\sigma$ and $\psi_L$ commute and that $\sigma$ is a homeomorphism). By the uniqueness statement of Thm.\ \ref{thm:Mpsiequal0}, it suffices to show that the new and original action coincide on $\Gamma\times (M_K)^{\psi=0}.$ We shall show that these actions coincide even as actions $\Gamma_L\times M_K\to M_K:$ For $\gamma\in\Gamma$, $f\in \cR_K(\mathfrak{Y})$ and $m\in M$ we calculate
\begin{align*}
  \sigma^{-1}\left(\sigma(\gamma)(\sigma(f\otimes m) ) \right) & =  \sigma^{-1}\left( \gamma(\sigma( f)\otimes m) \right) \\
  & =  \sigma^{-1}\left( \gamma(\sigma( f))\otimes\gamma( m) \right) \\
    & = \sigma^{-1}\left(  \sigma(\gamma( f) \right))\otimes \gamma( m)\\
    & =\gamma (f)\otimes \gamma( m)\\
    & =\gamma (f\otimes m).
\end{align*}
Here we used firstly that $\sigma$ acts trivially on $\gamma$ (or rather $\mathrm{ev}_\gamma$) as they are already defined over $L$ (via the Fourier transformation) and secondly, that the $G_L$- and $\Gamma_L$-actions commute. Since this equality holds for all $\sigma\in G_L,$ the claim follows.
\end{proof}

Taking $G_L$-invariants of \eqref{f:actionRonpsiequal0}
therefore induces - upon using \eqref{f:G_Linvariants} and \eqref{f:decentRGamma} - the following separately continuous action
\begin{align}
  \cR_L(\Gamma_L)\times M^{\psi=0} \to M^{\psi=0}
\end{align}
which   extends the $\Gamma_L$-action. We thus obtain the following

\begin{theorem}\phantomsection\label{thm:Mpsiequal0-decent}
\begin{enumerate}
\item The $\Gamma_L$-action on $M$ (in $\mathfrak{M}^{an}(\cR_{L}(\mathfrak{X}))$ or   $\mathfrak{M}^{an}(\cR_{L}(\mathbf{B}))$) extends   to a   separately continuous $\cR_L(\Gamma_L)$-action on $M^{\psi_L=0}$ (with respect to the $LF$-topology on $\cR_L(\Gamma_L) $ and the subspace topology on $M^{\psi_L=0} $). If $M\xrightarrow{f}N$ is a homomorphism of $L$-analytic $(\varphi_L,\Gamma_L)$-modules, then $M\xrightarrow{f^{\psi_L=0}}N$ is $\cR_L(\Gamma_L)$-equivariant with regard to this action.

\item If $\mathfrak{Y}=\mathfrak{X}$ then $M^{\psi_L=0} $ is a free $\cR_L(\Gamma_L)$-module of rank $\mathrm{rk}_{\cR}M.$ More precisely, if $e_1,\ldots, e_r$ is a basis of $M$ over $\cR_L(\mathfrak{X}),$ then an $\cR_L(\Gamma_L)$-basis of $M^{\psi_L=0}$ is given by $\mathrm{ev}_1\varphi_M(e_1),\ldots, \mathrm{ev}_1\varphi_M(e_r).$
\end{enumerate}
\end{theorem}

\begin{proof}
It is easy to check that also the $\cR_L(\Gamma_L)$-equivariance of $f^{\psi_L=0}$ follows by descent. Therefore only (ii) remains to be shown. But this is an immediate consequence of the fact noted above, that the family $(\mathrm{ev}_1\varphi_M(e_i))$ forms a $G_L$-invariant basis of $(M_K)^{\psi=0}$ as
$ \cR_K(\Gamma_L)$-module, by just taking $G_L$-invariants again.
\Footnote{For the uniqueness we would need Remark \ref{rem:denseX} for $L$.}
\end{proof}

\begin{remark}
\begin{enumerate}
\item For each complete intermediate field $L\subseteq K' \subseteq \mathbb{C}_p$ we obtain   an analogous structure theorem for $(M_{K'})^{\psi=0}$ over $\cR_{K'}(\Gamma_L)$  by   replacing $L$ by $K'$ everywhere in the above reasoning.
\item Since for $\mathfrak{Y}=\mathbf{B}$ the basis $(\eta(1,Z)\varphi_M(e_i))$ of $(M_K)^{\psi=0}$ as $\cR_K(\Gamma_L)$-module is visibly not $G_L$-invariant, we cannot conclude the analogue of Thm.\ \ref{thm:Mpsiequal0-decent}(ii) in this case.
\end{enumerate}
\end{remark}

\subsubsection{The Mellin transform and twists}

Extending the Mellin transform from Lemma \ref{Mellin} we introduce
the map
\begin{equation*}
  \mathfrak{M} : \cR_K(\Gamma_L) \xrightarrow{\;\cong\;} \cR_K(\mathfrak{X})^{\psi_L=0},\;\; \lambda\mapsto \lambda(\ev_1) \ ,
\end{equation*}
which is an isomorphism by Thm.\ \ref{thm:Mpsiequal0-decent}. If $\Omega\in K$, then its composite with the LT-isomorphism is the isomorphism
\begin{equation*}
  \mathfrak{M}_{LT} : \cR_K(\Gamma_L) \xrightarrow{\;\cong\;} \cR_K(\mathbf{B})^{\psi_L=0},\;\; \lambda\mapsto \lambda(\eta(1,Z)) \ .
\end{equation*}

%

Recall the twist operators $Tw_{\chi}$ from section \ref{subsec:twisting}.

\begin{lemma}\label{lem:LT-twisting}
  The diagram
\begin{equation}\label{f:twistX}
  \xymatrix{
    \cR_K(\Gamma_L) \ar[d]_{Tw_{\chi_{LT}}} \ar[r]^-{\mathfrak{M}} & \cR_K(\mathfrak{X})^{\psi_L=0} \ar[d]^{ \partial^{\mathfrak{X}}_{\mathrm{inv}}}_{\cong} \\
   \cR_K(\Gamma_L) \ar[r]^-{\mathfrak{M}} & \cR_K(\mathfrak{X})^{\psi_L=0}   }
\end{equation}
is commutative; in particular, the right hand vertical map is an isomorphism.
\end{lemma}
\begin{proof}
The commutativity can be checked after base change. Assuming $\Omega\in K$ the diagram corresponds by Remark \ref{rem:logXlogLT} to the diagram
\begin{equation}\label{f:twist}
  \xymatrix{
    \cR_K(\Gamma_L) \ar[d]_{Tw_{\chi_{LT}}} \ar[r]^-{\mathfrak{M}_{LT}} & \cR_K(\mathbf{B})^{\psi_L=0} \ar[d]^{\frac{1}{\Omega}\partial_{\mathrm{inv}}}_{\cong} \\
   \cR_K(\Gamma_L) \ar[r]^-{\mathfrak{M}_{LT}} & \cR_K(\mathbf{B})^{\psi_L=0} .  }
\end{equation}

Now, the corresponding result for $\cR_K(\Gamma_L)$ replaced by $D(\Gamma_L,K)$ is implicitly given in sections \ref{subsec:twisting} and \ref{subsec:LT}. \cite[\S 1.2.4]{Co2} establishes, for $\gamma \in \Gamma_L$, the relation $\partial_{\mathrm{inv}}\circ \gamma=\chi_{LT}(\gamma)\gamma \circ \partial_{\mathrm{inv}}$ as operators on $\cR_K(\mathbf{B})$. It follows by $K$-linearity and continuity that the relation of operators $\partial_{\mathrm{inv}}\circ \lambda=Tw_{\chi_{LT}}(\lambda)\circ \partial_{\mathrm{inv}}$ holds for all $\lambda\in D(\Gamma_L,K)$.
By continuity of the action of $\cR_K(\Gamma_L) = \mathbb{Z}[\Gamma_L] \otimes_{\mathbb{Z}[\Gamma_n]} \cR_K(\Gamma_n)$ on $\cR_K(\mathbf{B})^{\psi_L=0}$ it suffices to check the compatibility for the element $Y_n^{-1}$, where $Y_n \in D(\Gamma_n,K)$, for $n >> 0$, has been defined at the end of section \ref{sec:groupRobba}. Using that $Tw_{\chi_{LT}}$ is multiplicative and that $\partial_\mathrm{inv} (\eta(1,Z)) = \Omega \eta(1,Z)$ the claim follows from the relation
\begin{align*}
 Tw_{\chi_{LT}}(Y_n^{-1})\eta(1,Z)&=Tw_{\chi_{LT}}(Y_n)^{-1}\frac{1}{\Omega}\partial_{\mathrm{inv}}\left( Y_n Y_n^{-1}\eta(1,Z)\right)\\
&=Tw_{\chi_{LT}}(Y_n)^{-1}Tw_{\chi_{LT}}(Y_n)\frac{1}{\Omega}\partial_{\mathrm{inv}}\left( Y_n^{-1}\eta(1,Z)\right)\\
&=\frac{1}{\Omega}\partial_{\mathrm{inv}}\left( Y_n^{-1}\eta(1,Z)\right).
\end{align*}
\end{proof}

\begin{lemma}\label{lem:varphinR}
Assume $\Omega\in K$ and let $n_1$ be as in Lemma \ref{lem:Ked}. Then, for $n\geq n_1,$ the map $\mathfrak{M}_{LT}$ induces isomorphisms
\begin{equation*}
  \cR_K(\Gamma_n)\cong \varphi_L^n(\cR_K(\mathbf{B}))\eta(1,Z)\;\;(\subseteq \cR_K(\mathbf{B})^{\psi_L=0} )
\end{equation*}
of $\cR_K(\Gamma_n)$-modules and
\begin{equation*}
  D(\Gamma_n,K)\cong \varphi_L^n({\mathcal{O}_K(\mathbf{B})})\eta(1,Z)\;\;(\subseteq {\mathcal{O}_K(\mathbf{B})}^{\psi_L=0} )
\end{equation*}
of $D(\Gamma_n,K)$-modules.
\end{lemma}
\begin{proof}{  By taking limits the first isomorphism follows from Lemma \ref{lem:technical}(ii) in combination with Lemma \ref{lem:operatorH}(i), both applied to $M^I=\cR_K(\mathbf{B}_I)$. The isomorphism of the latter restricts visibly to the isomorphism ${\mathcal{O}_K(\mathbf{B}_I)}\cong \varphi_L^n({\mathcal{O}_K(\mathbf{B}_I)})\eta(1,Z)$ while ${\mathcal{O}_K(\mathbf{B}_I)}$ is a free $\mathcal{O}_K(\hat{\ell}_n^*\circ\kappa(\mathbf{B}_I))$-module with basis $1$ by an obvious analogue of the former reference. Hence we  obtain the  second isomorphism by the same reasoning.}
\end{proof}

\newpage

\subsection{Explicit elements}\label{sec:explicit}

 There are two sources for explicit elements in the distribution algebras $D(o_L^\times,L)$ and $D(U_n,L)$, where in this section we fix an $n \geq n_1$, i.e.,  $\log : U_n \xrightarrow{\cong} \pi_L^n o_L$ is an isomorphism. First of all we have, for any group element $u \in o_L^\times$, resp.\ $u \in U_n$, the Dirac distribution $\delta_u$ in $D(o_L^\times,L)$, resp.\ in $D(U_n,L)$. As in section \ref{sec:add-Robba} the corresponding holomorphic function $F_{\delta_u} = \ev_u$ is the function of evaluation in $u$.

\begin{lemma}\phantomsection\label{Dirac-zeros}
\begin{itemize}
  \item[i.] Let $u \in o_L^\times$ be any element not of finite order; then the zeros of the function $\ev_u - 1$ on $\mathfrak{X}^\times$ are exactly the characters $\chi$ of finite order such that $\chi(u) = 1$.
  \item[ii.] For any $1 \neq u \in U_n$ the zeros  of the function $\ev_u - 1$ on $\mathfrak{X}_n^\times$ all have multiplicity one.
\end{itemize}
\end{lemma}
\begin{proof}
i. Obviously the zeros of $\ev_u - 1$ are the characters $\chi$ such that $\chi(u) = 1$. On the other hand consider any locally $L$-analytic character $\chi : o_L^\times \rightarrow \mathbb{C}_p^\times$. Its kernel $H := \ker(\chi)$ is a closed locally $L$-analytic subgroup of $o_L^\times$. Hence its Lie algebra $\Lie(H)$ is an $L$-subspace of $\Lie(o_L^\times) = L$. We see that either $\Lie(H) = L$, in which case $H$ is open in $o_L^\times$ and hence $\chi$ is a character of finite order, or $\Lie(H) = 0$, in which case $H$ is zero dimensional and hence is a finite subgroup of $o_L^\times$. If $\chi(u) = 1$ then, by our assumption on $u$, the second case cannot happen.

ii. (We will recall the concept of multiplicity further below.) Because of the isomorphism $\mathfrak{X}_n^\times \cong \mathfrak{X}$ it suffices to prove the corresponding assertion in the additive case. Let $0 \neq a \in o_L$ and let $\chi \in \mathfrak{X}(\mathbb{C}_p)$ be a character of finite order such that $\chi(a) = 1$. By \cite{ST2} we have an isomorphism between $\mathfrak{X}_{/\mathbb{C}_p}$ and the open unit disk $\mathbf{B}_{/\mathbb{C}_p}$. Let $z \in \mathbf{B}(\mathbb{C}_p)$ denote the image of $\chi$ under this isomorphism. By \cite[Prop.\ 3.1]{ST2}  and formula $(\diamond\diamond)$ on p.\ 458, the function $\ev_a - 1$ corresponds under this isomorphism to the holomorphic function on $\mathbf{B}(\mathbb{C}_p)$ given by the formal power series
\begin{equation*}
  F_{at'_0}(Z) = \exp(g \Omega \log_{LT}(Z)) - 1 \ ,
\end{equation*}
where $\Omega \neq 0$ is a certain period. By assumption we have $F_{at'_0}(z) = 0$. On the other hand the formal derivative of this power series is
\begin{equation*}
  \frac{d}{dZ} F_{gt'_0}(Z) = g \Omega g_{LT}(Z) (F_{gt'_0}(Z) + 1) \ .
\end{equation*}
Since $g_{LT}(Z)$ is a unit in $o_L[[Z]]$ we see that $z$ is not a zero of this derivative. It follows that $z$ has multiplicity one as a zero of $F_{at'_0}(Z)$.
\end{proof}

The other source comes from the Lie algebra $\Lie(U_n) = \Lie(o_L^\times) = L$. We have the element
\begin{equation*}
  \nabla := 1 \in \Lie(o_L^\times) = L \ .
\end{equation*}
On the other hand there is the $L$-linear embedding (\cite[\S2]{ST1})
\begin{align*}
  \Lie(U_n) & \longrightarrow D(U_n,L) \\
  \mathfrak{x} & \longmapsto  [f \mapsto \frac{d}{dt} f(\exp_{U_n} (t\mathfrak{x}))_{| t=0} ]  \ ,
\end{align*}
which composed with the Fourier isomorphism becomes the map
\begin{align*}
  \Lie(U_n) & \longrightarrow \mathcal{O}_L(\mathfrak{X}_n^\times) \\
  \mathfrak{x} & \longmapsto  [\chi \mapsto d\chi(\mathfrak{x}) ]  \ .
\end{align*}
On the one hand we therefore may and will view $\nabla$ always as a distribution on $U_n$ or $o_L^\times$. On the other hand, using the formula before \cite[Lem.\ 1.28]{BSX}, one checks that the function $F_\nabla$ (corresponding to $\nabla$ via the Fourier isomorphism) on $\mathfrak{X}_n^\times$ is explicitly given by
\begin{equation}\label{f:Fnabla}
  F_\nabla(\chi) = \pi_L^{-n} \log(\chi(\exp(\pi_L^n))) \ .
\end{equation}

\begin{lemma}\label{nabla-zeros}
The zeros of the function $F_\nabla$ on $\mathfrak{X}_n^\times$ are precisely the characters of finite order each with multiplicity one.
\end{lemma}
\begin{proof}
Once again because of the isomorphism $\mathfrak{X}_n^\times \cong \mathfrak{X}$ it suffices to prove the corresponding assertion in the additive case. This is done in \cite[Lem.\ 1.28]{BSX}.
\end{proof}

To recall from \cite[\S1.1]{BSX} the concept of multiplicity used above and to explain a divisibility criterion in these rings of holomorphic functions we let $\mathfrak{Y}$ be any one dimensional smooth rigid analytic quasi-Stein space over $L$ such that $\mathcal{O}_L(\mathfrak{Y})$ is an integral domain. Under these assumptions the local ring in a point $y$ of the structure sheaf $\mathcal{O}_{\mathfrak{Y}}$ is a discrete valuation ring. Let $\mathfrak{m}_y$ denote its maximal ideal. The divisor $\divi(f)$ of any nonzero function $f \in \mathcal{O}_L(\mathfrak{Y})$ is defined to be the function $\divi(f) : \mathfrak{Y} \rightarrow \mathbb{Z}_{\geq 0}$ given by $\divi(f)(y) = n$ if and only if the germ of $f$ in $y$ lies in $\mathfrak{m}_y^n \setminus \mathfrak{m}_y^{n+1}$. By Lemma 1.1 in (loc.\ cit.) for any affinoid subdomain $\mathfrak{Z}\subseteq \mathfrak{Y}$ the set
\begin{equation}\label{f:divisorfinite}
  \{x\in \mathfrak{Z}| \divi(f)>0\} \mbox{ is finite.}
\end{equation}

\begin{lemma}\label{divisibility}
For any two nonzero functions $f_1, f_2 \in \mathcal{O}_L(\mathfrak{Y})$ we have $f_2 \in f_1 \mathcal{O}_L(\mathfrak{Y})$ if and only if $\divi(f_2) \geq \divi(f_1)$.
\end{lemma}
\begin{proof}
We consider the principal ideal $f_1 \mathcal{O}_L(\mathfrak{Y})$. As a consequence of \cite[Prop.\ 1.6 and Prop.\ 1.4]{BSX}  we have
\begin{equation*}
  f_1 \mathcal{O}_L(\mathfrak{Y}) = \{ f \in \mathcal{O}_L(\mathfrak{Y}) \setminus \{0\} : \divi(f) \geq \divi(f_1) \} \cup \{0\}.
\end{equation*}
\end{proof}

We now apply these results to exhibit a few more explicit elements in the distribution algebra $D(U_n,L)$, which will be used later on.

\begin{lemma}\label{fraction}
For any $1 \neq u \in U_n$ the fraction $\frac{\nabla}{\delta_u - 1}$ is a well defined element in the integral domain $D(U_n,L)$.
\end{lemma}
\begin{proof}
By the Fourier isomorphism we may equivalently establish that the fraction $\frac{F_\nabla}{\ev_u - 1}$ exists in $\mathcal{O}_L(\mathfrak{X}_u^\times)$. But for this we only need to combine the Lemmas \ref{Dirac-zeros}, \ref{nabla-zeros}, and \ref{divisibility}.
\end{proof}

The next elements will only lie in the Robba ring of $U_n$. Since $\mathfrak{X}_n^\times \cong \mathfrak{X}$ we deduce from Prop.\ \ref{quasi-Stein} and the subsequent discussion that there is an admissible covering $\mathfrak{X}_n^\times = \bigcup_{j \geq 1} \mathfrak{V}_{n,j}$ by an increasing sequence $\mathfrak{V}_{n,1} \subseteq \ldots \subseteq \mathfrak{V}_{n,j} \subseteq \ldots$ of affinoid subdomains $\mathfrak{V}_{n,j}$  with the following properties:
\begin{itemize}
  \item[--] The system $(\mathfrak{X}_n^\times \setminus \mathfrak{V}_{n,j})_{/\mathbb{C}_p}$ is isomorphic to an increasing system of one dimensional annuli. This implies:
      \begin{itemize}
        \item[--] $\mathcal{R}_L(\mathfrak{X}_n^\times)$ is the increasing union of the rings $\mathcal{O}_L(\mathfrak{X}_n^\times \setminus \mathfrak{V}_{n,j})$ and contains $\mathcal{O}_L(\mathfrak{X}_n^\times)$;
        \item[--] each $\mathcal{O}_L(\mathfrak{X}_n^\times \setminus \mathfrak{V}_{n,j})$ as well as $\mathcal{R}_L(\mathfrak{X}_n^\times)$ are integral domains.
      \end{itemize}
  \item[--] Each $\mathfrak{X}_n^\times \setminus \mathfrak{V}_{n,j}$ is a one dimensional smooth quasi-Stein space.
\end{itemize}
In particular, the $\mathcal{O}_L(\mathfrak{X}_n^\times \setminus \mathfrak{V}_{n,j})$ are naturally Fr\'echet algebras, and we may view $\mathcal{R}_L(\mathfrak{X_n^\times})$ as their locally convex inductive limit. We also conclude that Lemma \ref{divisibility} applies to each $\mathfrak{X}_n^\times \setminus \mathfrak{V}_{n,j}$.

We now fix a basis $b = (b_1, \ldots, b_d)$ of $U_n$ as a $\mathbb{Z}_p$-module such that $b_i \neq 1$ for any $1 \leq i \leq d$.

\begin{proposition}\label{Xi}
The fraction
\begin{equation*}
  \Xi_b :=  \frac{F_\nabla^{d-1}}{\prod_{i=1}^d (\ev_{b_i} - 1)}
\end{equation*}
is well defined in the Robba ring $\mathcal{R}_L(\mathfrak{X}_n^\times)$.
\end{proposition}
\begin{proof}
The zeros of the fraction $\frac{F_\nabla}{\ev_{b_i} - 1} \in \mathcal{O}_L(\mathfrak{X}_n^\times)$ are precisely those finite order characters which are nontrivial on $b_i$. Hence, if we fix a $1 \leq j \leq d$, then the product $\prod_{i \neq j} \frac{F_\nabla}{\ev_{b_i} - 1}$ still has a zero in any finite order character which is nontrivial on $b_i$ for at least one $i \neq j$. On the other hand the zeros of $\ev_{b_j} - 1$ are those finite order characters which are trivial on $b_j$ (and they have multiplicity one). Since only the trivial character is trivial on all $b_1, \ldots, b_d$ we see that all zeros of $\ev_{b_j} - 1$ with the exception of the trivial character occur also as zeros of the product $\frac{F_\nabla^{d-1}}{\prod_{i \neq j} (\ev_{b_i} - 1)}$. It follows that the asserted fraction $\Xi_b$ exists in $\mathcal{O}_L(\mathfrak{X}_n^\times \setminus \mathfrak{V}_{n,j})$ provided $j$ is large enough so that the trivial character is a point in $\mathfrak{V}_{n,j}$. Since $(\prod_{i=1}^d (\ev_{b_i} - 1))\Xi_b = F_\nabla^{d-1}$ and $\mathcal{R}_L(\mathfrak{X}_\Gamma^\times)$ is an integral domain, we see the independence of $j$.
\end{proof}

In fact, the proof of Prop.\ \ref{Xi} shows that $\Xi_b$ is a meromorphic function on $\mathfrak{X}_n^\times$ with a single pole at the trivial character, which moreover is a simple pole.  We abbreviate $\ell(b) := \prod_{i=1}^d \log(b_i)$.

\begin{proposition}\label{dep-b}
  For any other basis $b' = (b_1', \ldots, b_d')$ of $U_n$ as a $\mathbb{Z}_p$-module with $b_i' \neq 1$ we have
\begin{equation*}
  \ell(b') \Xi_{b'} - \ell(b) \Xi_b \in \mathcal{O}_L(\mathfrak{X}_n^\times) \ .
\end{equation*}
\end{proposition}
\begin{proof}
We only have to check that the asserted difference does not any longer have a pole at the trivial character. Both, $\Xi_{b'}^{-1}$ and $\Xi_b^{-1}$, are uniformizers in the local ring $\mathcal{O}_1$ of $\mathfrak{X}_n^\times$ in the trivial character. Hence we have in $\mathcal{O}_1$ an equality of the form
\begin{equation*}
  \frac{\Xi_{b'}}{\Xi_b} = x + \Xi_b^{-1} \cdot G
\end{equation*}
with some $x \in L$ and $G \in \mathcal{O}_1$. Our assertion amounts to the claim that $x = \prod_i \frac{\log(b_i)}{\log(b_i')}$. To compute $x$ we use \eqref{f:small-disk} which leads to the open embedding
\begin{equation*}
  \mathbf{B}(r_0)_{/L} \xrightarrow{\cong} \mathfrak{X}(r_0) \xrightarrow{\subseteq} \mathfrak{X} \xrightarrow{\ell_n^*} \mathfrak{X}_n^\times
\end{equation*}
which maps $y$ to the character $\chi_y(u) := \exp(\pi_L^{-n} \log(u)y)$ ( and, in particular, $0$ to the trivial character). Using \eqref{f:Fnabla} we see that $F_\nabla$ pulls back to the function $y \mapsto \pi_L^{-n} y$ on $\mathbf{B}(r_0)$. On the other hand $\ev_{b_i} - 1$ pulls back to $y \mapsto \exp(\pi_L^{-n} \log(b_i)y) - 1$. Hence $\Xi_b$ pulls back to the meromorphic function
\begin{equation*}
  y \longmapsto \frac{\pi_L^{-n(d-1)} y^{d-1}}{\prod_i (\exp(\pi_L^{-n}\log(b_i)y) - 1)} \ .
\end{equation*}
Its germ at zero lies in $\frac{1}{\pi_L^{-n} (\prod_i \log(b_i))y} (1 + y \mathcal{O}_0)$, where $\mathcal{O}_0$ denotes the local ring of $\mathbf{B}(r_0)_{/L}$ in zero. It follows that the germ of the pull back of $\frac{\Xi_{b'}}{\Xi_b}$ lies in $(\prod_i \frac{\log(b_i)}{\log(b_i')}) (1 + y \mathcal{O}_0)$.
\end{proof}

By Lemma \ref{nabla-zeros} the function $F_\nabla \Xi_b$ is holomorphic on $\mathfrak{X}_n^\times$ and has no zero in the trivial character.

\begin{lemma}\label{val-b}
The value of $F_\nabla \Xi_b$ at the trivial character is $\ell(b)^{-1}$.
\end{lemma}
\begin{proof}
We use the same strategy as in the previous proof. The function $\Theta_b$ pulls back to the function
\begin{equation*}
  y \longmapsto \frac{\pi_L^{-nd} y^d}{\prod_i (\exp(\pi_L^{-n}\log(b_i)y) - 1)}
\end{equation*}
on $\mathbf{B}(r_0)_{/L}$, and we have to compute its value at $0$. But visibly the above right hand side is a power series in $y$ with constant term $\frac{1}{\prod_{i=1}^d \log(b_i)}$.
\end{proof}

These last two facts suggest to renormalize our functions by setting
\begin{equation*}
  \overline{\Xi}_b := \ell(b) \Xi_b \qquad\text{and}\qquad  \Theta_b := F_\nabla \overline{\Xi}_b \ .
\end{equation*}
Choosing a field $K$ containing $\Omega$ we also let $\widetilde{\Xi_b}$ denote the image of $\overline{\Xi}_b$ under the composite map
\begin{equation}\label{f:tildeXi}
  \mathcal{R}_L(\mathfrak{X}_n^\times) \xrightarrow{\ell_{n*}} \mathcal{R}_L(\mathfrak{X}) \subseteq \mathcal{R}_K(\mathfrak{X}) \xrightarrow{\kappa^*} \mathcal{R}_K(\mathbf{B}) \ .
\end{equation}

\begin{remark}\label{rem:Xi}
Suppose that $K$ contains $\Omega$. We have
\begin{equation*}
  \widetilde{\Xi_b} = \frac{\ell(b) (\frac{\Omega}{\pi_L^n}\log_{LT}(Z))^{d-1}}{\prod_j (\exp(\log(b_j)\frac{\Omega}{\pi_L^n}\log_{LT}(Z))-1)} \ ,
\end{equation*}
and it follows from the proof of Prop.\  \ref{Xi} that $Z\widetilde{\Xi_b}$ belongs to ${\mathcal{O}_K(\mathbf{B})}$ with constant term $(\frac{\Omega}{\pi_L^n} )^{-1}$, whence
\begin{equation*}
  \widetilde{\Xi_b} \equiv \frac{\pi_L^n}{\Omega   Z} \bmod {\mathcal{O}_K(\mathbf{B})}.
\end{equation*}
\end{remark}
\begin{proof}
One checks that the map \eqref{f:tildeXi} sends a distribution $\mu$ to the map
\begin{equation*}
  g_\mu(z)=\mu(\exp(\Omega\frac{\log(-)}{\pi_L^n}\log_{LT}(z))) \ .
\end{equation*}
In particular, a Dirac distribution $\delta_a$ is sent to $\exp(\Omega\frac{\log(a)}{\pi_L^n}\log_{LT}(z)).$ Recall that the action of $\nabla$ as distribution sends a locally $L$-analytic function $f$ to $-\left(\frac{d}{dt}f(\exp(-t))\right)_{|t=0}, $ whence $\nabla$ is sent to $$\nabla\left( \exp(\Omega\frac{\log(-)}{\pi_L^n}\log_{LT}(z))\right)=-\left(\frac{d}{dt}\exp(\Omega\frac{\log(\exp(-t))}{\pi_L^n}\log_{LT}(z))\right)_{|t=0}=\frac{\Omega}{\pi_L^n}\log_{LT}(z).$$
\end{proof}

\begin{remark}\label{rem:ThetaMellin}
Recall that $\Theta_b$ lies in $\mathcal{O}_L(\mathfrak{X}_n^\times)$ and therefore can be viewed, via the Fourier transform, as a distribution in $D(U_n,L) \subseteq D(o_L^\times,L)$. If $K$ contains $\Omega$ and for sufficiently large $n$ the Mellin transform $\mathfrak{M}$ in Lemma \ref{Mellin} then satisfies
\begin{equation*}
  \kappa^* \circ \mathfrak{M}({\Theta_b} )= \varphi_L^n(\xi_b )\eta(1,Z)
\end{equation*}
with
\begin{equation*}
  \xi_b \equiv \frac{\log_{LT}(Z)}{Z} \bmod \log_{LT}(Z){\mathcal{O}_K(\mathbf{B})} \ .
\end{equation*}
\end{remark}
\begin{proof}
Consider the the element
\[F(X)=\frac{X}{\exp(X)-1}=1+XQ(X)\] with $Q(X)\in\qp[[X]]$  and let $r>0$ be such that $Q(X)$ converges on $|X|\leq r.$ We shall proof the claim within the Banach algebra $\cR^I_K(\mathbf{B})$ for $I=[0,r]$ (which contains ${\mathcal{O}_K(\mathbf{B})}$ and using that the actions on both rings are compatible). We assume for the operator norm that $\|\delta_{b_i}-1\|_I<\min(p^{-\frac{1}{p-1}},r)$ for all $i$ (otherwise we enlarge $n$ according to Lemma \ref{lem:Ked}). From \cite[Cor.\ 2.3.2, proof of Lem.\ 2.3.1]{BSX} it follows that $\nabla=\frac{\log(\delta_{b_i})}{\log(b_i)}$ as operators in the Banach algebra $A$ of continuous linear endomorphisms of $\cR^I_K(\mathbf{B})$ and
\begin{equation}\label{f:expnabla}
  \exp(\log(b_i)\nabla)=\exp(\log(\delta_{b_i}))=\delta_{b_i}
\end{equation}
in $A.$ Moreover,
\begin{equation}\label{f:opnorm}
 \|\log(\delta_{b_i})\|_I<\min(p^{-\frac{1}{p-1}},r)
\end{equation} for all $i$, whence $\|\nabla\|_I<\min(p^{-\frac{1}{p-1}},r)|\log(b_i)|.$ Then, as operators in $A$ we have
\begin{equation}
  \log(b_i)^{-1}+\nabla Q(\log(b_i)\nabla)  = \log(b_i)^{-1}F(\log(b_i)\nabla)=\frac{\nabla}{\exp(\log(b_i)\nabla)-1}=\frac{\nabla}{\delta_{b_i}-1}.
\end{equation}
Hence
\[\Theta_b = \frac{\ell(b) \nabla^{d}}{\prod (\exp(\log (b_j)\nabla)-1)}= 1 + \ell(b) \nabla g(\log (b_j)\nabla)\]
for some power series $g\in \cR^I_K(\mathbf{B})$. It  follows that
\begin{equation}\label{f:MTheta}
   \kappa^* \circ \mathfrak{M}({\Theta_b}) = \big(1 + \ell(b) \Omega \log_{LT}(Z)f(Z) \big)\eta(1,Z).
\end{equation}
for some $f(Z)\in  \cR^I_K(\mathbf{B})$ \Footnote{we know a priori that the eigenvalue lies in ${\mathcal{O}_K(\mathbf{B})}$}. Indeed, concerning the derived action we have
$$\nabla\left( \eta(1,Z)\right)=\left(\frac{d}{dt}\exp(\Omega\exp(t)\log_{LT}(Z))\right)_{|t=0}=\Omega \log_{LT}(Z)\eta(1,Z)$$ (cf.\ also \cite[end of \S 2.3]{BSX} for the fact that
\begin{align}\label{f:nablaonR}
  \nabla \mbox{ acts as }\log_{LT}(Z)\partial_\mathrm{inv} \mbox{  on } {\mathcal{O}_K(\mathbf{B})} \ )
\end{align}
and
\[\nabla(\Omega \log_{LT}(Z))=\Omega \log_{LT}(Z).\] Furthermore, we obtain inductively that
\[\nabla^i\eta(1,Z)=\left(\prod_{k=0}^{i-1}(\Omega \log_{LT}(Z)+k)\right)\eta(1,Z)\]
for all  $i\geq 0.$ The convergence of $f(Z)$ can be deduced using the operator norm \eqref{f:opnorm}.

On the other hand, according to \cite[Lem.\ 2.4.2]{BF} we have
\[\Theta_b\eta(1,Z) = \frac{\ell(b) \log_{LT}(Z)}{\varphi_L^n(Z)}g(Z)\] for some $g(Z)\in{\mathcal{O}_K(\mathbf{B})}.$ Since the element $\Theta_b\eta(1,Z)$ lies in $({\mathcal{O}_K(\mathbf{B})})^{\psi_L=0}$, we conclude from
\begin{align*}
0=\psi_L(\frac{\pi_L^{-1}\varphi_L(\log_{LT}(Z))}{\varphi_L^n(Z)}g(Z))=\frac{\pi_L^{-1}\log_{LT}(Z)}{\varphi_L^{n-1}(Z)}\psi_L(g(Z))
\end{align*}
that $g(Z)$ belongs to $({\mathcal{O}_K(\mathbf{B})})^{\psi_L=0}$, whence it is of the form $\sum_{a\in(o_L/\pi_L)^\times} \varphi_L(g_a(Z))\eta(a,Z) $ for some $g_a\in{\mathcal{O}_K(\mathbf{B})}$ by the analogue of \eqref{f:Mspiequal0} for ${\mathcal{O}_K(\mathbf{B})}.$ From Lemma \ref{lem:varphinR} we derive that, for some $a(Z)\in{\mathcal{O}_K(\mathbf{B})},$ we have
\[\Theta_b\eta(1,Z) = \ell(b) \varphi_L^n(a(Z))\eta(1,Z). \]
Since the decomposition in \eqref{f:Mspiequal0} is direct, we conclude that $g(Z)=\varphi_L(g_1(Z))\eta(1,Z)$ and $\frac{\log_{LT}(Z)}{\varphi_L^n(Z)}\varphi_L(g_1(Z))=\varphi_L^n(a(Z)),$ whence $\log_{LT}(Z)$ divides $\varphi_L^n(a(Z)Z)$. Since $\varphi_L^n$ sends the zeroes of $\log_{LT}(Z)$, i.e., the points in $LT(\pi_L)=\bigcup_{k} LT[\pi_L^k]$, surjectively onto itself, we conclude by Lemma \ref{divisibility}   that $\log_{LT}(Z)$ divides also $a(Z)Z $ in ${\mathcal{O}_K(\mathbf{B})}$  and that there exists $c(Z)\in {\mathcal{O}_K(\mathbf{B})}$ such that
\begin{equation}\label{f:Thetaeta}
  \kappa^* \circ \mathfrak{M}({\Theta_b}) = \ell(b) \varphi_L^n\left(\frac{\log_{LT}(Z)}{Z}c(Z)\right)\eta(1,Z).
\end{equation}
 Comparing \eqref{f:Thetaeta} with the first description \eqref{f:MTheta} gives the claim as $c(0) = \ell(b)^{-1}$ because evaluation at $0$ is compatible with the embedding ${\mathcal{O}_K(\mathbf{B})}\subseteq \cR^I_K(\mathbf{B})$ and $\frac{\log_{LT}(Z)}{Z}(0) = 1$ by \eqref{f:tLT}.
\end{proof}

\newpage

\subsection{Pairings}\label{sec:pairings}

In this section we discuss various kinds of pairings. 
The starting point is Serre duality on $\mathfrak{X}$ which induces a (residue) pairing
\[  < \ ,\ >_{\mathfrak{X}}: \cR_{L}(\mathfrak{X}) \times \cR_{L}(\mathfrak{X}) \to L,\]
as we have seen in \eqref{f:pairing-RobbaX}.
Similarly, Serre duality on $\mathfrak{X}^\times $
 induces  a pairing
\begin{equation}\label{f:bothpair}
   < \ ,\ >_{\Gamma_L}: \cR_L(\Gamma_L) \times \cR_L(\Gamma_L) \to L
\end{equation}
for the Robba ring of $\Gamma_L$, which by definition is the Robba ring of its character variety $\mathfrak{X}_{\Gamma_L}\cong\mathfrak{X}^\times $ (induced by  the isomorphism $\chi_{LT}:\Gamma_L\to o_L^\times$) as constructed in \eqref{f:pairing-mult-charvar}. This pairing, as defined in subsection \ref{subsec:res-group}, is actually already characterized by its restriction to $\cR_L(\Gamma_n)$ for any $n\geq n_0$ and thus is by construction  and the functoriality properties of section \ref{sec:duality} closely related to the pairing $ < \ ,\ >_{\mathfrak{X}}$ using the 'logarithm' $\cR_L(\Gamma_n)\xrightarrow{(\ell_n)_*}\cR_{L}(\mathfrak{X})$, see diagram \eqref{f:XmultXn}.

In contrast, the commutative diagram
\[ \xymatrix{
   \cR_L(\mathfrak{X}^\times) \ar[dd]_{(-)(\ev_1d\log_\mathfrak{X})} \ar[r]^(0.6){\cdot d\log_{\mathfrak{X}^\times }} & \Omega^1_{\mathfrak{X}^\times} \ar[d]^{\mathrm{res}_{\mathfrak{X}^\times}} \\
       & L  \\
   (\Omega^1_{\mathfrak{X}})^{\psi=0} \ar[r]^(0.6){\ev_{-1}\cdot} & \Omega^1_{\mathfrak{X} } \ar[u]_{\mathrm{res}_{\mathfrak{X}}}.  }\]
from Thm.\  \ref{thm:residuumidentity} in subsection \ref{subsec:res-alternative} relates  the pairing $  < \ ,\ >_{\Gamma_L} $
to the pairing $< \ ,\ >_{\mathfrak{X}}$ 
{\it in a highly non-trivial, non-obvious way}. The resulting description of $ < \ ,\ >_{\Gamma_L} $ in \eqref{f:defpairing<<<}  forms one main ingredient in the proof of the abstract reciprocity formula \ref{thm-recproclawKKK}   below in subsection \ref{subsec:abstract-rec}.

Based on the (generalized) residue pairings \eqref{f:res-pair-general} in subsection \ref{subsec:res-modules}
\begin{equation*}
\{\;,\;\}: \check{M} \times M \to L,
\end{equation*}
with $\check{M}:=Hom_{\cR}(M,\Omega^1_{\cR})$
 the pairing \eqref{f:bothpair} induces for any (analytic) $(\varphi_L,\Gamma_L)$-module $M$ over $\cR_{L}(\mathfrak{X})$ an  Iwasawa pairing \eqref{f:IWpairing}
\begin{equation*}
\KKl\;,\;\KKr_{Iw}: \check{M}^{\psi_L=\frac{q}{\pi_L}} \times M^{\psi_L=1} \to D(\Gamma_L,L)
\end{equation*}
in subsection \ref{subsec:IwasawaPairing}, which   behaves  well with twisting (cf.\ Lemma \ref{lem:twist}).


By construction and the comparison isomorphism \eqref{f:compM} for Kisin-Ren modules - the second main ingredient -
the pairing $ \KKl\;,\;\KKr_{Iw}$ is closely related to a pairing
\[[\;,\;]:{\cR_{L}(\mathfrak{X})^{\psi_L=0}\otimes_L D_{cris,L}(V^*(1))} \times  {\cR_{L}(\mathfrak{X})^{\psi_L=0}\otimes_L D_{cris,L}(V(\tau^{-1}))}  \to \cR_L(\Gamma_L)\]
induced from the natural pairing for $D_{cris,L}$. The precise relationship is the content of an abstract form of a reciprocity formula, see Thm.\  \ref{thm-recproclawKKK}. As a consequence we shall later derive a concrete reciprocity formula, i.e.,  the adjointness of Berger's and Fourquaux' big exponential map with our regulator map, see Thm.\  \ref{thm:adjointness}.

\subsubsection{The residuum pairing for modules}\label{subsec:res-modules}

Throughout our coefficient field $K$ is a complete intermediate extension $L \subseteq K \subseteq \mathbb{C}_p$. Let $\mathfrak{Y}$ be either $\mathfrak{X} $ or $\mathbf{B}$ and $\cR=\cR_{K}(\mathfrak{Y})$. Consider the residuum map $\mathrm{res}_\mathfrak{X}$ defined after \eqref{f:RobbaX-selfdual} and the residuum map
\[\mathrm{res}_\mathbf{B}: \Omega^1_{\cR} \to K,\;\;\; \sum_i a_iZ^idZ\mapsto a_{-1}.\]
Recall that we are using the operator $\psi_L:=\frac{q}{\pi}\psi_L^\mathfrak{Y}$ on $\cR.$

Moreover, we define $\upiota_*: \cR_K(\Gamma_L)\to  \cR_K(\Gamma_L)$ to be the map which is induced by sending $\gamma\in\Gamma_L$ to its inverse $\gamma^{-1},$ i.e.,  the involution of the group induces an isomorphism on the multiplicative character variety, which in turn gives rise to $\upiota_*$. The corresponding involution on $\cR_K(\Gamma_{n_0})$, also denoted by $\upiota_*$, satisfies the commutative diagram
\begin{equation}\label{f:involution}
\xymatrix{
  \cR_K(\Gamma_{n_0}) \ar[d]_{\upiota_*}^{\cong} \ar[r]^{\hat{\ell}_{{n_0}*}}_{\cong} & \cR_K(\mathfrak{X}) \ar[d]^{\upiota}_{\cong} \\
   \cR_K(\Gamma_{n_0}) \ar[r]^{\hat{\ell}_{{n_0}*}}_{\cong} & \cR_K(\mathfrak{X})   }
\end{equation}
where the involution $\upiota$ on $\cR_K(\mathfrak{X})$ sends $\ev_{x}$ to $\ev_{-x}$.

Setting $\check{M}:= \Hom_{\cR}(M,\Omega^1_{\cR}) \cong \Hom_{\cR}(M,{\cR})(\chi_{LT})$, for any finitely generated projective $\cR$-module $M$, we obtain more generally the pairing
\begin{equation}\label{f:res-pair-general}
   \{\;,\;\}:=\{\;,\;\}_M: \check{M} \times M \to K,\;\;\; (g,f)\mapsto \mathrm{res}_\mathfrak{Y}(g(f)),
\end{equation}
\Footnote{Note that Colmez defines $\Omega \mathrm{res}_\mathbf{B}(\sigma_{-1}(g)(f ) )$ instead.} which satisfies the following properties:

\begin{lemma}\phantomsection\label{lem-pair} For $M$ in $\mathcal{M}(\cR)$ we have
\begin{enumerate}
\item $ \{\;,\;\}$ identifies $M$ and $\check{M}$ with the   (strong)   topological duals of $\check{M}$ and $M$, respectively.
\item $\{\varphi_L(g),\varphi_L(f)\}=\frac{q}{\pi_L}\{g,f\}$ for all $g\in\check{M}$ and $f\in M,$
\item $\{ \sigma(g),\sigma(f)\}=\{g,f\}$ for all $g\in\check{M}$, $f\in M,$ and $\sigma\in \Gamma_L,$
\item  $\{\psi_L(g),f\}=\{g,\varphi_L(f)\}$ and $\{\varphi_L(g),f\}=\{g,\psi_L(f)\}$ for all $g\in\check{M}$ and $f\in M.$
\end{enumerate}
\end{lemma}
\begin{proof}
(i) follows from the discussion in subsection \ref{subsec:boundary}. (ii) is a purely formal consequence from (iv). (iii) follows as in \eqref{f:vdP} with $\sigma_*$ instead of $\kappa_*.$  For (iv) we refer to Lemma \ref{RobbaX-pairing-phi-psi}.  For $\mathfrak{Y}=\mathbf{B}$  see also \cite[Prop.\ 3.17, Cor.\ 3.18, Prop.\ 3.19]{SV15}.
\end{proof}

\noindent
{\bf Convention:} For coherence of our notation we set $\log_\mathbf{B}:=\log_{LT}$ although in general this is {\it not} the standard logarithm!\\

\begin{proposition}\label{prop:residuepairing}
The pairing $< \ ,\ >_{\mathfrak{Y}}: \cR \times \cR \to K,$ $(f,g)\mapsto \mathrm{res}_{\mathfrak{Y}}(fgd\log_{\mathfrak{Y}})$, induces topological isomorphisms
\begin{equation*}
  \Hom_{K,cts}(\cR,K)\cong \cR \qquad\text{and}\qquad \Hom_{K,cts}(\cR/{\mathcal{O}_K(\mathfrak{Y})},K)\cong {\mathcal{O}_K(\mathfrak{Y})}.
\end{equation*}
\end{proposition}
\begin{proof}
See subsection \ref{subsec:boundary}.
\end{proof}

\begin{remark}\label{rem:compPairinigXB}
If we assume $\Omega\in K$, then these pairings can be compared via the LT-iso\-morphism $\kappa.$ By   \eqref{f:res-XB} we have
\[\Omega< \kappa^*(f) ,\kappa^*(g) >_{\mathbf{B}}=< f ,g>_{\mathfrak{X}}\]
for $f,g\in \cR_{K}(\mathfrak{X}).$
\end{remark}

{\bf Assume} henceforth that $M$ is an analytic $(\varphi_L,\Gamma_L)$-module over $\cR$ and recall from Proposition  \ref{prop:distributionaction} that the $\Gamma_L$-action on $M$ extends continuously to a $D(\Gamma_L,K)$-module structure.

\begin{corollary}\phantomsection\label{cor:D-equiv}
 The isomorphism $\check{M} \cong \Hom_{K,cts}(M,K)$ (induced by $\{\;,\;\}$) is $D(\Gamma_L,K)$-linear.
\end{corollary}
\begin{proof}
 This follows from Lemma \ref{lem-pair}(iii) since $\Gamma_L$ generates a dense subspcae of $D(\Gamma_L,K)$.
\end{proof}

Since $\frac{\pi_L}{q}\psi_L\circ \varphi_L= \id_M$ we have a canonical decomposition $M=\varphi_L(M)\oplus M^{\psi_L=0}.$ By Lemma \ref{lem-pair} we see that $M^{\psi_L=0}$ is the exact orthogonal complement of $\varphi_L(\check{M}),$ i.e., we obtain a canonical isomorphism
\begin{equation}\label{f:psizero}
   \check{M}^{\psi_L=0} \cong \Hom_{K, cts}(M^{\psi_L=0},K) .
\end{equation}

\begin{lemma}\label{lem:sesquilinear}
The isomorphism \eqref{f:psizero} is $\cR_K(\Gamma_L)$-equivariant, i.e., we have for all $\check{m}\in \check{M}^{\psi_L=0}$, $m\in M^{\psi_L=0}$, and $\lambda\in \cR_K(\Gamma_L)$ that
\begin{equation*}
  \{\lambda\check{m},m\}=\{\check{m},\upiota_*(\lambda)m\} \ .
\end{equation*}
\end{lemma}
\begin{proof}
This is clear for $D(\Gamma_L,K)$ by Cor.\ \ref{cor:D-equiv}. Without loss of generality we may and do assume that $\Omega$ belongs to $K.$ It then follows for the localization $D(\Gamma_L,K)_{Y_{n_1}^{\mathbb{N}}}$, where we use the notation and considerations from subsection \ref{subsec:phiGamma}, especially Lemma \ref{lem:technical} and its proof. Since $D(\Gamma_L,K)_{Y_{n_1}^{\mathbb{N}}}$ is dense in $\cR_K(\Gamma_L)$ by \ref{rem:denseX} the assertion now is a consequence of the continuity property in Thm.\  \ref{thm:Mpsiequal0}.
\end{proof}

\subsubsection{The duality pairing \texorpdfstring{$< , >_{\Gamma_L}$}{<,>} for the group Robba ring}\label{subsec:res-group}

Using the isomorphisms \eqref{f:chi-iso} induced by the Lubin-Tate character $\chi_{LT}$ we now carry over structures concerning the (multiplicative) character varieties $\mathfrak{X}^\times, \mathfrak{X}^\times_n$ to those of the groups $\Gamma_L, \Gamma_n$. In particular, we use analogous  notation $\mathrm{res}_{\Gamma_L}, \mathrm{res}_{\Gamma_n},\log_{\Gamma_L},\log_{\Gamma_n}$ for corresponding objects. In this sense we introduce and recall from \eqref{f:pairing-mult-charvar} the pairing
\begin{equation}
\label{f:defpairinig} <\;,\;>_{\Gamma_L}: \cR_K(\Gamma_L) \times \cR_K(\Gamma_L) \xrightarrow{ }  K
\end{equation}
and similarly $<\;,\;>_{\Gamma_n} $ from \eqref{f:XmultXn}.
 This pairing is of the form
\begin{equation}
\label{f:defpairinig2}<\;,\;>_{\Gamma_L}: \cR_K(\Gamma_L) \times \cR_K(\Gamma_L) \xrightarrow{mult}\cR_K(\Gamma_L) \xrightarrow{ \varrho} K,
\end{equation}
where\Footnote{In \cite[Rem.\ 2.3.12]{KPX} the factor $-\frac{p-1}{p}$ shows up for the correct normaliszation of the Residuum map to be compatible with Tate's trace map on $H^2$ in Galois cohomology! Are these occurrences related? Probably this factor will effect the   Tate duality pairing and interpolation property. } 
\begin{align*}
 \varrho=<1,->_{\Gamma_L}:\cR_K(\Gamma_L) & \to \Omega^1_{\cR_K(\Gamma_L)} \xrightarrow{\mathrm{res}_{\Gamma_L}} K \\
  f & \mapsto f d\log_{\Gamma_L}\mapsto \mathrm{res}_{\Gamma_L}(f d\log_{\Gamma_L})
\end{align*}
has also the following description - writing  $pr_{n,m} $  and similarly $pr_{L,m} $ for the projection maps induced by \eqref{f:cross},\eqref{f:Rcrossproduct} -
\begin{align}\label{f:groupRobba-pairing-explizit}
\notag \varrho:\cR_K(\Gamma_L) & =\mathbb{Z}[\Gamma_L] \otimes_{\mathbb{Z}[\Gamma_{n_0}]}\cR(\Gamma_{n_0}) \xrightarrow{} K \\
  f & \mapsto \frac{q- 1}{q} (\frac{q}{\pi_L})^{n_0}\mathrm{res}_{\mathfrak{X}}(\hat{\ell}_{n_0*}\circ pr_{L,n_0}(f)d\log_\mathfrak{X})
\end{align}
with $n_0$ as defined at the beginning of subsection \ref{sec:groupRobba}.
Indeed, using \eqref{f:res-mult-charvar}, \eqref{f:XmultXn} we obain
\begin{align*}
<1,f>_{\Gamma_L}&= \mathrm{res}_{\Gamma_L}(fd\log_{\Gamma_L})\\
&=\frac{q- 1}{q}q^{n_0}\mathrm{res}_{\Gamma_{n_0}}(pr_{L,n_0}(f)d\log_{\Gamma_{n_0}})\\
&=\frac{q- 1}{q}(\frac{q}{\pi_L})^{n_0}\mathrm{res}_{\mathfrak{X} }(\hat{\ell}_{n_0*}\circ pr_{L,n_0}(f)d\log_{\mathfrak{X} })
\end{align*}
because $(\ell_n^*)^*(1)=1$.

The following properties follow immediately from the definition:
\begin{lemma}\label{lem:nabla}
We have for all $f,\lambda,\mu\in \cR_K(\Gamma_L)$ that
\begin{enumerate}
\item $< \lambda,f\mu>_{\Gamma_L}= < f\lambda,\mu >_{\Gamma_L},$
\item $<\lambda,\mu>_{\Gamma_L}=<\mu,\lambda >_{\Gamma_L}$.
\end{enumerate}
\end{lemma}

\begin{remark} For $n\geq n_0$ we have the projection formula $pr_{ L,n}(\iota_{n,*}(x)y)=xpr_{L,n}(y)$ and \eqref{f:RobbaX-pairing-Ln} translates into
\begin{equation}
\label{f:projformula} <\iota_{n,*}(x),y>_{\Gamma_L}\, = \,(q-1) q^{n-1}<x,pr_{L,n}(y)>_{\Gamma_n}
\end{equation}
for $x\in \cR(\Gamma_n)$, $y\in \cR(\Gamma_L)$ and the canonical inclusion $\cR(\Gamma_n)\xrightarrow{\iota_{ n,*} }\cR(\Gamma_L)$. Analogous formulae hold for $\Gamma_m$ with $n\geq m\geq n_0$instead of $ \Gamma_L$ by \ref{RobbaX-pairing-phi-psi} (ii).
\end{remark}

%

\begin{remark}[Frobenius reciprocity]\label{rem:Frobrec}
The projection map $pr_{\Gamma_L,\Gamma_n} : \cR_K(\Gamma_L) \rightarrow \cR_K(\Gamma_n)$ induces an isomorphism
\begin{equation*}
  \Hom_{\cR_K(\Gamma_L)}(N,\cR_K(\Gamma_L)) \cong \Hom_{\cR_K(\Gamma_n)}(N,\cR_K(\Gamma_n))
\end{equation*}
for any $\cR_K(\Gamma_L) $-module $N$; the inverse sends $f$ to the homomorphism $x \mapsto \sum_{g\in\Gamma_L/U} g\iota_{n,*}\circ f(g^{-1}x)$.
\end{remark}

The following proposition    translates the results at the end of subsection \ref{subsec:boundary} into the present setting.

\begin{proposition}\label{prop:residuepairingRGamma}
The pairing $<\;,\;>_{\Gamma_L}: \cR_K(\Gamma_L) \times \cR_K(\Gamma_L) \to K$ induces topological isomorphisms
\[ \Hom_{K,cts}(\cR_K(\Gamma_L),K)\cong \cR_K(\Gamma_L) \mbox{ and } \Hom_{K,cts}(\cR_K(\Gamma_L)/D(\Gamma_L,K),K)\cong D(\Gamma_L,K).\]
\end{proposition}

\begin{proposition}\label{prop:dualitygroupRobba}
Assume $\Omega\in K$ and    $M$ in $\mathcal{M}(\cR)$. Then the map
\begin{align}\label{f:psizeroequiv}
  \Hom_{\cR_K(\Gamma_L)}(M^{\psi_L=0},\cR_K(\Gamma_L))^\upiota & \xrightarrow{\;\cong\;} \Hom_{K, cts}(M^{\psi_L=0},K) \xrightarrow[\eqref{f:psizero}]{\cong} \check{M}^{\psi_L=0}  \\
  F & \longmapsto \rho \circ F             \nonumber
\end{align}
is an isomorphism of $\cR_K(\Gamma_L)$-modules, where the superscript ``$\upiota$'' on the left hand side indicates that $\cR_K(\Gamma_L)$ acts through the involution $\upiota_*$.
\end{proposition}
\begin{proof} According to Thm.\ \ref{thm:Mpsiequal0} the $\cR_K(\Gamma_L)$-module $ M^{\psi_L=0}$ is finitely generated free. Hence it suffices to show that the map
\begin{align*}
  \Hom_{\cR_K(\Gamma_L)}(\cR_K(\Gamma_L),\cR_K(\Gamma_L)) & \longrightarrow \Hom_{K, cts}(\cR_K(\Gamma_L),K)   \\
  F & \longmapsto \rho \circ F
\end{align*}
is bijective. But this map is nothing else than the duality isomorphism in Prop.\ \ref{prop:residuepairingRGamma}.
\end{proof}

%

\Footnote{
\begin{remark}
We rather can start with  the identification
\[\mathcal{R}(\Gamma_n) \cong \varphi_L^n(\mathcal{R}),\]
extending the isomorphism  $D(\Gamma_n,K)\xrightarrow{\ell_{*}} D(\pi_L^no_L,K).$
Then, for $m\geq n\geq n_0$ and $k:=m-n$, we have the commutative diagram
\begin{equation*}
\xymatrix{
  \cR(\Gamma_m) \ar[d]_{\iota'_{m,n}}^{injective} \ar[r]^{\ell_{*}}_{\cong} & \varphi_L^m(\cR) \ar[d]^{injective} \\
   \cR(\Gamma_n) \ar[r]^{\ell_{*}}_{\cong} & \varphi_L^n(\cR)   }
\end{equation*}
and we may also describe the pairing as
\begin{equation*}
<\;,\;>: \cR_K(\Gamma_L) \times \cR_K(\Gamma_L) \xrightarrow{mult}\cR_K(\Gamma_L)=\cR(\Gamma_{n_0})\rtimes_{\Gamma_{n_0}}\Gamma_L\xrightarrow{pr_{L,n_0}}\cR(\Gamma_{n_0})\xrightarrow{\ell}\varphi_L^{n_0}(\cR)\subseteq \cR\xrightarrow{ res } K,
\end{equation*}
 which obviously coincides with the one in the previous subsection (and, in particular, is independent of the choice of $n_0$).
\end{remark}}


\Footnote{
\paragraph{Twists}

One checks that the isomorphism
\[D(\Gamma_n,K)\xrightarrow{(\ell_n)_*}D(o_L,K)\cong\mathcal{O}_K(\mathfrak{X})\cong\mathcal{O}_K(\mathbf{B})\]
sends a distribution $\mu$ to the map $g_\mu(z)=\mu(\exp(\Omega\ell_n(-)\log_{LT}(z)))$. In particular, a Dirac distribution $\delta_\gamma$ is sent to $\exp(\Omega\ell_n(\gamma)\log_{LT}(z)).$

Let $z_{\chi}\in \mathbf{B}(K)$ be the point which corresponds via the LT-isomorphism to the character $\chi:o_L\to L^\times,$ $a\mapsto \exp(\pi_L^na),$ i.e., satisfying
\[\exp(\pi_L^na)=\exp(a\Omega\log_{LT}(z_{\chi}))\]
or
\[\pi_L^n=\Omega\log_{LT}(z_{\chi}).\]

\begin{lemma} \label{lem:twist-additive}
We have the commutative diagram
\[\xymatrix{
  D(\Gamma_n,K) \ar[d]_{Tw_{\chi_{LT}}} \ar[rrr]^-{LT\circ \mathrm{Fourier} \circ (\ell_n)_*} &&& \mathcal{O}_K(\mathbf{B}) \ar[d]^{(z_{\chi}+_{LT}Z)^*} \\
   D(\Gamma_n,K) \ar[rrr]^-{LT\circ \mathrm{Fourier} \circ (\ell_n)_*} &&& \mathcal{O}_K(\mathbf{B})   }\]
\end{lemma}
\begin{proof}
This follows from properties 1. and 3. in section \ref{subsec:twisting} together with the fact that the LT-isomorphism $\kappa : \mathfrak{X}_K \xrightarrow{\cong} \mathbf{B}_K$ is an isomorphism of group varieties.
\end{proof}

Since, with respect to the maximal ideal $\mathfrak{m}_K$ of $o_K$,
\[z_{\chi}+_{LT}Z\equiv 0+_{LT}Z =Z \mod \mathfrak{m}_K\]
we can write
\begin{eqnarray}
\label{f:Z'}
Z':=z_{\chi}+_{LT}Z=\alpha Z(1+\frac{\beta}{Z})
\end{eqnarray}
with $\alpha\in 1+\mathfrak{m}_K,$ $\beta\in z_\chi o_K[[Z]].$

Sending $Z$ to $z_{\chi}+_{LT}Z$ defines a continuous $K$-linear ring automorphism $\eta:\cR\to\cR$,  { What is the relation to $Tw_{\chi_{LT}}$ on $\cR(\Gamma_n)?$} which in turn induces an automorphism of $\Omega^1_\cR$ sending $f(Z)dZ$ to $\eta(F)d\eta(Z).$

\begin{lemma} \label{lem:res-teist} For all $\omega\in \Omega^1_\cR$ we have
\[Res(\eta(\omega))=Res(\omega).\]
\end{lemma}

\begin{proof}
By the same reasoning as in the proof of   \cite[Lem.\ 2.1.19]{KPX} this is reduced to the statement and proof of \cite[Rem.\ 3.4 ii.]{SV}. {
This is contained in \cite[Lem.\ 2.1.19]{KPX} for $\qp$-affinoid algebras instead of $K$, which does not apply here directly, but its strategy of proof.
For the convenience of the reader we sketch a proof. Since $d\cR$ is dense in  the kernel of $\mathrm{Res}_Z$ and $\mathrm{Res}_{Z'}$, these maps have the same kernel. It therefore suffices to show that $\mathrm{Res}_Z(\frac{dZ'}{Z'}) = 1$. From \eqref{f:Z'} we derive
\begin{equation*}
  \tfrac{dZ'}{Z'} = \tfrac{dZ}{Z} +  \tfrac{d(1+\frac{\beta}{Z})}{1+\frac{\beta}{Z}} \ .
\end{equation*}

We have $(1+\frac{\beta}{Z})^{p^m} = 1+\frac{ \beta_m'}{Z^{p^m}}$ for some $\beta_m' =\sum (\beta_m')_iZ^i\in o_K[[Z]]$ with $|(\beta_m')_i|\leq\max_{0\leq i\leq p^m}|\begin{pmatrix}
                                                                        p^m \\
                                                                        i \\
                                                                      \end{pmatrix}z_\chi^i
 |$. Therefore, if $m \geq 1$ is sufficiently big, then $\log(1+\frac{ \beta'_m}{Z^{p^m}})$ converges $p$-adically, which suffices by the reasoning of (loc.\ cit.) (but it is not clear whether it converges also in $\cR$).
It follows that $p^m\frac{d(1+\frac{\beta}{Z})}{1+\frac{\beta}{Z}} = \frac{d(1+\frac{ \beta_m'}{Z^{p^m}})}{1+\frac{ \beta_m'}{Z^{p^m}}} =  d(\log (1+\frac{ \beta_m'}{Z^{p^m}}))$ and therefore that $\mathrm{Res}_Z(\tfrac{d(1+\frac{\beta}{Z})}{1+\frac{\beta}{Z}}) = 0$.}
\end{proof}

\begin{lemma}\label{lem:diff-twist}
$\eta(d\log_{LT})=d\log_{LT}$
\end{lemma}

\begin{proof}
\begin{align*}
 d \eta(\log_{LT}(Z)) & =d\log_{LT}(z_{\chi}+_{LT}Z)\\
  &=d\log_{LT}(z_{\chi})+d\log_{LT}(Z)\\
  &=d\log_{LT}(Z)
\end{align*}
\end{proof}

From the definitions and Lemmata \ref{lem:twist-additive},  \ref{lem:res-teist} and \ref{lem:diff-twist} we conclude the following proposition.}
The following twist invariance is   just Lemma \ref{RobbaXtimes-pairing-twist}
\begin{proposition}\label{cor:twist<} Let $U$ be $\Gamma_L$ or $\Gamma_n$ for $n\geq n_0$. Then, for all $\lambda,\mu \in \cR(U)$ we have
 \[<Tw_{\chi_{LT}}(\mu),Tw_{\chi_{LT}}(\lambda)>_{U}=<\mu,\lambda>_{U}.\]
\end{proposition}

\subsubsection{A residuum identity and an alternative description of \texorpdfstring{$<\; ,\;>_{\Gamma_L}$}{<,>}}\label{subsec:res-alternative}
Let $\sigma_{-1}\in \Gamma_L$ be the element with $\chi_{LT}(\sigma_{-1})=-1.$ Consider the continuous $K$-linear map
 \begin{align*}
 \varsigma:\cR_K(\Gamma_L) &\to K, \\\lambda&\mapsto  \mathrm{res}_{\mathfrak{X} }(\mathfrak{M}(\sigma_{-1})\mathfrak{M}^{\Omega^1}(\lambda))
 \end{align*}
 where $\mathfrak{M}^{\Omega^1}:\cR_K(\Gamma_L) \xrightarrow{\cong}{\Omega^1_{\cR(\mathfrak{X})}}^{\psi=0} \subseteq\Omega^1_{\cR(\mathfrak{X})}$ sends $\lambda$ to \begin{align}\label{f:Twdlog}
  \lambda (\ev_1 d\log_{\mathfrak{X} })&=(Tw_{\chi_{LT}}(\lambda)(\ev_1))d\log_{\mathfrak{X}},
 \end{align}
 whence we also have
 \begin{align}
 \label{f:?}\mathrm{res}_{\mathfrak{X} }(\mathfrak{M}(\sigma_{-1})\mathfrak{M}^{\Omega^1}(\lambda))=\mathrm{res}_{\mathfrak{X} }(\mathfrak{M}(\sigma_{-1})\mathfrak{M} (Tw_{\chi_{LT}}(\lambda))d\log_{\mathfrak{X} }).
 \end{align}
Recall the definition from $\varrho$ from \eqref{f:defpairinig}.
\begin{theorem}\label{thm:residuumidentity} We have
  \[\varsigma= {\frac{q}{q-1} }\varrho  ,\] i.e., the following identity for the residue map holds
  \[(\frac{q}{\pi_L})^{n_0}\mathrm{res}_{\mathfrak{X}}\bigg(\hat{\ell}_{n_0*}\circ pr_{L,n_0}(\lambda)d\log_\mathfrak{X}\bigg)= \mathrm{res}_{\mathfrak{X} }\bigg(\ev_{-1}\lambda\bigg(\ev_1 d\log_{\mathfrak{X}}\bigg) \bigg)   \]
  for all $\lambda\in \cR_K(\Gamma_L) ,$ i.e., the following diagram commutes
 \[ \xymatrix{
   \cR_K(\mathfrak{X}^\times) \ar[dd]_{(-)(\ev_1d\log_\mathfrak{X})} \ar[r]^(0.6){\cdot d\log_{\mathfrak{X}^\times }} & \Omega^1_{\mathfrak{X}^\times} \ar[d]^{\mathrm{res}_{\mathfrak{X}^\times}} \\
       & K  \\
   (\Omega^1_{\mathfrak{X}})^{\psi=0} \ar[r]^(0.6){\ev_{-1}\cdot} & \Omega^1_{\mathfrak{X} } \ar[u]_{\mathrm{res}_{\mathfrak{X}}}.  }\]
\end{theorem}

\begin{remark}
Compare with  \cite[Prop.\ 2.2.1, 3.2.1]{benois} where also residue identities play a crucial role in the proof of his reciprocity formula.
\end{remark}

The proof of this theorem requires some preparation.

\begin{lemma}\label{lem:Tw<<<}
For all $\lambda\in \cR_K(\Gamma_L)$  and $j\in \mathbb{Z}$  we have
\[\varsigma(Tw_{\chi_{LT}^j}(\lambda))=\varsigma(\lambda).\]
\end{lemma}

\begin{proof} For the proof we may and do assume that $\Omega$ belongs to $K.$  Since then
 \[\mathrm{res}_\mathfrak{X}( \partial_\mathrm{inv}^\mathfrak{X}(f)d\log_{\mathfrak{X}})=\Omega\mathrm{res}_\mathbf{B}( \frac{1}{\Omega}\partial_\mathrm{inv}(\kappa^*(f))d\log_{LT})=  \mathrm{res}_\mathbf{B}( d\kappa^*(f))=0\] for any $f$ by Remark \ref{rem:logXlogLT}, \eqref{f:res-XB} and \cite[Prop.\ 2.12]{FX}, the case $j=1$ follows directly from the relation \eqref{f:?}
  using with $g:=Tw_{\chi_{LT}}(\lambda)$ that
\begin{align*}
\partial_\mathrm{inv}^\mathfrak{X}(\mathfrak{M}(\sigma_{-1})\mathfrak{M}(g))&=\partial_\mathrm{inv}^\mathfrak{X}(\mathfrak{M}(\sigma_{-1}))\mathfrak{M}(g)+\mathfrak{M}(\sigma_{-1})\partial_\mathrm{inv}^\mathfrak{X}(\mathfrak{M}(g))\\
&\stackrel{\text{\tiny   \eqref{f:twistX}}}{=} \mathfrak{M}(Tw_{\chi_{LT}}(\sigma_{-1}))\mathfrak{M}(g)+\mathfrak{M}(\sigma_{-1})\mathfrak{M}(Tw_{\chi_{LT}}(g))\\
&= -\mathfrak{M}(\sigma_{-1})\mathfrak{M}(g)+\mathfrak{M}(\sigma_{-1})\mathfrak{M}(Tw_{\chi_{LT}}(g)).
\end{align*}
From this  the general case is immediate.
\end{proof}

\begin{lemma}\label{lem:interpolation}
Let $\lambda\in D(\Gamma_L,\mathbb{C}_p)$ with $\ev_{\chi_{LT}^j}(\lambda)=0$ for infinitely many $j,$ then $\lambda=0.$
\end{lemma}

\begin{proof}
 On the character variety the characters $\chi_{LT}^j$ corresponds to points which converge to the trivial character.  It follows  that $\lambda$ corresponds to the trivial function, since otherwise its divisor of zeroes would have only finitely many zeroes in any disk with fixed radius strictly smaller than $1$ by \eqref{f:divisorfinite}, which would contradict the assumptions.
\end{proof}

Now fix a $\mathbb{Z}_p$-basis $b=(b_1,\ldots,b_d)$ of $U_{n_0}$ with all $b_i \neq 1$ and set $\ell^*(b):=\ell^*_\Gamma(b):=q^{-n_0}\ell(b)\in o_L^\times$ with $\Gamma:=\Gamma_{n_0}$.
According to section \ref{sec:explicit} we  may define the operator
\begin{align*}
\widehat{\Xi_b}:=&q^{-n_0}\chi_{LT}^*({\overline{\Xi}_b})=\ell^*(b)\chi_{LT}^*({\Xi_b})
\end{align*}
in $\cR_K(\Gamma)$. Let $\mathrm{aug}: D(\Gamma,K)\to K$ denote the augmentation map, induced by the trivial map $\Gamma\to \{1\}.$

\begin{lemma}  The element $\widehat{\Xi_b}$   induces - up to the constant $q^{-n_0}$ - the augmentation map
\begin{equation}
\label{f:Xi} <\widehat{\Xi_b}, ->_{\Gamma_{n_0}}= q^{-n_0}\mathrm{aug} : D(\Gamma_{n_0},K)\to K.
\end{equation}
Moreover, we have
\begin{equation}\label{f:varrhovalue}
  \varsigma( \widehat{\Xi_b})=1=\frac{q}{q-1}
  \varrho(\widehat{\Xi_b}).
\end{equation}
\end{lemma}

\begin{proof} We may and do assume $\Omega\in K$   by compatibility of   $\mathrm{res}$ with respect to (complete) base change \eqref{f:basechangeres}.
 Since $\kappa^*(\hat{\ell}_{n_0*}(\widehat{\Xi_b}))=q^{-n_0}\widetilde{\Xi}_{b}\equiv \frac{\pi_L^{n_0}}{q^{n_0}\Omega }\frac{1}{ Z}\mod {\mathcal{O}_K(\mathbf{B})} $ by Remark \ref{rem:Xi}, one has for every $\lambda\in D(\Gamma,K)$
 \begin{align*}
   <\widehat{\Xi_b}, \lambda>_\Gamma &\stackrel{\text{\tiny   \eqref{f:projformula}}}{=}\frac{1}{(q-1) q^{n_0-1}}<\widehat{\Xi_b}, \lambda>_{\Gamma_L}\\
   & \stackrel{\text{\tiny   \eqref{f:groupRobba-pairing-explizit}}}{=}  q^{-n_0} \mathrm{res}_{\mathfrak{X} }((\frac{q}{\pi_L})^{n_0}\kappa^*(\hat{\ell}_{n_0*}(\widehat{\Xi_b} \lambda))d\log_{\mathfrak{X} })\\ &\stackrel{\text{\tiny   \eqref{f:res-XB}}}{=}q^{-n_0} \mathrm{res}_\mathbf{B}(\Omega(\frac{q}{\pi_L})^{n_0}\kappa^*(\hat{\ell}_{n_0*}(\widehat{\Xi_b} \lambda))g_{LT}dZ)\\
   &=q^{-n_0}\mathrm{res}_\mathbf{B}(\frac{1}{ Z}\kappa^*(\hat{\ell}_{n_0*}( \lambda))g_{LT}dZ)\\
   &=q^{-n_0}\mathrm{aug}(\lambda),
 \end{align*}
 where we use for the last equation that $g_{LT}(Z)$ has constant term $1$  and the fact that the augmentation map corresponds via Fourier theory and the LT-isomorphism to the `evaluation at $Z=0$' map. Taking $\lambda=1$ we see that $\varrho(\widehat{\Xi_b})=<\widehat{\Xi_b}, 1>_{\Gamma_L}=\frac{q-1}{q}.$

 For the other equation of the second claim one has   by definition of $\varsigma$
 \begin{align*}
  \varsigma(\widehat{\Xi_b})& \stackrel{\text{\tiny   \eqref{f:res-XB}}}{=} \Omega\ell^*(b)\mathrm{res}_\mathbf{B}(\kappa^*\Big(\mathfrak{M}(\sigma_{-1}) \mathfrak{M}(Tw_{\chi_{LT}}(\chi_{LT}^*(\Xi_b)))\Big)d\log_{LT})\\
  &=\Omega\ell^*(b)\mathrm{res}_\mathbf{B}( \mathfrak{M}_{LT}(\sigma_{-1}) \mathfrak{M}_{LT}(Tw_{\chi_{LT}}(\chi_{LT}^*(\Xi_b)))d\log_{LT})\\
  & =\ell^*(b)\mathrm{res}_\mathbf{B}(\mathfrak{M}_{LT}(\sigma_{-1}) \log_{LT}(Z)\partial_\mathrm{inv}\mathfrak{M}_{LT}(\chi_{LT}^*(\Xi_b))\frac{d\log_{LT}}{\log_{LT}(Z)})\\
  & =\ell^*(b)\mathrm{res}_\mathbf{B}(\mathfrak{M}_{LT}(\sigma_{-1}) \mathfrak{M}_{LT}(\nabla\chi_{LT}^*(\Xi_b))\frac{d\log_{LT}}{\log_{LT}(Z)})\\
  & =\ell^*(b)\mathrm{res}_\mathbf{B}(\mathfrak{M}_{LT}(\sigma_{-1}) \frac{\pi_L^{n_0}\log_{LT}(Z)}{\varphi_L^{n_0}(Z\ell(b))}\eta(1,Z)\frac{d\log_{LT}}{\log_{LT}(Z)})\\
  & =\ell^*(b)\pi_L^{n_0} \mathrm{res}_\mathbf{B}(\eta(-1,Z ) \frac{1}{\varphi_L^{n_0}(Z\ell(b))}\eta(1,Z)d\log_{LT})\\
  & =\ell^*(b)\pi_L^{n_0} \mathrm{res}_\mathbf{B}(\eta(1-1,Z) \frac{1}{\varphi_L^{n_0}(Z\ell(b))}d\log_{LT})\\
   & =\ell^*(b)\pi_L^{n_0} \mathrm{res}_\mathbf{B}(\varphi_L^{n_0}\left(\eta(0,Z) \frac{1}{Z\ell(b)}d\log_{LT}\right))\\
  &=\frac{\ell^*(b)}{\ell(b)}\pi_L^{n_0}(\frac{q}{\pi_L})^{{n_0}}\mathrm{res}_\mathbf{B}( \frac{1}{Z}g_{LT}dZ)\\
  &=1,
 \end{align*}
where we use \eqref{f:twist} in the third equation, the fact that $\nabla$ acts on $\cR$ as $\log_{LT}(Z)\partial_\mathrm{inv}$ (cp.\  \eqref{f:nablaonR}) in the fourth equation,  Remark \ref{rem:ThetaMellin} for the fifth equation, Lemma \ref{lem-pair} (iv) with $\psi_L(1)=\frac{q}{\pi_L}$ for the penultimate equation
and finally for the last equation that $g_{LT}(Z)$ has constant term $1$.
\end{proof}

\begin{proof}[Proof of Thm.\  \ref{thm:residuumidentity}]
Since the equality can also be checked after base change by \eqref{f:basechangeres} we may and do assume
that $\Omega$ belongs to $K.$
 Due to Prop.\  \ref{prop:residuepairingRGamma} there exists $g\in D(\Gamma_L,K)$ such that $\varsigma(\lambda)=<g,\lambda>_{\Gamma_L}$ for all $\lambda\in \cR_K(\Gamma_L),$ because $\varsigma$ sends $D(\Gamma_L,K)$ to zero. We claim that
 \begin{equation}\label{f:twistvarsigma}
   Tw_{\chi_{LT}^j}(g)=g
 \end{equation}
 for all $j\in\mathbb{Z}:$  By Prop.\  \ref{cor:twist<} and Lemma \ref{lem:Tw<<<}  we have
 \begin{align*}
   <Tw_{\chi_{LT}^j}(g),f>_{\Gamma_L} & =<g, Tw_{\chi_{LT}^{-j}}(f)>_{\Gamma_L}\\
    & =\varsigma( Tw_{\chi_{LT}^{-j}}(f)) \\
     & =\varsigma(f)\\
     &= <g,f>_{\Gamma_L}
 \end{align*}
 for all $f\in \cR_K(\Gamma_L).$

 Now it follows from \eqref{f:twistvarsigma} combined with Lemma \ref{lem:interpolation}   that $g$ is constant (and equal to $\ev_{\chi_{LT}^0}(g)$), i.e., $\varsigma(-)= g <1,->_{\Gamma_L}=g\varrho(-).$ Finally, it follows from \eqref{f:varrhovalue} that $g={\frac{q}{q-1} }.$
\end{proof}

\begin{corollary}
The pairing $\frac{q}{q-1} \kkl , \kkr_{\Gamma_L}$     makes the following diagram commutative
\begin{equation}\label{f:defpairing<<<}
\xymatrix{
       {\cR_K(\mathfrak{X})^{\psi_L=0}}  \ar@{}[r]|{\times} &  {(\Omega_{\cR(\mathfrak{X})}^1)^{\psi_L=0}}  \ar[r]^{\phantom{mm}mult} &  \Omega_{\cR(\mathfrak{X})}^1\ar[r]^{\mathrm{res}_{\mathfrak{X} }} & {\phantom{.}K\phantom{.}}\ar@{=}[d]\\
\frac{q}{q-1}  \kkl ,\kkr_{\Gamma_L}: \cR_K(\Gamma_L)\ar[u]_{\sigma_{-1}\mathfrak{M}\circ\upiota_*} \phantom{ < , >: } \ar@{}[r]|{\times} &  \cR_K(\Gamma_L) \ar[u]_{  \mathfrak{M}^{\Omega^1}}   \ar@{->}[rr] &&
 K,}
\end{equation}
i.e.,  we have
\begin{align}
\label{f:<<<}
\frac{q}{q-1} \notag\kkl\mu, \lambda\kkr_{\Gamma_L}&=\{  \mathfrak{M}(\sigma_{-1}\upiota_*(\mu))  , \mathfrak{M}^{\Omega^1}(\lambda)   \}_{\Omega^1}\\
&=\mathrm{res}_{\mathfrak{X} }(\sigma_{-1}\mathfrak{M}(\upiota_*(\mu))\mathfrak{M}^{\Omega^1}(\lambda))\\
&=\mathrm{res}_{\mathfrak{X} }(\mathfrak{M}(\sigma_{-1}\upiota_*(\mu))\mathfrak{M}(Tw_{\chi_{LT}}(\lambda))d\log_{\mathfrak{X}})\notag\\
&=\mathrm{res}_{\mathfrak{X} }(\mathfrak{M}(\upiota_*(\mu))\mathfrak{M}(Tw_{\chi_{LT}}(\sigma_{-1}\lambda))d\log_{\mathfrak{X}}).\notag
\end{align}
\end{corollary}

\begin{proof}
By Thm.\  \ref{thm:residuumidentity}, the definition of $\varsigma$ and of $<\; ,\;>_{\Gamma_L}$ we have
\begin{align}\label{f:mulambda}
{\frac{q}{q-1} }\kkl\mu, \lambda\kkr_{\Gamma_L}&=\frac{q}{q-1} <1,\mu\lambda>_{\Gamma_L}\notag\\
&=\{  \mathfrak{M}(\sigma_{-1})  , \mathfrak{M}^{\Omega^1}(\mu\lambda)   \}_{\Omega^1}\\
&=\{  \mathfrak{M}(\sigma_{-1}\upiota_*(\mu))  , \mathfrak{M}^{\Omega^1}(\lambda)   \}_{\Omega^1}\notag
\end{align}
where we use Lemma \ref{lem:sesquilinear} for the last equation.
\end{proof}

\begin{lemma}\label{lem:nabla<<<}
We have for all $\lambda,\mu\in \cR_K(\Gamma_L)$ that
$\kkl\lambda,\mu\kkr_{\Gamma_L}=-\kkl\upiota_*(\lambda), \upiota_*(\mu) \kkr_{\Gamma_L}.$
\end{lemma}

\begin{proof}
 Using \eqref{f:<<<} for the first and third equation, property 3. in Subsection \ref{subsec:twisting} applied to $\upiota$  and the fact that $Tw_{\chi_{LT}}(\sigma_{-1})=-\sigma_{-1} $ for the second equation  and Prop.\  \ref{cor:twist<} for the last one, we see that
\begin{align*}
\frac{q}{q-1}  \kkl\mu,\lambda\kkr_{\Gamma_L}
   &=\mathrm{res}_{\mathfrak{X} }\big(\mathfrak{M}(Tw_{\chi_{LT}}(\sigma_{-1}\lambda))\mathfrak{M}(\upiota_*(\mu))d\log_{\mathfrak{X}}\big)\\
 &=-\mathrm{res}_{\mathfrak{X} }\big(\mathfrak{M}(\sigma_{-1}\upiota_*(Tw_{\chi_{LT}^{-1}}(\upiota_*(\lambda))))\mathfrak{M}(\upiota_*(\mu))d\log_{\mathfrak{X}}\big)\\
   &=-\frac{q}{q-1}\kkl Tw_{\chi_{LT}^{-1}}(\upiota_*(\lambda)),Tw_{\chi_{LT}^{-1}}(\upiota_*(\mu))\kkr_{\Gamma_L}\\
   &=  -\frac{q}{q-1}\kkl\upiota_*(\lambda) ,\upiota_*(\mu)\kkr_{\Gamma_L}.
\end{align*}
\end{proof}

\subsubsection{The Iwasawa pairing for \texorpdfstring{$(\varphi_L,\Gamma_L)$}{(phi,Gamma)}-modules over the Robba ring}\label{subsec:IwasawaPairing}

As before let $\mathfrak{Y}$ be either $\mathfrak{X} $ or $\mathbf{B}$ and $\cR=\cR_{K}(\mathfrak{Y})$ and $M$ be in $\mathcal{M}^{an}(\cR),$ where $K$ is any complete intermediate extension $L \subseteq K \subseteq \mathbb{C}_p$.
Using Prop.\  \ref{prop:residuepairingRGamma} we define the pairing
\[
\KKl\;,\;\KKr_{Iw}^0:=\KKl\;,\;\KKr_{M,Iw}^0:\check{M}^{\psi_L=0}\times M^{\psi_L=0} \to \cR_K(\Gamma_L) \ ,
\]
which is $\cR_K(\Gamma_L)$-$\upiota_*$-sesquilinear in the sense that
\begin{equation}\label{sesquilinear}
  \lambda\KKl \check{m},m\KKr_{Iw}^0=\KKl\lambda \check{m},m\KKr_{Iw}^0=\KKl\check{m},\upiota_*(\lambda) m\KKr_{Iw}^0
\end{equation}
for all $\lambda\in\cR_K(\Gamma_L)$ and $\check{m}\in\check{M}^{\psi_L=0},\; m\in M^{\psi_L=0}$. This requires the commutativity of the diagram
\begin{equation*}
  \xymatrix{
   \cR_K(\Gamma_L) \ar@{=}[d]_{} \ar@{}[r]|{\times} & {\phantom{mmm} \check{M}^{\psi_L=0}\times M^{\psi_L=0}}  \ar@{-->}[d]_{\KKl\;,\;\KKr_{Iw}^0} \ar[r]^{} & K\phantom{,} \ar@{=}[d]^{} \\
     \cR_K(\Gamma_L) \ar@{}[r]|{\times}&  \cR_K(\Gamma_L) \ar[r]^-{\kkl\,,\,\kkr_{\Gamma_L}} & K,   }
\end{equation*}
in which the upper line sends $(f,x,y)$ to $\{f(x), y\}_M,$ where the latter pairing is \eqref{f:res-pair-general}.  Indeed, the property
\[\KKl \lambda\check{m},m\KKr_{Iw}^0=\KKl \check{m},\upiota_*(\lambda)m\KKr_{Iw}^0\]
 follows from the corresponding property for $\{ ,\}_M$ by Lemma \ref{lem:sesquilinear}, while with regard to the second one
\[\lambda\KKl \check{m},m\KKr_{Iw}^0=\KKl \lambda\check{m},m\KKr_{Iw}^0    \]
we have for all $f\in \cR_K(\Gamma_L)$
 \begin{align*}
  \kkl f,\KKl\lambda\check{m},m\KKr_{Iw}^0\kkr_{\Gamma_L}&=\{f\cdot \lambda\check{m},m\}\\
&=\kkl \lambda f,\KKl\check{m},m\KKr_{Iw}^0\kkr_{\Gamma_L}\\
&=\kkl f,\lambda\KKl\check{m},m\KKr_{Iw}^0\kkr_{\Gamma_L}
 \end{align*}
 by Lemma \ref{lem:nabla}. Note that the pairing $\KKl\;,\;\KKr_{Iw}^0$ induces the
isomorphism \eqref{f:psizeroequiv}.

We  set
\[\mathcal{C}:=(\frac{\pi_L}{q}\varphi_L -1){M}^{\psi_L=1}\;\; \mbox{ and } \check{\mathcal{C}}:= ( \varphi_L -1)\check{M}^{\psi_L=\frac{q}{\pi_L}} \]
and we shall need the following
\begin{lemma}\label{lem:fbalanced}
 For $f\in D(\Gamma_L,K)$ we have $\{f\cdot(\varphi_L-1)x,(\frac{\pi_L}{q}\varphi_L -1) y\}=0$ for all $x\in \check{M}^{\psi_L=\frac{q}{\pi_L}}$ and $y\in M^{\psi_L=1}.$
\end{lemma}
\begin{proof}
Straightforward calculation using Lemma \ref{lem-pair} above, cp.\ \cite[Lem.\ 4.2.7]{KPX}.
\end{proof}
This Lemma combined with the second statement of Prop.\  \ref{prop:residuepairingRGamma} implies that the restriction of $\KKl\;,\;\KKr_{Iw}^0$ to $\check{\mathcal{C}}\times \mathcal{C}$, which by abuse of notation we denote   by the same symbol, is characterized by the commutativity of the diagram
\begin{equation*}
  \xymatrix{
     \check{\mathcal{C}}\phantom{m}\times \phantom{m}\mathcal{C}  \ar@{-->}[d]_{\KKl\;,\;\KKr_{Iw}^0} \ar@{}[r]|{\times} & \cR_K(\Gamma_L)/D(\Gamma_L,K) \ar@{=}[d]_{} \ar[r]^{} & K\phantom{,} \ar@{=}[d]^{} \\
     D(\Gamma_L,K) \ar@{}[r]|{\times}&  \cR_K(\Gamma_L)/D(\Gamma_L,K) \ar[r]^-{\kkl\,,\,\kkr_{\Gamma_L}} & K,   }
\end{equation*}
in which the upper line sends $( x,y, f)$ to $\{f(x), y\}_M.$ In particular, it takes values in $D(\Gamma_L,K).$

Finally, we obtain a $D(\Gamma_L,K)$-$\upiota_*$-sesquilinear pairing $\KKl\;,\;\KKr_{Iw} :=\KKl\;,\;\KKr_{M,Iw}$ which by definition fits into the following commutative diagram
\begin{equation*}
\xymatrix{
     {\KKl\;,\;\KKr_{M,Iw}:  {\check{M}^{\psi_L=\frac{q}{\pi_L}}}}\phantom{ [ , ]: }  \ar[d]_{ \varphi_L -1}\ar@{}[r]|{\times} &  {M^{\psi_L=1}} \ar[d]_{ \frac{\pi_L}{q}\varphi_L -1} \ar@{-->}[r] & D(\Gamma_L,K)\ar@{=}[d]\\
    {\KKl\;,\;\KKr_{M,Iw}^0: \check{\mathcal{C}}\phantom{\KKl\;,\;\KKr_{Iw}^0dd }} \ar@{}[r]|{\times} &  \mathcal{C}   \ar[r]^{} &D(\Gamma_L,K). }
\end{equation*}
Altogether we obtain the following
\begin{theorem}
 There is a  $D(\Gamma_L,K)$--$\upiota_*$-sesquilinear  pairing
\begin{equation}\label{f:IWpairing}
\{\;,\;\}_{Iw}: \check{M}^{\psi_L=\frac{q}{\pi_L}} \times M^{\psi_L=1} \to D(\Gamma_L,K).
\end{equation}
It is characterized by the following property
\begin{equation}
\label{f:definingproperty} <f,\{\check{m},m\}_{Iw}>_{\Gamma_L}=\{f\cdot(\varphi_L-1)\check{m},(\frac{\pi_L}{q}\varphi_L -1)m\} \mbox{ for all } f\in \cR_K(\Gamma_L),\; \check{m}\in\check{M},\; m\in M.
\end{equation}
\end{theorem}

\begin{remark}\label{rem:IWpairU}
For any $n\geq n_0,$ we obtain similarly as in \eqref{f:IWpairing} $D(\Gamma_n,K)$-$\upiota_*$-sesquilinear pairings
\begin{equation*}
\KKl\; ,\;\KKr_{Iw,\Gamma_n}: \check{M}^{\psi_L=\frac{q}{\pi_L}} \times M^{\psi_L=1} \to D(\Gamma_n,K).
\end{equation*}
It follows immediately from the definitions, the projection formulae \eqref{f:projformula} 
and Frobenius reciprocity \ref{rem:Frobrec} that
\[\KKl\;,\;\KKr_{Iw,\Gamma_n}:=(q-1) q^{n-1}pr_{L,n}\circ \KKl\;,\;\KKr_{Iw}.  \]
\end{remark}

If $\chi : \Gamma_L \longrightarrow o_L^\times$ is any continuous character with representation module $W_\chi=o_Lt_\chi$ then, for any $(\varphi_L,\Gamma_L)$-module $M$ over $\cR$, we have the twisted $(\varphi_L,\Gamma_L)$-module $M(\chi)$   where $M(\chi) := M\otimes_{o_L}W_\chi$ as $\cR$-module, $\varphi_{M(\chi)}(m\otimes w) := \varphi_M(m)\otimes w$, and $\gamma | M(\chi)(m\otimes w) :=  \gamma|M (m)\otimes  \gamma|W_\chi(w)=\chi(\gamma) \cdot\gamma|M(m)\otimes w$ for $\gamma \in \Gamma_L$.  It follows that $\psi_{M(\chi)}(m\otimes w) = \psi_M(m)\otimes w$. For the character $\chi_{LT}$ we take $W_{\chi_{LT}}=T=o_L \eta$ and $W_{\chi_{LT}^{-1}}=T^*=o_L\eta^*$ as representation module, where $T^*$ denotes the $o_L$-dual with dual basis $\eta^*$ of $\eta$.

Consider the $\cR_K$-linear (but of course not $\cR_K(\Gamma_L)$-linear) map
\[tw_\chi:M \to M(\chi), \; m\mapsto m\otimes t_\chi.\]

\begin{lemma}\label{lem:twist}
There is a commutative diagram
\begin{equation*}
\xymatrix{
 {\check{M}(\chi_{LT}^{-j})^{\psi_L=\frac{q}{\pi_L}}}  \ar@{}[r]|{\times} &  {{M}(\chi_{LT}^{j})^{\psi_L=1}}\ar[r]^-{\KKl,\KKr_{Iw}}  & D(\Gamma_L,\Cp) \\
   {\check{M}^{\psi_L=\frac{q}{\pi_L}}}\ar[u]_{ tw_{\chi_{LT}^{-j}} }  \ar@{}[r]|{\times} &  {{M}^{\psi_L=1}}\ar[r]^-{\KKl,\KKr_{Iw}}\ar[u]_{ tw_{\chi_{LT}^{j}}}  & D(\Gamma_L,\Cp).\ar[u]_{Tw_{\chi_{LT}^j}} }
\end{equation*}
\end{lemma}

\begin{proof}
We have for all $f\in\cR_K(\Gamma_L),$
\begin{align*}
  \kkl f,\KKl tw_{\chi_{LT}^{-j}}(\check{m}),tw_{\chi_{LT}^{j}}(m)\KKr_{Iw} \kkr_{\Gamma_L} &=\{f\cdot\left((\varphi_L-1)\check{m}\otimes \eta^{\otimes -j}\right),(\frac{\pi_L}{q}\varphi_L -1)m\otimes \eta^{\otimes j}\} \\
   &=\{\left(Tw_{\chi_{LT}^{-j}}(f)\cdot(\varphi_L-1)\check{m} \right)\otimes \eta^{\otimes -j},(\frac{\pi_L}{q}\varphi_L -1)m\otimes \eta^{\otimes j}\} \\
   &=\{\left(Tw_{\chi_{LT}^{-j}}(f)\cdot (\varphi_L-1)\check{m}\right) ,(\frac{\pi_L}{q}\varphi_L -1)m\} \\
   &= \kkl Tw_{\chi_{LT}^{-j}}(f),\KKl\check{m},m\KKr_{Iw}\kkr_{\Gamma_L}\\
   &= \kkl f,Tw_{\chi_{LT}^{j}}(\KKl\check{m},m\KKr_{Iw})\kkr_{\Gamma_L}
\end{align*}
where we used Corollary \ref{cor:twist<} for the last equation. The second equation is clear for $\delta$-distributions and hence extends by the uniqueness result of Thm.\  \ref{thm:Mpsiequal0}, cf.\ the proof of Thm.\  \ref{Rrs-extension}.
\end{proof}

\subsubsection{The abstract reciprocity formula}\label{subsec:abstract-rec}

We keep the notation from the preceding subsection and set $t_\mathfrak{Y}:=\log_\mathfrak{Y}$.

\paragraph{Compatibility of the Iwasawa pairing under comparison isomorphisms}

Let $M,N$ be (not necessarily \'{e}tale) $L$-analytic $(\varphi_L,\Gamma_L)$-modules over $\cR$. We extend the  action of $\Gamma_L,$ $\varphi_L$ and $\psi_L$ to the $\cR[\frac{1}{t_{\mathfrak{Y}}}]$-module $M[\frac{1}{t_{\mathfrak{Y}}}]$ (and in the same way to $N[\frac{1}{t_{\mathfrak{Y}}}]$) as follows:\footnote{Since $t_{\mathfrak{Y}}^k=\varphi_L(\frac{t_{\mathfrak{Y}}^k}{\pi_L^k})$ one checks that $\psi_L(t_{\mathfrak{Y}}^km)= \frac{t_{\mathfrak{Y}}^k}{\pi_L^k}\psi_L(m)$  by the projection formula. In particular, the definition is independent of the chosen denominator.}
\begin{align*}
  \gamma \frac{m}{t_{\mathfrak{Y}}^k}&:=\frac{\gamma m}{\gamma t_{\mathfrak{Y}}^k}=\frac{\frac{\gamma m}{\chi_{LT}^k(\gamma)}}{t_{\mathfrak{Y}}^k},  \\ \varphi_L(\frac{m}{t_{\mathfrak{Y}}^k})&:=\frac{\varphi_L(m)}{\varphi_L(t_{\mathfrak{Y}}^k)}=\frac{\frac{\varphi_L(m)}{\pi_L^k}}{t_{\mathfrak{Y}}^k} \mbox{ and}\\
  \psi_L(\frac{m}{t_{\mathfrak{Y}}^k})&:=\frac{\pi_L^k\psi_L(m)}{t_{\mathfrak{Y}}^k}.
\end{align*}

Now we assume that there is an isomorphism
\begin{align*}
  c:\cR[\frac{1}{t_{\mathfrak{Y}}}]\otimes_{\cR} M & \xrightarrow{\cong}\cR[\frac{1}{t_{\mathfrak{Y}}}]\otimes_{\cR} N
\end{align*}
of $(\varphi_L,\Gamma_L)$-modules over $\cR[\frac{1}{t_{\mathfrak{Y}}}].$

\begin{lemma}
\begin{enumerate}
\item $(M[\frac{1}{t_{\mathfrak{Y}}}])^{\psi_L=0}=(M^{\psi_L=0})[\frac{1}{t_{\mathfrak{Y}}}]:=\{\frac{m}{t_{\mathfrak{Y}}^k}|m\in M^{\psi_L=0},\, k\geq 0\} .$ \Footnote{We give this formal definition as $M^{\psi_L=0} $ does not seem to be a module over a ring containing $t_{\mathfrak{Y}}$ such that localisation with respect to $t_{\mathfrak{Y}}$ makes sense a priori.   $\varphi(\cR)$ does work! }
\item The (separatedly continuous) $\cR_K(\Gamma_L)$-action on $ M^{\psi_L=0}$ extends to a (separatedly continuous  with respect to direct limit topology) action of  $\cR_K(\Gamma_L)$ on  $(M[\frac{1}{t_{\mathfrak{Y}}}])^{\psi_L=0}.$
\end{enumerate}
\end{lemma}

\begin{proof} For (i) note that
$0=\psi_L(\frac{m}{t_{\mathfrak{Y}}^k})=\frac{\pi_L^k\psi_L(m)}{{t_{\mathfrak{Y}}^k}}$ if and only if ${\psi_L(m)}=0. $ For (ii) take for any $f\in \cR_K(\Gamma_L)$  the direct limit of the following commutative diagram
\[\xymatrix{
  M^{\psi_L=0} \ar[d]_{f} \ar[r]^{t_\mathfrak{Y}} & M^{\psi_L=0}\ar[d]_{Tw_{\chi_{LT}^{-1}}(f)} \ar[r]^{t_{\mathfrak{Y}}} & \cdots \ar[r]^{t_{\mathfrak{Y}}} & M^{\psi_L=0} \ar[d]_{Tw_{\chi_{LT}^{-i}}(f)} \ar[r]^{t_{\mathfrak{Y}}} & \cdots  \\
  M^{\psi_L=0} \ar[r]^{t_\mathfrak{Y}} & M^{\psi_L=0} \ar[r]^{t_{\mathfrak{Y}}} & \cdots \ar[r]^{t_{\mathfrak{Y}}} & M^{\psi_L=0} \ar[r]^{t_{\mathfrak{Y}}} & \cdots   }\]
This defines a (separatedly continuous) action.
\end{proof}

Consider the composite map
\begin{align*}
 \check{c}:\cR[\frac{1}{t_{\mathfrak{Y}}}]\otimes_{\cR} \check{M}&\cong \Hom_{\cR[\frac{1}{t_{\mathfrak{Y}}}]}(\cR[\frac{1}{t_{\mathfrak{Y}}}]\otimes_{\cR} M,\cR[\frac{1}{t_{\mathfrak{Y}}}]\otimes_{\cR}\Omega^1_{\cR})\\ &\cong \Hom_{\cR[\frac{1}{t_\mathfrak{Y}}]}(\cR[\frac{1}{t_{\mathfrak{Y}}}]\otimes_{\cR} N,\cR[\frac{1}{t_{\mathfrak{Y}}}]\otimes_{\cR} \Omega^1_{\cR})\\ &\cong \cR[\frac{1}{t_{\mathfrak{Y}}}]\otimes_{\cR} \check{N}
\end{align*}
where the second isomorphism is $ {(c^{-1})^*}.$

\begin{lemma}
$c^{\psi_L=0}$ and $\check{c}^{\psi_L=0}$ are $\cR_K(\Gamma_L)$-equivariant.
\end{lemma}

\com{
\begin{proof}
 Consider, for $n\in \mathbb{Z}$, the $(\varphi_L,\Gamma_L)$-modules (!)      $M_n:=t_{\mathfrak{Y}}^{-n}M$ over $\cR$ and note that the inclusion $(M_n)^{\psi_L=0}\subseteq (M[\frac{1}{t_{\mathfrak{Y}}}])^{\psi_L=0} $ is $\cR_K(\Gamma_L)$-equivariant by construction of the action. Now, since $M,N$ are finitely generated over $\cR,$ there exists $n_0\geq 0$ such that $c$ restricts to a homomorphism $c_0:M\to N_{n_0}$ of $(\varphi_L,\Gamma_L)$-modules over $\cR$, whence $c_0^{\psi_L=0}:M^{\psi_L=0}\to N^{\psi_L=0}_{n_0}\subseteq (N[\frac{1}{t_{\mathfrak{Y}}}])^{\psi_L=0}$ is $\cR_K(\Gamma_L)$-equivariant by the functoriality of Thm.\ \ref{thm:Mpsiequal0-decent} and similarly for the induced maps $c_n:M_n\to N_{n_0+n}$ for all $n\geq 0.$ The equivariance for $c^{\psi_L=0}$ follows by taking direct limits. \\ Similarly, for some $n_0\geq 0$, the inverse $b$ of $c$ induces homomorphisms $b_n:N_{-n_0-n}\to M_{-n}$ of $(\varphi_L,\Gamma_L)$-modules over $\cR$  all $n\in \mathbb{Z}.$ We obtain  homomorphisms of $(\varphi_L,\Gamma_L)$-modules over $\cR$
\begin{align*}
 \check{c}_n: (\check{M})_{ n}&=\Hom_{\cR}(  M, t_{\mathfrak{Y}}^{ -n}\Omega^1_{\cR})\\
&\cong \Hom_{\cR}(   M_{ -n}, \Omega^1_{\cR})\\
&\cong \Hom_{\cR}(   N_{-n_0 -n}, \Omega^1_{\cR})\\
&\cong (\check{N})_{n_0 +n}
\end{align*}
where the third isomorphism is $(b_n)^*.$ As above $(\check{c}_n)^{\psi_L=0}$ is $\cR_K(\Gamma_L)$-equivariant and the claim follows by taking direct limits.

\end{proof}}

\begin{lemma}\label{lem:checkIW}
The following diagram commutes on the vertical intersections
\[\xymatrix{
  \check{M}^{\psi_L=0} \ar@{^(->}[d]_{}  \ar@{}[r]|{\times} & M^{\psi_L=0} \ar@{^(->}[d]_{} \ar[r]^-{\{,\}^0_{M,Iw}} &\cR_K(\Gamma_L) \ar@{=}[ddd]^{} \\
  (\cR[\frac{1}{t_{\mathfrak{Y}}}]\otimes_{\cR} \check{M})^{\psi_L=0} \ar[d]_{\check{c}}^{\cong} \ar@{}[r]|{\times} & (\cR[\frac{1}{t_{\mathfrak{Y}}}]\otimes_{\cR} {M})^{\psi_L=0} \ar[d]_{c}^{\cong} \\
 (\cR[\frac{1}{t_{\mathfrak{Y}}}]\otimes_{\cR} \check{N})^{\psi_L=0}   \ar@{}[r]|{\times} & (\cR[\frac{1}{t_{\mathfrak{Y}}}]\otimes_{\cR} N)^{\psi_L=0}    \\
 \check{N}^{\psi_L=0}\ar@{^(->}[u] \ar@{}[r]|{\times} & N^{\psi_L=0}\ar@{^(->}[u]\ar[r]^-{\{,\}^0_{N,Iw}} &\cR_K(\Gamma_L),  }\]
 i.e., if $\check{m}\in \check{M}, m\in M,\check{n}\in\check{N}, n\in N$ with $\check{c}(\check{m})=\check{n}$ and $c(m)=n,$ then
\[\{\check{m},m\}_{M,Iw}^0=\{\check{n},n\}_{N,Iw}^0.\]
\end{lemma}

\begin{proof}
 By definition of the Iwasawa pairings we have for all $f\in\cR_K(\Gamma_L)$
\begin{align*}
  <f,\{\check{n},n\}_{N,Iw}^0>_{\Gamma_L}&=\{f\cdot\check{n},n\}_N  \\
    &= \{f\cdot\check{c}(\check{m}),c(m)\}_N \\
    &= \{\check{c}(f\cdot\check{m}),c(m)\}_N \\
    &= \mathrm{res}_{\mathfrak{Y}}(\check{c}(f\cdot\check{m})(c(m))\\
    &= \mathrm{res}_{\mathfrak{Y}}(((f\cdot\check{m})\circ c^{-1})(c(m))\\
    &= \mathrm{res}_{\mathfrak{Y}}((f\cdot\check{m})(m))\\
    &=\{f\cdot\check{m},m\}_M\\
    &=<f,\{\check{m},m\}_{M,Iw}^0>_{\Gamma_L}
\end{align*}
whence the claim. Here we use the $\cR_K(\Gamma_L)$-equivariance of $\check{c}$ in the third equality.
\end{proof}

Now let $D$ be any  $\varphi_L$-module over $L$ of finite dimension, say  $d$, (with trivial $\Gamma_L$-action) and consider the $(\varphi_L,\Gamma_L)$-module $N:=\cR\otimes_L D$ over $\cR$ (with diagonal actions)\Footnote{The isomorphism $\cR \otimes_{\varphi_L,\cR}N\cong N$ follows from the fact that it is true for $\cR$ itself (i.e., $D=L$) and the bijectivity of $\varphi_L$ on $D.$}  Since $N \cong \cR^d$ as $\Gamma_L$-module, it is $L$-analytic. Moreover, we have $\check{N}\cong\Omega^1_\cR\otimes D^*$ with $D^*=\Hom_L(D,L)$ being the dual $\varphi_L$-module.      We set
\[\tilde{\Omega}:=\left\{
                    \begin{array}{ll}
                      1, & \hbox{if $\mathfrak{Y}=\mathfrak{X}$;} \\
                      {\Omega}, & \hbox{if $\mathfrak{Y}=\mathbf{B}$ (and $\Omega\in K$).}
                    \end{array}
                  \right.
\]

\begin{lemma}\label{lem:Iwasawacheck} {  If $\mathfrak{Y}=\mathbf{B},$ we assume $\Omega\in K.$  }
There is a commutative diagram
\[\xymatrix{
   (\Omega^1_\cR\otimes D^*)^{\psi_L=0}   \ar@{}[r]|{\times} & (\cR\otimes_L D)^{\psi_L=0}  \ar[rr]^-{\tilde{\Omega}\frac{q-1}{q}\{,\}^0_{N,Iw}} && \cR_K(\Gamma_L) \ar@{=}[d]^{} \\
 \cR_K(\Gamma_L)\otimes_L D^*\ar[u]_{\cong}^{\mathfrak{M}^{\Omega^1}\otimes \id} \ar@{}[r]|{\times} & \cR_K(\Gamma_L)\otimes_L D\ar[u]_{\cong}^{\sigma_{-1}\mathfrak{M}\circ\upiota_*\otimes\id}\ar[rr]^-{} && \cR_K(\Gamma_L),   }\]
 where the bottom line is the $\cR_K(\Gamma_L)$-linear extension of the canonical pairing between $D^*$ and $D$, i.e., it maps $(\lambda\otimes l,\mu\otimes d)$ to $\lambda\mu l(d).$
\end{lemma}

\begin{proof}
Let $\check{d}_j$ and $d_i$ be a basis of $D^*$ and $D, $ respectively, and
 $ x=\sum_j \lambda_j\cdot \check{d}_j   $ and $ y=\sum_i \mu_i\cdot d_i.$ Then, by definition of $\KKl,\KKr_{Iw}^0$ we have for all $\lambda\in \cR_K(\Gamma_L)$
\begin{align*}
 \kkl\lambda,\KKl (\mathfrak{M}^{\Omega^1}\otimes & \id)(x),(\sigma_{-1}\mathfrak{M}\circ\upiota_*\otimes\id)(y)\KKr^0_{Iw}\kkr_{\Gamma_L}\\
  &= \KKl (\lambda\mathfrak{M}^{\Omega^1}\otimes \id)(x),(\sigma_{-1}\mathfrak{M}\circ\upiota_*\otimes\id)(y)\KKr\\
   &=\{\sum_j (\lambda\lambda_j)\cdot (\ev_1d\log_{\mathfrak{Y}}\otimes\check{d}_j ) ,\sum_i \upiota_*(\mu_i)\cdot \ev_{-1}\otimes d_i\}\\
  &=\sum_{i,j} \{(\lambda\lambda_j\mu_i)\cdot (\ev_1d\log_{\mathfrak{Y}})\otimes\check{d}_j,  \ev_{-1}\otimes d_i\}\\
    &=\sum_{i,j} \mathrm{res}_{\mathfrak{Y}}\bigg(\check{d}_j\big(d_i\big)\ev_{-1}(\lambda\lambda_j\mu_i)\cdot (\ev_1d\log_{\mathfrak{Y}}) \bigg). \\
    &=\sum_{i,j}\check{d}_j\big(d_i\big) \mathrm{res}_{\mathfrak{Y}}\bigg(\ev_{-1}(\lambda\lambda_j\mu_i)\cdot (\ev_1d\log_{\mathfrak{Y}}) \bigg). \\
\end{align*}
Here, for the third equation we used property (iii) in Lemma \ref{lem-pair}.
On the other hand we can pair the image $ \sum_{i,j}\lambda_j\mu_i\check{d}_j(d_i))$ of $(x,y)$ under the bottom pairing with $\lambda$ using the description \eqref{f:mulambda}
\begin{align*}
{ \frac{q}{q-1} } \kkl\lambda,\sum_{i,j}\lambda_j\mu_i\check{d}_j(d_i))\kkr_{\Gamma_L}&=\sum_{i,j}\check{d}_j(d_i)\{  \mathfrak{M}(\sigma_{-1})  , \mathfrak{M}^{\Omega^1}(\lambda\lambda_j\mu_i)   \}\\
 &=\sum_{i,j}\check{d}_j(d_i)\mathrm{res}_{\mathfrak{X}}\bigg(\ev_{-1}(\lambda\lambda_j\mu_i)\cdot (\ev_1d\log_{\mathfrak{X}}) \bigg),
 %
\end{align*}
whence comparing with the above gives the result for  $\mathfrak{Y}=\mathfrak{X},$using Prop.\  \ref{prop:residuepairingRGamma}. If $\mathfrak{Y}=\mathbf{B},$ we obtain the factor $\Omega$ due to Remark \ref{rem:compPairinigXB}.
\end{proof}

\begin{definition}\label{def:etaleana}
An $L$-analytic $(\varphi_L,\Gamma_L)$-module $M$ over $\mathcal{R}$ is called \'{e}tale,
if it is semistable and of slope $0.$ We write $\mathfrak{M}^{an,\acute{e}t}(\cR)$
for the category of \'{e}tale, analytic $(\varphi_L,\Gamma_L)$-modules over
$\mathcal{R}.$
\end{definition}

Crucial is the following

\begin{theorem}\label{thm:BergerEquiv}{
  There are equivalences of categories
  \begin{align*}
    Rep_L^{an}(G_L) & \longleftrightarrow \mathfrak{M}^{an,\acute{e}t}(\cR_L(\mathbf{B})) \\
    V & \mapsto D_{\mathrm{rig}}^\dagger(V).
  \end{align*}
and
\begin{align*}
    Rep_L^{an}(G_L) & \longleftrightarrow \mathfrak{M}^{an,\acute{e}t}(\cR_L(\mathfrak{X})) \\
    V & \mapsto D_{\mathrm{rig}}^\dagger(V)_\mathfrak{X},
  \end{align*}
 where the functor is defined in the proof below.
  }
\end{theorem}

\begin{proof}{
Thm.\  D in \cite{Be16} and Thm.\  3.27 in \cite{BSX}, which states an equivalence of categories
 \begin{align}\label{f:equivalenceBSX}
    \mathfrak{M}^{an,\acute{e}t}(\cR_L(\mathbf{B})) & \longleftrightarrow \mathfrak{M}^{an,\acute{e}t}(\cR_L(\mathfrak{X})) \\
   M & \mapsto M_\mathfrak{X}.\notag
  \end{align}}
\end{proof}

We recall  the
definition of the subring $ \mathbf{B}^\dagger_L $ of $\cR_L(\mathbf{B})$ by  defining  first $\tilde{\mathbf{A}}:=W(\mathbb{C}_p^\flat)_L$ and
\[
\tilde{\mathbf{A}}^\dagger:=\{x=\sum_{n\geq 0}\pi_L^n[x_n]\in\tilde{\mathbf{A}}: |\pi_L^{n}\|x_n|_\flat^r \xrightarrow{n \rightarrow \infty} 0\mbox{ for some } r>0\}.
\]
Then we set  ${\bf A}^\dagger:=\tilde{\bf A}^{\dagger}\cap {\bf A}$,  ${\bf B}^\dagger:={\bf A}^\dagger[\frac{1}{\pi_L}]$ as well as
${\bf A}^\dagger_L:=({\bf A}^\dagger)^{H_L} $ and
${\bf B}^\dagger_L:=({\bf B}^\dagger)^{H_L}.$

It follows from the proof of \cite[Thm.\ 10.1]{Be16} that for $V\in Rep_L^{an}(G_L) $
we have $D_{\mathrm{rig}}^\dagger(V)=\cR_L(\mathbf{B})\otimes_{\mathbf{B}^\dagger_L}D^\dagger(V)$,
where $D^\dagger(V) $ belongs to $\mathfrak{M}^{\acute{e}t}(\mathbf{B}^\dagger_L)$.
From the theory of Wach modules we actually know that $D_{LT}(V)$ is even of finite
height, if $V$ is crystallin in addition:
\begin{align*}
 D^\dagger(V) & =\mathbf{B}^\dagger_L\otimes_{\mathbf{A}_L^+}N(T)=\mathbf{B}^\dagger_L\otimes_{\mathbf{B}_L^+}N(V)
\end{align*}
for any Galois stable $o_L$-lattice $T\subseteq V.$ From the big diagram in section \ref{sec:KR} we thus obtain the following diagram, in which the horizontal maps are equivalences of categories.

\begin{equation*}
  \xymatrix{
     \mathrm{Mod}_{\mathbf{B}_L^+}^{\varphi_L,\Gamma_L,an} \ar[d]_{\mathcal{O} \otimes_{\mathbf{B}_L^+} -} \ar[d]  \ar[r]^{\mathbf{B}_L \otimes_{\mathbf{B}_L^+} -}_-{\simeq} &   \mathfrak{M}^{et,cris}(\mathbf{B}_L)   \\
    \mathrm{Mod}_{\mathcal{O}}^{\varphi_L,\Gamma_L,0} \ar[d]_{\cR_L(\mathbf{B})\otimes_\cO -} \ar@<1ex>[r]^-{V_L\circ D}_-{\simeq}  & \ar[d]_{\subseteq}   \Rep_L^{cris,an}(G_L) \ar@<1ex>[l]^-{\mathcal{M}\circ D_{cris,L}} \ar[u]_{D_{LT}(-)}^{\simeq} \ar[ul]_{N(-)}\\
    \mathfrak{M}({\cR_L(\mathbf{B})})^{an,\acute{e}t}  & \Rep_L^{an}(G_L) \ar[l]^-{\simeq}_-{ D^\dagger(V)}  }
\end{equation*}

Here $ \mathfrak{M}^{et,cris}(\mathbf{B}_L)$ denotes the essential image of $ \Rep_L^{cris,an}(G_L) $ under $D_{LT}(-)$ in $\mathfrak{M}^{et}(\mathbf{B}_L)$ with $\mathbf{B}_L:=\mathbf{A}_L[\frac{1}{\pi_L}].$

Now let $T$ be an $o_L$-lattice in an $L$-linear continuous representation of $G_L$ such that $V^*(1)$ (and hence $V(\tau^{-1})$) is $L$-analytic and crystalline: Then it follows from \cite{KR} and the discussion above that
   \[ M:=D_{\mathrm{rig}}^\dagger(V(\tau^{-1}))=\cR_L(\mathbf{B})\otimes_{{\mathcal{O}_K(\mathbf{B})}}\mathcal{M}(D_{cris,L}(V(\tau^{-1})))=\cR_L(\mathbf{B})\otimes_{\mathbb{A}_L^+}N(T(\tau^{-1}))                                      \] as well as
  \[\check{M}=D_{\mathrm{rig}}^\dagger(V^*(1))=\cR_L(\mathbf{B})\otimes_{{\mathcal{O}_K(\mathbf{B})}}\mathcal{M}(D_{cris,L}(V^*(1)))=\cR_L(\mathbf{B})\otimes_{\mathbb{A}_L^+}N(T^*(1)) \]
and the comparison isomorphism \eqref{f:comp-iso} induces isomorphisms
\begin{align*}
\mathrm{comp}_M:M[\frac{1}{t_{\mathfrak{Y}}}]\cong  \cR_L(\mathbf{B})[\frac{1}{t_{\mathfrak{Y}}}]\otimes_L D_{cris,L}(V(\tau^{-1}))
\end{align*}
and
\begin{align*}
\mathrm{comp}_{\check{M}}:\check{M}[\frac{1}{t_{\mathfrak{Y}}}]\cong  \cR_L(\mathbf{B})[\frac{1}{t_{\mathfrak{Y}}}]\otimes_L D_{cris,L}(V^*(1)).
\end{align*}
 By \cite[§3.4/5]{BSX} an analogue of Kisin-Ren modules exists for $\mathfrak{Y}=\mathfrak{X}$, i.e., if we take $M:=D_{\mathrm{rig}}^\dagger(V(\tau^{-1}))_\mathfrak{X}$ and $\check{M}=D_{\mathrm{rig}}^\dagger(V^*(1)) _\mathfrak{X}$ we obtain
analogous comparison isomorphisms
\begin{align}\label{f:compM}
\mathrm{comp}_M:M[\frac{1}{t_{\mathfrak{Y}}}]\cong  \cR_L(\mathfrak{X})[\frac{1}{t_{\mathfrak{Y}}}]\otimes_L D_{cris,L}(V(\tau^{-1}))
\end{align}
and
\begin{align*}
\mathrm{comp}_{\check{M}}:\check{M}[\frac{1}{t_{\mathfrak{Y}}}]\cong  \cR_L(\mathfrak{X})[\frac{1}{t_{\mathfrak{Y}}}]\otimes_L D_{cris,L}(V^*(1)).
\end{align*}
which this time stem from \cite[Prop.\ 3.42]{BSX} by base change $\cR_L(\mathfrak{X})\otimes_{{\mathcal{O}_K(\mathfrak{X})}}-$ using the inclusion $\mathcal{O}_L(\mathfrak{X})[Z^{-1}]\subseteq \cR_L(\mathfrak{X})[\frac{1}{t_{\mathfrak{Y}}}].$
Moreover, these comparison isomorphism for $\mathbf{B}$ and $\mathfrak{X}$ are compatible with regard to the equivalence of categories \eqref{f:equivalenceBSX} by \cite[Thm.\ 3.48]{BSX}.
Note that for $c=\mathrm{comp}_M$ and $D= D_{cris,L}(V(\tau^{-1})) $ we have
\begin{align}\label{f:checkcomp}
  \mathrm{comp}_{\check{M}} & =    (\mathrm{comp}_{\Omega^1_\cR}\otimes_L\id_{D^*})\circ\check{c}
\end{align}
using the identifications $\Omega^1_\cR\cong\cR(\chi_{LT})$ and
\begin{align*}
  D_{cris,L}(V^*(1)) & \cong D^*\otimes D_{cris,L}(L(\chi_{LT})).
\end{align*}
We set  $b:=\mathrm{comp}_{\Omega^1_\cR}(t_{\mathfrak{Y}}^{-1}d\log_{LT})=\frac{1}{t_{\mathfrak{Y}}}\otimes \eta\in D_0:=D_{cris,L}(L(\chi_{LT}))$  {  and
\[\tilde{\nabla}:=\left\{
                    \begin{array}{ll}
                      \nabla, & \hbox{if $\mathfrak{Y}=\mathfrak{X}$;} \\
                      \frac{\nabla}{\Omega}, & \hbox{if $\mathfrak{Y}=\mathbf{B}$ (and $\Omega\in K$).}
                    \end{array}
                  \right.
\]}
\begin{remark}\label{rem:nablaR}
As operators on $\cR$ we have the equalities
\[ \nabla=t_\mathfrak{Y}\partial^\mathfrak{Y}_{inv}\mbox{ and } \tilde{\nabla}=t_\mathfrak{Y}\tilde{\partial}^\mathfrak{Y}_{inv},\]
where we define $\partial^\mathbf{B}_{inv}:= \partial_{inv}$ and \[\tilde{\partial}^\mathfrak{Y}_{inv}:=\left\{
                    \begin{array}{ll}
                      \partial^\mathfrak{X}_{inv}, & \hbox{if $\mathfrak{Y}=\mathfrak{X}$;} \\
                      \frac{\partial_{inv}}{\Omega}, & \hbox{if $\mathfrak{Y}=\mathbf{B}$ (and $\Omega\in K$).}
                    \end{array}
                  \right.\]
 Indeed, for $\mathfrak{Y}=\mathbf{B}$  the fact \eqref{f:nablaonR} grants these equalities of operators  on the subring  ${\mathcal{O}_K(\mathbf{B})}$. Concerning the ring  $\cR_K(\mathbf{B})$ we note that $\nabla$ is acting as a continuous derivation as can be shown similarly as in \cite[Lem.\ 2.1.2]{KR}, while for the operator $t_{\mathbf{B}}\partial_{\mathrm{inv}}$ this is clear anyway. Thus the same equalities hold for $\cR$ by \ref{rem:denseB}. Indeed, on the localisation ${\mathcal{O}_K(\mathbf{B})}_{Z^{\mathbb{N}}}$ it extends uniquely by the derivation property and then it extends uniquely by continuity to $\cR_K(\mathbf{B})$. Regarding $\mathfrak{Y}=\mathfrak{X}$ note that all operators are defined over $K$. Since the equality can be checked over $\mathbb{C}_p$, the claim follows from Remark \ref{rem:logXlogLT} and the previous case $\mathfrak{Y}=\mathbf{B}$.
\end{remark}
\begin{lemma}\label{lem:compOmega1}  The following diagram commutes
\[\xymatrix{
  \Omega^1_\cR[\frac{1}{t_{\mathfrak{Y}}}]\otimes D^*  \ar[rrr]^-{\mathrm{comp}_{\Omega^1_\cR}\otimes_L\id_{D^*}} & && \cR[\frac{1}{t_{\mathfrak{Y}}}]\otimes D^*\otimes D_0 \\
  (\Omega^1_\cR\otimes_L D^* )^{\psi_L=0} \ar@{^(->}[u]  &&& \cR^{\psi_L=0}\otimes D^*\otimes D_0\ar@{^(->}[u]^{\tilde{\nabla}}  \\
   \cR_K(\Gamma_L)\otimes_L D^*\ar[u]_{\cong}^{\mathfrak{M}^{\Omega^1}\otimes \id_{D^*}} \ar[rrr]^-{\id_{  \cR_K(\Gamma_L)\otimes_L D^*}\otimes b} & && \cR_K(\Gamma_L)\otimes_L D^*\otimes D_0\ar[u]_{\cong}^{\mathfrak{M}\otimes \id_{D^*\otimes D_0}}    }\]
{  assuming $\Omega\in K,$ if $\mathfrak{Y}=\mathbf{B}$.   }
\end{lemma}

\begin{proof}   We first give the proof for $\mathfrak{Y}=\mathbf{B}$.
Observe, since on $D^*$ we have the identity throughout, that the commutativity of the above diagram follows from the commutativity of
\begin{equation}\label{f:nablacomp}
 \xymatrix{
    \cR_K(\Gamma_L) \ar@{=}[dd]_{} \ar[r]^{\mathfrak{M}^{\Omega^1}} & (\Omega^1_\cR)^{\psi_L=0} \ar@{^(->}[r]^{} & \Omega^1_\cR[\frac{1}{t_{\mathfrak{Y}}}]^{\psi_L=0} \ar[ddd]^{\mathrm{comp}_{\Omega^1_\cR}}_{\cong} \\
    &\cR(\chi_{LT})^{\psi_L=0}\ar[u]_{\cong}& \\
 \cR_K(\Gamma_L) \ar[d]_{\tilde{\nabla}} \ar[r]^(0.3){ \mathfrak{M}\otimes b} & \cR^{\psi_L=0} \otimes_L D_{cris,L}(L(\chi_{LT}))\ar[u]^{\tilde{\partial}^\mathfrak{Y}_{inv}\otimes t_{\mathfrak{Y}}} \ar[d]_{t_{\mathfrak{Y}}\tilde{\partial}^\mathfrak{Y}_{inv}\otimes\id }  &  \\
   \cR_K(\Gamma_L) \ar[r]^(0.3){\mathfrak{M}\otimes b} & \cR^{\psi_L=0} \otimes_L D_{cris,L}(L(\chi_{LT})) \ar@{^(->}[r]^{} & (\cR[\frac{1}{t_{\mathfrak{Y}}}]\otimes_L   D_{cris,L}(L(\chi_{LT})))^{\psi_L=0}  }
\end{equation}
where the map $\tilde{\partial}^\mathfrak{Y}_{inv}\otimes t_{\mathfrak{Y}}:\cR \otimes_L D_{cris,L}(L(\chi_{LT}))\to  \cR(\chi_{LT}) $ sends $f\otimes\frac{1}{t_{\mathfrak{Y}}}\otimes\eta$ to $\tilde{\partial}^\mathfrak{Y}_{inv} f\otimes \eta$ and the composite with the natural identification $\cR(\chi_{LT})\cong \Omega^1$, which sends $\eta$ to $d\log_{LT},$ is the map $\frac{d}{\Omega}:\cR\to \Omega^1_{\cR}$ upon identifying $\cR \otimes_L D_{cris,L}(L(\chi_{LT}))$ with $\cR$ by sending $f\otimes\frac{1}{t_{\mathfrak{Y}}}\otimes\eta$ to $f.$
Remark \ref{rem:nablaR} implies the commutativity of the left lower corner while for the upper left corner it follows from  \eqref{f:twist}, the easily checked identity $\tilde{\partial}^\mathfrak{Y}_{inv}\eta(1,Z)= \eta(1,Z)$ and \eqref{f:Twdlog}
\begin{align*}
\mathfrak{M}^{\Omega^1}(\lambda)&=\big(Tw_{\chi_{LT}}(\lambda)\cdot\eta(1,Z)\big)d\log_{LT}\notag\\
&=\big(Tw_{\chi_{LT}}(\lambda)\cdot\tilde{\partial}^\mathfrak{Y}_{inv}\eta(1,Z)\big)d\log_{LT}\notag\\
&= \tilde{\partial}^\mathfrak{Y}_{inv}(\lambda\cdot \eta(1,Z))d\log_{LT}\notag\\
\end{align*}
Finally, since
$\eta(1,Z)\otimes b\in\cR^{\psi_L=0} \otimes_L D_{cris,L}(L(\chi_{LT}))$ is sent up to $\eta(1,Z)d\log_{LT}$ and down to $t_{\mathfrak{Y}}\eta(1,Z)\otimes b$, the compatibility with $\mathrm{comp}_{\Omega^1_\cR}$ is easily checked.
The same proof works for   $\mathfrak{Y}=\mathfrak{X}$ by using \eqref{f:twistX} instead of \eqref{f:twist} and replacing $\eta(1,Z)$ and $\frac{d}{\Omega}$ by $\mathrm{ev}_1$ and $d$, respectively.
\end{proof}

Now we introduce a pairing  - if $\mathfrak{Y}=\mathbf{B}$  assuming $\Omega\in K$ as usual -
\[[\;,\;]:=[\;,\;]_{D_{cris,L}(V(\tau^{-1}))}:{\cR^{\psi_L=0}\otimes_L D_{cris,L}(V^*(1))} \times  {\cR^{\psi_L=0}\otimes_L D_{cris,L}(V(\tau^{-1}))}  \to \cR_K(\Gamma_L)\]
by requiring that the following diagram becomes commutative
\begin{equation}\label{f:def[]}
   \xymatrix{
   \cR^{\psi_L=0}\otimes D^*\otimes D_0   \ar@{}[r]|{\times} & \cR^{\psi_L=0} \otimes_L D \ar[r]^-{[\;,\;]} &\cR_K(\Gamma_L) \ar@{=}[d]^{} \\
  \cR_K(\Gamma_L)\otimes_L D^*\otimes D_0\ar[u]_{\cong}^{\mathfrak{M}\otimes \id_{D^*\otimes D_0}} \ar@{}[r]|{\times} & \cR_K(\Gamma_L)\otimes_L D\ar[u]_{\cong}^{\sigma_{-1}\mathfrak{M}\circ\upiota_*\otimes\id}\ar[r]^(0.5){} &\cR_K(\Gamma_L),   }
\end{equation}
 where the bottom line sends $(\lambda\otimes l\otimes\beta b,\mu\otimes d)$ to $\lambda\mu \beta l(d).$

Combining the Lemmata \ref{lem:Iwasawacheck} and \ref{lem:compOmega1} we obtain for $N=\cR\otimes_L {D_{cris,L}(V(\tau^{-1}))}$
\begin{lemma} \label{lem:[]versusIW} $[-,-]_{D_{cris,L}(V(\tau^{-1}))}=\frac{q-1}{q}\{ {\nabla}(\mathrm{comp}_{\Omega^1_\cR}\otimes_L\id_{D^*})^{-1}(-),-\}^0_{N,Iw}.$

\end{lemma}

Setting
\begin{align*}
 M' & :=\mathrm{comp}^{-1}(\cR^{\psi_L=0}\otimes_L D_{cris,L}(V(\tau^{-1})) ) \mbox{  and}\\
  \check{M}' & :=\mathrm{comp}^{-1}(\cR^{\psi_L=0}\otimes_L D_{cris,L}(V^*(1)) )
\end{align*}  we obtain

\begin{theorem}\label{thm-recproclawKKK} Assume $\Omega\in K$, if $\mathfrak{Y}=\mathbf{B}$.
For all $x\in \check{M}'\cap (\check{M}^{\psi_L=0}) $ and $y\in M'\cap ({M}^{\psi_L=0})$ it holds
\[\frac{q-1}{q}\KKl {\nabla}x,y\KKr^0_{Iw}=[x,y],\]
i.e., the following diagram commutes on the vertical intersections
\[\xymatrix{
  \check{M}^{\psi_L=0} \ar@{^(->}[d]_{}  \ar@{}[r]|{\times} & M^{\psi_L=0} \ar@{^(->}[d]_{} \ar[rr]^-{\frac{q-1}{q}{{\nabla}}\{,\}^0_{M,Iw}} && \cR_K(\Gamma_L) \ar@{=}[ddd]^{} \\
  (\cR[\frac{1}{t_{\mathfrak{Y}}}]\otimes_{\cR} \check{M})^{\psi_L=0} \ar[d]_{\mathrm{comp}_{\check{M}}}^{\cong} \ar@{}[r]|{\times} & (\cR[\frac{1}{t_{\mathfrak{Y}}}]\otimes_{\cR} {M})^{\psi_L=0} \ar[d]_{\mathrm{comp}_M}^{\cong} \\
 (\cR[\frac{1}{t_{\mathfrak{Y}}}]\otimes_L D_{cris,L}(V^*(1)))^{\psi_L=0}   \ar@{}[r]|{\times} & (\cR[\frac{1}{t_{\mathfrak{Y}}}]\otimes_L D_{cris,L}(V(\tau^{-1})))^{\psi_L=0}    \\
 {\cR^{\psi_L=0}\otimes_L D_{cris,L}(V^*(1))}\ar@{^(->}[u] \ar@{}[r]|{\times} & {\cR^{\psi_L=0}\otimes_L D_{cris,L}(V(\tau^{-1}))} \ar@{^(->}[u]\ar[rr]^-{[,]_{{D_{cris,L}(V(\tau^{-1}))}}} && \cR_K(\Gamma_L).}\]
\end{theorem}

\begin{proof}
  Combine Lemmata \ref{lem:[]versusIW} and \ref{lem:checkIW} using \eqref{f:checkcomp}.
\end{proof}
\paragraph{Interpretation of the abstract reciprocity formula in terms of the $D_{cris,L}$-pairing}

The canonical pairing $D_{cris,L}(V^*(1))\times D_{cris,L}(V(\tau^{-1}))\to D_{cris,L}(L(\chi_{LT}))$ extends to a pairing {  - if $\mathfrak{Y}=\mathbf{B}$  assuming $\Omega\in K$ as usual - }
\begin{equation*}
\xymatrix{
 \cR^{\psi_L=0}\otimes_L D_{cris,L}(V^*(1)) \ar@{}[r]|{\times} &  \cR^{\psi_L=0}\otimes_L D_{cris,L}(V(\tau^{-1})) \ar[r]^(0.5){[,]_{cris}} & \cR^{\psi_L=0}\otimes_L D_{cris,L}(L(\chi_{LT}))}
\end{equation*}
by requiring that the following diagram is commutative (in which the lower one is induced by multiplication within $\cR_K(\Gamma_L)$ and the natural duality paring on $D_{cris,L}$) {\scriptsize
\begin{equation}\label{f:def[]cris}
\xymatrix{
       {\cR^{\psi_L=0}\otimes_L D_{cris,L}(V^*(1))} \ar@{}[r]|{\times} &  {\cR^{\psi_L=0}\otimes_L D_{cris,L}(V(\tau^{-1}))}  \ar[r]^(0.5){[,]_{cris}} & {\cR^{\psi_L=0}}\otimes_L D_{cris,L}(L(\chi_{LT}))\\
  \cR_K(\Gamma_L)\otimes_L D_{cris,L}(V^*(1))  \ar[u]_{  \mathfrak{M}\otimes \id}\ar@{}[r]|{\times} &  \cR_K(\Gamma_L)\otimes_L D_{cris,L}(V(\tau^{-1}))  \ar[u]_{ \sigma_{-1}\mathfrak{M}\circ \upiota_*\otimes \id}  \ar[r] & \cR_K(\Gamma_L)\otimes_L D_{cris,L}(L(\chi_{LT}))\ar[u]_{  \mathfrak{M}\otimes \id}}
\end{equation}}

Note that \[\mathrm{comp}\left([x,y]\cdot { \ev_1}\otimes(t_{\mathfrak{Y}}^{-1}\otimes \eta)\right)=[x  ,y]_{cris}.\] Hence using the diagram \eqref{f:nablacomp} Thm.\  \ref{thm-recproclawKKK} is also equivalent to
\[\mathrm{comp}\circ \mathfrak{M}^{\Omega^1}\circ \tilde{\Omega}\frac{q-1}{q}\KKl x,y\KKr^0_{Iw}=[\mathrm{comp}(x),\mathrm{comp}(y)]_{cris},\]
i.e., the 'commutativity' (whenever it makes sense) of the following diagram {\tiny
\begin{equation}
\xymatrix{
  D_{\mathrm{rig}}^\dagger(V^*(1))^{\psi_L=\frac{q}{\pi_L}}[\frac{1}{t_{\mathfrak{Y}}}] \ar[d]_{1-\varphi_L} \ar@{}[r]|{\times}& D_{\mathrm{rig}}^\dagger(V(\tau^{-1}))^{\psi_L=1} [\frac{1}{t_{\mathfrak{Y}}}] \ar[d]_{1-\frac{\pi_L}{q}\varphi_L} \ar@{-->}[r]^(0.65){\tilde{\Omega}\frac{q-1}{q}\{,\}_{Iw}} & \cR_K(\Gamma_L) \ar@{=}[d] &\\
   D_{\mathrm{rig}}^\dagger(V^*(1))^{\psi_L=0}[\frac{1}{t_{\mathfrak{Y}}}]\ar[d]_{\mathrm{comp}}^{\cong} \ar@{}[r]|{\times} & D_{\mathrm{rig}}^\dagger(V(\tau^{-1}))^{\psi_L=0} [\frac{1}{t_{\mathfrak{Y}}}] \ar[d]_{\mathrm{comp}}^{\cong} \ar@{-->}[r]^(0.65){\tilde{\Omega}\frac{q-1}{q}\{,\}_{Iw}^0} &  \cR_K(\Gamma_L) \ar[r]^{\mathfrak{M}^{\Omega^1}} & \cR(\chi_{LT})[\frac{1}{t_{\mathfrak{Y}}}]^{\psi_L=0} \ar[d]^{\mathrm{comp}}_{\cong} \\
  {\cR^{\psi_L=0}[\frac{1}{t_{\mathfrak{Y}}}]\otimes_L D_{cris,L}(V^*(1))} \ar@{}[r]|{\times} &  {\cR^{\psi_L=0}[\frac{1}{t_{\mathfrak{Y}}}]\otimes_L D_{cris,L}(V(\tau^{-1}))}    \ar@{-->}[rr]^(0.5){[,]_{cris}} && {\cR^{\psi_L=0}[\frac{1}{t_{\mathfrak{Y}}}]}\otimes_L D_{cris,L}(L(\chi_{LT}))   }
\end{equation}}
for $\mathfrak{Y}=\mathbf{B}$ while for $\mathfrak{Y}=\mathfrak{X}$ one has to decorate the $D_{\mathrm{rig}}^\dagger$s with index $\mathfrak{X}  .$\\

\noindent
{\bf Question:} Can one extend the definition of $[\;,\;]$ and $\KKl\;,\;\KKr$ to $(\check{M}[\frac{1}{t_{\mathfrak{Y}}}])^{\psi_L=0}\times (M[\frac{1}{t_{\mathfrak{Y}}}])^{\psi_L=0}$ by perhaps enlarging the target $\cR_K(\Gamma_L)$ by an appropriate localisation, which reflects the inversion of $t_{\mathfrak{Y}}$ somehow?
\newpage

\section{Application}\label{sec:application}

\subsection{The regulator map}

Recall that we write $\tau^{-1}=\chi_{LT}\chi_{cyc}^{-1}.$ Let   $T$ be in $\Rep_{o_L,f}^{cris}(G_L)$    such that $T(\tau^{-1})$ belongs to  $\Rep_{o_L,f}^{cris,an}(G_L)$ with all   Hodge-Tate weights   in $[0,r],$ and such that $V:=L\otimes_{o_L} T$ does not have any quotient isomorphic to $L(\tau).$ Then we define the regulator maps
\begin{align*}
\mathbf{L}_V: & H^1_{Iw}(L_\infty/L,T)\to  D(\Gamma_L,\Cp)\otimes_L D_{cris,L}(V(\tau^{-1})), \\
  \mathcal{L}_V^0: & H^1_{Iw}(L_\infty/L,T)\to   \mathcal{O}_L(\mathbf{B})^{\psi_L=0}\otimes_L D_{cris,L}(V(\tau^{-1})), \\
  \mathcal{L}_V: & H^1_{Iw}(L_\infty/L,T)\to D(\Gamma_L,\Cp)\otimes_L D_{cris,L}(V)
\end{align*}
as (part of) the composite \small
\begin{align}\label{f:defregulator}
 \notag H^1_{Iw}(L_\infty{/L},T)&\cong D_{LT}(T(\tau^{-1}))^{\psi_L=1}=N(T(\tau^{-1}))^{\psi_{D_{LT}(T(\tau^{-1}))}=1} \xrightarrow{(1-\frac{\pi_L}{q}\varphi_L)}\varphi^*_L(N(V(\tau^{-1})))^{\psi_L=0}\\
  &\hookrightarrow \mathcal{O}^{\psi_L=0} \otimes_L D_{cris,L}(V(\tau^{-1})) \subseteq {\mathcal{O}_{\Cp}(\mathbf{B})}^{\psi_L=0}\otimes_L D_{cris,L}(V(\tau^{-1}))\\
\notag &\xrightarrow{\mathfrak{M}^{-1}\otimes \id}D(\Gamma_L,\Cp)\otimes_L D_{cris,L}(V(\tau^{-1}))\to D(\Gamma_L,\Cp)\otimes_L D_{cris,L}(V)
\end{align}\normalsize
using  \cite[Thm.\ 5.13]{SV15}, Lemma \ref{notrivialquot}, the inclusion \eqref{f:phistarcomp} and where the last map sends $\mu\otimes d \in D(\Gamma_L,\Cp)\otimes_L D_{cris,L}(V(\tau^{-1}))$ to $\mu\otimes d \otimes \mathbf{d}_1\in D(\Gamma_L,\Cp)\otimes_L D_{cris,L}(V(\tau^{-1}))\otimes_L D_{cris,L}(L(\tau))\cong D(\Gamma_L,\Cp)\otimes_L D_{cris,L}(V). $ Note that $D:=D_{cris,L}(L(\tau))=D^0_{dR,L}(L(\tau)) = L \mathbf{d}_1$ with $\mathbf{d}_1= t_{LT} t_{\mathbb{Q}_p}^{-1} \otimes (\eta^{\otimes -1} \otimes \eta^{cyc})$, where $L(\chi_{LT})=L\eta$ and $L(\chi_{cyc})=L\eta^{cyc}.$

 Alternatively, in order to stress that the regulator is essentially the map $1- \varphi_L, $ one can rewrite this as \small
 \begin{align}
 & H^1_{Iw}(L_\infty/L,T)\cong D_{LT}(V(\tau^{-1}))^{\psi_L=1}=N(T(\tau^{-1}))^{\psi_{D_{LT}(T(\tau^{-1}))}=1} \hookrightarrow N(V(\tau^{-1}))^{\psi_{D_{LT}(V(\tau^{-1}))}=1}\otimes_L D\\
  \notag &\xrightarrow{1- \varphi_L}\varphi^*_L(N(V(\tau^{-1})))^{\psi_L=0}\otimes_L D
  \hookrightarrow \mathcal{O}^{\psi_L=0} \otimes_L D_{cris,L}(V(\tau^{-1}))\otimes_L D \subseteq {\mathcal{O}_{\Cp}(\mathbf{B})}^{\psi_L=0}\otimes_L D_{cris,L}(V)\\
\notag &\xrightarrow{\mathfrak{M}^{-1}\otimes \id}D(\Gamma_L,\Cp)\otimes_L D_{cris,L}(V)
\end{align}\normalsize
 where the $ \hookrightarrow $ in the first line sends $n$ to $n\otimes \mathbf{d}_1$ and the $\varphi_L$ now acts diagonally.
By construction, this regulator map $\mathcal{L}_V$ takes values in $ D(\Gamma_L,K)^{G_L,*}\otimes_L D_{cris,L}(V),$ where the twisted action of $G_L$ on the distribution algebra is induced by the Mellin-transform as in (ii) of Prop.\  \ref{prop:twistinvariance}.\\

 We write $\nabla_{\mathrm{Lie}}\in \mathrm{Lie}(\Gamma_L)$ for the element in the Lie algebra of $\Gamma_L$ corresponding to $1$ under the identification $\mathrm{Lie}(\Gamma_L)=L $.
\begin{proposition}
   \label{prop:regulatortwisting}
   The regulator maps for $V$ and $V(\chi_{LT})$ - assuming that both representations satisfy the conditions above - are related by
   \[ \mathcal{L}_{V(\chi_{LT})}(x \otimes \eta) = \nabla_{\mathrm{Lie}} \cdot \left( \frac{1}{\Omega} \mathrm{Tw}_{\chi_{LT}^{-1}}(\mathcal{L}_{V}(x))\otimes t_{LT}^{-1} \eta \right),\]
   i.e., the following $\Gamma_L$-equivariant diagram commutes:
\begin{equation*}
   \xymatrix{
     H^1_{Iw}(L_\infty/L,V(\chi_{LT}))  \ar[d]_{\cong } \ar[rr]^-{ \mathcal{L}_{V(\chi_{LT})}} && D(\Gamma_L,\mathbb{C}_p)\otimes_L D_{cris,L} (V(\chi_{LT})) \\
      H^1_{Iw}(L_\infty/L,V)\otimes_{o_L} L(\chi_{LT}) \ar[rr]^-{\mathcal{L}_{V}\otimes  L(\chi_{LT}) } && D(\Gamma_L,\mathbb{C}_p)\otimes_L D_{cris,L}(V)\otimes_L  L(\chi_{LT})\ar[u]^{\frac{\nabla_{\mathrm{Lie}}Tw_{\chi^{-1}}}{\Omega} \otimes t_{LT}^{-1} },  }
 \end{equation*}
  \end{proposition}
\begin{proof}
  Analogous to \cite[Pro.\ 3.1.4]{LVZ15}. Note that the period $\Omega$ enters due to \eqref{f:twist}.
\end{proof}
This twisting property can be used to drop the condition concerning the Hodge-Tate weights in the definition of the regulator map, i.e., upon  replacing $D(\Gamma_L,\mathbb{C}_p)$ in the target by its  total ring of quotients one can extend the regulator map as usual to all  $T$   in $\Rep_{o_L,f}^{cris}(G_L)$    such that $T(\tau^{-1})$ belongs to  $\Rep_{o_L,f}^{cris,an}(G_L)$.

In order to better understand the effect of twisting we have the following
\begin{lemma}\label{lem:nablatwisting}
For $\mu\in D(\Gamma_L,K)$ we have
\[\frac{1}{\Omega}(\nabla_{\mathrm{Lie}}Tw_{\chi^{-1}})(\mu)= \mathfrak{M}^{-1}(t_{LT}\mathfrak{M}(\mu))\]
and for all $n\geq 1$
\[(\nabla_{\mathrm{Lie}}Tw_{\chi^{-1}}(\mu))(\chi_{LT}^n)=n\mu(\chi_{LT}^{n-1}).\]
\end{lemma}
\begin{proof}
  The first claim follows by combining \eqref{f:nablaMellin} with \eqref{f:twist}, while the second claim is just Lemma \ref{log} applied to the first.
\end{proof}

One significance of regulator maps is that it should interpolate (dual) Bloch-Kato exponential maps. We shall prove such interpolation formulae in subsection \ref{subsec:interpolation} by means of a reciprocity formula.

\subsubsection{The basic example}\label{sec:reg-ex}

Setting $U:= \varprojlim_n  o_{L_n}^\times   $ with   transition maps given by the norm we are looking for a map
\[\mathcal{L}: U\otimes_{\mathbb{Z}}\TLT^* \to D(\Gamma_L,\Cp)\otimes_L D_{cris,L}(L(\tau) )\]

such that
\begin{equation}\label{f:CWinterpolation}
  \frac{ \Omega^r}{r!} \frac{1-\pi_L^{-r}}{1-\frac{\pi_L^r}{q}} \mathcal{L}(u\otimes a \eta^*)(\chi_{LT}^{ r})\otimes (t_{LT}^{r-1} \otimes \eta^{\otimes -r+1} )=CW(u\otimes a \eta^{\otimes -r})
\end{equation}
for all $r\geq 1,u\in U,a\in o_L,$ where $CW$ denotes the diagonal map in

\begin{theorem}[A special case of Kato's explicit reciprocity law, {\cite[Cor.\ 8.7]{SV15}}] \label{thm:Kato}
For $r\geq 1$ the diagram
\begin{equation*}
  \xymatrix{
  {U} \otimes_{\mathbb{Z}} \TLT^{\otimes -r} \ar[d]_{-\kappa\otimes \id} \ar[rrdd]^{\qquad ''(1-\pi_L^{-r})r\psi_{CW}^r({_-}) \mathbf{d}_r\, ''} &   \\
  H^1_{Iw}(L_\infty{/L},\TLT^{\otimes -r}(1)) \ar[d]_{\mathrm{\mathrm{cor}}}  &   \\
  H^1(L, \TLT^{\otimes -r}(1)) \ar[rr]^-{\exp^*} & & D^0_{dR,L}(V_\pi^{\otimes -r}(1)) = L \mathbf{d}_r ,  }
\end{equation*}
commutes, i.e., the diagonal map sends $u\otimes a \eta^{\otimes -r}$ to
\begin{equation*}
  a(1-\pi_L^{-r})r\psi_{CW}^r(u) \mathbf{d}_r = a \frac{1-\pi_L^{-r}}{(r-1)!} \partial_{\mathrm{inv}}^r \log g_{u,\eta}(Z)_{| Z=0} \mathbf{d}_r \
\end{equation*}
with $\mathbf{d}_r := t_{LT}^rt_{\mathbb{Q}_p}^{-1} \otimes (\eta^{\otimes -r} \otimes \eta^{cyc})$.
\end{theorem}

We set $\mathcal{L}=\frak{L}\otimes \mathbf{d}_1$ with $\frak{L}$ given as follows
\[\mathfrak{L}: U\otimes \TLT^*\xrightarrow{\nabla} o_L[[\omega_{LT}]]^{\psi_L={1}}\xrightarrow{(1-\frac{\pi_L}{q}\varphi)}{\mathcal{O}_{\Cp}(\mathbf{B})}^{\psi_L=0}\xrightarrow{\log_{LT}\cdot}{\mathcal{O}_{\Cp}(\mathbf{B})}^{\psi_L=0}\xrightarrow{\mathfrak{M}^{-1}}D(\Gamma_L,\Cp),\]
where the map $\nabla$ has been defined in \cite[\S 6]{SV15} as
the homomorphism
\begin{align*}
  \nabla : U \otimes_\mathbb{Z} T^* & \longrightarrow  o_L[[\omega_{LT}]]^{\psi=1} \\
       u \otimes a\eta^* & \longmapsto a \frac{\partial_\mathrm{inv}(g_{u,\eta})}{g_{u,\eta}}(\omega_{LT}) \ .
\end{align*}
 Note that due to the multiplication by $\log_{LT}$ the maps $\mathcal{L},\; \frak{L}$ are not $\Gamma_L$-equivariant.
Using Lemmata \ref{log}, \ref{phi} we obtain

\begin{align}\label{f:decentKato}\notag
\frak{L}(u\otimes a\eta^*)(\chi_{LT}^r)&= a \mathfrak{M}^{-1}(\log_{LT}(1-\frac{\pi_L}{q}\varphi)\partial_\mathrm{inv} \log g_{u,\eta})(\chi_{LT}^r)\\
&= a \Omega^{-1} r \mathfrak{M}^{-1}(1-\frac{\pi_L}{q}\varphi)\partial_\mathrm{inv} \log g_{u,\eta})(\chi_{LT}^{r-1})\\ \notag
&= ar  \Omega^{-r} (1-\frac{\pi_L}{q}\pi_L^{r-1})(\partial_\mathrm{inv}^{r-1}\partial_\mathrm{inv} \log g_{u,\eta})_{|Z=0}\\\notag
&= ar  \Omega^{-r} (1-\frac{\pi_L^r}{q})(\partial_\mathrm{inv}^{r-1}\partial_\mathrm{inv} \log g_{u,\eta})_{|Z=0},
\end{align}
i.e., $\mathcal{L}$ satisfies \eqref{f:CWinterpolation}, indeed.
By construction and Proposition \ref{prop:twistinvariance} the image of $\mathcal{L}$ actually lies in the ${G_L}$-invariants:
\[\mathcal{L}: U\otimes_{\mathbb{Z}}\TLT^* \to D(\Gamma_L,K)^{G_L}\otimes_L D_{cris,L}(L(\tau) ).\]

 We claim that
\be\label{f:claim}U\otimes_{\mathbb{Z}} \TLT^*  \xrightarrow{-\kappa \otimes \TLT^*}  H^1_{Iw}(L_\infty/L,o_L(\tau)) \xrightarrow{\mathcal{L}_{L(\tau\chi_{LT})}\otimes o_L(\chi_{LT}^{-1}) \otimes t_{LT}}D(\Gamma_L,\Cp)\otimes_L D_{cris,L}(L(\tau))
\ee
  coincides with
 \[\mathcal{L}: U\otimes_{\mathbb{Z}}\TLT^* \to D(\Gamma_L,\Cp)\otimes_L D_{cris,L}(L(\tau) ).\]

Indeed,  from to the commutativity of  the following diagram (cp.\ with  \cite[Appendix C]{LVZ15} for $L=\Qp$), in which $\mathcal{L}_{L(\tau\chi_{LT})}\otimes \mathbf{d}_1^{\otimes -1}$ or more generally  $\mathcal{L}_{L(\tau\chi_{LT}^r)}\otimes \mathbf{d}_1^{\otimes -1}$, $r\geq 1$, shows up at $\diamondsuit$, the above claim immediately follows by tensoring the diagram for $r=1$ with $o_L(\chi_{LT}^{-1})$ and then composing with the multiplication by $t_{LT}$; we set $e_r := t_{LT}^{-r}\otimes \eta^{\otimes r}\in D_{cris,L}(L(\chi_{LT}^r))$: \Footnote{Also for $r=0$ the diagram commutes up to the upper diagram perhaps.}

\[\footnotesize
\xymatrix{
  {{U}\otimes\TLT^{\otimes(r-1)}} \ar[d]_{\nabla\otimes\eta^{\otimes r}} \ar[rr]^{-\kappa\otimes\TLT^{\otimes(r-1)}} && {\H^1_{Iw}(L_\infty/L,o_L(\tau\chi_{LT}^r))} \ar[d]^{\cong}\ar@<0ex>`r[d]`_ld[ddddd]^{ \diamondsuit}[dddd] \\
   {(o_L[[\omega_{LT}]]\otimes \eta^{\otimes r})^{\psi_L=1}} \ar[d]_{(1-\frac{\pi_L}{q}\vp_L)\otimes \id} \ar@{^(->}[r]^{ } & (\omega_{LT}^{-r}o_L[[\omega_{LT}]]\otimes \eta^{\otimes r})^{\psi_L=1}\ar@{=}[r] \ar[d]^{(1-\frac{\pi_L}{q}\vp_L)\otimes \id} & N(o_L(\chi_{LT}^r))^{\psi_L=1} \ar[d]^{1-\frac{\pi_L}{q}\vp_L}\\
  {o_L[[\omega_{LT}]][\frac{1}{p}]^{\psi_L=0}\otimes \eta^{\otimes r}} \ar[d]_{\partial_\mathrm{inv}^{-r}\otimes t_{LT}^{-r}} \ar@{^(->}[r]^{ } & {\varphi(\omega_{LT})^{-r}o_L[[\omega_{LT}]][\frac{1}{p}]^{\psi_L=0}\otimes \eta^{\otimes r}} \ar[d]^{t_{LT}^r\otimes t_{LT}^{-r}}\ar@{=}[r] & {\varphi^*(N(L(\chi_{LT}^r)))^{\psi_L=0}\ar[d]^{comp}}\\
    {\cO^{\psi_L=0}\otimes e_r}\ar[dd]_{\mathfrak{M}^{-1}\otimes\id} \ar[r]_(0.45){=t_{LT}^r \partial_\mathrm{inv}^r\otimes\id }^(0.45){\mathfrak{l}_0 \cdots \mathfrak{l}_{r-1}\otimes \id} & {\phantom{m}\cO^{\psi_L=0}\otimes e_r}\ar@{=}[r] & {\cO^{\psi_L=0}\otimes_L D_{cris,L}(L(\chi_{LT}^r))} \ar[d]^{ \mathfrak{M}^{-1}\otimes\id}\\
      & &{D(\Gamma_L,K)^{G_L}\otimes_{L} D_{cris,L}(L(\chi_{LT}^r))} \\
  {D(\Gamma_L,K)^{G_L}\otimes e_r} \ar[d]_{\id\otimes t_{LT}^r} \ar@{=}[rr]^{ } && {D(\Gamma_L,K)^{G_L}\otimes_{L} D_{cris,L}(L(\chi_{LT}^r))} \ar[d]_{\id\otimes t_{LT}^r}\ar[u]_{ \mathfrak{l}_{L(\chi_{LT}^r)}} \\
   {D(\Gamma_L,K)^{G_L}\otimes_L \eta^r} \ar@{=}[rr]^{ } & &{{D(\Gamma_L,K)^{G_L}\otimes_L L(\chi_{LT}^r)} },
   }
   \]
where $\mathfrak{l}_i := t_{LT} \partial_\mathrm{inv} - i$, $\partial_\mathrm{inv}= \frac{d}{dt_{LT}}$. Note that we have
  \be\label{f:nablaMellin}
  \mathfrak{M}^{-1}(\mathfrak{l}_0 f)=\lim_{\gamma \to 1} \frac{\delta_\gamma(\mathfrak{M}^{-1}(f))-\mathfrak{M}^{-1}(f)}{\ell(\gamma))}=\nabla_{\mathrm{Lie}}\mathfrak{M}^{-1}(f),
  \ee
  see \cite[Lem.\ 2.1.4]{KR} for the fact that $\nabla_{\mathrm{Lie}}=t_{LT} \partial_\mathrm{inv} $ as operators on $\mathcal{O}.$
  By abuse of notation we thus also write $\mathfrak{l}_i=\nabla_{\mathrm{Lie}}-i$ for the corresponding element in  $D(\Gamma_L,K),$ compare \cite[\S 2.3]{ST1} for the action of $\mathrm{Lie}(\Gamma_L)$ on and its embedding into $D(\Gamma_L,K).$ Moreover we set $\mathfrak{l}_{L(\chi_{LT}^r)}= \prod_{i=0}^{r-1} \mathfrak{l}_i.$
Note that  $\partial_\mathrm{inv}$ is invertible on $\cO^{\psi_L=0}$ by \cite[Prop. 3.12]{FX}. Finally the map
\[comp: \varphi^*(N(o_L(\chi_{LT}^r)))^{\psi_L=0}\to \cO^{\psi_L=0}\otimes_L D_{cris,L}(L(\chi_{LT}^r))\]
is  \eqref{f:phistarcomp}.

Inspired by Proposition \ref{prop:regulatortwisting} - we define  $\mathcal{L}_{L(\tau)}$ - since $L(\tau)$ does not satisfy the conditions from the beginning of this chapter while $ L(\tau\chi_{LT}) $ does - as a  twist of $\mathcal{L}_{L(\tau\chi_{LT})}$  by requiring the commutativity of the following diagram:
 \begin{equation*}
   \xymatrix{
     H^1_{Iw}(L_\infty/L,o_L(\tau))  \ar[d]_{\cong } \ar[rr]^-{ \mathcal{L}_{L(\tau)}} && D(\Gamma_L,\mathbb{C}_p)\otimes_L D_{cris,L} (L(\tau)) \ar[d]^{\frac{\nabla_{\mathrm{Lie}}Tw_{\chi^{-1}}}{\Omega} \otimes t_{LT}^{-1} }\\
      H^1_{Iw}(L_\infty/L,o_L(\tau\chi_{LT}))\otimes_{o_L} o_L(\chi_{LT}^{-1}) \ar[rr]^-{\mathcal{L}_{L(\tau\chi_{LT})}\otimes o_L(\chi_{LT}^{-1}) } && D(\Gamma_L,\mathbb{C}_p)\otimes_L D_{cris,L}(L(\tau\chi_{LT}))\otimes_L  L(\chi_{LT}^{-1})  }
 \end{equation*}
 which is possible due to the commutativity of  the above diagram. Then
\[\mathcal{L}: U\otimes_{\mathbb{Z}}\TLT^* \to D(\Gamma_L,\Cp)\otimes_L D_{cris,L}(L(\tau) )\] also coincides with
\begin{equation}
\label{f:regulatortau}
U\otimes_{\mathbb{Z}} \TLT^*  \xrightarrow{-\kappa \otimes \TLT^*}  H^1_{Iw}(L_\infty/L,o_L(\tau)) \xrightarrow{(\frac{1}{\Omega } \nabla_{\mathrm{Lie}}Tw_{\chi^{-1}}\otimes\id)\circ\mathcal{L}_{L(\tau)}}D(\Gamma_L,\Cp)\otimes_L D_{cris,L}(L(\tau))
\end{equation}
 by Proposition \ref{prop:regulatortwisting}.

 We refer the interested reader to  \S 5 of \cite{ST2} for an example of a CM-elliptic curve $E$ with supersingular reduction  at $p$ in which they attach to a norm-compatible sequence of elliptic units $e(\mathfrak{a})$ (in the notation of \cite[II 4.9]{dS}) a   distribution $\mu(\mathfrak{a})\in D(\Gamma_L,K)$ in \cite[Prop.\ 5.2]{ST2} satisfying a certain interpolation property with respect to the values of the attached (partial) Hecke-$L$-function. Without going into any detail concerning their setting and instead referring the reader to the notation in (loc.\ cit.) we just want to point out that up to twisting this distribution is the image of $\kappa(e(\mathfrak{a}))\otimes \eta^{-1}$ under the regulator map $\mathcal{L}_{L(\tau)}$:
 \[\mathcal{L}_{L(\tau)}(\kappa(e(\mathfrak{a}))\otimes \eta^{-1})=\Omega Tw_{\chi_{LT}}(\mu(\mathfrak{a}))\otimes \mathbf{d}_1.\]
 Here, $L=\mathbf{K}_p=\mathbf{F}_\wp$ (in their notation) is the unique unramified extension of $\mathbb{Q}_p$ of degree $2,$ $\pi_L=p,$ $q=p^2,$ and the Lubin-Tate formal group is $\hat{E}_\wp$ while $K=\widehat{L_\infty}$.

 Indeed, we have a commutative diagram
 \begin{equation}\label{f:Col}\xymatrix{
   U \ar[d]_{-\kappa(-)\otimes \eta^{-1}} \ar[rr]^-{Col} && D(\Gamma_L,K) \ar[d]^{\Omega Tw_{\chi_{LT}}\otimes \mathbf{d}_1} \\
   H^1_{Iw}(L_\infty/L, o_L(\tau)) \ar[rr]^-{\mathcal{L}_{L(\tau)}} &&  D(\Gamma_L,K)\otimes_L D_{cris,L}(L(\tau) ),   }
 \end{equation}
 where the Coleman map $Col$ is given as the composite in the upper line of the following commutative diagram
 \begin{equation}\label{f:ColL}
   \xymatrix{
     U \ar[d]_{} \ar[r]^{\log g_{-}} & \cO^{\psi_L=\frac{1}{\pi_L} } \ar[d]_{\partial_\mathrm{inv}} \ar[r]^{1-\frac{\varphi_L}{p^2} } & \cO^{\psi_L=0 } \ar[d]_{\partial_\mathrm{inv}} \ar@{=}[r]^{} & \cO^{\psi_L=0 } \ar[d]_{\mathfrak{l}_0} \ar[r]^{\mathfrak{M}^{-1}} & D(\Gamma_L,K) \ar[d]^{\nabla_\mathrm{Lie}} \\
     U\otimes T^*_p \ar[r]^{\nabla} &  {\cO^{\psi_L=1 }} \ar[r]^{1-\frac{\pi_L}{q}\varphi_L } & {\cO^{\psi_L=0 } } \ar[r]^{  \log_{LT}\cdot} & {\cO^{\psi_L=0 }}  \ar[r]^{\mathfrak{M}^{-1}} & D(\Gamma_L,K), }
 \end{equation}
 in which the second line is just $\mathfrak{L}.$ Then the commutativity of \eqref{f:Col} follows by comparing \eqref{f:ColL} with \eqref{f:regulatortau}. Finally, $Col(e(\mathfrak{a}))=\mu(\mathfrak{a})(=\mathfrak{M}^{-1}(g_\mathfrak{a}(Z))\mbox{ in their notation)}$ holds by construction in (loc.\ cit.) upon noting that  on ${{\mathcal{O}}^{\psi_L=\frac{1}{\pi_L}} }$ the operator $1-\frac{\pi}{p^2}\varphi_L\circ\psi_L$, which is used implicitly to define $g_\mathfrak{a}(Z)(= (1-\frac{\pi}{p^2}\varphi_L\circ\psi_L )\log Q_\mathfrak{a}(Z))$, equals $1-\frac{\varphi_L}{p^2} .$

 %
%
%


\subsection{Relation to Berger's and Fourquaux' big exponential map}

 Let  $V$ denote a $L$-analytic representation of $G_L$ and take an integer $h\geq 1$   such that \linebreak $\mathrm{Fil}^{-h}D_{cris,L}(V)=D_{cris,L}(V)$ and such that $D_{cris,L}(V)^{\varphi_L=\pi_L^{-h}}=0$ holds. Under these conditions in \cite{BF} a big exponential map \`{a} la Perrin-Riou
\begin{align*}
\Omega_{V,h}:\left(\cO^{\psi_L=0}\otimes_L D_{cris,L}(V)\right)^{\Delta=0}\to D_{\mathrm{rig}}^\dagger(V)^{\psi_L=\frac{q}{\pi_L}}
\end{align*}
is constructed as follows: According to \cite[Lem.\ 3.5.1]{BF}  there is an exact sequence
\begin{align*}
0\to \bigoplus_{k=0}^h t_{LT}^kD_{cris,L}(V)&^{\varphi_L=\pi_L^{-k}}\to \left( \cO\otimes_{o_L} D_{cris,L}(V)\right)^{\psi_L=\frac{q}{\pi_L}}\xrightarrow{1-\varphi_L} \\ & \cO^{\psi_L=0}\otimes_{L} D_{cris,L}(V)\xrightarrow{\Delta}\bigoplus_{k=0}^hD_{cris,L}(V)/(1-\pi_L^k\varphi_L)D_{cris,L}(V)\to 0,
\end{align*}
where, for $f\in \cO\otimes_L D_{cris,L}(V)$,  $\Delta(f)$ denotes the image of $\bigoplus_{k=0}^h(\partial_\mathrm{inv}^k\otimes \id_{D_{cris,L}(V)})(f)(0)$ in $\bigoplus_{k=0}^hD_{cris,L}(V)/(1-\pi_L^k\varphi_L)D_{cris,L}(V).$ Hence, if $f\in \left(\cO^{\psi_L=0}\otimes_L D_{cris,L}(V)\right)^{\Delta=0}$ there exists $y\in  \left( \cO\otimes_{o_L} D_{cris,L}(V)\right)^{\psi_L=\frac{q}{\pi_L}}$ such that $f=(1-\varphi_L)y.$ Setting $\nabla_i:=\nabla-i$ for any integer $i$, one observes that $\nabla_{h-1}\circ \ldots \circ \nabla_0$ annihilates $ \bigoplus_{k=0}^{h-1} t_{LT}^kD_{cris,L}(V)^{\varphi_L=\pi_L^{-k}}$ whence $\Omega_{V,h}(f):=\nabla_{h-1}\circ \ldots \circ \nabla_0(y)$ is well-defined and belongs under the comparison isomorphism \eqref{f:comp-iso} to $D_{\mathrm{rig}}^\dagger(V)^{\psi_L=\frac{q}{\pi_L}}$ by Prop.\  \ref{elementarydiv}.

Note that $\left(\cO^{\psi_L=0}\otimes_L D_{cris,L}(V)\right)^{\Delta=0}=\cO^{\psi_L=0}\otimes_L D_{cris,L}(V) $ if $D_{cris,L}(V)^{\varphi_L=\pi_L^{-k}}=0$ for all $0\leq k\leq h.$ If this does not hold for $V$ itself, it does hold for $V(\chi_{LT}^{-r})$ for $r$ sufficiently large (with respect to the same $h$).

 In the case $L=\Qp$ the above map specialises to the exponential map due to Perrin-Riou and satisfies the following adjointness property with Loeffer's and Zerbes' regulator map, see \cite[A.2.2]{LVZ15}, where the upper pairing and notation are introduced:\Footnote{{ In (loc.\ cit.) are probably the left and right side interchanged, does this cause a sign \com{??} Ich denke, ja: die obere Paarung liefert beim Verauschen ein Vorzeichen, da dies f\"{u}r das cup-Produkt so ist, w\"{a}hrend die untere nicht! Andererseits scheint mir in (loc.\ cit) von 3.3.7 nach A.2.2 ein Vorzeichen verlorengegangen zu sein, so dass sich beide kompensieren. }}
\begin{equation*}
\xymatrix{
   {D(\Gamma,{\mathbb{Q}_p})\otimes_{\Lambda_L} H_{Iw}(\mathbb{Q}_p,V^*(1)) } \ar@{}[r]|{\times} &  {D(\Gamma,{\mathbb{Q}_p})\otimes_{\Lambda_{\mathbb{Q}_p}} H_{Iw}(\mathbb{Q}_p,V)}\ar[d]_{\gamma_{-1}\mathcal{L}_V}  \ar[r] & D(\Gamma,{\mathbb{Q}_p}) \ar@{=}[d]_{} \\
      D(\Gamma,{\mathbb{Q}_p})\otimes_{{\mathbb{Q}_p}} D_{cris,{\mathbb{Q}_p}}(V^*(1)) \ar[u]_{\Omega_{V^*(1),1}}  \ar@{}[r]|{\times} &  D(\Gamma,{\mathbb{Q}_p})\otimes_{\mathbb{Q}_p} D_{cris,{\mathbb{Q}_p}}(V)\ar[r] & D(\Gamma,{\mathbb{Q}_p})}
\end{equation*}
In fact this is a variant of Perrin-Riou's reciprocity law comparing $ \Omega_{V,h}$  with $ \Omega_{V^*(1),1-h}$ .

For $L\neq \Qp$ the issue of $L$-analyticity requires that $V^*(1)$ is $L$-analytic for the construction of $ \Omega_{V^*(1),1-h}$, which then implies that $V$ is not $L$-analytic. Instead our regulator map is available and
the purpose of this subsection is to prove an analogue of the above adjointness for arbitrary $L$.

\begin{theorem}[Reciprocity formula/Adjointness of Big exponential and regulator map] \label{thm:adjointness}\phantom{mmmmmm}  Assume that $V^*(1)$ is $L$-analytic with $\mathrm{Fil}^{-1}D_{cris,L}(V^*(1))=D_{cris,L}(V^*(1))$ and \linebreak $D_{cris,L}(V^*(1))^{\varphi_L=\pi_L^{-1}}=D_{cris,L}(V^*(1))^{\varphi_L=1}=0$. Then the following diagram consisting of $D(\Gamma_L,K)$-$\upiota_*$-sesquilinear pairings (in the sense of \eqref{sesquilinear}) commutes:
\begin{equation}\label{f:adjoint}
\xymatrix{
    {D_{\mathrm{rig}}^\dagger(V^*(1))^{\psi_L=\frac{q}{\pi_L}}}\ar@{}[r]|{\times} &  {D(V(\tau^{-1}))^{\psi_L=1}}\ar[d]_{\mathcal{L}_V^0}  \ar[r]^-{\frac{q-1}{q}\{,\}_{Iw}} & D(\Gamma_L,\Cp) \ar@{=}[d] \\
    \cO^{\psi_L=0}\otimes_L D_{cris,L}(V^*(1)) \ar[u]_{\Omega_{V^*(1),1}}   \ar@{}[r]|{\times} &  \cO^{\psi_L=0}\otimes_L D_{cris,L}(V(\tau^{-1})) \ar[r]^-{[,]_{}} & D(\Gamma_L,\Cp).}
\end{equation} Note that the terms on the right hand side of the pairings are all defined over $L!$
\end{theorem}

\begin{proof}
This follows from the abstract reciprocity formula \ref{thm-recproclawKKK} (with $M:=D_{\mathrm{rig}}^\dagger (V(\tau^{-1}))$ as before) by construction.   Indeed, assuming that $z\in \cO^{\psi_L=0}\otimes_L D_{cris,L}(V^*(1))$ and $y\in D(V(\tau^{-1}))^{\psi_L=1} $ we have that  $(1-\frac{\pi_L}{q}\varphi_L)y\in M'\cap (M^{\psi_L=0})$ (see \eqref{f:defregulator}) and $\mathrm{comp}^{-1}((1-\varphi_L)x)\in\check{M}'$ for $x\in  \left( \cO\otimes_{L} D_{cris,L}(V^*(1))\right)^{\psi_L=\frac{q}{\pi_L}}$ such that $z=(1-\varphi_L)x.$   Moreover, $\mathrm{comp}^{-1}((1-\varphi_L)x)\in \check{M}^{\psi_L=0}$
   by  Prop.\  \ref{elementarydiv} as $V^*(1)$ is positive by assumption.
Recall that $ \mathrm{comp}^{-1}(\nabla x)$ is an element in $D_{\mathrm{rig}}^\dagger(V^*(1))^{\psi_L=\frac{q}{\pi_L}}$  again by  Prop.\  \ref{elementarydiv}.  We thus obtain
\begin{align*}
 \frac{q-1}{q} \KKl \mathrm{comp}^{-1}(\nabla x),y \KKr_{Iw}&= \frac{q-1}{q}\KKl \nabla\mathrm{comp}^{-1}((1-\varphi_L)x),(1-\frac{\pi_L}{q}\varphi_L)y \KKr^0_{Iw}\\
  &=[(1-\varphi_L)x,\mathrm{comp}((1-\frac{\pi_L}{q}\varphi_L)y)].
\end{align*}
By definition of the big exponential and regulator map the latter is equivalent to
\begin{align*}
  \KKl \Omega_{V^*(1),1}(z),y \KKr_{Iw}&= [z,\mathcal{L}_V^0(y)].
\end{align*}
\end{proof}

We also could consider the following variant of the big exponential map (under the assumptions of the theorem)
\begin{align*}
\mathbf{\Omega}_{V,h}: D(\Gamma_L,\mathbb{C}_p)\otimes_L D_{cris,L}(V^*(1))\to D_{\mathrm{rig}}^\dagger(V)^{\psi_L=\frac{q}{\pi_L}}
\end{align*}
by extending scalars from $L$ to $\mathbb{C}_p$ and composing the original one with $\Omega^{-h}$\footnote{This means to replace $\nabla$ by $\frac{\nabla}{\Omega}$ in order to achieve twist invariance of the big exponential map, see the remark below.} times \[
 D(\Gamma_L,\mathbb{C}_p)\otimes_L D_{cris,L}(V^*(1))\xrightarrow{ \mathfrak{M}\otimes \id}({\mathcal{O}_K(\mathbf{B})})^{\psi_L=0}\otimes_L D_{cris,L}(V^*(1)).\]

 \begin{corollary}[Reciprocity formula/Adjointness of Big exponential and regulator map] \label{cor:adjointness}
\phantom{mmm} Under the assumptions of the theorem the following diagram  of $D(\Gamma_L,K)$-$\upiota_*$-sesquilinear pairings  commutes:\Footnote{According to \cite[Prop.\ A.2.2]{LVZ15} in the cyclotomic case we have the relation
\[<\mathcal{L}_V(x),\omega>_{Iw,cris}=\gamma_{-1}<x,\Omega_{V^*(1),1,\xi}(\omega)>_{Iw}\]
for $x\in H^1_{Iw}(\qp.V)$ and $\omega\in D(\Gamma,\qp)\otimes D_{cris}(V^*(1))$ (partly in the notation of that article). Note that the factor $\gamma_{-1}$ is already contained in the definition \eqref{f:def[]cris} of the pairing $[,]_{cris}.$ It is derived from the equation (12) in (loc.\ cit.)
\[\Omega_{V,1,\xi}(\mathcal{L}_V(x))=\nabla(x)  \mbox{ modulo torsion}\]
for all $x\in D(\Gamma,\qp)\otimes_{\Lambda(\Gamma)} H^1_{Iw}(\qp,V)$ and the reciprocity formula from Thm.\ 3.3.7 in (loc.\ cit.)
\[<\mathcal{L}_V(x),\mathcal{L}_{V^*(1)}(y)>_{Iw,cris}=\gamma_{-1}\nabla<x,y>_{Iw},\]
(In (loc.\ cit.) there appears a sign on the right hand side, which  I do not understand and which according to David Loeffler is a typo!) {We assume that the HT weights of $V$ are in $[0,r].$ Then $\mathcal{L}_{V^*(1)} $ is defined as follows:  choose $h$ sufficiently large that $V^*(1+h)$ has HT weights $\geq0$ and set \[\mathcal{L}_{V^*(1)}(y)=(\nabla_{-1}\cdots \nabla_{-h})^{-1}Tw_{\chi_{cyc}^{-h}}(\mathcal{L}_{V^*(1+h)}(tw_{\chi_{cyc}^{-h}}(y)))\otimes t^he_{-h}\in FracD(\Gamma,\qp)\otimes D_{cris}.\] }}
\begin{equation}\label{f:adjointcor}
\xymatrix{
    {D_{\mathrm{rig}}^\dagger(V^*(1))^{\psi_L=\frac{q}{\pi_L}}}\ar@{}[r]|{\times} &  {D(V(\tau^{-1}))^{\psi_L=1}}\ar[d]_{\frac{\sigma_{-1}\mathbf{L}_V}{{{\Omega}}} } \ar[r]^-{\frac{q-1}{q}\{,\}_{Iw}} & D(\Gamma_L,\Cp) \ar@{=}[d] \\
    D(\Gamma_L,\mathbb{C}_p)\otimes_L D_{cris,L}(V^*(1)) \ar[u]_{\mathbf{\Omega}_{V^*(1),1}}   \ar@{}[r]|{\times} & D(\Gamma_L,\mathbb{C}_p)\otimes_L D_{cris,L}(V(\tau^{-1})) \ar[r]^-{[,]^0} & D(\Gamma_L,\Cp),}
\end{equation}
where $[-,-]^0=[ \mathfrak{M}\otimes \id (-)  , \sigma_{-1}\mathfrak{M}\otimes \id (-)   ]   ,$ i.e.,
\begin{equation}\label{f:def[]0}
[\lambda\otimes\check{d},\mu \otimes d]^0\cdot \eta(1,Z)\otimes (t_{LT}^{-1}\otimes \eta)=\lambda\iota_*(\mu)\cdot  \eta(1,Z)\otimes[\check{d}, d]_{cris},
\end{equation}  where $D_{cris,L}(V^*(1))\times D_{cris,L}(V(\tau^{-1}))\xrightarrow{[\;,\;]_{cris}} D_{cris,L}(L(\chi_{LT}))$ is the canonical pairing.
\end{corollary}

\begin{remark}\label{rem:BigExpTwist}
By \cite[Cor.\ 3.5.4]{BF} we have $\Omega_{V,h}(x)\otimes \eta^{\otimes j}=\Omega_{V(\chi_{LT}^j), h+j}(\partial_\mathrm{inv}^{-j}x\otimes t_{LT}^{-j}\eta^{\otimes j})$ and $\mathfrak{l}_h\circ \Omega_{V,h}=\Omega_{V,h+1},$ whence we obtain $\mathbf{\Omega}_{V,h}(x)\otimes \eta^{\otimes j}=\mathbf{\Omega}_{V(\chi_{LT}^j),h+j}(Tw_{\chi_{LT}^{-j}}(x)\otimes t_{LT}^{-j}\eta^{\otimes j})$ and $\mathfrak{l}_h\circ \mathbf{\Omega}_{V,h}=\mathbf{\Omega}_{V,h+1}.$
\end{remark}

\subsubsection{Some homological algebra}\label{sec:homol}

 Let  $X\xrightarrow{f}Y$ be a morphism of cochain complexes. Its mapping cone $ \mathrm{cone}(f)$  is defined as $X[1]\bigoplus Y$ with differential
 $d^i_{\mathrm{cone}(f)}:=\begin{pmatrix}
 d_{X[1]}^i & 0 \\
 f[1]^i& d_Y^i \\
\end{pmatrix}$
 (using column notation)  and we define the mapping fibre of $f$ as $\mathrm{Fib}(f):=\mathrm{cone}(f)[-1]$. Here the translation  $X[n]$ of a complex $X$ is given by $X[n]^i:=X^{i+n}$ and $d^i_{X[n]}:=(-1)^nd_X^{i+n}.$ Alternatively,  we may consider $f$ as a double cochain complex concentrated horizontally in degree $0$ and $1$
 and form the total complex (as in \cite[Def.\ 18.3/tag 012Z]{stacks-project}).
  Then the associated total complex coincides with $\mathrm{Fib}(-f)$.

For a complex $(X^\bullet, d_X )$ of topological $L$-vector spaces we define its $L$-dual $((X^*)^\bullet, d_{X^*})$ to be the complex with \[(X^*)^i:=\mathrm{Hom}_{L,cts}(X^{-i},L)\] and \[d_{X^*}(f):=(-1)^{\mathrm{deg}(f)-1} f\circ d_X.\]

More generally, for two complexes $(X^\bullet, d_X )$ and $(Y^\bullet, d_Y )$ of topological $L$-vector spaces we define the complex $\Hom^\bullet_{L,cts}(X^\bullet,Y^\bullet)$ by
\[\Hom^n_{L,cts}(X^\bullet,Y^\bullet)=\prod_{i\in\mathbb{Z}}\Hom_{L,cts}(X^i,Y^{i+n})\]
with differentials $df=d\circ f +(-1)^{\mathrm{deg}(f)-1} f\circ d.$ Note that the canonical isomorphism
\[\Hom^\bullet(X^\bullet,Y^\bullet)[n]\xrightarrow{\cong}\Hom^\bullet(X^\bullet,Y^\bullet[n])\]
does not involve any sign, i.e., it is given by the identity map in all degrees.

Also we recall that the tensor product of two complexes $X^\bullet$ and $Y^\bullet$ is given by
\[(X^\bullet \otimes_L Y^\bullet)^i:=\bigoplus_{n} X^n\otimes_L Y^{i-n}\] and
\[d(x\otimes y)=dx \otimes y + (-1)^{\mathrm{deg}(x)}x\otimes d y.\]
The adjunction morphism on the level of complexes
\[\mathrm{adj}:\Hom^\bullet_{L,cts}(X^\bullet\otimes_LY^\bullet,Z^\bullet)\to \Hom^\bullet_{L,cts}(Y^\bullet,\Hom^\bullet_{L, cts}(X^\bullet,Z^\bullet))\] sends $u$ to $(y\mapsto(x\mapsto (-1)^{\deg(x)\deg(y)}u(x\otimes y))).$ It is well-defined and continuous with respect to the projective tensor product topology and the strong topology for the Homs. Furthermore, by definition we have the  following
  commutative diagram
\begin{equation}\label{f:adjointmorphism}
\xymatrix{
   X^\bullet \otimes_L Y^\bullet \ar[d]_{\id\otimes \mathrm{adj}(u)} \ar[r]^-{u} & L[-2] \ar@{=}[d] \\
  X^\bullet\otimes_L \Hom^\bullet_{L,cts}(X^\bullet,L[-2]) \ar[r]^-{ev_2} & L[-2]  },
\end{equation}
  where $ev_2$ sends $(x,f)$ to $(-1)^{\deg(x)\deg(f)}f(x).$

%

\begin{lemma}\label{lem:strict}
Let $(\cC^\bullet,d^\bullet)$ be a complex in the category of  locally convex topological $L$-vector spaces.
\begin{enumerate}
  \item   If $\cC$ consists of Fr\'{e}chet spaces and  $h^i(\cC^\bullet)$ is finite-dimensional over $L$, then  $d^{i-1}$ is strict and has closed image.
  \item  If  $d^{i}$ is strict, then $h^{-i}(\cC^*)\cong h^{i}(\cC)^*.$
\end{enumerate}
\end{lemma}

\begin{proof}
(i) Apply the argument from \cite[§ IX, Lem.\ 3.4]{BW} and use the open mapping theorem \cite[Prop.\ 8.8]{NFA}. (ii) If
\[\xymatrix@C=0.5cm{
    A \ar[rr]^{\alpha} && B \ar[rr]^{\beta} && C  }\]
forms part of the complex with $B$ in degree $i$, one immediately obtains a map
\[\ker(\alpha^*)/\mathrm{im}(\beta^*)\to \left(\ker(\beta)/\mathrm{im}(\alpha) \right)^*,\] where  $\ker(\beta) $ carries the subspace topology and $\ker(\beta)/\mathrm{im}(\alpha)  $ the quotient topology. Now use  the Hahn-Banach theorem \cite[Cor.\ 9.4]{NFA} for the strict maps $B/\ker(\beta)\hookrightarrow C$ (induced from $\beta$) and $\ker(\beta)\hookrightarrow B$ in order to show that this map is   an isomorphism.
\end{proof}

\begin{definition}
  A locally convex topological vector space is called an LF-space, if it is the direct limit of a countable
family of Fr\'{e}chet spaces, the limit being formed in the category of locally convex vector spaces.
\end{definition}

\begin{remark}\phantomsection\label{rem:nonstrict}
\begin{enumerate}
\item If $V\xrightarrow{\alpha} W$ is a continuous linear map of Hausdorff $LF$-spaces with finite dimensional cokernel, then $\alpha$ is strict and has closed image by the same argument used in (i) of the previous lemma. However, since a closed subspace of an LF-space need not be an LF-space, we cannot achieve the same conclusion for complexes by this argument as   $ \ker(d^{i})$ may fail to be an LF-space, whence one cannot apply the open mapping theorem, in general.   But consider the following special situation. Assume that the complex $\cC^\bullet$ consists of LF-spaces and  $h^i(\cC^\bullet)$ is finite-dimensional. If moreover $\cC^{i+1}=0,$ i.e., $\cC^{i }=\ker(d^{i}),$ then $d^{i-1}$ is strict and $h^{1-i}(\cC^*)\cong h^{i-1}(\cC)^*.$
\item
If $d^i$ is not strict, the above proof still shows that we obtain a surjection $h^{-i}(\cC^*) \twoheadrightarrow h^{i}(\cC)^*.$
\end{enumerate}
\end{remark}

However, for a special class of LF-spaces and under certain conditions we can say more about how forming duals and cohomology interacts.

\begin{lemma}
\label{lem:LFdual} Let $(\cC^\bullet,d^\bullet)=\varinjlim_r  (\cC^\bullet_r,d^\bullet_r)$ be a complex in the category of  locally convex topological $L$-vector spaces arising as
regular inductive limit of complexes of Fr\'{e}chet spaces, i.e., in each degree $i$  the transition maps in the countable sequence $(\cC^i_r)_r$ are injective and for each bounded subset
$B\subseteq \cC^i$ there exists an $r \geq 1$ such that $B$ is contained in $\cC^i_r$ and is
bounded as a subset of the Fr\'{e}chet space $\cC^i_r$. Then,
\begin{enumerate}
  \item we have topological isomorphisms $(\cC^\bullet)^*\cong\varprojlim_r (\cC^\bullet_r)^*,$
  \item  if, in addition,  $\varprojlim^1_{r\geq 0} h^{i}((\cC^\bullet_r)^*)=0$ for all $i$, we have a long exact sequence
  \begin{equation*}\label{f:lim1dual}
     \xymatrix@C=0.5cm{
       \ldots \ar[r] & h^i((\cC^\bullet)^*) \ar[r]^{} & \varprojlim_{r\geq 0} h^{i}((\cC^\bullet_r)^*) \ar[r]^{ } & h^{i-1}(\varprojlim^1_{r\geq 0}(\cC^\bullet_r)^* ) \ar[r] & h^{i+1}((\cC^\bullet)^*)\ar[r] & \ldots ,}
   \end{equation*}
  \item if, in addition to (ii), the differentials $d_r^\bullet$ are strict, e.g., if all $h^i(\cC^\bullet_r)$ have finite dimension over $L$, and $\varprojlim^1_{r\geq 0}(\cC^\bullet_r)^*=0,$ we have isomorphisms \[h^i((\cC^\bullet)^*)\cong\varprojlim_{r\geq 0} h^{-i}(\cC^\bullet_r)^*.\]

\end{enumerate}
\end{lemma}

\begin{proof}
(i) is \cite[Thm:\ 11.1.13]{PGS} while (ii), (iii) follows from (i) and \cite[Ch.\ 3, Prop.\ 1]{Lu} applied to the inverse system $ ((\cC^\bullet_r)^*)_r$ combined with Lemma \ref{lem:strict}. .
\end{proof}

\subsubsection{Koszul complexes}\label{sec:KoszulCom}

In this paragraph we restrict to the situation $U\cong \mathbb{Z}_p^d$ and fix topological generators $\gamma_1,\ldots \gamma_d$  of $U$ and we set $\Lambda:= \Lambda(U).$ Furthermore, let $M$ be any complete linearly topologized  $o_L$-module    
with a continuous $U$-action. Then by \cite[Thm.\ II.2.2.6]{La} this actions extends to continuous $\Lambda$-action and one has $\mathrm{Hom}_{\Lambda, cts}(\Lambda,M)=\mathrm{Hom}_{\Lambda}(\Lambda,M).$

Consider the (homological) complexes $K_\bullet(\gamma_i):=[\Lambda\xrightarrow{\gamma_i-1}\Lambda]$ concentrated in degrees $1$ and $0$ and define
\begin{align*}
K_\bullet:=& K_\bullet^U:= K_\bullet(\gamma):={\bigotimes_{{\substack{\phantom{i=1}{\Lambda}\\i=1}}}^d}K_\bullet(\gamma_i),\\
K^\bullet(M):=&K^\bullet_U(M):=\Hom_\Lambda^{\bullet} (K_\bullet,M)\cong \Hom_\Lambda^{\bullet }(K_\bullet,\Lambda)\otimes_\Lambda M=K^\bullet(\Lambda)\otimes_\Lambda M,\\
K_\bullet(M):=&K_\bullet\otimes_\Lambda M \mbox{ (homological complex)},\\
K_\bullet(M)^\bullet:=&(K_\bullet\otimes_\Lambda M)^\bullet \mbox{ (the associated cohomological complex)}.
\end{align*}\Footnote{ Hier lassen wir $\Lambda$ per Funktorialit\"{a}t auf $\Hom_\Lambda^{\bullet}(K_\bullet,M)$ operieren, d.h. $(\lambda f)(x)=f(x\lambda)$. Dann wird aus $Hom^{\bullet}_\Lambda(K_\bullet(\gamma_i),\Lambda)$ der Komplex $ [\Lambda\xrightarrow{-(\gamma_i-1)}\Lambda]$ in den Geraden $0$ and $1$.}
  If we want to indicate the dependence on $\mathbf{\bf\gamma}=(\gamma_1,\ldots \gamma_d)$ we also write $K^\bullet(\mathbf{\gamma},M)$ instead of $K^\bullet(M)$  and similarly for other notation; moreover, we shall   use the notation $\mathbf{\bf\gamma}^{-1}=(\gamma_1^{-1},\ldots,\gamma_d^{-1})$ and $\mathbf{\bf\gamma}^{p^n}=(\gamma_1^{p^n},\ldots \gamma_d^{p^n})$ .
Note that in each degree these complexes consists of a direct sum of finitely many copies of $M$ and will be equipped with the corresponding direct product topology.

The complex $K_\bullet$ will be identified with the exterior algebra complex $\bigwedge^\bullet_\Lambda \Lambda^d$ of the free $\Lambda$-module with basis $e_1,\ldots,e_d$, for which the differentials $d_q: \bigwedge^q_\Lambda \Lambda^d  \to{\bigwedge}^{q-1}_\Lambda \Lambda^d $ with respect to  the standard basis $e_{i_1,\ldots, i_q}=e_{i_1}\wedge \cdots \wedge e_{i_q}$, $1\leq i_1<\cdots <i_q\leq d,$ is given by the formula
\[d_q(a_{i_1,\ldots,i_q})=\sum_{k=1}^q (-1)^{k+1} (\gamma_{i_k}-1)a_{i_1,\ldots,\widehat{i_k},\ldots,i_q}.\]

Then the well-known selfduality (compare \cite[Prop.\ 17.15]{Ei} although the claim there is not precisely the same) of the Koszul complex, i.e., the isomorphism of complexes
 \begin{equation}\label{f:selfdual}
   K_\bullet(\Lambda)^\bullet \cong K^\bullet(\Lambda)[d]
\end{equation}
can be explicitly described in degree $-q$ as follows (by identifying ${\bigwedge}^d_\Lambda \Lambda^d=\Lambda e_1\wedge \cdots \wedge e_d=\Lambda$):
\begin{align*}
{\bigwedge}^q_\Lambda \Lambda^d&\xrightarrow{\alpha_{-q}}\Hom_\Lambda^{}({\bigwedge}^{d-q}_\Lambda \Lambda^d ,\Lambda)\\
e_{i_1,\ldots, i_q}&\mapsto \phantom{m}\mathrm{sign}(I,J)e^*_{j_1,\ldots, j_{d-q}},
\end{align*}
where $e^*_1,\ldots e^*_d$ denotes the dual basis of $e_1,\ldots,e_d,$ the elements $e^*_{j_1,\ldots, j_{d-q}}=e^*_{j_1}\wedge \cdots \wedge e^*_{j_{d-q}}$, $1\leq j_1<\cdots <j_{d-q}\leq d,$ form a (dual) basis of   $\Hom_\Lambda({\bigwedge}^{d-q}_\Lambda \Lambda^d ,\Lambda)$, the indices $J=(j_k)_k$ are complementary to $I=(i_n)_n$ in the following sense $\{i_1,\ldots,i_q\}\cup\{j_1,\ldots,j_{d-q}\}=\{1,\ldots,d\}$ and $\mathrm{sign}(I,J)$ denotes the sign of the permutation $[i_1,\ldots,i_q,j_1,\ldots,j_{d-q}].$ Indeed, the verification that the induced diagram involving the differentials from cohomological degree $-q$ to $-q+1$
\[\xymatrix{
  {\bigwedge}^q_\Lambda \Lambda^d \ar[d]_{d_q} \ar[r]^(0.35){\alpha_{-q}} & \Hom_\Lambda^{}({\bigwedge}^{d-q}_\Lambda \Lambda^d ,\Lambda) \ar[d]^{(-1)^d(-1)^{d-q-1}d^*_{d-q+1}} \\
  {\bigwedge}^{q-1}_\Lambda \Lambda^d \ar[r]^(0.35){\alpha_{-q+1}} & \Hom_\Lambda^{}({\bigwedge}^{d-q+1}_\Lambda \Lambda^d ,\Lambda)   }
\]\footnote{The signs $(-1)^d$ and $(-1)^{d-q-1}$ result from the shift by $d$ and the sign rule for complex-homomorphisms, respectively. }
commutes, relies on the observation that
\[\mathrm{sign}(I,J)\mathrm{sign}(I_{\widehat{k}},J_k)^{-1}=(-1)^{q-k+l-1},\]
where $I_{\widehat{k}}:=(i_1,\ldots,\widehat{i_k},\ldots,i_q)$ denotes the sequence which results from $I$ by omitting $i_k$ while $J_k=(j_1,\ldots,j_{l-1},i_k,j_l,\ldots i_{d-q})$ denotes the sequence which arises from $J$ by inserting $i_k$ at position $l$ with regard to the strict increasing ordering: The permutations  $[i_1,\ldots,i_q,j_1,\ldots,j_{d-q}]$ and $[i_1,\ldots,\widehat{i_k},\ldots,i_q,j_1,\ldots,j_{l-1},i_k,j_l,\ldots,j_{d-q}]$ differ visibly by $q-k+l-1$ transpositions.

Now we assume that $M$ is any complete locally convex $L$-vector space  with continuous $U$-action such that its strong dual is again complete with continuous $U$-action. Then we obtain isomorphisms of complexes
\begin{align}\label{f:selfualM}
K^\bullet(\gamma,M)^{*}&=\Hom_{L,cts}^{\bullet}(  \Hom_\Lambda^{\bullet}(K_\bullet(\gamma),\Lambda)\otimes_\Lambda M  ,L)\notag\\
&\cong \Hom_{\Lambda}^{\bullet}(\Hom_\Lambda^{\bullet}(K_\bullet(\gamma^{-1}),\Lambda),\Hom_{L,cts}(M,L))\notag\\
&\cong \Hom_{\Lambda}^{\bullet}(\Hom_\Lambda^{\bullet}(K_\bullet(\gamma^{-1}),\Lambda),\Lambda)\otimes_\Lambda \Hom_{L,cts}(M,L)\notag\\
&\cong K_\bullet(\gamma^{-1},\Lambda)^\bullet \otimes_\Lambda \Hom_{L,cts}(M,L)\\
&\cong K^\bullet(\gamma^{-1},\Lambda)[d]\otimes_\Lambda M^*\notag\\
&\cong K^\bullet(\gamma^{-1},M^*)[d],\notag
\end{align}
where in the second line we use the adjunction morphism;  the isomorphism in the  fourth line being the biduality morphism (according to \cite[(1.2.8)]{Ne})
\begin{align*}
 K_\bullet(\Lambda)^\bullet&\xrightarrow{\cong} \Hom_{\Lambda}^{\bullet}(\Hom_\Lambda^{\bullet}(K_\bullet,\Lambda),\Lambda)\\
 x&\mapsto (-1)^{i}x^{**}
\end{align*}\Footnote{$(-1)^{-i^2}=(-1)^{ik}$ for $k=-i$ is the sign recipe from Nekovar.}
with the usual biduality of modules
\begin{align*}
   K_\bullet(\Lambda)^i & \xrightarrow{\cong} \Hom_{\Lambda}(\Hom_\Lambda(K_{-i},\Lambda),\Lambda) \\
  x &\mapsto (x^{**}: f\mapsto f(x))
\end{align*}
involves a sign, while the isomorphism in the third last line stems from \eqref{f:selfdual} together with Lemma \ref{lem-pair} (i). Note that the isomorphism in the second last line does not involve any further signs by \cite[(1.2.15)]{Ne}.

We finish this subsection by introducing restriction and corestriction maps concerning the change of group for Koszul complexes. To this end let $U_1\subseteq U$ be the open subgroup generated by $\gamma_1^{p^n},\ldots \gamma_d^{p^n}$. Then $\Hom_\Lambda^{\bullet}(-,M)$ applied to the tensor product of the diagrams
\[\xymatrix{
  \Lambda(U) \ar@{=}[d]_{} \ar[r]^{\gamma_i^{p^n}-1} &  \Lambda(U)  \\
   \Lambda(U)\ar[r]^{\gamma_i -1} &  \Lambda(U)  \ar[u]_{\sum_{k=0}^{p^n-1}\gamma_i^k} }\]
gives a map $cor^{U_1}_U:K_{U_1}^\bullet(\gamma^{p^n})(M)\to K_U^\bullet(\gamma)(M)$ which we call corestriction map  and which is compatible under \eqref{f:KoszulCtsCoh} below with the corestriction map on cocylces (for appropriate choices of representatives in the definition of the latter). Using the diagram
\[\xymatrix{
  \Lambda(U) \ar[d]_{\sum_{k=0}^{p^n-1}\gamma_i^k} \ar[r]^{\gamma_i^{p^n}-1} &  \Lambda(U)\ar@{=}[d]_{}   \\
   \Lambda(U)\ar[r]^{\gamma_i -1} &  \Lambda(U)   }\]
instead, one obtains the restriction map $res_{U_1}^U:K_U^\bullet(\gamma)(M)\to K_{U_1}^\bullet(\gamma^{p^n})(M),$ again compatible under \eqref{f:KoszulCtsCoh} with the restriction map on cocycles.

\subsubsection{Continuous and analytic cohomology}\label{sec:cts-analytic}

For any profinite group $G$ and topological abelian group $M$ with continuous $G$-action we write $\cC^\bullet := \cC^\bullet(G,M)$  for  the continuous (inhomogeneous) cochain complex  of $G$ with coefficients in $M$ and $H^*(G,M) := h^*(\cC^\bullet(G,M))$ for continuous group cohomology. Note that $\cC^0(G,M)=M.$

If $G$ is moreover a $L$-analytic group and $M=\varinjlim_s\varprojlim_r M^{[r,s]}$ with Banach spaces $M^{[r,s]} $  a LF space  with a pro-$L$-analytic action of $G,$   i.e., a locally analytic action on each  $M^{[r,s]}$, which means that for all $m \in M^{[r,s]}$ there exist an open $L$-analytic subgroup $\Gamma_n\subseteq\Gamma$ in the notation of subsection \ref{sec:groupRobba} such that the orbit map of $m$ restricted to $\Gamma_n$ is a power series of the form $g(m)=\sum_{k\geq 0} \ell(g)^k m_k$ for a sequence $m_k$ of elements in $M^{[r,s]}$ with $\pi_L^{nk}m_k$ converging to zero. Following \cite[\S 5]{Co2} we write $\cC^\bullet_{an} := \cC^\bullet_{an}(G,M)$  for  the locally $L$-analytic cochain complex  of $G$ with coefficients in $M$ and $H^*_{an}(G,M) := h^*(\cC^\bullet_{an}(G,M))$ for  locally $L$-analytic group cohomology.  More precisely, if $\mathrm{Maps}_{loc L-an}(G,M^{[r,s]})$ denotes the space of locally $L$-analytic maps from $G$ to  $M^{[r,s]}$, then \[C_{an}^n(G,M)=\varinjlim_s\varprojlim_r\mathrm{Maps}_{loc L-an}(G,M^{[r,s]})\] is the space of locally $L$-analytic functions (locally with values in  $\varprojlim_r M^{[r,s]} $ for some $s$ and such that the composite with the projection onto $M^{[r,s]}$ is locally $L$-analytic for all $r$). Note that again $\cC^0_{an}(G,M)=M$ and that there are canonical homomorphisms
\begin{align}
\label{f:an-cts} &\cC^\bullet_{an}(G,M)\hookrightarrow \cC^\bullet(G,M),\\
 &H^\bullet_{an}(G,M)\to H^\bullet(G,M).
\end{align}

Let $f$ be any continuous endomorphism of $M$ which commutes with the $G$-action. We define
\begin{equation}\label{f:f-coh}
  H^0(f,M) := M^{f=1} \qquad\text{and} \qquad H^1(f,M) := M_{f=1}
\end{equation}
as the kernel and cokernel of the map $M \xrightarrow{f-1} M$, respectively.

The endomorphism $f$ induces an operator on $\cC^\bullet$ or $\cC^\bullet_{an}$ 
and  we denote by $\mathcal{T} := \mathcal{T}_{f,G}(M)$ and $\mathcal{T}^{an} := \mathcal{T}_{f,G}^{an}(M)$ the mapping fibre of  $\cC^\bullet(G,f)$ and $\cC^\bullet_{an}(G,f)$, respectively.

 Again there are canonical homomorphisms
\begin{align}
\label{f:an-ctsT} \mathcal{T}_{f,G}^{an}(M)&\hookrightarrow \mathcal{T}_{f,G}(M),\\\label{f:an-ctsh}
 h^\bullet(\mathcal{T}_{f,G}^{an}(M))&\to h^\bullet(\mathcal{T}_{f,G}(M)).
\end{align}

For $?$ either empty or $an$, one of the corresponding double complex spectral sequences is
\begin{equation}\label{f:spectralsequenceI}
  _{I}E_2^{i,j}=H^i(f,H^j_?(G,M)) \Longrightarrow h^{i+j}(\mathcal{T}^?)
\end{equation}
\footnote{Naively, one would expect that the second corresponding double complex spectral sequences looks like \begin{equation*}
  _{II}E_2^{i,j}=H^i_?(G,H^j(f,M)) \Longrightarrow h^{i+j}(\mathcal{T}^?) \ .
\end{equation*} But this would require to first of all give sense to the required structure of $ H^j(f,M)$ as topological/analytic $G$-module! In low degrees this can be achieved and we obtain an exact sequence
\begin{equation*}
  0 \longrightarrow H^1_?(G,M^{f=1}) \longrightarrow h^1(\mathcal{T}^?) \longrightarrow (M_{f=1})^G \xrightarrow{\ \delta\ } H^2_?(G,M^{f=1}) \ .
\end{equation*} See \cite{Th}. If $M^{f=1}$ is again a LF-space with pro-$L$-analytic $G$-operation, one might be able to interpret the second spectral sequence in low degrees.
}
It degenerates into the short exact sequences
\begin{equation*}
  0 \longrightarrow H^{i-1}_{?}(G,M)_{f=1} \longrightarrow h^i(\mathcal{T}_{f,G}^?(M) ) \longrightarrow H^i_?(G,M)^{f=1} \longrightarrow 0.
\end{equation*}

In (loc. cit.) as well as in \cite{BF} analytic cohomology is also defined for the semigroups $\Gamma_L\times \Phi$ and $\Gamma_L\times \Psi$ with $\Phi=\{\varphi_L^n|n\geq 0\}$ and $\Psi=\{(\frac{\pi}{q}\psi_L)^n|n\geq 0\},$ if $M$ denotes an  $L$-{\it analytic}   $(\varphi_L,\Gamma_L)$-module over the Robba ring $\cR.$

\begin{remark}
 Any $L$-{\it analytic}   $(\varphi_L,\Gamma_L)$-module $M$ over the Robba ring $\cR$ is a pro-L-analytic $\Gamma_L$-module by the discussion at the end of the proof of  \cite[Prop.\ 2.25]{BSX}, whence it is also an $L$-{\it analytic} $\Gamma_L\times \Phi$- and $\Gamma_L\times \Psi$-module as $\Phi$ and $\Psi$ possess the discrete structure as $L$-analytic manifolds.
\end{remark}

\begin{proposition}\label{Herr}
We have canonical isomorphisms
\begin{equation*}
 h^i(\mathcal{T}_{\varphi_L,\Gamma_L}^{an}(M))\cong H^i_{an}(\Gamma_L\times \Phi,M)\cong H^i_{an}(\Gamma_L\times \Psi,M)\cong  h^i(\mathcal{T}_{\frac{\pi}{q}\psi_L,\Gamma_L}^{an}(M)).
\end{equation*}
and an   exact sequence
\begin{equation}\footnotesize
 \label{f:SSandim1}
  0 \longrightarrow H^1_{an}(\Gamma_L,M^{\psi_L=\frac{q}{\pi}}) \longrightarrow h^i(\mathcal{T}^{an}_{\frac{\pi}{q}\psi_L,\Gamma_L}(M)) \longrightarrow (M_{\psi_L=\frac{q}{\pi}})^{\Gamma_L} \longrightarrow   H^2_{an}(\Gamma_L,M^{\psi_L=\frac{q}{\pi}})\longrightarrow h^2(\mathcal{T}^{an}_{\frac{\pi}{q}\psi_L,\Gamma_L}(M)) \ .
\end{equation}
\Footnote{{The second last group vanishes, if $\Gamma_L$ acts analytically on each $\varprojlim_s M^{[r,s]}.$ See also \cite[Remarks 5.1.6/7]{Th} for other criteria.}}
\end{proposition}

\begin{proof} The isomorphism in the middle is \cite[Cor.\ 2.2.3]{BF}.
For the two outer isomorphism we refer the reader to \cite[3.7.6]{Th}. The   exact sequence is the extension \cite[Thm.\ 5.1.5]{Th} of
\cite[Thm.\ 2.2.4]{BF}.
\end{proof}

Note that, for $U\subseteq U'$, the restriction and corestriction homomorphisms $ \cC^\bullet(U',M)\xrightarrow{\mathrm{res}} \cC^\bullet(U,M) $ and $ \cC^\bullet(U,M)\xrightarrow{\mathrm{cor}} \cC^\bullet(U',M) $ induce maps on $\mathcal{T}_{f,U'}(M)\xrightarrow{\mathrm{res}}\mathcal{T}_{f,U}(M)$ and \linebreak $\mathcal{T}_{f,U}(M)\xrightarrow{\mathrm{cor}}\mathcal{T}_{f,U'}(M)$, respectively.


We write $\mathrm{Ext}^1_{\mathfrak{C}}(A,B) $ for isomorphism classes of extensions of $B$ by $A$ in any abelian  category $\mathfrak{C}.$ Furthermore, we denote by $\mathfrak{M}_U(R) $ ($\mathfrak{M}^{\acute{e}t}_U(R)$, $\mathfrak{M}^\dagger_U(R) $ )  the category of all (\'{e}tale, overconvergent) $(\varphi_L,U)$-modules over $R$, respectively, and by  $\mathrm{Rep}_L^\dagger(G_{L_\infty^U})$ the category of overconvergent representations of $G_{L_\infty^U}$   consisting of those representations  $V$ of $ G_{L_\infty^U} $ such that  $\dim_{{\bf B}^\dagger_L} D^\dagger(V) = \dim_LV$ with  $D^\dagger(V):=(\mathbf{B}^\dagger\otimes_L V)^{H_L}.$

\begin{theorem}\label{thm:Herr}
Let $V$ be in $\Rep_{L}(G_L)$ and $U \subseteq \Gamma_L$ be any open subgroup.
\begin{enumerate}
\item For $  D(V)$ the corresponding $(\varphi_L,\Gamma_L)$-module over $\mathbf{B}_L$ we have canonical isomorphisms
\begin{equation}\label{f:cohLphiB}
  h^* = h^*_{U,V}: H^*(L_\infty^U,V) \xrightarrow{\ \cong\ }   h^*(\mathcal{T}_{\varphi_{L},U}(D(V)) )
\end{equation}
which are functorial in $V$ and compatible with restriction and corestriction.
\item If $V$ is in addition overconvergent  there are isomorphisms  
\begin{align}\label{f:cohLphiOVer}
h^0(\mathcal{T}_{\varphi_{L},U}(D_{rig}^\dagger(V)) )&\cong V^{G_{L_\infty^U}}, \\
\label{f:cohLphi1Ver}h^1(\mathcal{T}_{\varphi_{L},U}(D_{rig}^\dagger(V)) ) & \cong
    H^1_\dagger(L_\infty^U,V),
\end{align}
which are functorial in $V$ and compatible with restriction and corestriction and where by definition $H^1_{\dagger}(L_\infty^U,V)\subseteq H^1(L_\infty^U,V) $ classifies  the overconvergent extensions of $L$ by $V$. In particular, these $L$-vector spaces have finite dimension.
\item If $V$ is in addition $L$-analytic, then we have
\begin{equation}
H^1_{an}(L_\infty^U,V) \xrightarrow{\ \cong\ }   h^1(\mathcal{T}_{\varphi_{L},U}^{an}(D_{rig}^\dagger(V)) )
\end{equation}
where  by definition \footnote{Note that the absolute Galois group of $L_\infty^U $ is not $L$-analytic, so this group has not been defined earlier.} $H^1_{an}(L_\infty^U,V)\subseteq H^1_\dagger(L_\infty^U,V) \subseteq H^1(L_\infty^U,V)$ classifies  the $L$-analytic extensions of $L$ by $V$.
\end{enumerate}
\end{theorem}

\begin{proof}
(i) is \cite[Thm.\ 5.1.11.]{Ku} or \cite[Thm.\ 5.1.11.]{KV}.
The  statement (iii) is \cite[Prop.\ 2.2.1]{BF} combined with Prop.\  \ref{Herr} while (ii) follows from \cite{FX} (the reference literally only covers the case $U=\Gamma_L$, but the same  arguments allow to extend the result to general $U$) as follows: Firstly, by Lemma \ref{lem:ext} below one has an isomorphism $h^1(\mathcal{T}_{\varphi_{L},U}(D_{rig}^\dagger(V)) )\cong \mathrm{Ext}^1_{\mathfrak{M}_U(\mathcal{R}_L)}(\mathcal{R}_L,D_{rig}^\dagger(V)).$ Then use  the HN-filtration \`{a} la Kedlaya to see that any extension of \'{e}tale $(\varphi_L,U)$-modules is \'{e}tale again, whence \[\mathrm{Ext}^1_{\mathfrak{M}_U(\mathcal{R}_L)}(\mathcal{R}_L,D_{rig}^\dagger(V))=\mathrm{Ext}^1_{\mathfrak{M}^{\acute{e}t}_U(\mathcal{R}_L)}(\mathcal{R}_L,D_{rig}^\dagger(V))\] and the latter group equals \[\mathrm{Ext}^1_{\mathfrak{M}^\dagger_U(\mathcal{R}_L)}(\mathcal{R}_L,D_{rig}^\dagger(V))\cong\mathrm{Ext}^1_{\mathrm{Rep}_L^\dagger(G_{L_\infty^U})}(L,V)= H^1_\dagger(L_\infty^U,V)\] by  Prop.\ 1.5 and 1.6 in (loc.\ cit.). For the claim in degree $0$ one has to show that the inclusion $D^\dagger(V)\subseteq D_{rig}^\dagger(V)$ induces an isomorphism on $\varphi_L$-invariants, which follows from  \cite[Hypothesis 1.4.1, Prop.\ 1.2.6]{Ked08}.\footnote{Since the {\it strong hypothesis}  holds by \cite[Hypothesis 1.4.1, Prop.\ 1.2.6]{Ked08}   we also obtain an isomorphism on the $\varphi_L$-coinvariants $H^1(\varphi_L,-)$. Then the second spectral sequence above or a similar argument via the Koszul complexes as in Prop.\  \ref{cor:herr} implies that the canonical base change map induces an isomorphism  $h^*(\mathcal{T}_{\varphi_{L},U}(D^\dagger(V)) )  \cong h^*(\mathcal{T}_{\varphi_{L},U}(D_{rig}^\dagger(V)) ) . $ Cp.\ \cite[proof of Prop.\ 2.7]{Li}.}
\end{proof}

\begin{lemma}\label{lem:ext}    Let $M$ be   in $\mathfrak{M}_U(\mathcal{R})$. Then we have a canonical isomorphism
\[h^1(\mathcal{T}_{\varphi_{L},U}(M) )\cong \mathrm{Ext}^1_{\mathfrak{M}_U(\mathcal{R}_L)}(\mathcal{R}_L,M).\]
\end{lemma}

\begin{proof}
Starting with a class $z=[( c_1,-c_0)]$ in $h^1(\mathcal{T}_{\varphi_L,U}(M)  )$ with $c_1\in C^1(M)$ and $c_0\in C^0(M)=M$ (i.e., we work with {\it inhomogeneous } continuous cocycle)  satisfying the cocycle property
\begin{equation}\label{f:cocycle}
   c_1(\sigma\tau)=\sigma c_1(\tau)+c_1(\sigma) \mbox{ for all } \sigma,\tau\in U, \;\; \mbox{ and } \;\;
   (\varphi_L-1)c_1(\tau)=(\tau-1)c_0 \mbox{ for all } \tau\in U,
\end{equation}
we define an extension of $(\varphi_L,U)$-modules
\[\xymatrix@C=0.5cm{
  0 \ar[r] & M \ar[rr]^{ } && E_c \ar[rr]^{ } && {\cR_L} \ar[r] & 0 }\]
  with $E_c:=M\times {\cR_L}$ as ${\cR_L}$-module, $g(m,r):=(gm+gr\cdot c_1(g),gr)$ for $g\in U$ and $\varphi_{E_c} ((m,r)):=(\varphi_M(m)+\varphi_L(r) c_0,\varphi_L(r))$; note that  this defines a (continuous) group-action by the first identity in \eqref{f:cocycle}, while the $U$- and $\varphi_L$-action commute by the second identity in \eqref{f:cocycle}. If we change the representatives $(c_1,-c_0)$ by the coboundary induced by $m_0\in M$, then sending $(0,1)$ to $(-m_0,1)$ induces an isomorphism of extensions from the first to the second one, whence our map is well-defined.

  Conversely, if $E$ is any such extension, choose a lift $e\in E$ of $1\in\cR_L$ and define
  \[c_1(\tau):=(\tau-1)e \in M, \;\;\; c_0:=(\varphi_E-1)e,\]
  which evidently satisfy the cocycle conditions \eqref{f:cocycle}. Choosing another lift $\tilde{e}$ leads to a cocylce which differs from the previous one by the coboundary induced by $\tilde{e}-e\in M,$ whence the inverse map is well-defined.

  One easily verifies that  these maps are mutually inverse to each other.
\end{proof}

\begin{question}
Can one show that $h^2(\mathcal{T}_{\varphi_{L},U}(D_{rig}^\dagger(V)) ) $ is finite-dimensional (and related to $H^2(L_\infty^U,V)$) and that the groups $h^i(\mathcal{T}_{\varphi_{L},U}(D_{rig}^\dagger(V)) )$ vanish for $i\geq 3?$
\end{question}

\begin{remark}\label{rem:FX}
By \cite[Thm.\ 0.2, Rem.\ 5.21]{FX} it follows that the inclusions
\[H^1_{an}(L_\infty^U,V)\subseteq H^1_\dagger(L_\infty^U,V) \subseteq H^1(L_\infty^U,V)\] are in general strict. More precisely, the codimension for the left one equals  $([L_\infty^U:\mathbb{Q}_p]-1)\dim_L V^{G_{L_\infty^U}}.$
\end{remark}

Let us recall Tate's local duality in this context.

\begin{proposition}[Local Tate duality]\label{Tate-local}
Let $V$ be an object in $\Rep_{L}(G_L)$, and $K$ any finite extension of $L$. Then the cup product and the local invariant map induce perfect pairings of finite dimensional $L$-vector spaces
\begin{equation*}
  H^i(K,V) \times  H^{2-i}(K,\Hom_{\mathbb{Q}_p}(V,\mathbb{Q}_p(1)) ) \longrightarrow H^2(K, \mathbb{Q}_p(1)) = \mathbb{Q}_p
  \end{equation*}
and
\begin{equation*}
  H^i(K,V) \times H^{2-i}(K,\Hom_{L}(V,L(1)) ) \longrightarrow H^2(K, L(1)) = L
\end{equation*}
where $-(1)$ denotes the Galois twist by the cyclotomic character. In other words, there are canonical isomorphisms
\begin{equation*}
  H^i(K,V)\cong H^{2-i}(K,V^*(1))^* \ .
\end{equation*}
\end{proposition}

\begin{proof} This is well known. For lack of a reference (with proof) we sketch the second claim (the first being proved similarly).
Choose a Galois stable $o_L$-lattice $T\subseteq V$ and denote by ${_{\pi_L^n}}A$ the kernel of multiplication by $ {{\pi_L^n}}$ on any $o_L$-module $A.$ Observe that we have short exact sequences
\[\xymatrix@C=0.5cm{
  0 \ar[r] & H^i(K,T)/\pi_L^n \ar[rr]^{ } && H^i(K,T/\pi_L^nT) \ar[rr]^{ } && {_{\pi_L^n}} H^{i+1}(K,T) \ar[r] & 0 }\] for $i\geq 0$ and similarly for $T$ replaced by $T^*(1)=\Hom_{o_L}(T,o_L(1))$. By \cite[Prop.\ 5.7]{SV15} (remember the normalisation given there!) the cup product induces isomorphism
  \[H^i(K,T/\pi_L^nT)\cong \Hom_{o_L}(H^{2-i}(K,T^*(1)/\pi_L^nT^*(1)),o_L/\pi_L^n)\] such that we obtain altogether canonical maps
  \[H^i(K,T)/\pi_L^n\to \Hom_{o_L}(H^{2-i}(K,T^*(1))/\pi_L^n,o_L/\pi_L^n)\cong \Hom_{o_L}(H^{2-i}(K,T^*(1)),o_L)/\pi_L^n.\] Using that the cohomology groups are finitely generated $o_L$-modules and isomorphic to the inverse limits of the corresponding cohomology groups with coefficients modulo $\pi_L^n$ we see that the inverse limit of the above maps induces a surjective map
  \[H^i(K,T)\twoheadrightarrow \Hom_{o_L}(H^{2-i}(K,T^*(1)),o_L)\] with finite kernel, whence the claim after tensoring with $L$ over $o_L$ using the isomorphism $H^i(K,T)\otimes_{o_L}L\cong H^i(K,V)$ and analogously for $T^*(1).$
\end{proof}

Now let $W$ be a $L$-analytic representation of $G_L$  and set
\[ H^1_{/\dagger}(L_\infty^U,W^*(1)):=H^1_{\dagger}(L_\infty^U,W)^*,\]
which, by local Tate duality and Thm.\  \ref{thm:Herr}, is a quotient of $H^1(L_\infty^U,W^*(1)).$ By definition, the local Tate pairing induces a non-degenerate pairing
\begin{equation}\label{f:Tatedagger}
 <,>_{Tate,L,\dagger}\ : H^1_{\dagger}(L_\infty^U,W) \times H^1_{/\dagger}(L_\infty^U,W^*(1)) \longrightarrow H^2(L,L(1))\cong L.
\end{equation}

In order to compute this pairing more explicitly in certain situations we shall use Koszul-complexes. For this we have to assume first that $U$ is torsionfree. Following \cite[\S 4.2]{CoNi} we obtain for  any complete linearly topologised $o_L$-module $M$ with continuous $U$-action    a quasi-isomorphism \footnote{(unique up to homotopy, i.e., unique in the derived category of $o_L$-linear topological $U$-modules.) We have not yet defined any topology on the cocycles nor do we know whether the references says anything about it! $M$ is allowed to be any complete linearly topologised $o_L$-module with continuous $U$-action by \cite[V.1.2.6]{La}}
\begin{equation}\label{f:KoszulCtsCoh}
  K^\bullet_U(M)\xrightarrow{\simeq} \cC^\bullet(U,M)
\end{equation}
which arises as follows: Let $X_\bullet:=X_\bullet(U)$ and $Y_\bullet=Y_\bullet(U)$ denote the completed standard complex \cite[V.1.2.1]{La}, i.e., $X_n=\mathbb{Z}_p[[U]]^{\hat{\otimes}(n+1)}$, and the standard complex computing group cohomology, i.e., $Y_n=\mathbb{Z}_p[U]^{\otimes (n+1)}.$ Then, by \cite[Lem.\ V.1.1.5.1]{La} we obtain a diagram of complexes
\begin{equation}
\label{f:cupproductcompatible}
\xymatrix{
  Y_\bullet(U) \ar[d]_{ } \ar[r]^-{ \Delta} & Y_\bullet(U\times U)\cong Y_\bullet(U)\otimes_{\zp} Y_\bullet(U) \ar[d]^{} \\
  X_\bullet(U) \ar[d]_{ } \ar[r]^-{\Delta} & X_\bullet(U\times U)\cong X_\bullet(U)\widehat{\otimes}_{\zp}  X_\bullet(U) \ar[d]^{ } \\
  K_\bullet^U \ar[r]^-{\Delta} & K_\bullet^{U\times U} \cong K_\bullet^U \widehat{\otimes}_{\zp} K_\bullet^U , }
\end{equation}
which commutes up to homotopy (of filtered $\Lambda$-modules)
. Here the maps $\Delta$ are induced by the diagonal maps $U\to U\times U,$ e.g., $\zp[[U]]\to \zp[[U\times U]]\cong\zp[[U]]\widehat{\otimes}_{\zp}\zp[[U]].$  The first column induces a morphism
\[\Hom_{\Lambda}(K_\bullet^U,M)\to \Hom_{\Lambda, cts}(X_\bullet(U),M)\to \Hom_{\zp[U],cts}(Y_\bullet(U),M),\]
which is \eqref{f:KoszulCtsCoh}. The upper line induces as usual the cup product  on continuous group cohomology
\[H^r(U,M)\times H^s(U,N)\xrightarrow{\cup_U} H^{r+s}(U,M\otimes N)\]
via
\begin{multline*}
  \Hom_{\zp[U],cts}(Y_\bullet(U),M)\times \Hom_{\zp[U],cts}(Y_\bullet(U),N) \\  \xrightarrow{\times} \Hom_{\zp[U]\otimes\zp[U],cts}( Y_\bullet(U)\otimes_{\zp} Y_\bullet(U),M\otimes N)
  \xrightarrow{\Delta^*} \Hom_{\zp[U],cts}( Y_\bullet(U) ,M\otimes N).
\end{multline*}
The lower line induces analogously the Koszul-product
\[K^r_U(M)\times K^s_U(N)\xrightarrow{\cup_K} K^{r+s}_U(M\otimes N).\]
By diagram \eqref{f:cupproductcompatible} both products are compatible with each other.

Let $f$ be any continuous endomorphism of $M$ which commutes with the $U$-action; it  induces an operator on $K^\bullet(M)$
and  we denote by $K_{f,U}(M):=\mathrm{cone}\left( K^\bullet(M)\xrightarrow{f-\id}K^\bullet(M) \right)[-1]$ the mapping fibre of  $K^\bullet(f)$. Then  the quasi-isomorphism \eqref{f:KoszulCtsCoh} induces a quasi-isomorphism
\begin{equation}\label{f:KoszulCtsCohPhi}
  K_{\varphi,U}(M)\xrightarrow{\simeq} \mathcal{T}_{\varphi,U}(M).
\end{equation}

\begin{remark}
\label{rem:hyper-cup-product}
By a standard procedure   cup products   can be extended to hyper-cohomology (defined via total complexes), we follow \cite[(3.4.5.2)]{Ne}, but for the special case of a cone, see also \cite[Prop.\ 3.1]{niziol}. In particular, we obtain compatible cup products $\cup_K$   and $\cup_U$   for  $K_{\varphi,U}(M)$ and $\mathcal{T}_{\varphi,U}(M),$ respectively.
\end{remark}

Now we allow some arbitrary open subgroup $U\subseteq \Gamma_L$ and let $L'=L_\infty^U.$  Note that we obtain a decomposition $U\cong \Delta\times U'$ with
  a subgroup $U'\cong\mathbb{Z}_p^d$ of $U$ and  $\Delta$  the torsion subgroup of $U$. By  Lemma \ref{lem:HS} we obtain a canonical isomorphism
  \begin{equation}\label{f:KoszulCtsCohPhiDelta}
  K_{\varphi,U'}(M^\Delta)\xrightarrow{\simeq} \mathcal{T}_{\varphi,U}(M).
\end{equation}

Now let $M$ be a finitely generated projective $\cR$-module $M$  with continuous $U$-action. Then $M^*=\check{M}$ is again a finitely generated projective $\cR$-module $M$  with continuous $U$-action by Lemma \ref{lem-pair} (i). Hence $M$ as well as $M^\Delta$ satisfies the assumptions of \eqref{f:selfualM} and we have isomorphisms    \footnote{For $X\xrightarrow{f}Y$ we have $\mathrm{cone}(f)^*\cong    \mathrm{cone}(f^*)[-1] $    the isomorphism being realized by multiplying with $(-1)^i$ on $(X^*)^{i}$   and  $ \mathrm{cone}(f[n])= \mathrm{cone}((-1)^nf)[n]$.}
\begin{align}\notag
K_{\varphi,U}(M^\Delta)^*&\cong \mathrm{cone}\left( K^\bullet(M^\Delta)^* \xrightarrow{\varphi^*-1}K^\bullet(M^\Delta)^*   \right)\\ \label{f:KpsiU}
&= \mathrm{cone}\left( K^\bullet((M^\Delta)^*)[d]\xrightarrow{\varphi^*-1}K^\bullet((M^\Delta)^*)[d]  \right) \\
&= K_{\varphi^*,U}((M^*)_\Delta)[d+1]\notag\\
&=K_{\psi,U}(\check{M}_\Delta )[d+1]\notag\\
&=K_{\psi,U}(\check{M}^\Delta )[d+1]\notag.
\end{align}
The last isomorphism is induced by  the canonical isomorphism $\check{M}^\Delta\cong \check{M}_\Delta.$

Now note that \begin{equation}
\label{f:innerhom}
 {D_\rig^\dagger(W)}^\vee\cong D_\rig^\dagger(W^*(\chi_{LT}))
\end{equation}
for any $L$-analytic representation $W$ by the fact that the functor $D_\rig^\dagger$ respects inner homs, (cp. \cite[Remark 5.6]{SV} for the analogous case $D_{LT} $). Hence the tautological pairing $ev_2$ from \eqref{f:adjointmorphism} together with the above isomorphism \eqref{f:KpsiU} induces the following pairing:
\begin{equation}
\label{f:cupPsi}
\cup_{K,\psi}\ : h^{1}(K_{\varphi,U'}(D_{rig}^\dagger(W)^\Delta)) \quad\times\quad h^1(K_{\psi,U'}(D_{rig}^\dagger(W^*(\chi_{LT}))^\Delta) [d-1] ) \longrightarrow L
\end{equation}

\begin{remark}
For $U=U'$ and $M= D_{rig}^\dagger(W)   $, on the level of cochains
this pairing   is given as follows:\Footnote{  To see whether signs from the isomorphism \eqref{f:KpsiU} are  relevant, we specialize to the cyclotomic situation, i.e., $d=1$: On p. 1065 after \cite[Notation 2.3.13]{KPX} with trivial $\Delta$ one encounters an analogous diagram in the cyclotomic situation, which in contrast to ours does not arise by dualising the complex $K_{\varphi,U}(M)$. Instead they calculate the cup product of $K_{\varphi,\gamma}(M)$ with $K_{\varphi,\gamma}(\check{M})$  and identify the latter complex with $K_{\psi,\gamma^{-1}}(\check{M})$. Then their resulting pairing in degree $1$ looks like \[M\oplus M \times \check{M}\oplus \check{M}\to L,((x,y),(x',y'))\mapsto \{x, y'\}-\{y,x'\}.\] With our notation we obtain the following commutative diagram, which - up to shift by $2$ and maybe multiplication by $-1$ in all degrees - corresponds to the isomorphism
\begin{equation}
{\footnotesize \xymatrix{
{K_{\varphi,\gamma}(M )^*[-2]:} &  0 \ar[r]  & \check{M} \ar[d]_{=}  \ar[rr]^-{- \begin{pmatrix}
                           1-\psi \\
\gamma^{-1} - 1  \\
                               \end{pmatrix}
 } & & \check{M}\oplus \check{M}  \ar[d]^{\Xi}  \ar[rr]^{\begin{pmatrix}
                                        1-  \gamma^{-1} \  1-\psi\\                                       \end{pmatrix}
  } &&   \check{M} \ar[d]_{=}  \ar[r] & 0 \\
{K_{\psi,\gamma^{-1}}(\check{M})}: &  0 \ar[r]   & { \check{M}  }  \ar[rr]^-{\begin{pmatrix}
                                                    1-  \gamma^{-1} \\
                                                        {1-\psi} \\
                                                    \end{pmatrix}
    } && {\check{M}\oplus  \check{M} } \ar[rr]^-{\begin{pmatrix}
                                               {1-\psi} \  \gamma^{-1}-1 \\
                                             \end{pmatrix}
      } & & { \check{M} }  \ar[r] & 0   &
   }}\end{equation}
with $ \Xi(x,y)=(y,-x)$. This is the origin for the sign and both formulae are compatible.

}

\[\check{M}\oplus K^{d-1}(\check{M}) \times K^1(M)\oplus M \to L, ((x,y),(x',y'))\mapsto    \{y',x\} - y(x'),\]
where we again use  that $K^{d-1}(\check{M})\cong K^1(M)^*$   and where $\{\;,\;\}$ denotes the pairing \eqref{f:res-pair-general}. More generally, we have the following diagram

\begin{equation}
\label{f:duality}
{\scriptsize \xymatrix{
{K_{\varphi,U}(M ):} &  0 \ar[r]  & M   \ar[r]^-{ \left(\begin{smallmatrix}
                                 d^0_K \\
                                 {1-\varphi} \\
                               \end{smallmatrix}\right) }
  & K^1(M)\oplus M    \ar[r]^-{\left(\begin{smallmatrix}
                                          d^1_K & 0 \\
                                          {1-\varphi} & -d^0_K \\
                                        \end{smallmatrix}\right) }
 &   K^2(M)\oplus K^1(M)  \ar[r]^{ } & & \\
 &&\times &\times &\times \\
{K_{\psi,U}(\check{M})[d-1]}:   &  \ar[r]^{ }   & {K^{d-1}(\check{M})\oplus K^{d-2}(\check{M}) } \ar[d] \ar[r]^-{\left(\begin{smallmatrix}
                                                      d_K^{d-1} & 0 \\
                                                        {1-\psi} & -d_K^{d-2} \\
                                                    \end{smallmatrix}\right) }
     & {\check{M}\oplus K^{d-1}(\check{M})} \ar[d]\ar[r]^-{\left(\begin{smallmatrix}
                                               {1-\psi} & -d_K^{d-1} \\
                                             \end{smallmatrix}\right) }
       &  { \check{M} }\ar[d] \ar[r]^{ } & 0   &\\
&& L&L&L\\
\mbox{in degrees:} & &0&1&2
   }}\end{equation}
\end{remark}

Recall that $W=V^*(1)$   is $L$-analytic and set  $M=D^\dagger_{rig}(W)$ as well as  $\check{M}=D^\dagger_{rig}(V(\tau^{-1}))=D^\dagger_{rig}(W^*(\chi_{LT})).$ We obtain a Fontaine-style, explicit map
\begin{equation}\label{f:Fontaine-SV}
pr_U:D_{rig}^\dagger(V(\tau^{-1}))^{\psi=1} \to h^1(K_{\psi,U'}(\check{M}^\Delta)[d-1]  ),\; m\mapsto [(\bar{m},0)],
\end{equation}
where $\bar{m}=\frac{1}{\#\Delta}\sum_{\delta\in \Delta}\delta m$ denotes the image of $m$ under the map $\check{M}\twoheadrightarrow\check{M}_\Delta\cong\check{M}^\Delta.$

\begin{remark}\label{rem:res-cores}
Let $U_1\subseteq U$ an open subgroup with torsion subgroups $\Delta_1$ and $\Delta,$
respectively. Assume that the torsionfree parts $U_1'$ and $U'$ are generated by
$\gamma_1^{p^n},\ldots \gamma_d^{p^n}$ and $\gamma_1,\ldots \gamma_d $, respectively.
Then, for $M$ any complete locally convex $L$-vector space with continuous $U$-action,
the restriction and corestriction maps of Koszul-complexes from section
\ref{sec:KoszulCom} extend by functoriality to the mapping fibre
\begin{align}
cor^{U_1}_U:=&cor^{U_1'}_{U'}\circ
K_{\varphi,U_1'}(N_{\Delta/\Delta_1}):K_{\varphi,U_1'}(M^{\Delta_1})\to
K_{\varphi,U'}(M^\Delta)\notag\\ \notag
res^U_{U_1}:=&K_{\varphi,U_1'}(\iota) \circ
res^{U'}_{U_1'}:K_{\varphi,U'}(M^\Delta)\to K_{\varphi,U_1'}(M^{\Delta_1})
\end{align}
Here $N_{\Delta/\Delta_1} :M^{\Delta_1}\to M^\Delta$ denotes the norm/trace map
sending $m$ to $\sum_{\delta\in \Delta/\Delta_1}\delta m$ while $\iota: M^\Delta\to
M^{\Delta_1}$ is the inclusion. Taking duals as in \eqref{f:KpsiU} we also obtain
\begin{align}
cor^{U_1}_U:=&(res^U_{U_1})^*[1-d]:K_{\psi,U_1'}(M^{\Delta_1})\to
K_{\psi,U'}(M^\Delta)\notag\\\notag
res^U_{U_1}:=&(cor^{U_1}_U)^*[1-d]:K_{\psi,U'}(M^\Delta)\to K_{\psi,U_1'}(M^{\Delta_1})
\end{align}\Footnote{Do these map also arise directly from the Koszul complex variants
in the same way as for $K_{U,\varphi}?$}
(co)restriction maps for the $\psi$-Herr complexes.

Since inflation is compatible with restriction and corestriction one checks that the
above maps are compatible under the isomorphism \eqref{f:cohLphiB} with the usual maps
in Galois cohomology. Moreover, they define such maps on $H^1_\dagger$ and
$H^1_{/\dagger}$ via \eqref{f:cohLphi1Ver} and
$h^1(K_{\psi,U'}(D_{rig}^\dagger(W^*(\chi_{LT}))^\Delta[d-1])\cong
H^1_{/\dagger}(L',W^*(1)).$

By the discussion at the end of section \ref{sec:KoszulCom} the restriction map
\linebreak $K_{\varphi,U'}(M^\Delta)\xrightarrow{\mathrm{res}^{U}_{U_1}}
K_{\varphi,U_1'}(M^{\Delta_1}) $ and corestriction map $
K_{\varphi,U_1'}(M^{\Delta_1})\xrightarrow{\mathrm{cor}^{U_1}_{U}}
K_{\varphi,U'}(M^\Delta)$ in degree $0$
  are given as inclusion $M^\Delta\hookrightarrow M^{\Delta_1}$ and norm $M^{\Delta_1}
  \xrightarrow{N_{U',U_1'}\circ N_{\Delta/\Delta_1}} M^\Delta$, respectively, where
  \[N_{U',U_1'}:=\prod_{i=1}^d\sum_{k=0}^{p^n-1}\gamma_i^k\in \Lambda(U').\]
  Hence, by duality the restriction map
 $
 K_{\psi,U}(\check{M}^{\Delta})[d-1]^2\xrightarrow{\mathrm{res}^{U}_{U_1}}K_{\psi,U_1}(\check{M}^{\Delta_1})[d-1]^2$
  and corestriction map $
  K_{\psi,U_1}(\check{M}^{\Delta_1})[d-1]^2\xrightarrow{\mathrm{cor}^{U_1}_{U}}
  K_{\psi,U}(\check{M}^{\Delta})[d-1]^2 $ are given by the norm
  $\check{M}^{\Delta}\xrightarrow{(\Delta:\Delta_1)\upiota(N_{
  U',U_1'})}\check{M}^{\Delta_1}$ and projection map $
  \check{M}^{\Delta_1}\xrightarrow{\frac{1}{(\Delta:\Delta_1)}N_{\Delta/\Delta_1}}\check{M}^\Delta
  $, respectively. Here $ \upiota $ denotes the involution of $\Lambda(U)$ sending $u$
  to $u^{-1}$.
Note that the latter two descriptions also hold for the first components of $
K_{\psi,U}(\check{M}^{\Delta})[d-1]^1\xrightarrow{\mathrm{res}^{U}_{U_1}}K_{\psi,U_1}(\check{M}^{\Delta_1})[d-1]^1$
  and $ K_{\psi,U_1}(\check{M}^{\Delta_1})[d-1]^1\xrightarrow{\mathrm{cor}^{U_1}_{U}}
  K_{\psi,U}(\check{M}^{\Delta})[d-1]^1 $, respectively.
Hence, we obtain
\[\mathrm{cor}^{U_1}_{U}\circ pr_{U_1}=pr_{U} \mbox{ and } \mathrm{res}^{U}_{U_1}\circ
pr_{U}=pr_{U_1}\circ N_{\Delta/\Delta_1}\circ \upiota(N_{U',U_1'}). \]
\end{remark}

Berger and Fourquaux in contrast define a different Fontaine-style map in \cite[Thm.\ 2.5.8]{BF} for an $L$-analytic representation $Z$ and $N=D_\rig^\dagger(Z)$\footnote{We do not know whether this map coincides with the following composite we had used in older versions and which uses the shuffle maps from Propostion \ref{Herr}:
\begin{align}
\label{f:Fontaine-BF-thefootnote}h^1_{L_\infty^U,Z}:&D_\rig^\dagger(Z)^{\psi_L=\frac{q}{\pi_L}}\to
H^1_{an}(U,D_\rig^\dagger(Z)^{\psi_L=\frac{q}{\pi_L}})\to h^1(\mathcal{T}_{\frac{\pi}{q}\psi_L,U}^{an}(N))\cong h^1(\mathcal{T}_{\varphi_L,U}^{an}(M))\\\notag
&\phantom{mmmmmmmmmmmmmmmmmmmmmm}\to h^1(\mathcal{T}_{\varphi_L,U}(N)) \cong H^1_\dagger(L_\infty^U,Z) .
\end{align}
 Here the second map is stemming from the spectral sequence \eqref{f:SSandim1}, the third from Propostion \ref{Herr}, the fourth is the natural map \eqref{f:an-ctsh}, while the last one is \eqref{f:cohLphi1Ver}.
}
\begin{align}
\label{f:Fontaine-BF}h^1_{L_\infty^U,Z}:&D_\rig^\dagger(Z)^{\psi_L=\frac{q}{\pi_L}}\to
H^1_{an}(U,D_\rig^\dagger(Z)^{\psi_L=\frac{q}{\pi_L}})\to  h^1(\mathcal{T}_{\varphi_L,U}(N)) \cong h^1(K_{\varphi_L,U}(N^\Delta)),\\ \notag
&y\mapsto\phantom{mmmmmmmmmmm} [c_b(y)]\mapsto\phantom{mmm}[(c_b(y),-m_c)]\mapsto  [(\tilde{c_b}(y),-\tilde{m_c})],
\end{align}
  in which  the cocycle $h^1_{L_\infty^U,Z}(y)$ is given in terms of the pair $(c_b(y),-m_c)$ in the notation of Thm.\  2.5.8 in (loc.\ cit.): $m_c$ is the unique element in $D_\rig^\dagger(Z)^{\psi_L=0} $ such that \begin{equation}
  \label{f:mc}
  (\varphi_L-1)c_b(y)(\gamma)=(\gamma-1)m_c\end{equation} for all $\gamma\in U$ and this pair defines the extension class in the sense of Lemma \ref{lem:ext}. Here, the first map is implicity given by Prop.\  2.5.1 in (loc.\ cit.), the second one is
   the composite from maps arising in   Cor.\ 2.2.3, Thm.\ 2.2.4, of (loc.\ cit.) with the natural map from analytic to continuous cohomology
   \[H^1_{an}(U,D_\rig^\dagger(Z)^{\psi_L=\frac{q}{\pi_L}})\to  H^1_{an}(U\times \Psi,D_\rig^\dagger(Z))\cong H^1_{an}(U\times \Phi,D_\rig^\dagger(Z)) \to H^1(U\times \Phi,D_\rig^\dagger(Z))\]
   combined with the interpretation of extension classes (see \S 1.4 in (loc.\ cit.) and Lemma \ref{lem:ext}), while the last one is \eqref{f:KoszulCtsCohPhiDelta} (the concrete image $(\tilde{c_b}(y),-\tilde{m_c})$ will be of interest for us only in the situation where $\Delta$ is trivial, when $\tilde{m_c}={m_c}$).

According to \cite[Prop.\ 2.5.6, Rem.\ 2.5.7]{BF} this map also satisfies
\begin{equation}
\label{f:coresh1} cor^U_{U'}\circ h^1_{L_\infty^{U'},Z}=h^1_{L_\infty^U,Z}.
\end{equation}

Since ${D_\rig^\dagger(V(\tau^{-1}))})^\vee\cong D_\rig^\dagger(V^*(1))$ by \eqref{f:innerhom},  concerning the Iwasawa-pairing we have the following
\begin{proposition}\label{prop:Tate0}
For a $G_L$-representation $V$ such that $V^*(1)$ is $L$-analytic the following diagram consisting of $D(\Gamma_L,K)$-$\upiota_*$-sesquilinear pairings (in the sense of  \eqref{sesquilinear}) is commutative
\begin{equation*}
\xymatrix{
h^{1}(K_{\varphi,U'}(D_{rig}^\dagger(V^*(1)(\chi_{LT}^{j}))^\Delta))\ar@{}[r]|{\times} & h^1(K_{\psi,U'}(D_{rig}^\dagger(V(\chi_{LT}^{-j}))^\Delta) [d-1] )  \ar[r]^-{\cup_{K,\psi}} &   L\subseteq\Cp \\
    {D_{\mathrm{rig}}^\dagger(V^*(1))^{\psi_L=\frac{q}{\pi_L}}}\ar[u]_{h^1_{L,V^*(1)(\chi_{LT}^{j})}\circ tw_{\chi_{LT}^{j}}} \ar@{}[r]|{\times} &  {D_{\mathrm{rig}}^\dagger(V(\tau^{-1}))^{\psi_L=1}}\ar[u]_{pr_U\circ tw_{\chi_{LT}^{-j}} } \ar[r]^-{    \frac{q-1}{q} \{,\}_{Iw,\Gamma_L}}& D(\Gamma_L,\Cp)\ar[u]_{ ev_{\chi_{LT}^{-j}}} }
\end{equation*}
taking $U'\times \Delta=U=\Gamma_L.$
\end{proposition}

\begin{proof}
By Lemma \ref{lem:twist}, it suffices to show the case $j=0,$ i.e., the trivial character $\chi_{triv}.$ Furthermore,
 it suffices to show the following statement for any subgroup of the form $\Gamma_n$ without any $p$-torsion:
\begin{equation}\label{f:evLn}
 q^{-n} ev_{L_n,\chi_{triv|\Gamma_n}}\circ \{x,y\}_{Iw,\Gamma_n}=h^1_{L_n,V^*(1)}(x)\cup_{K,\psi} pr_{\Gamma_n}(y)
\end{equation}
for $ x\in {D_{\mathrm{rig}}^\dagger(V^*(1))^{\psi_L=\frac{q}{\pi_L}}}, y\in  {D_{\mathrm{rig}}^\dagger(V(\tau^{-1}))^{\psi_L=1}}.$

Indeed, by Remark \ref{rem:res-cores}, for every such $n$, we have the commutative diagram
\begin{equation*}
\xymatrix{
h^{1}(K_{\varphi,\Gamma_n}(D_{rig}^\dagger(V^*(1)))) \ar[d]^{\mathrm{cor}}\ar@{}[r]|{\times} & h^1(K_{\psi,\Gamma_n}(D_{rig}^\dagger(V) ) [d-1] )  \ar[r]^(0.8){\cup_{K,\psi}} &  L\ar@{=}[d] \\
h^{1}(K_{\varphi,U'}(D_{rig}^\dagger(V^*(1))^\Delta))\ar@{}[r]|{\times} &h^1(K_{\psi,U'}(D_{rig}^\dagger(V)^\Delta) [d-1] )  \ar[u]_{\mathrm{res}}  \ar[r]^(0.8){\cup_{K,\psi}} &  L.}
\end{equation*}
Hence we 
obtain using \eqref{f:coresh1}
\begin{align*}
 h^1_{L',V^*(1)}(x)\cup_{K,\psi} pr_{U}(y) &=(\mathrm{cor}\circ h^1_{L_n, V^*(1)}(x)) \cup_{K,\psi} pr_{U}(y) \\
 &=h^1_{L_n,V^*(1)}(x)\cup_{K,\psi}( \mathrm{res}\circ pr_{U}(y))\\
 &=h^1_{L_n,V^*(1)}(x)\cup_{K,\psi}( pr_{\Gamma_n}(N_{\Delta}\circ \upiota(N_{U',\Gamma_n})y)),
\end{align*}
where we use Remark \ref{rem:res-cores} for the last equality. On the other hand one easily checks that\footnote{This is obvious if you decompose $D(U,\Cp)=\bigoplus_{g\in U/\Gamma_n}D(\Gamma_n,\Cp)g$ with respect to the inverses of the representatives used in the definition of $ N_{\Delta}\circ \upiota(N_{U',\Gamma_n})$.}
\[ev_{L_n,\chi_{triv|\Gamma_n}}\circ pr_{U,\Gamma_n}\circ N_{\Delta}\circ \upiota(N_{U',\Gamma_n})=ev_{L,\chi_{triv}}:D(U,\Cp)\to\Cp, \]
whence
\begin{align*}
 ev_{L,\chi_{triv}}\circ \frac{q-1}{q}\{x,y\}_{Iw,U}&= \frac{q-1}{q} ev_{L_n,\chi_{triv|\Gamma_n}}\circ pr_{U,\Gamma_n}( N_{\Delta}\circ \upiota(N_{U',\Gamma_n})\{x,y\}_{Iw,U})\\
 &= \frac{q-1}{q}ev_{L_n,\chi_{triv|\Gamma_n}}\circ   pr_{U,\Gamma_n}( \{x,N_{\Delta}\circ \upiota(N_{U',\Gamma_n})y\}_{Iw,U})\\
 &= \frac{q-1}{q} [U:\Gamma_n]^{-1}ev_{L_n,\chi_{triv|\Gamma_n}}\circ \{x,N_{\Delta}\circ \upiota(N_{U',\Gamma_n})y\}_{Iw,\Gamma_n}\\
 &=  q^{-n}ev_{L_n,\chi_{triv|\Gamma_n}}\circ \{x,N_{\Delta}\circ \upiota(N_{U',\Gamma_n})y\}_{Iw,\Gamma_n}
\end{align*}
where we have used Remark \ref{rem:IWpairU}  for the last equation.

In order to prove \eqref{f:evLn} choose $n=n_0$ (see section \ref{sec:groupRobba}). As recalled  in \eqref{f:Fontaine-BF})  the map \[h^1_{L_{n_0},V^*(1)}:D_\rig^\dagger(V^*(1))^{\psi_L=\frac{q}{\pi_L}}\to h^1(K_{\varphi_L,\Gamma_{n_0}}(D_\rig^\dagger(V^*(1))))    \]
is given by the cocycle $h^1_{L_{n_0},V^*(1)}(x)$   in terms of the pair $(\tilde{c_b}(x),-m_c).$ Note that we have
\begin{equation*}
 m_c=\widehat{\Xi_b}(\varphi_L-1)x.
\end{equation*}
Indeed, by \cite[Thm.\ 2.5.8]{BF} we have $c_b(x)(b_j^k)=(b_j^k-1)\widehat{\Xi_b}x$ for all $j, k\geq 0$, which together with \eqref{f:mc} and the uniqueness of $m_c$ (loc.\ cit.) implies the claim.
On the other hand we have the map \eqref{f:Fontaine-SV}
\begin{equation*}
pr_{\Gamma_{n_0}}: D_{rig}^\dagger(V(\tau^{-1}))^{\psi_L=1} \to h^1(K_{\psi,\Gamma_{n_0}}(D_{rig}^\dagger(V(\tau^{-1})))[d-1]  ),\;\;\;y\to \mbox{ class of }(y,0).
\end{equation*}
Thus the pairing $ \cup_{K,\psi} $ sends by construction (see diagram \eqref{f:duality}) the above classes to
\begin{align*}
h^1_{L_n,V^*(1)}(x)\cup_{K,\psi} pr_{\Gamma_n}(y)&=0(\tilde{c_b}(x))+\{-\widehat{\Xi_b}(\varphi_L-1)x,y\}\\
&=\{\widehat{\Xi_b}(\varphi_L-1)x, (\frac{\pi_L}{q}\varphi_L-1)y\}\\
&=<\widehat{\Xi_b},\{x,y\}_{Iw,\Gamma_n}>_{\Gamma_n}\\
&=  q^{-n}\mathrm{aug}(\{x,y\}_{Iw,\Gamma_n}).
\end{align*}
Here the second equality holds  due to Lemma  \ref{lem-pair} because the left hand side belongs to $D_\rig^\dagger(V^*(1))^{\psi_L=0}$, the third one is \eqref{f:definingproperty} while the last one comes from \eqref{f:Xi}.
\end{proof}

\begin{proposition}\label{prop:Koszul-product}
For $W$ an $L$-analytic representation we have a canonical commutative diagram
\[\footnotesize\xymatrix{
\cup_{K,\psi}:h^{1}(K_{\varphi,U'}(D_{rig}^\dagger(W)^\Delta))\ar[d]^{\cong}_b \ar@{}[r]|(0.45){\times} & { \phantom{} }h^1(K_{\psi,U'}(D_{rig}^\dagger(W^*(\chi_{LT}))^\Delta) [d-1] )\ar[d]^{\cong}_a \ar[rr]^(0.45){} &&   L\ar@{=}[d]\\
  <,>_{Tate,L,\dagger}:
H^{1}_{\dagger}(L',W)\ar@{^(->}[d]\ar@{}[r]|(0.65){\times} & H^1_{/\dagger}(L',W^*(1)) \ar[r]^(0.45){} & H^2(L',L(1))\ar@{=}[d]\ar[r]^(0.65){\cong} & L\ar@{=}[d]\\
    <,>_{Tate,L'}:
H^{1} (L',W)\ar@{}[r]|(0.65){\times} & H^1 (L',W^*(1))\ar@{->>}[u]^{pr} \ar[r]^(0.45){} & H^2(L',L(1))\ar[r]^(0.65){\cong} & L.}
\]
Moreover, the isomorphism $a$ is compatible with the middle map of the diagram \eqref{f:diagcup3}.
 \end{proposition}

\begin{proof}
The lower square of pairings comes from Tate duality as in Prop.\ \ref{Tate-local} and \eqref{f:Tatedagger}. Its commutativity holds by definition.
  In the upper square of pairings the left upper vertical isomorphism $b$ arises from \eqref{f:cohLphi1Ver} combined with \eqref{f:KoszulCtsCohPhiDelta}, while the middle vertical isomorphism $a$ is uniquely determined as adjoint of the latter because both pairings are non-degenerate: The middle one by definition of  $H^1_{/\dagger}$ while the upper one due to Corollary \ref{cor:finitedim}(ii) with $W=V^*(1)$. Therefore one immediately checks that  $a^{-1}\circ pr$  is   induced by the cohomology of the middle map (going down) in diagram \eqref{f:diagcup3}   (being the same as the middle map (going from right to left) of  diagram \eqref{f:PsiSlashdagger} upon identifying $h^1\left( K_{\psi,U'}^\bullet({D}(V(\tau^{-1}))^\Delta)[d-1]\right)$ and $H^1(L',V)$ by the isomorphism described there).
\end{proof}

Combining the last two propositions we get the following result.
\begin{corollary}\label{prop:Tate}
For a $G_L$-representation $V$ such that $V^*(1)$ is $L$-analytic the following diagram is commutative
\begin{equation*}
\xymatrix{
H^1_{\dagger}(L,V^*(1)(\chi_{LT}^{j}))\ar@{}[r]|{\times} &H^1_{/\dagger}(L,V(\chi_{LT}^{-j}))  \ar[rr]^-{<,>_{Tate,L,\dagger}} && H^2(L,L(1))\cong L\subseteq\Cp \\
    {D_{\mathrm{rig}}^\dagger(V^*(1))^{\psi_L=\frac{q}{\pi_L}}}\ar[u]_{h^1_{L,V^*(1)}\circ tw_{\chi_{LT}^{j}}} \ar@{}[r]|{\times} &  {D_{\mathrm{rig}}^\dagger(V(\tau^{-1}))^{\psi_L=1}}\ar[u]_{pr_{L}\circ tw_{\chi_{LT}^{-j}} } \ar[rr]^-{ \frac{q-1}{q}\{,\}_{Iw}} && D(\Gamma_L,\Cp).\ar[u]_{ ev_{\chi_{LT}^{-j}}} }
\end{equation*}
\end{corollary}

\begin{remark}
For applications it might be useful to renormalize $\cup_{K,\psi}$ by the factor $\frac{q}{q-1}$, i.e., setting $\cup'_{K,\psi}:= \frac{q}{q-1} \cup_{K,\psi}.$ Then we would get rid of the factor $\frac{q-1}{q}$ in front of the Iwasawa pairing in the above results. Moreover and more important the new normalisation would be compatible with the cyclotomic situation taking $L=\mathbb{Q}_p$, $\pi_L=p=q,$ i.e., the upper pairing in Proposition \ref{prop:Tate0} would coincide - at least up to a sign - with the cup product pairing of Galois cohomology
\[\footnotesize\xymatrix{
H^{1} (\mathbb{Q}_p,V^*(j+1))\ar@{}[r]|(0.5){\times} & H^1 (\mathbb{Q}_p,V(-j)) \ar[r]^(0.45){} & H^2(\mathbb{Q}_p,\mathbb{Q}_p(1))\ar[r]^(0.65){inv}_(0.65){\cong} & \mathbb{Q}_p.}
\]
using Tate's trace map
\[inv:H^2(\mathbb{Q}_p,\mathbb{Q}_p(1))\cong\mathbb{Q}_p\]
given by class field theory, if one chooses $Z=\gamma-1$ for a topological generator $\gamma$ of $\Gamma_{\mathbb{Q}_p}$ satisfying $\log\chi_{cyc}(\gamma)=1$. This follows from \cite[Prop.\ 1.3.4, Thm.\ 2.2.6]{benois}, \cite[Thm.\ 5.2,Rem.\ 5.3]{Her} and \cite[Rem.\ 2.3.11/12]{KPX}: they claim that $-\frac{p-1}{p}inv$ corresponds to the trace map from the second cohomology group of the $\varphi$-Herr complex induced by sending $f\otimes\eta$ to $\frac{1}{\log\chi_{cyc}(\gamma)} \mathrm{res}_{\omega_{cyc}}(f\frac{d\omega_{cyc}}{1+\omega_{cyc}}  ).$
\Footnote{The remaining sign difference here may arise by interchanging the factors of the alternating cup-product or from different sign conventions.}
\end{remark}

With respect to evaluating at a character we have the following analogue of Corollary \ref{prop:Tate}.

\begin{proposition}\label{prop:cris}
For a $G_L$-representation $V$ such that $V^*(1)$ is $L$-analytic the following diagram is commutative \small
\begin{equation*}\scriptsize
\xymatrix{
\mathbb{C}_p\otimes_L D_{cris,L}(V^*(1)(\chi_{LT}^{j}))\ar@{}[r]|{\times} & \mathbb{C}_p\otimes_LD_{cris,L}(V(\tau^{-1})(\chi_{LT}^{-j}))  \ar[r]^(0.5){[,]_{cris}^{}} & \mathbb{C}_p\otimes_LD_{cris,L}(L(\chi_{LT}))\cong\Cp \\
    {D(\Gamma_L,\mathbb{C}_p)\otimes_L D_{cris,L}(V^*(1))}\ar[u]_{ev_{\chi_{LT}^{-j}}\otimes t_{LT}^{-j} \eta^{\otimes j} } \ar@{}[r]|{\times} &  {D(\Gamma_L,\mathbb{C}_p)\otimes_L D_{cris,L}(V(\tau^{-1})) }\ar[u]_{ev_{\chi_{LT}^{j}}\otimes t_{LT}^{j} \eta^{\otimes -j} } \ar[r]^(0.6){[,]^0}& D(\Gamma_L,\Cp),\ar[u]_{ev_{\chi_{LT}^{-j}}} }
\end{equation*}\normalsize
where, for the identification in the right upper corner we choose $t_{LT}^{-1}\otimes \eta$ as a basis.
\end{proposition}

\begin{proof}
Using Lemma \ref{lem:twistcris} below the statement is reduced to $j=0.$ Evaluation of \eqref{f:def[]0} implies the claim in this case.
\end{proof}

\begin{lemma}\label{lem:twistcris}
There is a commutative diagram \small
\begin{equation*}
\xymatrix{
 {D(\Gamma_L,\mathbb{C}_p)\otimes_L D_{cris,L}(V^*(1)(\chi_{LT}^{j}))}  \ar@{}[r]|{\times} &  {D(\Gamma_L,\mathbb{C}_p)\otimes_L D_{cris,L}(V(\tau^{-1})(\chi_{LT}^{-j}))}\ar[r]^(0.7){[,]^0}  & D(\Gamma_L,\Cp) \\
  {D(\Gamma_L,\mathbb{C}_p)\otimes_L D_{cris,L}(V^*(1))}\ar[u]_{ Tw_{\chi_{LT}^{-j}}\otimes t_{LT}^{-j} \eta^{\otimes j}  }  \ar@{}[r]|{\times} &  {D(\Gamma_L,\mathbb{C}_p)\otimes_L D_{cris,L}(V(\tau^{-1})) }\ar[r]^(0.7){[,]^0}\ar[u]_{Tw_{\chi_{LT}^{j}}\otimes t_{LT}^{j} \eta^{\otimes -j} }  & D(\Gamma_L,\Cp)\ar[u]_{Tw_{\chi_{LT}^{-j}}}. }
\end{equation*}\normalsize
\end{lemma}

\begin{proof}The claim follows immediately from \eqref{f:def[]0}, the compatibility of the usual $D_{cris}$-pairing with twists and the fact that $Tw_{\chi_{LT}^j}(\lambda\iota_*(\mu))=Tw_{\chi_{LT}^j}(\lambda)\iota_*(Tw_{\chi_{LT}^{-j}}(\mu))$ holds.
\end{proof}

\subsubsection{The interpolation formula for the regulator map}\label{subsec:interpolation}

 In this subsection we are going to prove the   interpolation property for $\mathcal{L}_V$.
 First recall that we introduced in section \ref{sec:KR} the notation $D_{dR,L'}(V) := ( B_{dR} \otimes_{\qp} V)^{G_{L'}}.$  Since $B_{dR}$ contains the algebraic closure $\overline{L}$ of $L$ we have the isomorphism
\begin{equation*}
  B_{dR} \otimes _{\mathbb{Q}_p} V = (B_{dR} \otimes _{\mathbb{Q}_p} L) \otimes_L V  \xrightarrow{\;\cong\;} \prod_{\sigma \in G_{\mathbb{Q}_p}/G_L} B_{dR} \otimes_{\sigma,L} V
\end{equation*}
which sends $b \otimes v$ to $(b \otimes v)_\sigma$. The tensor product in the factor $B_{dR} \otimes_{\sigma,L} V$ is formed with respect to $L$ acting on $B_{dR}$ through $\sigma$. With respect to the $G_L$-action the right hand side decomposes according to the double cosets in $G_L \backslash G_{\mathbb{Q}_p}/G_L$. It follows, in particular, that $D_{dR}^{\id}(V) := ( B_{dR} \otimes_L V)^{G_L}$ is a direct summand of $D_{dR,L}(V)$ and we denote by $pr^{\id}$ the corresponding projection. Similarly, $tan_{L,\id}(V) := ( B_{dR}/B^+_{dR} \otimes_L V)^{G_L}  $ is a direct summand of $tan_L(V):= ( B_{dR} \otimes_L V)^{G_L} $.
 More generally, also the filtration $D_{dR,L}^i(V)$ decomposes into direct summands.

According to \cite[Appendix A]{SV15} the dual Bloch-Kato exponential map is uniquely determined by the commutativity of the following diagram, in which all pairings are perfect:
\begin{equation}
\label{f:dualexp}
\xymatrix{
H^{1} (L',W)\ar[d]^{\exp_{L',W}^*}  \ar@{}[r]|(0.45){\times} & { \phantom{} }H^1 (L',W^*(1)) \ar[rr]^-{<,>_{Tate,L'}} &&   L\ar@{^(->}[d]\\
D_{dR,L'}^0(W)\ar@{^(->}[d]\ar@{}[r]|(0.45){\times} & tan_{L'}(W^*(1)) \ar[u]_{\exp_{L',W^*(1)}} \ar[r]^{} & D_{dR,L'}(\mathbb{Q}_p(1))\ar@{=}[d]\ar[r]^-{\cong} & L'\ar@{=}[d]\\
    D_{dR,L'}(W)
\ar@{}[r]|(0.45){\times} & D_{dR,L'}(W^*(1))\ar@{->>}[u]^{pr} \ar[r]^{} & D_{dR,L'}(\mathbb{Q}_p(1))\ar[r]^-{\cong} & L'.}
\end{equation}
In the Lubin-Tate setting we can also consider the dual of the identity component $\exp_{L',W^*(1),\id} $ of $\exp_{L',W^*(1)}: $
\footnote{
 We have the compatibility of the following pairings
\[
\footnotesize
\xymatrix{
{ D_{cris,L}(V^*(1) )}\ar[d]^{\cong} & \times &  { D_{cris,L}(V(\tau^{-1}))}\ar[d]^{\cong}   \ar[rr]^-{[,]^{}_{cris}}  && D_{cris,L} (L(\chi_{LT}))\ar[d]^{\cong} \ar[r]^(0.7){\cong} &  L \ar[d]_{=}\\
{ D_{dR}^{\id}(V^*(1) )} \ar[d]^{injective}  & \times &  { D_{dR}^{\id}(V(\tau^{-1}))} \ar[rr]^-{[,]^{id}_{dR}}  && D_{dR}^{\id}(L(\chi_{LT}))\ar[r]^-{\cong} & L\ar[dd]_{=} \ar[dl]_{ } \\
 { D_{dR}(V^*(1) )}  & \times &  {D_{dR}(V(\tau^{-1}))} \ar[u]^{pr^{\id}}_{surjective} \ar[r]^{}  & D_{dR}(L(\chi_{LT})) \ar[r]^-{\cong}& L\otimes_{\qp} L\ar[dr]^{pr^{\id}}\\
  {  D_{dR}(V^*(1))}\ar[u]^{=}    & \times &  { D_{dR}(V ) }\ar[u]_{\id\otimes t_{cyc}t_{LT}^{-1}\otimes\eta\otimes\eta_{cyc}^{-1} }\ar[rr]^-{[,]_{dR}}  && D_{dR}({\qp}(1))\ar[r]^-{\cong}& L }
\]
\comb{  cf.\ \cite[Appendix,(57)]{SV15}.}}
\begin{equation}
\label{f:dualexpId}
\xymatrix{
H^{1} (L',W)\ar[d]^{\widetilde{\exp}_{L',W,\id}^*}  \ar@{}[r]|(0.45){\times} & { \phantom{} }H^1 (L',W^*(1)) \ar[rr]^-{<,>_{Tate,L'}} &&   L\ar@{^(->}[d]\\
D_{dR,L'}^{\id,0}(W(\tau^{-1}))\ar@{^(->}[d]\ar@{}[r]|(0.45){\times} & tan_{L',\id}(W^*(1)) \ar[u]_{\exp_{L',W^*(1),\id}} \ar[r]^{} & D_{dR,L'}^{\id}(L(\chi_{LT}))\ar@{=}[d]\ar[r]^-{\cong} & L'\ar@{=}[d]\\
    D_{dR,L'}^{\id}(W(\tau^{-1}))
\ar@{}[r]|(0.45){\times} & D_{dR,L'}^{\id}(W^*(1))\ar@{->>}[u]^{pr} \ar[r]^{} & D_{dR,L'}^{\id}(L(\chi_{LT}))\ar[r]^-{\cong} & L'.}
\end{equation}
Upon noting that under the identifications $D_{dR,L'}(\mathbb{Q}_p(1)) \cong L'$ and $D_{dR,L'}^{\id}(\mathbb{Q}_p(1)){\cong}  L'$ the elements $t_{\mathbb{Q}_p}\otimes \eta_{cyc}$  and $t_{LT}\otimes \eta$  are sent to $1,$ one easily checks that, if $W^*(1)$ is $L$-analytic, whence the inclusion $tan_{L',\id}(W^*(1)) \subseteq tan_{L'}(W^*(1)) $ is an equality and  $\exp_{L',W^*(1),\id} =\exp_{L',W^*(1)} $, it holds
\begin{equation}\label{f:dualexpcomp}
  \mathbb{T}_{\tau^{-1}}\circ \exp_{L',W}^*=\widetilde{\exp}_{L',W,\id}^*,
\end{equation}
where $\mathbb{T}_{\tau^{-1}}: D_{dR,L'}^0(W) \to D_{dR,L'}^{\id,0}(W(\tau^{-1})) $ is the isomorphism, which sends $b\otimes v$ to $b\frac{t_{\mathbb{Q}_p}}{t_{LT}}\otimes v\otimes \eta\otimes\eta_{cyc}^{\otimes -1};$ note that $\frac{t_{\mathbb{Q}_p}}{t_{LT}}\in (B_{dR}^+)^\times,$ whence the filtration is preserved.

Now let $W$ be an $L$-analytic, crystalline  $L$-linear representation of $G_L.$ Recall that $\eta=(\eta_n)_n$ denotes a fixed generator of $T_\pi$ and that the map $tw_{\chi_{LT}^j}:D^\dagger_{rig}(W)\to D^\dagger_{rig}(W(\chi_{LT}^j))$ has been defined before Lemma \ref{lem:twist}. For $D_{cris}$ twisting $D_{cris,L}(W)\xrightarrow{-\otimes e_j} D_{cris,L}(W(\chi_{LT}^j)) $ maps $d$ to $d\otimes e_j$ with $e_j:=t_{LT}^{-j}\otimes \eta^{\otimes j}\in D_{cris,L}(L(\chi_{LT}^j)).$

If we assume, in addition, that
\begin{enumerate}
\item $W$ has Hodge Tate weights   $\leq 0$,   whence $W^*(1)$ has Hodge Tate weights $\geq 1$ and $D_{dR,L}^0(W^*(1))=0$, and
\item $D_{cris,L}(W^*(\chi_{LT}))^{\varphi_L=\frac{q}{\pi_L}}=0,$
\end{enumerate}
then  $exp_{L,W^*(1)}:D_{dR,L}(W^*(1))\hookrightarrow H^1(L,W^*(1))$ is   injective    with image $H^1_e(L,W^*(1))=H^1_f(L,W^*(1))$ by our assumption (see \cite[Cor.\ 3.8.4]{BK}). We denote its inverse by
\[\log_{L,W^*(1)}:H^1_f(L,W^*(1))\to D_{dR,L}(W^*(1))\] and define
\[\widetilde{\log}_{L,W^*(1)}:H^1_f(L,W^*(1))\to D_{dR,L}(W^*(1))\xrightarrow{\mathbb{T}_{\tau^{-1}}} D_{cris,L}(W^*(\chi_{LT}))\]
where  (by abuse of notation) we also write $\mathbb{T}_{\tau^{-1}}: D_{dR,L}(W^*(1)) \to D_{dR,L}^{\id}(W^*(\chi_{LT}))=D_{cris,L}(W^*(\chi_{LT})) $ for the isomorphism, which sends $b\otimes v$ to $b\frac{t_{\mathbb{Q}_p}}{t_{LT}}\otimes v\otimes \eta\otimes\eta_{cyc}^{\otimes -1}.$ We obtain the following commutative diagram, which defines the dual map $\log_{L,W}^*$ being inverse to $\exp_{L,W}^*$ (up to factorisation over $H^{1} (L,W)/H^1_f(L,W) $):
\begin{equation}
\label{f:duallog}
\xymatrix{
H^{1} (L,W)/H^1_f(L,W) \ar@{}[r]|(0.5){\times} & { \phantom{} }H^1_f (L,W^*(1)) \ar[d]_{\log_{L,W^*(1)}}\ar[rr]^-{<,>_{Tate,L}} &&   L\ar@{=}[d]\\
D_{dR,L}(W)
\ar@{}[r]|(0.45){\times}\ar[u]^{\log_{L,W}^*}  & D_{dR,L}(W^*(1))  \ar[r]^(0.5){} & D_{dR,L}(\mathbb{Q}_p(1))\ar[r]^-{\cong} & L.}
\end{equation}
Similarly as above we obtain a commutative diagram more convenient for the Lubin-Tate setting:
\begin{equation}
\label{f:duallogId}
\xymatrix{
H^{1} (L,W)/H^1_f(L,W) \ar@{}[r]|(0.45){\times} & { \phantom{} }H^1_f (L,W^*(1)) \ar[d]_{\widetilde{\log}_{L,W^*(1)}} \ar[rr]^-{<,>_{Tate,L}} &&   L\ar@{=}[d]\\
    D_{dR,L}^{\id}(W)\ar[u]^{ {\log}_{L,W,\id}^*}
\ar@{}[r]|(0.45){\times} & D_{dR,L}^{\id}(W^*(\chi_{LT}))  \ar[r]^(0.45){} & D_{dR,L}^{\id}(L(\chi_{LT}))\ar[r]^-{\cong} & L.}
\end{equation}
We write $\mathrm{Ev}_{W,n}:\mathcal{O}_L \otimes_L D_{cris,L}(W)\to L_n\otimes_L D_{cris,L}(W)$ for the composite $\partial_{D_{cris,L}(W)}\circ \varphi_q^{-n}$ from the introduction of \cite{BF}, which actually sends $f(Z)\otimes d$ to $f(\eta_n)\otimes \varphi_L^{-n}(d).$ By abuse of notation we also use $\mathrm{Ev}_{W,0}$ for the analogous map   $ \mathcal{O}_K \otimes_L D_{cris,L}(W)\to K\otimes_L D_{cris,L}(W)$. For $x\in  D(\Gamma_L,K)\otimes_L D_{cris,L}(W)$ we denote by $x(\chi_{LT}^{j}) $ the image under the map $ D(\Gamma_L,K)\otimes_L D_{cris,L}(W)\to K\otimes_L D_{cris,L}(W),$ $\lambda\otimes d\mapsto \lambda(\chi_{LT}^{j})\otimes d.$

\begin{lemma}\label{lem:Ev}
Assume that $\Omega$ is contained in $K$. Then there are commutative diagrams
\begin{equation*}
\xymatrix{
  D(\Gamma_L,K)\otimes_L D_{cris,L}(W) \ar[d]_{ev_{\mathrm{triv}}} \ar[r]^-{\mathfrak{M}\otimes \id} & \mathcal{O}_K \otimes_L D_{cris,L}(W)\ar[d]^{\mathrm{Ev}_{W,0}} &&  \mathcal{O}_K \otimes_L D_{cris,L}(W) \ar[ll]_-{1-\varphi_L\otimes \varphi_L} \ar[d]^{\mathrm{Ev}_{W,0}}\\
 \phantom{D(\Gamma_L,)} K\otimes_L D_{cris,L}(W) \ar@{=}[r]^{ } &  \phantom{O} K\otimes_L D_{cris,L}(W)   && K\otimes_L D_{cris,L}(W) \ar[ll]_-{1-\id\otimes\varphi_L}    }
\end{equation*}
and
\begin{equation*}\scriptsize
\xymatrix{
  D(\Gamma_L,K)\otimes_L D_{cris,L}(W) \ar[d]_{Tw_{\chi_{LT}^{-j}}\otimes e_j} \ar[r]^-{\mathfrak{M}\otimes \id} & \mathcal{O}_K \otimes_L D_{cris,L}(W)\ar[d]^{(\frac{\partial}{\Omega})^{-j}\otimes e_j} &  \mathcal{O}_K \otimes_L D_{cris,L}(W) \ar[l]_{1-\varphi_L\otimes \varphi_L}\ar[d]^{(\frac{\partial}{\Omega})^{-j}\otimes e_j}\\
D(\Gamma_L,K)\otimes_L D_{cris,L}(W(\chi_{LT}^j))  \ar[r]^-{\mathfrak{M}\otimes \id} & \mathcal{O}_K \otimes_L D_{cris,L}(W(\chi_{LT}^j)) &  \mathcal{O}_K \otimes_L D_{cris,L}(W(\chi_{LT}^j)) \ar[l]_{1-\varphi_L\otimes \varphi_L}  .   }
\end{equation*}
In the latter we follow the (for $j>0$) abusive   notation $\partial^{-j}$  from \cite[Rem.\ 3.5.5.]{BF}.
\end{lemma}

\begin{proof}
For the upper diagram note that $\eta_0=0$ and $\left(\delta_g\cdot\eta(1,Z)\right)_{|Z=0}=1,$ from which the claim follows for Dirac distributions, whence in general. For the right square we observe that $\varphi_L(g(Z))_{|Z=0}=g(0).$ Regarding the lower diagram we use   \ref{lem:LT-twisting} and the relation $\partial_{\mathrm{inv}}\circ \varphi_L=\pi_L\varphi_L\circ \partial.$
\end{proof}

With this notation Berger's and Fourquaux' interpolation property reads as follows:

\begin{theorem}[{Berger-Fourquaux \cite[Thm.\ 3.5.3]{BF}}]
 Let $W$ be $L$-analytic and $h\geq 1$ such that $\mathrm{Fil}^{-h}D_{cris,L}(W)=D_{cris,L}(W)$. For any  $f\in \left(\mathcal{O}^{\psi=0}\otimes_L D_{cris,L}(W)\right)^{\Delta=0}$ and $y\in \left(\mathcal{O} \otimes_L D_{cris,L}(W)\right)^{\psi=\frac{q}{\pi_L}}$ with $f=(1-\varphi_L)y$ we have: If $h+j\geq 1$, then
\begin{align}\label{f:interjgeq1}\notag
  &h^1_{L_n,W(\chi_{LT}^j)}  (tw_{\chi_{LT}^j}(\Omega_{W,h}(f)))= \\
   & (-1)^{h+j-1}(h+j-1)!\left\{
                           \begin{array}{ll}
                             \exp_{L_n,W(\chi_{LT}^j)}\Big(q^{-n}\mathrm{Ev}_{W(\chi_{LT}^j),n}( \partial_{\mathrm{inv}}^{-j}y\otimes e_j)\Big) & \hbox{if $n\geq 1$;} \\
                   \exp_{L,W(\chi_{LT}^j)}\Big((1-q^{-1}\varphi_L^{-1})\mathrm{Ev}_{W(\chi_{LT}^j),0}   ( \partial_{\mathrm{inv}}^{-j}y\otimes e_j)\Big)         , & \hbox{if $n=0$.}
                           \end{array}
                         \right.
\end{align}
If $h+j\leq  0$, then
\begin{align}\label{f:interjleq0}\notag
\exp_{L_n,W(\chi_{LT}^j)}^* &\Big(  h^1_{L_n,W(\chi_{LT}^j)}  (tw_{\chi_{LT}^j}(\Omega_{W,h}(f)))\Big)= \\
   & \frac{1}{(-h-j)!}\left\{
                           \begin{array}{ll}
                            q^{-n}\mathrm{Ev}_{W(\chi_{LT}^j),n}(\partial_{\mathrm{inv}}^{-j}y\otimes e_j) & \hbox{if $n\geq 1$;} \\
                    (1-q^{-1}\varphi_L^{-1})\mathrm{Ev}_{W(\chi_{LT}^j),0}   (\partial_{\mathrm{inv}}^{-j}y\otimes e_j)         , & \hbox{if $n=0$.}
                           \end{array}
                         \right.
\end{align}
\end{theorem}
By abuse of notation we shall denote the base change $K\otimes_L-$ of the (dual) Bloch-Kato exponential map by the same expression. Using Lemma  \ref{lem:Ev}  we deduce the following interpolation property for the modified big exponential map with $x\in D(\Gamma_L,K)\otimes_L D_{cris,L}(W):$
If $j\geq 0$, then
\begin{align}\label{f:interjgeq1mod}\notag
  h^1_{L,W(\chi_{LT}^j)} & (tw_{\chi_{LT}^j}(\mathbf{\Omega}_{W,1}(x)))= \\
   & (-1)^{j}j! \Omega^{-j-1}
 \exp_{L,W(\chi_{LT}^j)}\Big((1-q^{-1}\varphi_L^{-1})(1-\varphi_L)^{-1}\left(x(\chi_{LT}^{-j})  \otimes e_j\right)\Big);
\end{align}
if $ j< 0$, then
\begin{align}\label{f:interjleq0mod}\notag
  h^1_{L,W(\chi_{LT}^j)}  & (tw_{\chi_{LT}^j}(\mathbf{\Omega}_{W,1}(f)))= \\
   & \frac{1}{(-1-j)!} \Omega^{-j-1}\log_{L,W(\chi_{LT}^j)}^*\Big(
 (1-q^{-1}\varphi_L^{-1})(1-\varphi_L)^{-1}\left(x(\chi_{LT}^{-j})  \otimes e_j\right)\Big),
\end{align}
assuming in both cases that  the operator $1-\varphi_L$ is invertible on $D_{cris,L}(W(\chi_{LT}^j))$ and for $j<0$ also that the operator $ 1-q^{-1}\varphi_L^{-1}$ is invertible on $D_{cris,L}(W(\chi_{LT}^j))$ (in order to grant the existence of $\log_{L,W(\chi_{LT}^j)}$).

Recall that the generalized Iwasawa cohomology of $T\in Rep_{o_L}(G_L)$ is defined by
\begin{equation*}
    H^*_{Iw}(L_\infty/L,T) := \varprojlim_K H^*(K,T)
\end{equation*}
where $K$ runs through the finite Galois extensions of $L$ contained in $L_\infty$ and
the transition maps in the projective system are the cohomological corestriction maps.
For $V:=T\otimes_{o_L}L\in Rep_L(G_L)$ we define
\begin{equation*}
    H^*_{Iw}(L_\infty/L,V) := H^*_{Iw}(L_\infty/L,T)\otimes_{o_L} L,
\end{equation*}
which is independent of the choice of $T.$ As usual we denote by $cor:H^*_{Iw}(L_\infty/L,T) \to H^*(L',T)$ the projection map and
analogously for rational coefficients. Similarly as in \eqref{f:Fontaine-SV} we have a
map
\begin{equation}\label{f:Fontaine-SVII}
pr_U:D(V(\tau^{-1}))^{\psi=1} \to h^1(K_{\psi,U'}(D(V(\tau^{-1}))^\Delta)[d-1] )\cong
H^1(L',V),\; m\mapsto [(\bar{m},0)].
\end{equation}
$m$ under the map $\check{M}\twoheadrightarrow\check{M}_\Delta\cong\check{M}^\Delta.$
Note that under the assumptions of Lemma \ref{notrivialquot} for $V( \tau^{-1})$ there
is a commutative diagram
\begin{equation}
\label{f:corpr}\xymatrix{
  H^1_{Iw}(L_\infty/L,T) \ar[d]_{cor} \ar[r]^{\cong } &
  D_{LT}(T(\tau^{-1}))^{\psi=1}\ar[d]_{pr_U} \ar@{^(->}[r]^{ } &
  D_{rig}^\dagger(V(\tau^{-1}))^{\psi=1} \ar[d]^{pr_U} \\
  H^1(L',V) \ar@{=}[r]^{ } & H^1(L',V) \ar@{->>}[r]^{ } & H^1_{/\dagger}(L',V), }
\end{equation}
where the right vertical map is induced by \eqref{f:Fontaine-SV}.
Indeed, for the commutativity of the left rectangle and the right rectangle we refer
the reader to \eqref{f:decent} and \eqref{f:PsiSlashdagger}, respectively. Let $y_{\chi_{LT}^{-j}}$ denote the image of $y$ under the map
 \[ H^1_{Iw}(L_\infty/L,T)\xrightarrow{ \cdot \otimes \eta^{\otimes -j}
 }H^1_{Iw}(L_\infty/L,T(\chi_{LT}^{-j}))\xrightarrow{\mathrm{cor}}H^1 (
 L,T(\chi_{LT}^{-j}))\to H^1 ( L,V(\chi_{LT}^{-j})). \]
The following result generalizes \cite[Thm.~A.2.3]{LVZ15} and \cite[Thm.\ B.5]{LZ}  from the cyclotomic  case.
\begin{theorem}\label{thm:adjointformula} Assume that $V^*(1)$ is $L$-analytic with $\mathrm{Fil}^{-1}D_{cris,L}(V^*(1))=D_{cris,L}(V^*(1))$ \linebreak and $D_{cris,L}(V^*(1))^{\varphi_L=\pi_L^{-1}}=D_{cris,L}(V^*(1))^{\varphi_L=1}=0$. Then it holds that for   $j\geq 0$
\begin{align*}
 \Omega^{j} \mathbf{L}_V(y)(\chi_{LT}^j)&=  j!\Big((1-\pi_L^{-1}\varphi_L^{-1})^{-1}(1-\frac{\pi_L}{q}\varphi_L)\widetilde{\exp}^*_{L,V(\chi_{LT}^{-j}),\id}(y_{\chi_{LT}^{-j}})\Big)\otimes e_{ j} \\
    &=  j!(1-\pi_L^{-1-j}\varphi_L^{-1})^{-1}(1-\frac{\pi_L^{j+1}}{q}\varphi_L)\Big(\widetilde{\exp}^*_{L,V(\chi_{LT}^{-j}),\id}(y_{\chi_{LT}^{-j}})\otimes e_{ j}\Big)
\end{align*}
and for $j\leq -1:$
\begin{align*}
 \Omega^{{violet}j} \mathbf{L}_V(y)(\chi_{LT}^j)&=  \frac{(-1)^{j}}{(-1-j)!}\Big((1-\pi_L^{-1}\varphi_L^{-1})^{-1}(1-\frac{\pi_L}{q}\varphi_L)\widetilde{\log}_{L,V(\chi_{LT}^{-j}),\id}(y_{\chi_{LT}^{-j}})\Big)\otimes e_{ j} \\
    &=  \frac{(-1)^{j}}{(-1-j)!}(1-\pi_L^{-1-j}\varphi_L^{-1})^{-1}(1-\frac{\pi_L^{j+1}}{q}\varphi_L)\Big(\widetilde{\log}_{L,V(\chi_{LT}^{-j}),\id}(y_{\chi_{LT}^{-j}})\otimes e_{ j}\Big),
\end{align*}
if the operators $1-\pi_L^{-1-j}\varphi_L^{-1}, 1-\frac{\pi_L^{j+1}}{q}\varphi_L $ or equivalently $1-\pi_L^{-1}\varphi_L^{-1}, 1-\frac{\pi_L}{q}\varphi_L$ are invertible on $D_{cris,L}(V(\tau^{-1}))$  and $D_{cris,L}(V(\tau^{-1}\chi_{LT}^j))$, respectively.
\end{theorem}
\begin{proof}
From the reciprocity formula in Corollary \ref{cor:adjointness} and Propositions \ref{prop:Tate} and \ref{prop:cris} we obtain for $x\in D(\Gamma_L,\mathbb{C}_p)\otimes_L D_{cris,L}(V^*(1)),$ $y\in D(V(\tau^{-1}))^{\psi_L=1}  $ and $j\geq 0$ using \eqref{f:corpr}
{\small
\begin{align*}
\notag[x(\chi_{LT}^{-j})&\otimes e_{j} , (-1)^{j}\mathbf{L}_V(y)(\chi_{LT}^{j})\otimes e_{-j}]_{cris}^{}\\
&={\Omega}
[x,\frac{\sigma_{-1}\mathbf{L}_V(y)}{{violet}\Omega}]^0(\chi_{LT}^{-j})\\
&={\Omega}{ \frac{q-1}{q}}\{\mathbf{\Omega}_{{V^*(1)},1}(x),y\}_{Iw}(\chi_{LT}^{-j})\\
&={ \Omega}<h^1_{L}\circ \notag tw_{\chi_{LT}^{j}}\left(\mathbf{\Omega}_{{V^*(1)},1}(x)\right),y_{\chi_{LT}^{-j}}>_{Tate}\\
\notag&={ \Omega}<(-1)^{j}j!\Omega^{-j-1}\exp_{L,V^*(1)(\chi_{LT}^j)}((1-q^{-1}\varphi_L^{-1})(1-\varphi_L)^{-1}(x(\chi_{LT}^{-j})\otimes e_ j),y_{\chi_{LT}^{-j}}>_{Tate}\\
\notag&=(-1)^{j}\Omega^{ -j}j![ (1-q^{-1}\varphi_L^{-1})(1-\varphi_L)^{-1}(x(\chi_{LT}^{-j})\otimes e_{j}),\widetilde{\exp}^*_{L,V(\chi_{LT}^{-j}),\id}(y_{\chi_{LT}^{-j}})]_{cris}
\\
\notag&=[ x(\chi_{LT}^{-j})\otimes e_{j},(-1)^{j}\Omega^{ -j}j!(1-\pi_L^{-1}\varphi_L^{-1})^{-1}(1-\frac{\pi_L}{q}\varphi_L)\widetilde{\exp}^*_{L,V(\chi_{LT}^{-j}),\id}(y_{\chi_{LT}^{-j}})]_{cris}
\end{align*}}
Here we used \eqref{f:interjgeq1mod} in the fourth equation for the interpolation property of $\mathbf{\Omega}_{{V^*(1)},1}. $  The fifth equation is the defining equation for the dual exponential map resulting from \eqref{f:dualexpId}. Furthermore, for the last equality we use that $ \pi_L^{-1}\varphi_L^{-1}$ is adjoint to $\varphi_L$ under the lower pairing. The claim follows since the evaluation map is surjective and $[\;,\;]_{cris}$ is non-degenerated.
\Footnote{Implicity we have used in the calculation that the operators $(1-q^{-1}\varphi_L^{-1}), \; (1-\varphi_L) $ are invertible, this can perhaps be avoided by a better arrangement!? }
Now assume that $j<0:$
{\small
\begin{align*}
\notag[x(\chi_{LT}^{-j})&\otimes e_{j} , (-1)^{j}\mathbf{L}_V(y)(\chi_{LT}^{j})\otimes e_{-j}]_{cris}^{}\\
&={ \Omega}
[x,\frac{\sigma_{-1}\mathbf{L}_V(y)}{{ \Omega}}]^0(\chi_{LT}^{-j})\\
&={ \Omega}{ \frac{q-1}{q}}\{\mathbf{\Omega}_{{V^*(1)},1}(x),y\}_{Iw}(\chi_{LT}^{-j})\\
&={ \Omega}<h^1_{L}\circ \notag tw_{\chi_{LT}^{j}}\left(\mathbf{\Omega}_{{V^*(1)},1}(x)\right),y_{\chi_{LT}^{-j}}>_{Tate}\\
\notag&={ \Omega}<\frac{1}{(-1-j)!} \Omega^{ -j-1} \log_{L,W(\chi_{LT}^j)}^*\Big(
 (1-q^{-1}\varphi_L^{-1})(1-\varphi_L)^{-1}\left(x(\chi_{LT}^{-j})  \otimes e_j\right)\Big),y_{\chi_{LT}^{-j}}>_{Tate}\\
\notag&=\frac{\Omega^{ -j}}{(-1-j)!}[ (1-q^{-1}\varphi_L^{-1})(1-\varphi_L)^{-1}(x(\chi_{LT}^{-j})\otimes e_{j}),\widetilde{\log}_{L,V(\chi_{LT}^{-j}),\id}(y_{\chi_{LT}^{-j}})]_{cris}
\\
\notag&=[ x(\chi_{LT}^{-j})\otimes e_{j},\frac{\Omega^{ -j}}{(-1-j)!}(1-\pi_L^{-1}\varphi_L^{-1})^{-1}(1-\frac{\pi_L}{q}\varphi_L)\widetilde{\log}_{L,V(\chi_{LT}^{-j}),\id}(y_{\chi_{LT}^{-j}})]_{cris}
\end{align*}}
\end{proof}

Now consider $V=L(\tau\chi_{LT})$ and $W=V(\chi_{LT}).$ Then the latter satisfies the condition of the Theorem and using Propositon \ref{prop:regulatortwisting} and Lemma \ref{lem:nablatwisting} one easily derives the following interpolation property concerning the former for $y=-\kappa(u)$, $u\in U$ d for all $r\geq 1:$
\begin{align*}
   \mathbf{L}_V(y)(\chi_{LT}^r)&=  \frac{r!}{\Omega^{ r}}\Big((1-\pi_L^{-1}\varphi_L^{-1})^{-1}(1-\frac{\pi_L}{q}\varphi_L)\widetilde{\exp}^*_{L,V(\chi_{LT}^{-r}),\id}(y_{\chi_{LT}^{-r}})\Big)\otimes e_{ r}\\
   &=\frac{r!}{\Omega^{ r}}(1-\pi_L^{-r})^{-1}(1-\frac{\pi_L^r}{q})\widetilde{\exp}^*_{L,V(\chi_{LT}^{-r}),\id}(y_{\chi_{LT}^{-r}})\otimes e_{ r}.\\
\end{align*}
On the other hand we have $ \mathbf{L}_V(y)\otimes\mathbf{d}_1=\cL_V(y) $ and hence by the claim concerning \eqref{f:claim}
\begin{align*}
  \mathbf{L}_V(y)(\chi_{LT}^r)\otimes\mathbf{d}_1\otimes\eta^{-1}\otimes t_{LT}&=\cL_V(y)(\chi_{LT}^r)\otimes\eta^{-1}\otimes t_{LT}\\
  &=\cL(-\kappa(u)\otimes \eta^{-1})(\chi_{LT}^r),
\end{align*}
whence
\[\cL(-\kappa(u)\otimes \eta^{-1})(\chi_{LT}^r)\otimes e_{1-r}=\frac{r!}{\Omega^{ r}}(1-\pi_L^{-r})^{-1}(1-\frac{\pi_L^r}{q}){\exp}^*_{L,V(\chi_{LT}^{-r}),\id}(y_{\chi_{LT}^{-r}}).\]
This is \eqref{f:CWinterpolation}, i.e.,  together with \eqref{f:decentKato} we have just obtained a new proof of Kato's reciprocity law \ref{thm:Kato} and we may consider Theorem \ref{thm:adjointformula} as a vast generalisation of it.

\newpage

\appendix

\section{Cup products and local Tate duality}\label{sec:cupprod}
The aim of this subsection is to discuss cup products and to prove Prop.\  \ref{prop:Koszul-product}. We fix some open subgroup $U\subseteq \Gamma_L$ and let $L'=L_\infty^U.$  Note that we obtain a decomposition $U\cong \Delta\times U'$ with
  a subgroup $U'\cong\mathbb{Z}_p^d$ of $U$ and  $\Delta$  the torsion subgroup of $U$.

\begin{lemma}\label{lem:HS}
Let $M_0$ be a complete linearly topologised $o_L$-module  with continuous $U$-action and with a continuous $U$-equivariant endomorphism $f$. Then there is a canonical quasi-isomorphism \[\mathcal{T}_{f,U}(M_0)[\frac{1}{\pi_L}]\simeq K^\bullet_{f,U'}(M_0[\frac{1}{\pi_L}]^\Delta).\] If $M_0$ is an $L$-vector space, the inversion of $\pi_L$ can be omitted on both sides.
\end{lemma}

\begin{proof}
  Let $\cC^\bullet_n(U,M_0) \subseteq \cC^\bullet(U,M_0)$ denote the subcomplex of normalized cochains. Since $\Delta$ is finite, \cite[Thm.\ 3.7.6]{Th} gives a canonical quasi-isomorphism:
\begin{equation*}
  \cC^\bullet_n(U,M_0) = \cC^\bullet_n(\Delta \times U',M_0) \xrightarrow{\simeq} \cC^\bullet_n(\Delta,\cC^\bullet_n(U',M_0)).
\end{equation*}
Here we understand the above objects in the sense of hypercohomology as total complexes of the obvious double complexes involved. After inverting $\pi_L$ we may compute the right hand side further as
\begin{equation*}
  \cC^\bullet_n(\Delta,\cC^\bullet_n(U',M_0))[\frac{1}{\pi_L}] = \cC^\bullet_n(\Delta,\cC^\bullet_n(U',M_0)[\frac{1}{\pi_L}]) \xleftarrow{\simeq} \cC^\bullet_n(U',M_0)[\frac{1}{\pi_L}]^\Delta = \cC^\bullet_n(U',M_0^\Delta)[\frac{1}{\pi_L}] \ .
\end{equation*}
Here the middle quasi-isomorphism comes from the fact that a finite group has no cohomology in characteristic zero. The right hand equality is due to the fact that $\Delta$ acts on the cochains through its action on $M_0$. Altogether we obtain a natural quasi-isomorphism
\begin{equation*}
  \cC^\bullet_n(U,M_0){ [\frac{1}{\pi_L}] } \cong \cC^\bullet_n(U',M_0^\Delta)[\frac{1}{\pi_L}] \ .
\end{equation*}
By using \cite[Prop.\ 3.3.3]{Th} we may replace the normalized cochains again by general cochains obtaining the left hand quasi-isomorphism in
\begin{equation*}
  \cC^\bullet(U,M_0){ [\frac{1}{\pi_L}] } \cong \cC^\bullet(U',M_0^\Delta)[\frac{1}{\pi_L}] \cong K^\bullet_{{  U'}}(M_0^\Delta)[\frac{1}{\pi_L}] = K^\bullet_{{ U'}}(M_0[\frac{1}{\pi_L}]^\Delta) \ .
\end{equation*}
The middle quasi-isomorphism is \eqref{f:KoszulCtsCoh}. The claim follows by taking mapping fibres of the attached map $f-1$ of complexes.
%
\end{proof}

\begin{proposition}\label{prop:finitedim} Let   $M  $  be a $\varphi_L$-module over  $\cR=\cR_K$ (cf.\ \ref{def:phimodule}) and $c\in K^\times.$ Then   $M/(\psi-c)(M)$ is finite-dimensional over $K$.
\end{proposition}

\begin{proof}(The proof follows closely the proof of \cite[Prop.\ 3.3.2]{KPX} in the cyclotomic situation) We set $\psi_c:=c^{-1}\psi$ and show that $M/(\psi_c-1)(M)$ is finite-dimensional over $K$.

Choose a model $M^{[r_0,1)}$ of $M$ with $1>r_0>p^{\frac{-1}{(q-1)e}}$ and put $r=r_0^{\frac{1}{q^2}}.$ Recall that for all $1>s\geq r$ we have maps $M^{[s,1)}\xrightarrow{\psi_c-1}M^{[s,1)}$ (where strictly speaking we mean $\psi_c$ followed by the  corresponding restriction). We first show that it suffices to prove that $\mathrm{coker}\left( M^{[r,1)}\xrightarrow{\psi_c-1}M^{[r,1)}  \right)$ has finite dimension over $K$. Indeed, given any $m\in M$ we have $m\in M^{[s,1)}$ for some $1>s\geq r$. Then there exists $k\geq 0$ such that $r\geq s^{q^k}\geq r_0$, whence $\psi_c^k(m)$ belongs to $M^{[r,1)}$ and represents the same class in $M/(\psi_c-1)(M)$ as $m.$

Choose a basis $\mathbf{e}_1',\ldots,\mathbf{e}_n'$ of $M^{[r_0,1)}$  and take $\mathbf{e}_i:= \varphi(\mathbf{e}_i')\in M^{[r^q,1)}$; by the $\varphi$-module property the latter elements also form a basis of $M^{[r^q,1)}$ . Note that by base change these two bases also give rise to  bases in $M^{[s,1)}$ for $1>s\geq r^q$.  Thus we find a matrix $F'$ with entries in $\cR^{[r^q,1)}$ such that $\mathbf{e}_j=\sum_iF_{ij}'\mathbf{e}_i'$ and we put $F=\varphi(F')$ with entries in $\cR^{[r,1)}$, i.e., $\varphi(\mathbf{e}_j)=\sum F_{ij}\mathbf{e}_i$. Similarly let $G$ be a matrix with values in $\cR^{[r^q,1)}\subseteq \cR^{[r,1)}$ such that $\mathbf{e}_j'=\sum_i G_{ij}\mathbf{e}_i$ and hence $\mathbf{e}_j=\varphi\left(\sum_i G_{ij}\mathbf{e}_i\right).$

We identify $ M^{[r,1)} $ with $\left(\cR^{[r,1)}\right)^n$ by sending $(\lambda_i)_i$ to $\sum_i\lambda_i \mathbf{e}_i$ and endow it for each $r\leq s<1$ with the  norm given by $\max_i|\lambda_i |_s.$ Note that then the ''semi-linear'' map  $\psi_c$ (followed by the  corresponding restriction) on $\left(\cR^{[r,1)}\right)^n$ is given by the matrix  $G$ as follows from the projection formula \eqref{f:projectionrs}:
\begin{equation*}
  \psi_c(\sum_j \lambda_j \mathbf{e}_j )  = \sum_j\psi_c(\lambda_j\varphi(\sum_i G_{ij}\mathbf{e}_i) )
     = \sum_{i,j}\psi_c(\lambda_j)G_{ij} \mathbf{e}_i.
\end{equation*}
  Moreover, the restriction of $\varphi: M^{[r,1)} \to M^{[r^\frac{1}{q},1)}$ to $\sum_I {\mathcal{O}_K(\mathbf{B})} e_i$ becomes the semi-linear map $\left({\mathcal{O}_K(\mathbf{B})}\right)^n \to \left(\cR^{[r,1)}\right)^n$ given by the matrix  $F$.

Consider, for $I$ any subset of the reals $\mathbf{R}$, the $K$-linear map $P_I:\cR^{[r,1)}\to \cR^{[r,1)}, \sum a_iZ^i\mapsto \sum_{i\in \mathbf{Z}\cap I}a_iZ^i$.
We then introduce $K$-linear operators $P_{I}$   and $Q_k$, $k\geq 0,$ on $ M^{[r,1)} $ by
\begin{align*}
  P_{I}((\lambda)_i) &:=(P_{I}(\lambda_i))_i, \\
  Q_k&:= P_{(-\infty,-k)}-\frac{c\pi_L}{q}\varphi\circ P_{(k,\infty)}, i.e.,\\
 Q_k((\lambda)_i) & = (P_{(-\infty,-k)}(\lambda)_i)_i -\frac{c\pi_L}{q} F\cdot(\varphi(P_{(k,\infty)}(\lambda_i)))_i,
\end{align*}
because $P_{(k,\infty)} $ factorises through ${\mathcal{O}_K(\mathbf{B})}.$
Then the $K$-linear operator $\Psi_k:=\id- P_{[-k,k]} +(\psi_c-1)Q_k$ of  $ M^{[r,1)} $ satisfies
\begin{align*}
  \Psi_k & =\psi_c\circ P_{(-\infty,-k)} + \frac{c\pi_L}{q}\varphi\circ P_{(k,\infty)}, i.e., \\
 \Psi_k((\lambda_i)_i)&=G\cdot(\psi_c(P_{(-\infty,-k)}(\lambda_i)))_i+F\cdot(\frac{c\pi_L}{q}\varphi(P_{(k,\infty)}(\lambda_i)))_i,
\end{align*}
whence its operator norm satisfies
\[\|\Psi_k\|_s\leq \max \{\|G\|_s \|\psi_c\circ P_{(-\infty,-k)}\|_s, \frac{c\pi_L}{q} \|F\|_s \|\varphi\circ P_{(k,\infty)}\|_s\}  .\]
 It is easy to check that, for $1>s>q^{\frac{-1}{q-1}},$ we have $\|\varphi\circ P_{(k,\infty)}\|_s\leq |Z|_s^{(q-1)k}=s^{(q-1)k}$ (using the norm relation after \eqref{f:isometry})
 and $\|\psi_c\circ P_{(-\infty,-k)}\|_s \leq C_s s^{k(1-q^{-1})}$ for some constant $C_s > 0.$  E.g. for the latter we have for $\lambda=\sum_i a_i Z^i\in \cR^{[r,1)} $
\begin{align*}
  |\sum_{i<-k}\psi_c( a_iZ^i)|_s & \leq \sup_{i<-k} |a_i\|\psi_c(Z^i)|_s \\
    & \leq \sup_{i<-k} |a_i| C_s s^{\frac{i}{q}}\\
  & \leq  C_s \sup_{i<-k} |a_i\|Z|_s^{ i}  s^{i({q}^{-1}-1)}\\
  & \leq C_s |\lambda|_s s^{-k({q}^{-1}-1)},
\end{align*}
where we use that by continuity of $\psi_c$ there exists $C_s$ such that
\[|\psi_c(Z^i)|_s \leq C_s |Z^i|_{s^{\frac{1}{q}}} = C_s s^{\frac{i}{q}} .\]

Thus we may and do choose $k$ sufficiently big such that $\|\Psi_k\|_r \leq \frac{1}{2}.$ Given $m_0\in M^{[r,1)}$ we define inductively  $m_{i+1}:=\Psi_k(m_i).$ This sequence obviously converges to zero with respect to the $r$-Gauss-norm. We shall show below that also for all $s\in(r^\frac{1}{q},1)$  the series $(m_i)_i$ tends to zero with respect to the Gauss norm $|\;|_s$, i.e., by cofinality the sum $m:=\sum_{i\geq 0} m_i$ converges in $M^{[r,1)}$ for the Frech\'{e}t-topology and satisfies
\[m-m_0=m-P_{[-k,k]}(m)+(\psi_c-1)Q_k(m),\]
i.e. $ P_{[-k,k]}(m)$ represents the same class as $m_0$ in $M^{[r,1)}/(\psi_c-1)(M)$. Since the image of $ P_{[-k,k]}$ has finite dimension, the proposition follows, once we have shown the following

\noindent
{\it Claim:} For all $s\in(r^\frac{1}{q},1)$ we have
\[|\Psi_k(m)|_s\leq \max \{\frac{1}{2}|m|_s, C_s\|G\|_s\left(\frac{s^{\frac{1}{q}}}{r}\right)^{-k} |m|_r , |\frac{c\pi_L}{q}|\|F\|_s\left(\frac{s^{q}}{r}\right)^{k'} |m|_r \}.\]
Indeed, we fix such $s$ and may choose $k'\geq k$ such that $\|\Psi_{k'}\|_s\leq\frac{1}{2}.$ Then $\Psi_k=\Psi_{k'} -\psi_c\circ P_{[-k',-k)}-\frac{c\pi_L}{q} \varphi \circ P_{(k,k']},$ whence the claim as for $\lambda\in \cR^{[r,1)} $
\begin{align*}
  |\psi_c\circ P_{[-k',-k)}(\lambda)|_s & \leq C_s \left(\frac{s^{\frac{1}{q}}}{r}\right)^{-k} |\lambda|_r \\
   | \varphi \circ P_{(k,k']}(\lambda)|_s & \leq \left(\frac{s^{q}}{r}\right)^{k'} |\lambda|_r \\
\end{align*}
by similar estimations  as above.
\end{proof}

\begin{remark}
This result answers the expectation from  \cite[Remark 2.3.7.]{BF} positively.
\end{remark}

\begin{corollary}\label{cor:finitedim}Let $V^*(1)$ be $L$-analytic and $M:= {D}^\dagger_{rig}(V^*(1))  $.
\begin{enumerate}
\item  The cohomology group $h^2(K^\bullet_{\psi,U'}(({M} )^\Delta)[d-1])$ is finite dimensional over $L$.
\item We have isomorphisms
\begin{align*}
  h^1(K_{\psi,U'}^\bullet({D}^\dagger_{rig}(V(\tau^{-1}))^\Delta)[d-1]  )^{*} &\cong h^{1}(K^\bullet_{\varphi,U'}({M} ^\Delta)) \\
    & \cong H_\dagger^{1}(L',V^*(1)),
\end{align*}
and
\begin{align*}
  h^2(K_{\psi,U'}^\bullet({D}^\dagger_{rig}(V(\tau^{-1}))^\Delta)[d-1]  )^{*} &\cong h^{0}(K^\bullet_{\varphi,U'}({M} ^\Delta)) \\
    & =(V^*(1))^{G_{L'}}.
\end{align*}
\end{enumerate}
\end{corollary}

\begin{proof}
 (i) Since $h^2(K^\bullet_{\psi,U'}(({M} )^\Delta)[d-1])$ is a quotient of  $(M/(\psi-1)(M))^{\Delta} $ by \eqref{f:duality} this follows from the Proposition. (ii)  We are in the situation of Remark \ref{rem:nonstrict} (i) with regard to $\cC=K_{\psi,U'}^\bullet({D}^\dagger_{rig}(V(\tau^{-1}))^\Delta)[d-1]$ and $i=2,3$ in the notation of the remark. Indeed,  for $h^3(\cC)=\cC^4=0$ by construction and $\cC^3=0$ as well as $h^2(\cC)$ is finite by (i). Hence the first isomorphism follows in both cases from  \eqref{f:KpsiU} using the reflexivity of $M.$ The second isomorphisms arise by Lemma \ref{lem:HS} together with \eqref{f:cohLphi1Ver} and  \eqref{f:cohLphiOVer}, respectively.
 \end{proof}

{

We quickly discuss the analogues of some results of \S I.6 in \cite{ChCo2}. First we remind the reader of  the definition of  $\tilde{\mathbf{A}}:=W(\mathbb{C}_p^\flat)_L$,
\[
\tilde{\mathbf{A}}^\dagger:=\{x=\sum_{n\geq 0}\pi_L^n[x_n]\in\tilde{\mathbf{A}}: |\pi_L^{n}\|x_n|_\flat^r \xrightarrow{n \rightarrow \infty} 0\mbox{ for some } r>0\}
\]
 ${\bf A}^\dagger:=\tilde{\bf A}^{\dagger}\cap {\bf A}$\footnote{\label{foot:leq1}
In the literature one also finds the subring $ \tilde{\mathbf{A}}^\dagger_{\leq 1}:=\bigcup_{r>0}W^r_{\leq 1}(\mathbb{C}_p^\flat)_L$ of $ \tilde{\mathbf{A}}^\dagger $ where $W^r_{\leq 1}(\mathbb{C}_p^\flat)_L=\{ x\in W^r(\mathbb{C}_p^\flat)_L|\; |x|_{r}\leq 1 \}$ consists of those $x\in \tilde{\mathbf{A}}$ such that   $|\pi_L^{n}\|x_n|_\flat^r \xrightarrow{n \rightarrow \infty} 0$ and $|\pi_L^{n}\|x_n|_\flat^r\leq 1$ for all $n.$ Denoting by $\tilde{\mathbf{A}}^{\dagger,s}_{St}$ the ring defined in \cite[Def.\ 3.4]{Ste} we have the equality $W^r_{\leq 1}(\mathbb{C}_p^\flat)_L = \tilde{\bf A}^{\dagger,\frac{q-1}{qr}}_{St}$ corresponding to $\tilde{\bf A}^{\dagger}_{(\frac{q-1}{qr})^-}$ in the notation of  \cite[\S II.1]{ChCo} for $q=p$. For these relations use the following  normalisations compatible with \cite{Ste}:
 $|\pi_L|=\frac{1}{q}$, $v_{\mathbb{C}_p^\flat}(\omega)=\frac{q}{q-1},$ $v_{\pi_L} (\pi_{L})=1,$ $|x|_\flat=q^{-v_{\mathbb{C}_p^\flat}}, $ $|\omega|_\flat=q^{-\frac{q}{q-1}}=|\pi_L|^{\frac{q}{q-1}},$  where $\omega=\omega_{LT}\mod \pi_L$ as in section \ref{sec:KR}.
Furthermore, $|x|_{r}=|\pi_L|^{V(x,\frac{q-1}{rq})}$ and $|\omega_{LT}|_{r}=|\pi_L|^{\frac{rg}{q-1}}=|{\omega}|^r_\flat,$ where $V(x,r):=\inf_k\left(v_{\mathbb{C}_p^\flat}(x_k)\frac{q-1}{rq}+k\right)$ for $x=\sum_{k\geq 0} \pi_L^k[x_k]\in \tilde{\bf A}.$ For $x\in \tilde{\mathbf{A}}^\dagger$ we have $V(x,r)=\frac{q-1}{rq}V_{St}(x,r)$ where $V_{St}(x,r) $ uses the notation in \cite{Ste}.
Note also that $\omega_{LT}^{-1}$ is contained in $W^{\frac{q-1}{q}}_{\leq 1}(\widehat{L_\infty}^\flat)_L$ by  \cite[Lem.\ 3.10]{Ste} (in analogy with \cite[Cor. II.1.5]{ChCo}).} and ${\bf A}^\dagger_L:=({\bf A}^\dagger)^{H_L}.$

\begin{remark}\label{rem:overconvergent}
 There  is also the following more concrete description for ${\bf A}^\dagger_L$ in terms of Laurent series in $\omega_{LT}:$
\begin{align*}
{\bf A}^\dagger_L=\{F(\omega_{LT})\in \mathbf{A}_L|& F(Z)\mbox{ converges on } \rho\leq |Z|<1 \mbox{ for some } \rho\in (0,1)\}\subseteq \mathbf{A}_L.
\end{align*}
Indeed this follows from  the analogue of \cite[Lem.\ II.2.2]{ChCo} upon noting that the latter holds with and without the integrality condition:  ''$rv_p(a_n)+n\geq 0$ for all $n\in\mathbb{Z}$'' (for $r\in \overline{\mathbf{R}}\setminus \mathbf{R}$) in the notation of that article.
\footnote{\label{foot:leq1imperfect} This description does not require any completeness property! A similar result holds for $\mathbf{A}_{\leq 1,L}^\dagger$ when requiring for the Laurent series in $\omega_{LT}$ in addition that $F(Z)$ takes values on $ \rho\leq |Z|<1$ of norm at most $1.$ More precisely, for $r<1$ (or equivalently $s(r):=\frac{q-1}{rq}>\frac{q-1}{q}$) $W^r(\mathbf{C}_p^\flat)_L$ and $W^r_{\leq1}(\mathbf{C}_p^\flat)_L$ correspond to  $\{F(\omega_{LT})\in\mathbf{A}_L| \; F(Z) \mbox{ converges on } |\omega|^r\leq |Z|<1\}$ and $\{F(\omega_{LT})\in\mathbf{A}_L| \; F(Z) \mbox{ converges on } |\omega|^r=q^{-\frac{1}{s(r)}}\leq |Z|<1 \mbox{ with values } |F(z)|\leq1 \}$, respectively. The latter condition on the values can also be rephrased as $s(r)v_{\pi_L}(a_m)+m\geq 0$ for all $m\in \z$ corresponding to $V(x,s(r))\geq 0$ on the Witt vector side if $F(Z)=\sum_m a_mZ^m\in\mathbf{A}_L.$
}
In particular we obtain canonical embeddings ${\bf A}^\dagger_L\subseteq {\bf B}^\dagger_L \hookrightarrow \cR_L$ of rings.
\end{remark}
Now consider the subring $A=\mathbf{A}_L^+[[\frac{\pi_L}{Z^{q-1}}]]=\{x=\sum_k a_k Z^k\in \mathbf{A}_L| v_{\pi_L}(a_k)\geq-\frac{k}{q-1} \}\subseteq \mathbf{A}_L$. For $x\in \mathbf{A}_L$ and each inter $n\geq 0,$ we define $w_n(x)$ to be the smallest integer $k\geq 0$ such that $x\in Z^{-k}A+\pi_L^{n+1}\mathbf{A}_L.$ It satisfies $w_n(x+y)\leq \max\{w_n(x),w_n(y)\}$ and $w_n(xy)\leq w_n(x)+w_n(y) $ (since $A$ is a ring) and $w_n(\varphi(x))\leq qw_n(x)$ (use that $\frac{\varphi(Z)}{Z^q}\in A^\times,$ whence $\varphi(Z^{-k})A=Z^{-qk}A$). Set for $n\geq 2, m\geq 0$ the integers $r(n):_=(q-1)q^{n-1},$ $l(m,n)=m(q-1)(q^{n-1}-1)=m(r(n)-(q-1))$ and define $\mathbf{A}_L^{\dagger,n}=\{x=\sum_k a_k Z^k\in \mathbf{A}_L|v_{\pi_L}(a_k)+\frac{k}{r(n)} \rightarrow\infty \mbox{ for } k\mapsto -\infty\}$. Then, by Remark \ref{rem:overconvergent} and thefootnote  \ref{foot:leq1}, \ref{foot:leq1imperfect}   we obtain that $\mathbf{A}_L^{\dagger}= \bigcup_{n\geq 2} \mathbf{A}_L^{\dagger,n}.$
\begin{lemma}\label{lem:criterion} Let $x=\sum_k a_k Z^k\in \mathbf{A}_L$ and $l\geq 0, n\geq 2.$ Then
\begin{enumerate}
\item we have \begin{align}
 w_m(x)\leq l \Leftrightarrow v_{\pi_L}(a_k)\geq \min\{m+1,-\frac{k+l}{q-1}\} \mbox{ for } k<-l.
\end{align}
\item $x\in \mathbf{A}_L^{\dagger,n}$ if and only if $w_m(x)-l(m,n)$ goes to $-\infty$ when $m$ runs to $\infty.$
\end{enumerate}
\end{lemma}
Item (ii) of the Lemma is an analogue of \cite[Prop.\ III 2.1 (ii)]{ChCo2} for $\mathbf{A}_L^{\dagger,n}$ instead of $ \mathbf{A}_{L,\leq 1}^{\dagger,n}=\{x=\sum_k a_k Z^k\in \mathbf{A}_L^{\dagger,n}|v_{\pi_L}(a_k)+\frac{k}{r(n)}\geq 0 \mbox{ for all } k\leq 0\}$.
\begin{proof}
(i) follows from the fact that $x\in Z^{-l}A $ if and only if $v_{\pi_L}(a_k)\geq  -\frac{k+l}{q-1}  \mbox{ for } k<-l.$ (ii) Let $M,N=M(q-1)\gg0$ be   arbitrary huge integers  and assume first that  $x\in \mathbf{A}_L^{\dagger,n}$. Then
\begin{equation}\label{f:wn}
   w_m(x)-l(m,n)\leq -N
 \end{equation}   is equivalent to
\begin{equation}\label{f:wnequ}v_{\pi_L}(a_k)\geq \min\{m+1,-\frac{k+l(m,n)-N}{q-1}\} \mbox{ for } k<-l(m,n)+N.
\end{equation}
 by (i). To verify this relation for $m$ sufficiently huge, we choose a $k_0\in\z$ such that $v_{\pi_L}(a_k)+\frac{k}{r(n)}\geq N\geq 0 $ for all $k\leq k_0.$ Now choose $m_0$ with $-l(m_0,n)<k_0$ and fix $m\geq m_0$. For $\frac{-k}{r(n)}>m$ we obtain $v_{\pi_L}(a_k)\geq m+1$, because $k<k_0$ holds. For $k$ with
\begin{equation}
\label{f:krn}k\geq -r(n)m \Leftrightarrow\frac{k}{r(n)}-\frac{k+l(m,n)}{q-1}\leq 0
\end{equation}
we obtain $v_{\pi_L}(a_k)\geq -\frac{k+l(m,n)-N}{q-1} $. Thus the above relation holds true. \\
Vice versa choose $m_0 $ such that \eqref{f:wn} holds for all $m\geq m_0$,  and fix \[k\leq k_0:=-r(n)\max\{Mq^{n-1},m_0\}.\]   Let $m_1$ be the unique integer satisfying $r(n)M-k\geq r(n)m_1\geq r(n)M-k-r(n).$ In particular, we have $m_1+1+\frac{k}{r(n)}\geq M$ and $k\geq -r(n)m_1,$ which implies
$-\frac{k+l(m_1,n)-N}{q-1}+\frac{k}{r(n)}\geq M$   by \eqref{f:krn}. Moreover, it holds $m_1\geq m_0$ and $k<-l(m_1,n)+N$ (using $k\leq r(n)M-r(n)m_1=-l(m_1,n)+q^{n-1}N-(q-1)m_1$ and $m_1>(q^{n-1}-1)M$ by our assumption on $k$). Hence we can apply \eqref{f:wnequ} to conclude $v_{\pi_L}(a_k)+\frac{k}{r(n)}\geq M$ as desired.
\end{proof}

 The analogue of Lemma 6.2 in (loc.\ cit.) holds by the discussion in \cite{SV15} after Remark 2.1. This can be used to show the analogue of Corollary 6.3, viz $w_n(\psi(x))\leq 1+ \frac{w_n(x)}{q}.$ Now fix a basis $(e_1,\ldots, e_d)$ of $ D(T) $ over $\mathbf{A}_L$ and denote by  $\Phi=(a_{ij})$  the matrix defined by $e_j=\sum_{i=1}^{d} a_{ij}\varphi(e_i).$   The proof of Lemma 6.4 then carries over to show that for $x=\psi(y)-y$ with $x,y\in D(T)$ we have
\begin{equation}\label{f:wpsi}
  w_n(y)\leq \max\{w_n(x),\frac{q}{q-1}\left(w_n(\Phi)+1\right)\},
\end{equation}
where $w_n(\Phi)=\max_{ij} w_n(a_{ij})$ and $w_n(a)=\max_i w_n(a_i)$ for $a=\sum_{i=1}^d a_i\varphi(e_i)$ with $a_i\in \mathbf{A}_L.$

\begin{lemma} Let $T\in \Rep_{o_L}(G_L)$ such that $T$ is free over $o_L$ and $V=T\otimes_{o_L}L$ is overconvergent.
Then the canonical map $D^\dagger(T)\to D(T)$ induces an isomorphism $ D^\dagger(T)/(\psi-1)(D^\dagger(T))\cong D(T)(\psi-1)(D(T)).$
\end{lemma}

\begin{proof}  We follow closely the proof of \cite[Lem.\ 2.6]{Li}, but note that he claims the statement for $D^\dagger_{\leq 1}(T). $\Footnote{{ Problem: Firstly, Liu claims in  \cite[Lem.\ 3.6]{Li}  that in the situation of the Lemma $D^\dagger_{\leq 1}(T) $ is free over $\mathbf{A}_{L,\leq 1}^\dagger $, for which I neither find a reference nor proof by my own. Secondly, \cite[Lem.\ I.6.4]{ChCo2} uses a basis of the form $\varphi(e_i)$ which is also needed to apply the projection formula. I do not understand how this translates to the setting in his proof, where he just chooses any basis of $D^\dagger_{\leq 1}(T)  $. }}
Choose a basis $e_1,\ldots, e_d$ of the ${\bf A}_L^\dagger$-module $D^\dagger(T)$, which is free because ${\bf A}_L^\dagger$ is a henselian discrete valuation ring with respect to the uniformiser $\pi_L,$ compare with \cite[Def.\ 2.1.4]{KedNew}. Since
  $V$ is overconvergent it is also a basis   of $D(T).$
  Due to \'{e}taleness and since $(\mathbf{A}_L^\dagger) \cap \mathbf{A}_L^{\times} = (\mathbf{A}_L^\dagger)^{\times}$ also $\varphi(e_1),\ldots,\varphi(e_d)$   is a basis of all these modules. Given $x=\psi(y)-y$ with $x\in D^\dagger(T) $  and $y\in D(T)$
  there is an $m>0$ such that  all $x_i, a_{ij}$ lie in $\mathbf{A}_{L}^{\dagger,m}$ for some $m$. Since $q\geq 2$ it follows from the   criterion in Lemma \ref{lem:criterion} (ii) combined with \eqref{f:wpsi} that all $y_i$ belong to $\mathbf{A}_{L}^{\dagger,m+1}$, whence
  $y\in D^\dagger(T) $. This shows injectivity. In order to show surjectivity we apply Nakayamas Lemma with regard to the ring $o_L$ upon recalling that $D(T)/(\psi-1)$  is of finite type over it. Indeed, by left exactness of $D^\dagger$ we obtain $D^\dagger(T)/\pi_LD^\dagger(T)\subseteq D^\dagger(T/\pi_LT)=D (T/\pi_LT).$ Since these are vector spaces over $\mathbf{E}_L$ of the same dimension, they are equal, whence
\[(D^\dagger(T)/(\psi-1))/(\pi_L)=(D^\dagger(T)/(\pi_L))/(\psi-1)=(D (T)/(\pi_L))/(\psi-1)=(D (T)/(\psi-1))/(\pi_L).\]
\end{proof}

%

\begin{corollary}\label{cor:Psidaggerquasi}\label{cor:herr}
Under the assumption of Lemma \ref{notrivialquot} for $V(\tau^{-1}) $, the inclusion of complexes
\[K_{\psi,U'}^\bullet({D}^\dagger(V(\tau^{-1}))^\Delta)\subseteq K_{\psi,U'}^\bullet({D}(V(\tau^{-1}))^\Delta)\] is a quasi-isomorphism.
\end{corollary}

\begin{proof}
Forming Koszul complexes with regard to $U'$ we obtain  the following diagram of (double) complexes with exact columns
\[\xymatrix{
0\ar[d] & 0\ar[d]\\
  K^\bullet({D}^\dagger(V(\tau^{-1}))^\Delta) \ar[d]_{} \ar[r]^{\psi-1} &  K^\bullet({D}^\dagger(V(\tau^{-1}))^\Delta)\ar[d]^{} \\
   K^\bullet( {D}(V(\tau^{-1}))^\Delta)\ar[d]_{} \ar[r]^{\psi-1} &  K^\bullet({D}(V(\tau^{-1}))^\Delta) \ar[d]^{} \\
  K^\bullet(({D}(V(\tau^{-1}))/{D}^\dagger(V(\tau^{-1})))^\Delta)\ar[d] \ar[r]^{\psi-1}_{\cong} & K^\bullet( ({D}(V(\tau^{-1}))/{D}^\dagger(V(\tau^{-1})))^\Delta)\ar[d] \\
  0& 0  }\]
in which the bottom line is an isomorphism of complexes because  under the assumptions $ \psi-1$ induces an automorphism of $ {D}(V(\tau^{-1}))/{D}^\dagger(V(\tau^{-1})) $  and as the action of $\Delta$ commutes with $\psi$. Hence,
going over to total complexes gives  an exact sequence
\[ 0\to K_{\psi,U'}^\bullet({D}^\dagger(V(\tau^{-1}))^\Delta) \to K_{\psi,U'}^\bullet({D}(V(\tau^{-1}))^\Delta)\to K_{\psi,U'}^\bullet(({D}(V(\tau^{-1}))/{D}^\dagger(V(\tau^{-1})))^\Delta)\to 0, \]
in which $K_{\psi,U'}^\bullet(({D}(V(\tau^{-1}))/{D}^\dagger(V(\tau^{-1})))^\Delta) $ is acyclic, whence  the statement follows.
\end{proof}

\begin{remark}
Instead of using Lemma \ref{notrivialquot} (for crystalline, analytic representations) one can probably show  by the same techniques as in \cite[Prop.\ III.3.2(ii)]{ChCo2} that for any overconvergent representation $V$ we have $D^\dagger(V)^{\psi=1}=D(V)^{\psi=1}.$
\end{remark}
The interest in the following diagram, the     commutativity of which is shown  before Lemma \ref{f:decent}, stems from the   discrepancy that the reciprocity law has been formulated and proved in the setting of $K_{\psi,U'}^\bullet({D}^\dagger_{rig}(V(\tau^{-1}))^\Delta)[d-1] $ while the regulator map originally lives in the setting of $K_{\psi,U'}^\bullet({D} (V(\tau^{-1}))^\Delta)[d-1] $:
\begin{equation}\label{f:diagcup3}
\xymatrix{
\cC^\bullet(G_{L'},V^*(1))\ar@{-->}[d]^{\simeq}\ar@{}[r]|(0.35){\times} & \cC^\bullet(G_{L'},V) \ar[r]^-{\cup_{G_{L'}}}\ar@{-->}[d]^{\simeq}_{e} & \cC^\bullet(L',L(1))
\ar@{-->}[r]^-{\mathrm{tr}_{\cC}} & L[-2]\ar@{=}[d]\\
K^\bullet_{\varphi,U'}(M^\Delta)
\ar@{}[r]|(0.35){\times} & K^\bullet_{\psi,U'}(D(V(\tau^{-1}))^\Delta)[d-1] 
\ar[rr]^-{\cup_{K,\psi}} & & L[-2]\ar@{=}[d]\\
K^\bullet_{\varphi,U'}(({M}^\dagger)^\Delta)\ar@{^(->}[d]^{\simeq} \ar@{^(->}[u]\ar@{}[r]|(0.35){\times} & K^\bullet_{\psi,U'}(D^\dagger(V(\tau^{-1}))^\Delta)[d-1] \ar@{^(->}[u]\ar@{^(->}[d] \ar[rr]^-{\cup_{K,\psi}} &  & L[-2]\ar@{=}[d] \\  
K^\bullet_{\varphi,U'}(({M}^\dagger_{rig})^\Delta) \ar@{}[r]|(0.35){\times} & K_{\psi,U'}^\bullet({D}^\dagger_{rig}(V(\tau^{-1}))^\Delta)[d-1] \ar[rr]^-{\cup_{K,\psi}} & & L[-2] \\}
\end{equation}
which in turn induces the commutativity of the lower rectangle in
the following  diagram (the upper rectangles commute obviously)
\begin{equation}\label{f:PsiSlashdagger}
\tiny
\xymatrix{
  {D}^\dagger_{rig}(V(\tau^{-1}))^{\psi=1} \ar[d]_{pr_U }  &  {D}^\dagger (V(\tau^{-1}))^{\psi=1} \ar[d]_{pr_U }  \ar@{_(->}[l] \ar@{^(->}[r]^-{a}  &  {D}(V(\tau^{-1}))^{\psi=1} \ar[d]^{pr_U } \\
 h^1\left( K_{\psi,U'}^\bullet({D}^\dagger_{rig}(V(\tau^{-1}))^\Delta)[d-1]\right) \ar[d]_{\cong }^{a}   &  h^1\left( K_{\psi,U'}^\bullet({D}^\dagger(V(\tau^{-1}))^\Delta)[d-1]\right)  \ar[l]^{ } \ar[r]^-{b } &  h^1\left( K_{\psi,U'}^\bullet({D}(V(\tau^{-1}))^\Delta)[d-1]\right)  \ar[d]^{\cong }_c \\
  H^1_{/\dagger}(L',V)   &  &  H^1(L',V)    \ar@{->>}[ll]_-{pr}}
\end{equation}
Here the vertical maps $pr_U$ are defined as in \eqref{f:Fontaine-SV}, $a$ and $pr$ are taken from Prop.\  \ref{prop:Koszul-product} while the isomorphism $c$ stems from \eqref{f:KpsiKupferer}.
The map $a$ is bijective under the assumption of Lemma \ref{notrivialquot}, which extends to the map $b$ by Corollary \ref{cor:Psidaggerquasi}.

} 

\section{Iwasawa cohomology and descent}

In this subsection we recall a crucial observation from \cite{Ku,KV}, which is based on  \cite{Ne} and  generalizes \cite[Thm.\ 5.13 ]{SV15}. As before let $U$ be an open subgroup of $\Gamma_L.$
We set $\mathbb{T}:=\Lambda(U)\otimes_{o_L}T$ with actions by $\Lambda(U):= o_L[[U]] $ via left multiplication on the left factor and by $g\in G_{L'}$ given as $\lambda\otimes t\mapsto \lambda \bar{g}^{-1}\otimes g(t)$, where $\bar{g}$ denotes the image of $g$ in  $U.$  We write
  $R\Gamma_{Iw}(L_\infty/L',T) $  for the continuous cochain complex $C^\bullet(U, \mathbb{T} )$ and recall that its cohomology identifies with $H^\bullet_{Iw}(L_\infty/L',T)$ by \cite[Lem.\ 5.8]{SV15}. For  any continuous endomorphism $f$ of $M$, we set  $\mathcal{T}_{f}(M):=[M\xrightarrow{f-1}M],$ a complex concentrated in degree $0$ and $1$.

The map $p:\mathbb{T}\to o_L\otimes_{\Lambda(U)} \mathbb{T}\cong T,\; t\mapsto 1\otimes t,$ and its dual $i:T^\vee(1)\to \mathbb{T}^\vee(1)$ induce  on cohomology the corestriction and restriction map, respectively, and they are linked by the following commutative diagram
\begin{equation}\label{f:Iwdual}
\xymatrix{
{\cC^\bullet(G_{L'},\mathbb{T})}\ar[d]^{p_*}\ar@{}[r]|(0.45){\times} & \cC^\bullet(G_{L'},\mathbb{T}^\vee(1)) \ar[r]^-{\cup_{G_{L'}}} & \cC^\bullet(L',L/o_L(1)) \ar@{-->}[r]^-{\mathrm{tr}_{\cC}} & L/o_L[-2]\ar@{=}[d]\\
\cC^\bullet(G_{L'},{T})\ar@{}[r]|(0.45){\times} & \cC^\bullet(G_{L'},\mathbb{T}^\vee(1)) \ar[u]^{i_*} \ar[r]^-{\cup_{G_{L'}}} & \cC^\bullet(L',L/o_L(1)) \ar@{-->}[r]^-{\mathrm{tr}_{\cC}} & L/o_L[-2]}
\end{equation}

By \cite[Prop.\ 1.6.5 (3)]{FK} (see also \cite[(8.4.8.1)]{Ne}) we have a canonical isomorphism
\begin{equation}\label{f:Iwdescent}
o_L\otimes^\mathbb{L}_{\Lambda(U)} R\Gamma(L',\mathbb{T})\cong R\Gamma({L'},o_L\otimes_{\Lambda(U)}\mathbb{T})\cong R\Gamma({L'},{T})
\end{equation}
where we denote by $R\Gamma(L',-) $ the complex $\cC^\bullet(G_{L'},-)$ regarded as an object of the derived category. Dually, by a version of Hochschild-Serre, there is a canonical isomorphism
\begin{equation}\label{f:discretedescent}
R\Hom_{\Lambda}(o_L,R\Gamma(L',\mathbb{T}^\vee(1)) )\cong R\Gamma(L',T^\vee(1)).
\end{equation}
It follows that the isomorphism
\[R\Gamma_{Iw}(L_\infty/L',T)\cong R\Hom_{o_L}(R\Gamma(L',\mathbb{T}^\vee(1)),L/o_L)[-2]\]
induced by the upper line of \eqref{f:Iwdual} induces an isomorphism
\begin{equation}
\label{f:Iwdualdiscrete}
o_L\otimes^\mathbb{L}_{\Lambda(U)} R\Gamma_{Iw}(L_\infty/L',T)\cong R\Hom_{o_L}(R\Hom_{\Lambda}(o_L,R\Gamma(L',\mathbb{T}^\vee(1)) ),L/o_L)[-2] ,\end{equation}
which is compatible  with the lower cup product pairing in \eqref{f:Iwdual} via the canonical identifications \eqref{f:Iwdescent} and \eqref{f:discretedescent}.

\begin{lemma}
There is a canonical isomorphism $R\Gamma(L',\mathbb{T}^\vee(1))\cong \mathcal{T}_\varphi(D(T^\vee(1)))$ in the derived category.
\end{lemma}

\begin{proof}
See \cite[Thm.\ 5.1.11]{KV}.
\end{proof}

For the rest of this section we assume that $U\subseteq \Gamma_L$ is an open \emph{torsionfree} subgroup.

\begin{lemma}\label{lem:TauphiTaupsi}
Let $T$ be in $\Rep_{o_L}(G_L)$ of finite length. Set $\Lambda:=\Lambda(U)$ and let $\gamma_1,\ldots, \gamma_d$ be topological generators of $U.$ Then we have a up to signs  canonical isomorphism of complexes
\begin{align*}
 \Hom_{\Lambda}^\bullet(K_\bullet(\gamma),\mathcal{T}_\varphi(D(T^\vee(1))))^\vee[-2]&\cong \mathrm{tot}\left(\mathcal{T}_\psi(D(T (\tau^{-1})))[-1]\otimes_\Lambda K_\bullet(\gamma^{-1})(\Lambda)^\bullet\right)
\end{align*}
where $-^\vee$ denotes forming the Pontrjagin dual.
\end{lemma}

\begin{proof}
  Upon noting that $\mathcal{T}_\varphi(D(T^\vee(1)))^\vee[-2]\cong \mathcal{T}_\psi(D(T (\tau^{-1})))[-1]$ (canonically up to a sign!) this is easily reduced to the following statement
\begin{align*}
 \Hom_{\Lambda}^\bullet(K_\bullet(\gamma),M)^\vee&\cong M^\vee \otimes_\Lambda K_\bullet(\gamma^{-1})(\Lambda)^\bullet,
\end{align*}
which can be proved in the same formal way as \eqref{f:selfualM}, and a consideration of signs.
\end{proof}

\begin{remark}\label{rem:KphiKpsi}
For every $M\in  \mathfrak{M}(\mathbf{A}_L)$ we have  a canonical isomorphism \[ \Hom_{\Lambda}^\bullet(K_\bullet^U,\mathcal{T}_\varphi(M))\cong K_{\varphi,U}(M)\] up to the sign $(-1)^n$ in degree $n$ and a non-canonical isomorphism  \[\mathrm{tot}\left(\mathcal{T}_\psi(M)[-1]\otimes_\Lambda K_\bullet^U(\Lambda)^\bullet\right)\cong  K_{\psi,U}(M)[d-1] \] (involving the self-duality of the Koszul complex). Here, the right hand sides are formed with respect to the same sequence of topological generators as the left hand sides.
\end{remark}

\begin{proof} By our conventions in section \ref{sec:homol}
$ K_{\varphi,U}(M)$ is the total complex of the double complex $\Hom^\bullet(K_\bullet(\Lambda)^\bullet,M)\xrightarrow{1-\varphi_*}\Hom^\bullet(K_\bullet(\Lambda)^\bullet,M)$. A comparison with the total Hom-complex (with the same sign rules as in secton \ref{sec:homol}) shows the first claim. For the second statement we have
\begin{align*}
  \mathrm{tot}\left(\mathcal{T}_\psi(M)[-1]\otimes_\Lambda K_\bullet(\Lambda)^\bullet\right) & \cong \mathrm{tot}\left(\mathcal{T}_\psi(M)\otimes_\Lambda K_\bullet(\Lambda)^\bullet\right)[-1] \\
&=\mathrm{tot}\left( \mathcal{T}_{\psi}\left( M\otimes_{\Lambda} K_\bullet(\Lambda)^\bullet \right) \right)[-1]\\
&\cong \mathrm{tot}\left( \mathcal{T}_{\psi}\left( M\otimes_{\Lambda} K^\bullet(\Lambda)[d] \right) \right)[-1]\\
&=\mathrm{tot}\left( \mathcal{T}_{\psi}\left( K^\bullet(M)[d]\right) \right)[-1]\\
&=\mathrm{cone}\left( K^\bullet_{U}(M)[d]\xrightarrow{1-\psi}K^\bullet_{U}(M)[d] \right)[-2]\\
&\cong K_{\psi,U}(M)[d-1].
\end{align*}
The first isomorphism involves a sign on $\mathcal{T}^1_\psi(M).$ The third isomorphisms stems from \eqref{f:selfdual} while the last isomorphism again involves signs.
\end{proof}

\begin{theorem}
There are canonical isomorphisms
\begin{align}
R\Gamma_{Iw}(L_\infty/L,T)&\cong \mathcal{T}_{\psi}\left( D(T(\tau^{-1}))\right)[-1]\\
\label{f:KpsiKupferer}
 K_{\psi,U}(D(T(\tau^{-1})))[d-1]&\xrightarrow{\simeq} R\Gamma({L'},T).
\end{align}
in the derived category $D_{perf}(\Lambda_{o_L}(\Gamma_L))$ of perfect complexes and in
the derived category $D^+(o_L)$ of bounded below cochain complexes of
${o_L}$-modules, respectively.
\end{theorem}

\begin{proof}
The first isomorphism is  \cite[Thm.\ 5.2.54 ]{KV} while the second one follows from this and \eqref{f:Iwdescent} as \Footnote{Do we need  the characterization $D_{perf}(\Lambda_{o_L}(U))=D^b_{fg}(\Lambda_{o_L}(U))$ from \cite[(4.2.8)]{Ne} as derived category of bounded complexes of $\Lambda_{o_L}(U)$-modules with finitely generated cohomology over $\Lambda_{o_L}(U)$?}
\begin{align*}
R\Gamma_{Iw}(L_\infty/L,T)\otimes^\mathbb{L}_{\Lambda_{o_L}(U)} o_L&\cong \mathcal{T}_{\psi}\left( D(T(\tau^{-1}))\right)[-1]\otimes^\mathbb{L}_{\Lambda} K_\bullet(\Lambda)^\bullet\\
&=\mathrm{tot}\left( \mathcal{T}_{\psi}\left( D(T(\tau^{-1}))\right)[-1]\otimes_{\Lambda}  K_\bullet(\Lambda)^\bullet \right)\\
&=K_{\psi,U}(D(T(\tau^{-1})))[d-1].
\end{align*}
by Remark \ref{rem:KphiKpsi}.
\end{proof}

By Lemma \ref{lem:TauphiTaupsi} and Remark \ref{rem:KphiKpsi} we see that, for $T$ be in $\Rep_{o_L}(G_L)$ of finite length,
\begin{equation}\label{f:0}
 K_{\varphi,U}(D(T^\vee(1)))= R\Hom_{\Lambda}(o_L,\mathcal{T}_\varphi(D(T^\vee(1)))))[2]
\end{equation}
is dual to
\begin{equation}\label{f:00}
 K_{\psi,U}(D(T (\tau^{-1})))= o_L\otimes^\mathbb{L}_{\Lambda(U)}\mathcal{T}_\psi(D(T (\tau^{-1})))[-1],
\end{equation}
such that the upper rectangle in the diagram \eqref{f:diagcup3} commutes by \eqref{f:Iwdualdiscrete}, taking inverse limits and inverting $\pi_L$.

\begin{lemma}\label{f:decent}Let $T$ be in $\Rep_{o_L}(G_L).$
 Then the left rectangle in \eqref{f:corpr} is commutative.
\end{lemma}

\begin{proof} (Sketch)
 By an obvious analogue of Remark \ref{rem:res-cores} it suffices to show the statement
 for $U=\Gamma_n\cong \mathbb{Z}_p^d.$ In this situation we have a homological spectral
 sequence
 \[H_{i,cts}(U, H^{-j}_{Iw}(L_\infty/L,T )) \Longrightarrow
 H_{cts}^{-i-j}(L',T)\]
 which is induced by \eqref{f:Iwdescent},
 see \cite[(8.4.8.1)]{Ne} for the statement and missing
notation. We may and do assume that $T$ is of finite length. Then, on the one hand, the map $H^1_{Iw}(L_\infty/L,T)\xrightarrow{cor}H^1({L'},T)$ is dual to $H^1({L'},T^\vee(1))\xrightarrow{res} H^1(L_\infty, T^\vee(1)),$ which sits in the five term exact sequence of lower degrees associated with the Hochschild-Serre spectral sequence. As explained just before this lemma the above homological spectral sequence arises by dualizing from the latter. Hence $cor$ shows up in the five term exact sequence of lower degrees associated with this homological spectral sequence.  On the other hand via the isomorphisms \eqref{f:Iwdescent} and \eqref{f:KpsiKupferer} the latter spectral sequence is isomorphic to
\[H_{i,cts}(U, h^{-j}(\mathcal{T}_{\psi}\left( D(T(\tau^{-1}))\right)[-1]) \Longrightarrow
 h^{-i-j}(K_{\psi,U}(D(T(\tau^{-1})))[d-1] )\] and one checks by inspection that $cor$ corresponds to $pr_{U}.$
\end{proof}

%
%
%
%

%

\Footnote{Since $\Omega^1\xrightarrow{Res}o_L $ is obviously surjective, the map $Res$ induces an isomorphism $\Omega^1/(\psi_L-1)\cong o_L.$ Hence, for $T=o_L(1),$ we obtain an isomorphism $H^2(K_{\psi,U}[d-1])\cong \left(\Omega^1/(\psi_L-1)\right)_{U}=\Omega^1/(\psi_L-1)\cong o_L,$ which up to a scalar has to coincide with Tate's invariant/trace map. One could try to use the methods from \cite[\S 2.2]{benois} to determine this scalar. There the $H^2$ of the $\varphi_L$-Herr complex is considered instead! So probably one would have to calculate the cup-product for the $\psi_L$-Herr complex to apply his methods!?}

\noindent
Peter Schneider,\\
Universit\"{a}t M\"{u}nster,  Mathematisches Institut,\\  Einsteinstr. 62,
48291 M\"{u}nster,  Germany,\\
 http://www.uni-muenster.de/math/u/schneider/  \\
pschnei@uni-muenster.de \\ \\

\noindent
Otmar Venjakob\\
Universit\"{a}t Heidelberg,  Mathematisches Institut,\\  Im Neuenheimer Feld 288,  69120
Heidelberg,  Germany,\\
 http://www.mathi.uni-heidelberg.de/$\,\tilde{}\,$venjakob/\\
venjakob@mathi.uni-heidelberg.de

\printendnotes[custom]

\end{document}